Alexander A. Ermolitski

# RIEMANNIAN MANIFOLDS WITH GEOMETRIC STRUCTURES

# PREFACE

The theory of structures on manifolds is a very interesting topic of modern differential geometry and its applications.

There are many results concerning various differential geometric structures on Riemannian manifolds.

The main aim of this book is to get a way of a union of such results in one scheme. It seems that introduced by the author a notion of the canonical connection $\overline{\nabla}$ and the second fundamental tensor field *h* adjoint to a structure is very useful for this purpose and, in many cases, it is more effective than the Riemannian connection $\nabla$. Especially, we pay attention to use of *h* to obtain classifications of structures and to the case of so-called quasi homogeneous structures.

Projections of structures on submanifolds are also considered in the book.

The introduced by the author class of Riemannian (locally) regular $\sigma-$manifolds is studied here too.

Further, the book is not a survey of what has been done in the theory of structures on manifolds, it is concerning only to restricted subjects.

The list of references is not intended to be a collection of prestigious papers and famous names. The references are simply limited to those closely connected with the topic. Thus, some works quoted here are important and some others may be not.

Now, let us sum up briefly the contents:

In ***Chapter 1***, we have given a short survey of *G*-structures, associated Riemannian metrics and have introduced the so-called canonical connection and the second fundamental tensor field $h = \nabla - \overline{\nabla}$ of a fixed pair *(P(G), g)*, where $\nabla$ is the Riemannian connection of a Riemannian metric *g* associated to a *G*-structure *P(G)*. The integrability of *G*-structures and the polar decomposition of an *O*-deformable *(1,1)* tensor field are considered here too.

In ***Chapter 2***, using $\overline{\nabla}$ and *h* we have got some results which are generalisations of those obtained by numerous authors for concrete structures. We consider torsion and curvature of $\overline{\nabla}$, so-called quasi homogeneous structures, isometries, affine transformations and holonomy fibre bundles of $\overline{\nabla}$. Homogeneous structures and projections of structures on submanifolds and foliations are also discussed.

In ***Chapter 3***, we have obtained a classification $\boldsymbol{T}=\boldsymbol{T}_1\oplus\boldsymbol{T}_2\oplus\boldsymbol{T}_3$ of *G*-structures over Riemannian manifolds as a decomposition on invariant irreducible subspaces of tensors of type *h* in $\otimes^3 T^*$ under the natural action of the orthogonal group. Structures of types $\boldsymbol{T}_1$, $\boldsymbol{T}_3$ have been studied more explicitly. Some algebraic construction has been considered for quasi homogeneous structures of type $\boldsymbol{T}_3$.

Naturally reductive homogeneous structures and nearly Kaehlerian manifolds are examples of such structures.

In ***Chapter 4***, we have introduced so-called Riemannian (locally) regular $\sigma$–manifolds which generalise on the one hand the spaces with reflections of O. Loos and on the other hand the Riemannian regular *s*-manifolds and have proved that every regular $\sigma$-manifold can be described as a fibre bundle over a regular *s*-manifold. Conversely, such fibre bundles and regular *s*-manifolds of order *2k* give examples of regular $\sigma$–manifolds.

We also consider the Lie algebra of infinitesimal automorphisms of Riemannian regular $\sigma$–manifolds, orbits under the action of the structural group *G*, the structure of locally regular $\sigma$–manifolds and submanifolds.

In ***Chapter 5***, we have our methods illustrated for the following structures:

a.p.R.s., i.e., almost-product Riemannian structure $(P, g)$, where $P^2 = I$;

a.H.s., that is, almost Hermitian structure $(J, g)$, where $J^2 = -I$;

f-s., i.e., a structure defined by an affinor $F$, $F^3 + F = O$, and an associated metric $g$.

For each of them the canonical connection $\overline{\nabla}$ and the second fundamental tensor field *h* have been computed and applied to study of geometry of manifolds with a.p.R.s., a.H.s., f-s. Integrability conditions, conformal changes of the Riemannian metric *g*, the parallel translation of structures along curves etc. are discussed here too.

We have the classification of A.Gray and L.M.Hervella rewritten in terms of tensor field *h*, i.e., the canonical connection have been obtained for each from 16 classes. We can apply these connections for example to construct characteristic classes.

A.p.R.s., a.H.s., f-s. can be appeared with help of an affinor *S* defined on a Riemannian locally regular *s* (or $\sigma$) –manifold. Conversely, we consider some conditions, when a.p.R.s., a.H.m., f-s. given on a Riemannian manifold *M* are indused such a tensor field *S* , i.e., M has a structure of a Riemannian regular s (or $\sigma$) –manifold.

In ***Chapter 6***, using the classification of A.Gray and L.M.Hervella and constructed the second fundamental tensor field *h* of almost contact metric structure we have obtained in terms of *h* a classification of such structures. There are $2^{12}$ classes. A similar classification was got by D.Chinea and C.Gonzalez (A.A.Alexiev and G.Ganchev) by a different method. Good relations have been found between both the classifications and this allows to apply various tensor characteristics of the canonical connection $\overline{\nabla}$ to study of every class.

Examples of all the classes adduced in classification Table are given in ***Chapter 7***. In particular, $\alpha$-Sasakian and $\alpha$-Kenmotsu structures are identified and it is shown that, when $\alpha = const$, they are quasi homogeneous. Further, the conditions of integrability, normality and the fundamental tensor fields $N^{(1)}, N^{(2)}, N^{(3)}, N^{(4)}$ are considered. We identify some of the classes studied by various authors with those obtained from classification Table. Riemannian locally regular $\sigma$–manifolds (R.l.r. $\sigma$–m.) with one-dimensional foliations of mirrors are discussed here too. We consider necessary and sufficient conditions for $M$ to be a R.l.r. $\sigma$–m., and, also, the induced almost contact metric structures (a.c.m.s.). In this case the canonical connection $\tilde{\nabla}$ of R.l.r. $\sigma$–m. and that $\overline{\nabla}$ of the induced a.c.m.s. are the same. R.l.r. $\sigma$–m. of order 3,4 are studied more explicitly.

Furher, we want to give the diagram of dependence of the Chapters

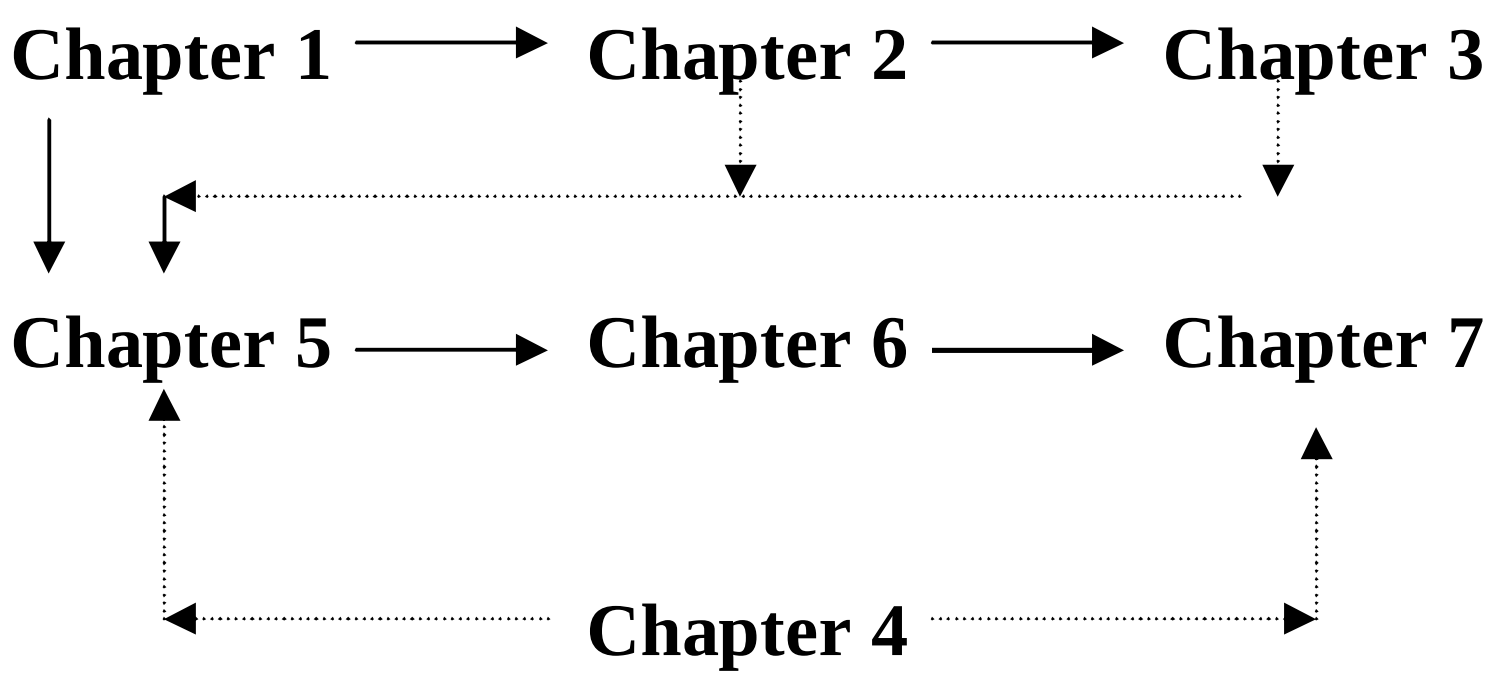

We remark that all the theorems, propositions, formulas and definitions are numbered for each chapter. For example, in each chapter, say, Chapter 5, Theorem 1.2 is in Chapter 1.

All the manifolds, maps, tensor fields etc. are supposed to be of class $C^{\infty}$. We consider connected Riemannian manifolds if not otherwise stated.

The book assumes a basic knowledge of the modern differential geometry (cf. Foundations of Differential Geometry by S.Kobayashi and K.Nomizu, Vol. I, II).

The author wish to express here the gratitude to his brother G.A.Ermolitski for support of this work.

Minsk, 1998 Alexander A. Ermolitski

# TABLE OF CONTENTS

## CHAPTER 5
## CLASSICAL STRUCTURES

## CHAPTER 6
## CLASSIFICATION OF ALMOST CONTACT METRIC STRUCTURES

## CHAPTER 7
## REMARKS ON GEOMETRY OF ALMOST CONTACT METRIC MANIFOLDS

# CHAPTER 1

# SECOND FUNDAMENTAL TENSOR FIELD OF G – STRUCTURE

In this chapter, we have given a brief survey of $G$ - structures and close notions. In §1, we give the basic definitions of fibre bundles, $G$ - structures, tensor $G$ - structures defined by $O$-deformable tensor fields and connections. In §2, associated Riemannian metrics are discussed and we define introduced by the author a notion of the second fundamental tensor field $h$ of a Riemannian $G$ - structure. The integrability of $G$ - structures and the polar decomposition of $O$ - deformable (1,1) tensor field are considered in §3.

We follow fairly closely [45], [46], [73].

## *§1. G - STRUCTURES AND CONNECTIONS*

**$1^0$** . Let $M$ be a differentiable manifold and $G$ be a Lie group.

**DEFINITION 1.1**. A principal fibre bundle over $M$ with group $G$ consists of a manifold $P$ and an action of $G$ on $P$ satisfying the following conditions:

(1) $G$ acts freely on $P$ on the right: $(u,a)\in P\times G \longrightarrow ua = R_a u \in P$;

(2) $M$ is the quotient space of $P$ by the equivalence relation induced by $G$, $M = P/G$, and the canonical projection $\pi : P \to M$ is differentiable;

(3) $P$ is locally trivial, that is, every point $x$ of $M$ has a neighbourhood U such that $\pi^{-1}(U)$ is isomorphic with $U\times G$ in the sense that there is a diffeomorphism $\Psi : \pi^{-1}(U) \to U\times G$ such that $\Psi(u) = (\pi(u), \Phi(u))$, where $\Phi$ is a mapping of $\pi^{-1}(U)$ into $G$ satisfying $\Phi(ua) = (\Phi(u))\, a$ for all $u\in\pi^{-1}(U)$ and $a\in G$.

A principal fibre bundle will be denoted by $P(M,G,\pi)$, $P(G)$ or simply $P$. We call $P$ the total space or the bundle space, $M$ the base space, $G$ the structure group and $\pi$ the projection. For each point $x$ of $M, \pi^{-1}(x)$ is a closed submanifold of $P$, called the fibre over $x$. If $u$ is a point of $\pi^{-1}(x)$, then $\pi^{-1}(x)$ is the set of points $ua$, $a\in G$, and is called the fibre through $u$. Every fibre is diffeomorphic to $G$.

From local triviality of *P(M,G)* we see that if *W* is a submanifold of *M*, then $\pi^{-1}(W)(W,G)$ is a principal fibre bundle. We call it the restriction of *P* to *W* and denote it by $P_W$ .

Given a principal fibre bundle *P(G)*, the action of *G* on *P* induces a homomorphism $\sigma$ of the Lie algebra $\underline{\boldsymbol{g}}$ of *G* into the Lie algebra *x(P)* of vector fields on *P*. $\sigma$ can be define as follows: for every *u*, let $\sigma_u$ be the mapping $a \in G \to ua \in P$. Then $(\sigma_u)_* A_e = (\sigma A)_u$. For each $A \in \underline{\boldsymbol{g}}$ $A^* \in \sigma(A)$ is called the fundamental vector field corresponding to *A*. Since the action of *G* sends each fibre into itself, $A_u^*$ is tangent to the fibre at each $u \in P$. As *G* acts freely on *P*, $A^*$ never vanishes on *P* if $A \neq 0$. The dimension of each fibre being equal to that of $\underline{\boldsymbol{g}}$, the mapping $A \to (A^*)_u$ of $\underline{\boldsymbol{g}}$ into $T_u(P)$ is a linear isomorphism of $\underline{\boldsymbol{g}}$ onto the tangent space at *u* of the fibre through *u*. We also see that for each $a \in G$, $(R_a)_* A^*$ is the fundamental vector field corresponding to $(ad(a^{-1}))A \in \underline{\boldsymbol{g}}$.

A homomorphism *f* of a principal fibre bundle *P'(M',G')* into another principal fibre bundle *P(M,G)* consists of a mapping $f': P' \to P$ and a homomorphism $f'': G' \to G$ such that *f'(u'a') = f'(u') f''(a)* for all $u' \in P'$ and $a' \in G'$. For the sake of simplicity, we shall denote *f'* and *f''* by the same letter *f*. Every homomorphism $f : P' \to P$ maps each fibre of *P'* into fibre of *P* and hence induces a mapping of *M'* into *M*, which will be also denoted by *f*. A homomorphism $f : P'(M',G') \to P(M,G)$ is called an imbedding or injection if the induced mapping $f : M' \to M$ is an imbedding and if $f : G' \to G$ is a monomorphism. By identifying *P'* with *f(P')*, *G'* with *f(G')* and *M'* with *f(M')*, we say that *P'(M',G')* is a subbundle of *P(M,G)*. If, moreover, *M' = M* and the induced mapping $f : M' \to M$ is the identity transformation of *M*, $f : P'(M',G') \to P(M,G)$ is called a reduction of the structure group *G* of *P(M,G)* to *G'*. The subbundle *P'(M,G')* is called a reduced bundle. Given *P(G)* and a Lie subgroup *G'* of *G*, we say that *G* is reducible to *G'* if there is a reduced bundle *P'(G')*.

**$2^0$**. Let *P(M,G)* be a principal fibre bundle and $\Lambda$ a manifold on which *G* acts on the left: $(a,\lambda) \in G \times \Lambda \to a\lambda \in \Lambda$. A construction of a fibre bundle *E(M,Λ,G,P)* associated with *P* with standard fibre $\Lambda$ is considered below. On the product manifold $P \times \Lambda$ the group *G* acts on the right as follows: $a \in G$ maps $(u,\lambda) \in P \times \Lambda$ into $(ua,a^{-1}\lambda) \in P \times \Lambda$. The quotient space of $P \times \Lambda$ by this group action is denoted by $E = P \times_G \Lambda$. The mapping $P \times \Lambda \to M$ which maps *(u,λ)* into $\pi(u)$ induces a mapping $\pi_E$, called the projection, of *E* onto *M*. For each $x \in M, \pi_E^{-1}(x)$ is called the fibre of $E$ over $x$. Every point $x$ of $M$ has a neighbourhood $U$ such that $\pi^{-1}(U)$ is isomorphic to $U \times G$. Identifying $\pi^{-1}(U)$ with $U \times G$, one can see that the action of $G$ on $\pi^{-1}(U) \times \Lambda$ on the right is given by $(x,a,\lambda) \to (x,ab,b^{-1}\lambda)$ for $(x,a,\lambda) \in U \times G \times \Lambda$ and $b \in G$.

The isomorphism $\pi^{-1}(U) \cong U \times G$ induces an isomorphism $\pi_E^{-1}(U) \cong U \times \Lambda$. Therefore a differentiable structure in $E$ can be introduced by the requirement that $\pi_E^{-1}(U)$ is an open submanifold of $E$ which is diffeomorphic with $U \times \Lambda$ under the isomorphism $\pi_E^{-1}(U) \cong U \times \Lambda$. The projection $\pi_E$ is then a differentiable mapping of $E$ onto $M$. Thus $E(M,\Lambda,G,P)$ or simply $E$ is called the fibre bundle over the base space $M$, with standard fibre $\Lambda$ and structure group $G$, which is associated with the principal fibre bundle $P$.

**3⁰** . We consider now $G$-structures. Let $M$ be an $n$-dimensional manifold. A linear frame $u$ at a point $x$ of $M$ is an ordered basis $X_1, \ldots, X_n$ of $T_x(M)$. Let $L(M)$ be the set of all linear frames $u$ at all points of $M$ and let $\pi$ be the mapping of $L(M)$ onto $M$ which maps a linear frame $u$ at $x$ into $x$. The general linear group $GL(n,R)$ acts on $L(M)$ on the right as follows. Let $a = (a_j^i) \in GL(n,R)$ and $u = (X_1, \ldots, X_n)$ be a linear frame at $x$. Then $ua$ is the frame $(Y_1, \ldots, Y_n)$ at $x$ defined by $Y_i = \sum_j a_j^i X_j$. It is well known that $L(M)$ $(M, GL(n,R))$ is a principal fibre bundle. We call $L(M)$ the bundle of linear frames over $M$. Let $e_1, \ldots, e_n$ be the natural basis for $R^n : e_1 = (1,0,\ldots,0), \ldots, e_n = (0,0,\ldots,1)$. A linear frame $u = (X_1, \ldots, X_n)$ at $x$ can be given as a linear mapping $u : R^n \to T(M)$ such that $ue_i = X_i$ for $i = 1,\ldots,n$.

**DEFINITION 1.2**. Let $G$ be a Lie subgroup of $GL(n,R)$. A reduction $P(G)$ of $L(M)$ to $G$ is called $G$ - structure over a manifold $M$.

We say that $G$ - structure $P(G)$ is conjugate with $G'$ – structure $P'(G')$ if there exists an element $a \in GL(n,R)$ that $P' = Pa$. In this case we have $G' = a^{-1}Ga$. $P'(G')$ will be denoted by $P_a(a^{-1}Ga)$. Let for $G$ - structure $P(G)$ $G$ be a Lie subgroup of a Lie group $G'$. Then it is easy to build the unique $G'$- structure $P'(G') \supset P(G)$ which is called the extension of $P(G)$ to $G'$.

**4⁰**. If $p$ is a point in $M$, we define $T_s^r(p)$ as the set of all $R$-multilinear

mappings of

$$T_p^*(M)\times\ldots\times T_p^*(M)\times T_p(M)\times\ldots\times T_p(M)$$

($T_p^*(M)$ - $r$ times, $T_p(M)$ - $s$ times)

into *R*. Otherwise $T_s^r(p)$ is the tensor product

(1.1) $T_s^r(p)=\otimes^r T_p(M)\otimes^s T_p^*(M)$.

If $K\in T_s^r$ is a tensor field on *M* of type $(r,s)$, $e_1,\ldots..,e_n$ is a basis of $T_p(M)$, $e^1,\ldots..,e^n$, is the dual basis of $T_p^*(M)$, then we have

(1.2) $K_p=k_{\mu_1\ldots\mu_s}^{\lambda_1\ldots\lambda_r}e_{\lambda_1}\otimes\ldots\otimes e_{\lambda_r}\otimes e^{\mu_1}\otimes\ldots\otimes e^{\mu_s}$

If we consider $T_s^r(R^n)$, then one can define an equivariant with respect to the group *GL(n,R)* map

$$K: L(M)\to T_s^r(R^n)$$

by the formula

$$K(n)(v_1,\ldots,v_s,\omega^1,\ldots,\omega^r)=K(uv_1,\ldots,uv_s,u^*\omega^1,\ldots,u^*\omega^r)$$

where $u\in L(M)$.

The action of *GL(n,R)* on *R* induces that on $T_s^r(R^n)$.

**DEFINITION 1.3**. A tensor field *K* is said to be *O*-deformable if *K(L(M))* belongs to one single orbit in $T_s^r(R^n)$ with respect to *GL(n,R)*. Let ***O*** be such an orbit and $K_0\in\boldsymbol{O}$. Then $K^{-1}(K_0)\subset L(M)$ is a *G*-structure, where *G* is the automorphism group of $K_0$. Selecting $K_0$. is the same as selecting $n^{r+s}$ components $k_{\mu_1\ldots\mu_s}^{\lambda_1\ldots\lambda_r}$ for *K*.

For any frame $u'\in L(M)$

$$K(u')=K_0'\in O\subset T_s^r(R_n)$$

and $K_o'$ is in the same *GL(n,R)* orbit ***O*** as $K_0$. So, there is $a \in GL(n,R)$ such that $K_0 = K_0' a$, then $K(u'a) = K_0$ and at any point $x \in M$ there is at least one frame *u'* belonging to $K^{-1}(K_0)$. For $b \in G$ we have $K(ub) = K(u)b = K_0 b = K_0$, because *b* is an element of the automorphism group *G* . Let $(a^i_j) \in GL(n,R)$ and $e_1,...,e_n$ be a basis of $T_p(M)$. We define the group *G* by

$$(1.3)\quad G = \{ (a^i_j) \in GL(n,R) : k^{j_1 \dots j_r}_{\mu_1 \dots \mu_s} a^{\lambda_1}_{j_1} \cdot \dots \cdot a^{\lambda_r}_{j_r} = k^{\lambda_1 \dots \lambda_r}_{k_1 \dots k_s} a^{k_1 \dots k_s}_{\mu_1 \dots \mu_s} \}$$

The group *G* is an algebraic group therefore *G* is a closed Lie group.

**DEFINITION 1.4**. A *G* - structure *P(G)* with a structure group of type (1.3) is called a tensor *G* - structure. The frames of *P(G)* are the same for which *K* has the components $k^{j_1 \dots j_r}_{\mu_1 \dots \mu_s}$ on *M*.

We remark that a *O*-deformable tensor field *K* on *M* defines the whole class of conjugated each other *G*-structures. It is clear that for $a \in GL(n,R)$ the tensor $K_0 a$ still belongs to the orbit $O \subset T^r_s(R^n)$ but the automorphism group can change *G'* instead of *G*, where $G' = a^{-1} G a$. So, the choice of structure from the class depends on the choice of a point $p \in M$ and a frame $(e_1, e_2, ..., e_n)$ at *p*.

For instance, if *F* is an *O*-deformable (1,1) tensor field, then at each point $p \in M$ there exists a frame $(e_1,...,e_n)$ for which $F_p$ has the Jordan normal form. We can take a *G*-structure consisting of all such frames from *L(M)*. Later on we shall denote by *P(K)* a *G*-structure defined by an *O*-deformable tensor field *K* if the choise of frames is clear from context.

In 196O I.M.Singer put up a question if there exist locally nonhomogeneous Riemannian manifolds for which the Riemannian metric *g* and the curvature tensor field *R* are simultaneously *O*-deformable. Many papers have been published in this direction; for a survey see [12].

Examples of various *G*-structures one can find in [6], [13], [45], [64].

**5⁰**. We define a connection in a principal fibre bundle (see [46],[73]). For each $u \in P$, let $T_u(P)$ be the tangent space of *P* at *u* and $G_u$ the subspace of $T_u(P)$ consisting of vectors tangent to the fibre through *u*.

**DEFINITION 1.5**. A connection $\Gamma$ in $P$ is an assignment of a subspace $Q_u$ of $T_u(P)$ to each $u \in P$ such that

(a) $T_u(P) = G_u \oplus Q_u$ ;

(b) $Q_{ua} = (R_a)_* Q_u$ for every $u \in P$ and $a \in G$ ;

(c) $Q_u$ depends differentiably on *u*.

Condition (b) means that the distribution $u \to Q_u$ is invariant by *G*. We call $G_u$ the vertical subspace and $Q_u$ the horizontal subspace of $T_u(P)$. A vector $X \in T_u(P)$ is called vertical (respectively horizontal) if it lies in $G_u$ (respectively $Q_u$). By (a), any vector *X* of *P* at *u* can be uniquely written as

$$X = Y + Z, \quad Y \in G_u, Z \in Q_u .$$

We call *Y* (respectively *Z*) the vertical (respectively horizontal) component of *X*.

Given a connection $\Gamma$ in *P*, we define a 1-form $\omega$ on *P* with values in the Lie algebra $\underline{\mathfrak{g}}$ of *G* as follows. For each $X \in T_u(P)$ we define $\omega(X)$ to be the unique $A \in \underline{\mathfrak{g}}$ such that $(A^*)_u$ is equal to the vertical component of *X*, where $A^*$ is the fundamental vector field corresponding to *A*. It is clear that $\omega(X) = O$ if and only if *X* is horizontal. The form $\omega$ is called the connection form of the given connection $\Gamma$. The projection $\pi : P \to M$ induces a linear mapping $\pi : T_u(P) \to T_x(M)$ for each $u \in P$, where $x = \pi(u)$. When a connection is given, $\pi$ maps $Q_u$ isomorphically onto $T_x(M)$. The horizontal lift or simply, lift of a vector field *X* on *M* is a unique vector field $X^*$ on *P* which is horizontal and $\pi(X_u^*) = X_{\pi(u)}$ for every $u \in P$. Given a connection in *P* and a vector field *X* on *M*, there is a unique horizontal lift $X^*$ of *X* which is invariant by $R_a$ for every $a \in G$. Conversely, every horizontal vector field $X^*$ on *P* invariant by *G* is the lift of a vector field *X* on *M*.

**$6^0$.** Let *L(M)* be the bundle of linear frames over a manifold *M*. The canonical form $\Theta$ of *L(M)* is the $R^n$-valued 1-form on *L(M)* defined by

$$(1.4) \quad \Theta(X) = u^{-1}(\pi(X)) \quad \text{for} \quad X \in T_u(P)$$

where *u* is considered as a linear mapping of $R^n$ onto $T_{\pi(u)}(M)$.

**DEFINITION 1.6**. A connection in *L(M)* over *M* is called a linear connection of *M*.

Condition *(b)* of *Definition 1.5* means that any connection in *G*-structure $P \subset L(M)$ can be extended to linear connection of *M*.

We suppose that *P(H)* is a subbundle of *P(G)* and *H* is a reductive Lie subgroup of *G*, that is

(a) $\underline{\mathfrak{g}} = \underline{\mathfrak{h}} \oplus \underline{\mathfrak{m}}$;

(b) $ad_G H(\underline{\mathfrak{m}}) = \underline{\mathfrak{m}}$.

If $\omega$ is a $\underline{\mathfrak{g}}$ - valued form of a connection in *P(G)* then we have

(1.5) $\omega(X) = \omega(X)_{|\underline{\mathfrak{h}}} + \omega(X)_{|\underline{\mathfrak{m}}}$ for $X \in T(P(G))$.

It is easy to see that $\overline{\omega}(X) = \omega(X)_{|\underline{\mathfrak{h}}}$ is the $\underline{\mathfrak{h}}$ - valued form of a connection in *P(H)* (see [46]).

Let $\varphi$ be the cross section of *L(M)* over the neighbourhood *U* which assigns to each $x \in U$ the linear frame $((X_1)_x, ..., (X_n)_x)$ and let $X, Y = \sum_k f^k X_k$ be vector fields on *M*. Then, according to [7], the form $\omega$ defines an affine connection $\nabla$ by

(1.6) $\nabla_{X_x} Y = \varphi(x)\omega(\varphi_* X_x)\varphi(x)^{-1} Y_x + \sum_k (Xf^k)(x)(X_k)_x$

where $\varphi(x)$ is considered as a mapping of $R^n$ onto $T_x(M)$.

## *§2. ASSOCIATED RIEMANNIAN METRICS AND SECOND FUNDAMENTAL TENSOR FIELD OF G – STRUCTURE*

**$1^0$.** Let *M* be a Riemannian manifold with a Riemannian metric $g = \langle\, ,\, \rangle$. The Riemannian metric g defines the subbundle *O(M)* = *P(O(n))* of *L(M)* containing all the orthonormal frames over *M*.

**DEFINITION 1.7** [15]. Let *P(G)* be a *G*-structure over *M*. The Riemannian metric *g* is called an associated one if there exists such an element $P_{a_0}(a_0^{-1}O(n)a_0)$ in the class

$$\{ P_a(a^{-1}O(n)a) \}_{a\in GL(n,R)}$$

that the structures *P(G)* and $P_{a_0}(a_0^{-1}O(n)a_0)$ have a common subbundle of the frames over *M.*

It is easy to verify that the given *Definition 1.7* is equivalent to the following

**DEFINITION 1.8.** A Riemannian metric *g* is said to be associated to a *G*-structure *P(G)* if there exists such an element $a \in GL(n,R)$ such that the intersection of the Riemannian G-structure *P(O(n))* with the subbundle *P(G)a* is non-empty. In this case the intersection is a principal subbundle *P(H)*, where $H = O(n) \cap a^{-1}Ga$.

We remark that if such a structure *P(H)* exists then any element of the class $\{ P_a(a^{-1}Ha) \}_{a\in O(n)}$ satisfies *Definition 1.8.*

The following theorems describe the existence of associated Riemannian metrics.

**THEOREM 1.1** [39]. Let *G* be a real Lie group with a finite number of connected components. Then for each *G*-structure *P(G)* over *M* there exists a reduction of *G* to the maximal compact subgroup *H* of *G.*

In the class of *G*-structures conjugated with *P(G)* there exists such a structure, let it be *P(G)* itself, that its maximal compact subgroup *H* is a subgroup of *O(n)*. We extend *P(H)* to the group *O(n)* and get the structure *P(O(n))* = *O(M)* over *M,* which defines the Riemannian metric *g* = < , > on *M*. We see from *Definition 1.8* that the metric g is an associated one to *P(G).*

**THEOREM 1.2.** Let *K* be an *O*-deformable tensor field on *M* and let *P(G)* be a corresponding *G*-structure. Then

a) there exists an associated metric *g* to *P(G)*

b) Riemmannian metric *g* is an associated one if and only if there exists such an affine connection $\nabla$ on *M* that $\nabla K = O$ and $\nabla g = O$.

**Proof**. The first conclusion of the theorem follows from *Theorem 1.1* because the Lie group $G$ is an algebraic one, see (1.3), and therefore $G$ has a finite number of connected components.

We consider the second conclusion. It is well known, see [46], [73], that any principal fibre bundle $P(M,G)$ admits a connection, if $M$ is paracompact. If $\omega$ is a connection in $P(H)=P(G)\cap P(O(n))$ with corresponding covariant derivative $\nabla$, see (1.6), then the parallel translation $\tau$ along a curve segment $\gamma$ in $M$ transfers the frames of $P(H)$ onto the frames, [46]. As $K$ and $g$ have the constant components on $M$ in these frames we see that $K$ and $g$ are invariant with respect to $\tau$, that is, $\nabla K = O$, $\nabla g=O$.

Conversely, if $\nabla K = O$, $\nabla g = O$ and $\omega$ is the corresponding to $\nabla$ linear connection, then we can construct the holonomy fibre bundle $P(u_0)$ of $\omega$ passing through the fixed orthonormal frame $u_0 \in O(M)$, [46]. It is clear that $K$ has constant on $M$ components in the frames of $P(u_0)$, hence $P(u_0)$ can be extended to the tensor $G$-structure $P(K)$. $P(u_0)$ is the common subbundle of $P(K)$ and $O(M)$ and, using *Definition 1.8,* we obtain that the Riemannian metric $g$ is an associated one to the tensor structure $P(G)$.

**QED**.

We remark that a searching for "the best associated metric" to the given $G$-structure $P(G)$ is one of the interesting subjects of modern differential geometry, see [9], [1O].

**2⁰.** Later on we shall consider only the fixed pair $(P(G), g)$ or $(P(H), g)$, where $H=G\cap O(n)$.

To construct the second fundamental field h of the structure $(P(H), g)$ we consider the Lie algebras $\underline{\boldsymbol{O}}$, $\underline{\boldsymbol{h}}$ of the Lie groups $O(n)$, $H$.

We assume that $\underline{\boldsymbol{O}} \neq \underline{\boldsymbol{h}}$

**EXAMPLE**. If $P(G)=P(GL^+(n,R))$, where

$GL^+(n,R)=\{a\in GL(n,R): \det a>0\}$, then $\underline{\boldsymbol{O}} = \underline{\boldsymbol{h}}$.

Let $\tilde{g}$ be a biinvariant Killing form on $O(n)$ and let $\underline{\boldsymbol{m}}$ be the orthogonal complement of $\underline{\boldsymbol{h}}$ to $\underline{\boldsymbol{O}}$

$$\underline{\boldsymbol{m}} = \underline{\boldsymbol{h}}^{\perp} = \{\tilde{x} \in \underline{\boldsymbol{O}} : \tilde{g}(\tilde{x},\tilde{y})=0, \forall \tilde{y} \in \underline{\boldsymbol{h}}\}.$$

We have

(1.7) $\underline{\boldsymbol{O}} = \underline{\boldsymbol{h}} \oplus \underline{\boldsymbol{m}}$;.

and from the biinvariance of $\tilde{g}$

(1.8) $ad(H)\underline{\boldsymbol{m}} = \underline{\boldsymbol{m}}$.

If $\omega$ is the $\underline{\boldsymbol{O}}$ - valued form of the Riemannian connection in *O(M)*, then from (1.7), (1.8) we get that $\overline{\omega} = \omega_{|\underline{\mathbf{h}}}$ defines some connection in *P(H)*, see (1.5). The connections $\omega$, $\overline{\omega}$ can be extended to linear connections with corresponding covariant derivatives $\nabla$, $\overline{\nabla}$.

**DEFINITION 1.9** [17]. The tensor field $h = \nabla - \overline{\nabla}$ is called the second fundamental tensor field of the structure *(P(G), g)*; the connection $\overline{\nabla}$ is called the canonical connection of the pair *(P(G), g)*.

To be correct we must verify an independence of $\overline{\omega}$ from the choice of *P(H)* from the class $\{P_a(a^{-1}Ha)\}_{a \in O(n)}$.

Let $P' = P_a(a^{-1}Ha)$ and $H' = a^{-1}Ha$, then $ad(a^{-1})(\underline{\mathbf{h}}) = \underline{\boldsymbol{h}}'$. From the biinvariance of the Killing metric $\tilde{g}$ it follows that $ad(a^{-1})(\underline{\mathbf{m}}) = \underline{\boldsymbol{m}}'$, where $\underline{\boldsymbol{m}}'$ is the orthogonal complement $\underline{\boldsymbol{m}}' = \underline{\boldsymbol{h}}'^{\perp}$ of $\underline{\boldsymbol{h}}'$ to $\underline{\boldsymbol{O}}$. We have $ad(a^{-1})(X_{\underline{\mathbf{h}}}) = (ad(a^{-1})X)_{|\underline{\mathbf{m}}'}$ and $ad(a^{-1})(X_{\underline{\mathbf{m}}}) = (ad(a^{-1})X)_{|\underline{\mathbf{m}}'}$. Let $\tilde{\omega}'$ be the $\underline{\boldsymbol{h}}'$ - valued form of the connection in *P'(H')*, then

$$\overline{\omega}'((R_a)_* X) = [\omega((R_a)_* X)]_{|\underline{\mathbf{h}}'} = [ad(a^{-1})\omega(X)]_{|\underline{\mathbf{h}}'} = ad(a^{-1})(\omega(X)_{|\underline{\mathbf{h}}})$$

$$= ad(a^{-1})(\overline{\omega}(X))$$

The condition $\overline{\omega}_{ua}((R_a)_* X) = ad(a^{-1})(\overline{\omega}_u(X))$, see [46], is equivalent to the item (b) of *Definition 1.5* hence the linear connections defined by $\overline{\omega}$, $\overline{\omega}'$ coincide.

**3⁰**. *(P(G), g)* is called the particular structure if $h=0$ and the nearly particular one if $h^{+} = 0$, where $h^{+}$ denotes the symmetric part of the tensor field *h*, i.e.,

$$h^{+}_{X}Y = \frac{1}{2}(h_X Y + h_Y X), \quad X, Y \in x(M).$$

It is evident that *(P(G), g)* is the nearly particular structure if and only if the connections $\nabla, \overline{\nabla}$ have the same geodesics. We constructed $\overline{\nabla}$ as a connection in *P(H)*, where *P(H)* is a subbundle of *O(M)*, hence we have got that

(1.9) $\overline{\nabla} g = O,$

that is, the affine connection $\overline{\nabla}$ is a metric connection.

Let $\varphi$ be a cross section of *L(M)* over a neighbourhood *U* which assigns to each $x \in U$ the linear frame $((X_1)_x, \ldots, (X_n)_x)$ and *X, Y* are vector fields on *M*. From (1.6) it follows that

(1.10) $(h_X Y)_x = (\nabla_X Y - \overline{\nabla}_X Y)_x = \varphi(x)\omega|_{\mathbf{m}}(\varphi_* X_x)\varphi(x)^{-1} Y_x,$

where $\varphi(x)$ is considered as the mapping of $R^n$ onto $T_x(M)$.

The following theorem shows that the pair *(g, h)* is defining for the structure *P(G)*.

**THEOREM 1.3.** Structures *P(G)* and $P'(a^{-1}Ga)$ are conjugate if and only if there exists their common associated Riemannian metric *g* and the second fundamental tensor field *h* of *P* coincides with one of *P'*.

**Proof**. If *P* and *P'* are conjugate then from *Definition 1.8* it follows that they have a common associated Riemannian metric *g* and from the construction of *h* it is evident that the tensor field *h* of *P* is equal to that *h'* of *P'*, therefore *h = h'*.

Conversely, if *g* is the common associated Riemannian metric and *h = h'* then $\overline{\nabla} = \overline{\nabla}'$. We consider the holonomy fibre bundles *P(u)* and *P'(u')* of the linear connection $\varpi$ corresponding to $\overline{\nabla}$ such that $P(u) \subset P(G)$ and $P'(u') \subset P'(a^{-1}Ga)$. The structures *P(u)* and *P'(u')* are conjugate, see [46], therefore their extensions *P(G)* and $P'(a^{-1}Ga)$ are also conjugate.

**QED.**

Thus, the pair *(g, h)* defines the class of *G* - structures, where *g* plays a role of the "first fundamental tensor field".

## §3. INTEGRABILITY OF G - STRUCTURES AND POLAR DECOMPOSITION OF O - DEFORMABLE (1,1) TENSOR FIELDS

**$1^0$**. We follow [45], [64].

**DEFINITION 1.1O.** A $G$ - structure $P(G)$ is called an integrable one if for each $p \in M$ there exists a coordinate neighbourhood $U$ with coordinates $(x_1,...,x_n)$ such that for every point $x \in U$ the frame $(x; \frac{\partial}{\partial x_1},...,\frac{\partial}{\partial x_n})$ belongs to $P(G)$.

We give some results about the integrability of tensor $G$ - structures.

**THEOREM 1.4** [55]. Let $P(G)$ be a structure defined by $O$ – deformable (1,1) tensor field $F$ and there exists a connection $\omega$ in $P(G)$ such that corresponding linear connection $\omega$ is without torsion. Then $P(G)$ is an integrable structure.

The following tensor field $N(F)$ on $M$ is called the Nijenhuis tensor of $F$

(1.11) $$N(F)(X,Y) = [FX,FY] - F[FX,Y] - F[X,FY] + F^2[X,Y],$$

$X,Y \in x(M)$

**THEOREM 1.5** [44]. Let $P(G)$ be a structure defined by $O$ – deformable (1,1) tensor field $F$ and for every eigenvalue of $F$ there exists one (generalized) Jordan box or several one-dimensional (generalized) Jordan boxes. Then $P(G)$ is an integrable structure if and only if the Nijenhuis tensor field $N(F)$ vanishes on $M$.

**THEOREM 1.6** [45]. $P(O(n))$ is an integrable structure if and only if the corresponding Riemannian curvature tensor field vanishes on $M$.

Let $O$ - deformable (1,1) tensor fields $F_1,...,F_k$ define the structure $P(G)$, that is, the frames of $P(G)$ are the same for which every $F_i$ has the constant coordinates on $M$ and let $h$ be the second fundamental tensor field of the pair $(P(G),g)$. Then we have the following

**THEOREM 1.7.** If *(P(G),g)* is the particular structure ($h=0$) , then every $F_i$ defines an integrable structure on *M*.

**Proof.** Since $h=0$ hence $\overline{\nabla} = \nabla$ is a connection without torsion. $\overline{\nabla}$ defines a connection $\overline{\omega}$ in *P(G)* and the rest follows from *Theorem 1.4*.

**QED.**

**2⁰**. It is well known, see [31], a polar decomposition of matrix *(F)* in Euclidian space: *(F) = (S)(P)* , where *(S)* is an unitary matrix and *(P)* is a symmetric positively semi-indefinite one. This decomposition is unique if *(F)* is a nonsingular marix.

We consider the polar decomposition of an *O*-deformable (1,1) tensor field *F*.

**DEFINITION 1.11.** Further, for simplicity, an *O*-deformable (1,1) tensor field *F* will be called an affinor.

Let *F* be an affinor, which defines such a structure *P(G)* that $P(G)\cap O(M)=P(H)$, $H = G\cap O(n)$, where *O(M)* is a subbundle of all the orthonormal frames over *M* determined by an associated to *P(G)* Riemannian metric *g* and *G* is the group of invariance of some fixed matrix *(F)* of *F*.

**PROPOSITION 1.8.** There exists an affinor $F^*$ on *M* such that $<F^*X,Y>=<X,FY>$ for any $X,Y\in x(M)$, the metric *g* is an associated one to $P(F^*)$ and the second fundamental tensor fields of *P(F)* and $P(F^*)$ coincide.

**Proof.** If $a\in H$, then $a^T(F)a=(F)$ and $a^T(F)^T a=(F)^T$, therefore the matrix $(F)^T$ is invariant under *H* and *P(H)* can be extended to the group $G^*$ of invariance of $(F)^T$. For any $b\in G^*\cap O(n)$ we have $b^T(F)^T b=(F)^T$, hence $b^T(F)b=(F)$ and $G^*\cap O(n)=H$. Thus *g* is an associated one to $P(G^*)$ and from construction of *h*, see *Definition 1.9*, it is evident that the second fundamental tensor fields of *P(G)* and $P(G^*)$ coincide.

**QED**.

We consider now the matrix $(P^2)=(F)^T(F)$. Let the matrix *(P)* be the nonnegative square root of $P^2$, [29], $G_2$ be the group of invariance of $P^2$ and $G_1$ be that of *(P)*.

**LEMMA 1.9.** We have $\boldsymbol{G_1 = G_2}$ and $\boldsymbol{H \subset G_2}$.

**Proof.** If $a^{-1}(P)a=(\overline{P})$, then $a^{-1}(P^2)a=(\overline{P}^2)$, where

$$(\overline{P}^2)=\begin{bmatrix} \boxed{0} & & & & & 0 \\ & \ddots & & & & \\ & & \boxed{\begin{matrix}\lambda_1^2 & & 0\\ & \ddots & \\ 0 & & \lambda_1^2\end{matrix}} & & & \\ & & & \ddots & & \\ 0 & & & & \boxed{\begin{matrix}\lambda_k^2 & & 0\\ & \ddots & \\ 0 & & \lambda_k^2\end{matrix}} \end{bmatrix},\quad (\overline{P})=\begin{bmatrix} \boxed{0} & & & & & 0 \\ & \ddots & & & & \\ & & \boxed{\begin{matrix}\lambda_1 & & 0\\ & \ddots & \\ 0 & & \lambda_1\end{matrix}} & & & \\ & & & \ddots & & \\ 0 & & & & \boxed{\begin{matrix}\lambda_k & & 0\\ & \ddots & \\ 0 & & \lambda_k\end{matrix}} \end{bmatrix},\ \lambda_j\geq 0$$

The groups $\overline{G_1}=a^{-1}G_1a,\ \overline{G_2}=a^{-1}G_2a$ are the groups of invariance of $(\overline{P}),\ (\overline{P}^2)$ respectively and it is evident that $\overline{G_1}=\overline{G_2}$, hence $G_1=G_2$. For $a\in H$ we have

$$a^T(P^2)a=(a^T(F)^T a)(a^T(F)a)=(F)^T(F)=(P^2)$$

and $a\in G_2$.

**QED.**

**THEOREM 1.1O.** Let *F* be an affinor, let *g* be an associated to *P(F)* Riemannian metric and let $\overline{\nabla}$ be the canonical connection of the pair *(P(F),g)*. Then we have

a) *F = SP*, where *S* and *P* are affinors, *S* is an unitary affinor and *P* is a symmetric positively semi-indefinite one with respect to *g*.

b) $\overline{\nabla}F=\overline{\nabla}S=\overline{\nabla}P=0$.

**Proof**. It follows from *Lemma 1.9* that *P(H)* can be extended to the group $G_1$ and $P(G_1)$ determines the affinor *P*. From the identity $(F)=(S)(P)$ for matrixes we see that *H* preserves the matrix *(S)*, hence an affinor *S* is defined by extension of *P(H)*. Since $\overline{\omega}$ is a connection in *P(H)*, where *P(H)⊂P(G)*, *P(H)⊂P(*$G_1$*)*, *P(H)⊂P(S)*, it is clear that $\overline{\nabla}F=\overline{\nabla}S=\overline{\nabla}P=0$. The rest is obvious.

**QED.**

REMARK. In general, the second fundamental tensor fields of *(P(G),g)*, *(P(*$G_1$*),g)*, *(P(S),g)* are not the same. For example, if *<FX,FY>=<X,Y>*, where *X,Y∈x(M)*, then *F=S* and *P=I*, therefore the structure *(I,g)* is a particular one, i.e. *h=0* , but the structure *(P(F),g)* can be various.

—–

# CHAPTER 2

# RIEMANNIAN GEOMETRY OF G – STRUCRURES

In this chapter, we study a pair *(P(G),g)* where *P(G)* is a *G*-structure over *M* and *g* is an associated with it Riemannian metric. There are many results about different concrete structures on manifolds considered with Riemannian metrics. Using the introduced notions of the second fundamentical tensor field h and of the canonical connection $\overline{\nabla}$ of the structure we want to unite such results in on scheme to get a general theory. In §1 we consider torsion and curvature of $\overline{\nabla}$, so-called quasi homogeneous structures and dependent notions. §2 is devoted to the study of isometries, affine transformations and holonomy fibre bundles of $\overline{\nabla}$. Homogeneous structures are discussed in §3 and projections on submanifolds and foliations in the last §4.

We follow especially closely [47], [67].

## *§1. TORSION AND CURVATURE, QUASI HOMOGENEOUS STRUCTURES*

**1⁰**. We consider a fixed pair *(P(G),g)* or *(P(H),g)*, where $H=G\cap O(n)$, and its second fundamental tensor field $h=\nabla-\overline{\nabla}$, where $\nabla$ is the Riemannian connection, $\overline{\nabla}$ is the canonical that of the structure. Tensor fields $h^+$, $h^-$ denote the symmetric and skew-symmetric parts of *h* respectively:

$$(2.1)\quad h^+_X Y=\frac{1}{2}(h_X Y+h_Y X),\ \ h^-_X Y=\frac{1}{2}(h_X Y-h_Y X),\ \ X,Y\in \mathfrak{X}(M).$$

**LEMMA 2.1**. Let $h_{XYZ}=<h_X Y,Z>$ for any $X,Y,Z\in\mathfrak{X}(M)$. Then we have

$$(2.2)\quad h_{XYZ}=-h_{XZY}.$$

**Proof**. Since $\nabla$ is a metric connection, then

$$X<Y,Z>=<\nabla_X Y,Z>+<Y,\nabla_X Z>$$

and it is analogously for $\overline{\nabla}$ from (1.9). Subtracting one equality from the other we have got the lemma.

**QED.**

Let $\overline{T}$ be the torsion tensor field for $\overline{\nabla}$, such a tensor field is equal to zero for the connection $\nabla$.

**LEMMA 2.2**. We have

(2.3) $\overline{T} = -2h^{-}$.

**Proof**. Subtracting from the equality

$$\overline{T}_X Y = \overline{\nabla}_X Y - \overline{\nabla}_Y X - [X,Y]$$

the equality

$$O = \nabla_X Y - \nabla_Y X - [X,Y]$$

we have obtained our lemma.

**QED.**

**PROPOSITION 2.3**. The pair *(P(G),g)* is a particular structure if and only

$h^{-} = 0$.

**Proof** obviously follows from (2.3) and from uniqueness of metric connection without torsion.

**QED.**

The Riemannian curvature tensor field *R* of the connection $\nabla$ is defined by

(2.4) $R_{XY} = [\nabla_X, \nabla_Y] - \nabla_{[X,Y]}, \quad X, Y \in \mathcal{X}(M)$

The curvature tensor field $\overline{R}$ of the canonical connection $\overline{\nabla}$ is defined analogously by

(2.5) $\overline{R}_{XY} = [\overline{\nabla}_X, \overline{\nabla}_Y] - \overline{\nabla}_{[X,Y]}, \quad X, Y \in \mathcal{X}(M)$

**THEOREM 2.4**. For any $X,Y,Z \in \mathcal{X}(M)$ we have

(2.6) $R_{XY}Z = \overline{R}_{XY}Z + (\overline{\nabla}_X h)(Y,Z) - (\overline{\nabla}_Y h)(X,Z) + [h_X, h_Y]Z + h_{\overline{T}_X Y}Z$

**Proof**. Since $\nabla = \overline{\nabla} + h$ we have

$$R_{XY} = [\nabla_X, \nabla_Y] - \nabla_{[X,Y]} = [\overline{\nabla}_X + h_X, \overline{\nabla}_Y + h_Y] - \overline{\nabla}_{[X,Y]} - h_{[X,Y]}$$
$$= \overline{R}_{XY} + [\overline{\nabla}_X, h_Y] + [h_X, \overline{\nabla}_Y] + [h_X, h_Y] - h_{[X,Y]}.$$

As

$$[\overline{\nabla}_X, h_Y] = \overline{\nabla}_X h_Y - h_Y \overline{\nabla}_X = (\overline{\nabla}_X h)(Y, .\ ) + h_{\overline{\nabla}_X Y}$$

and

$$[h_X, \overline{\nabla}_Y] = h_X \overline{\nabla}_Y - \overline{\nabla}_Y h_X = -((\overline{\nabla}_Y h)(X, .\ ) + h_{\overline{\nabla}_Y X}),$$

therefore

$$R_{XY} = \overline{R}_{XY} + (\overline{\nabla}_X h)(Y, .) - (\overline{\nabla}_Y h)(X, .) + [h_X, h_Y] + h_{\overline{\nabla}_X Y} - h_{\overline{\nabla}_Y X} - h_{[X,Y]}$$
$$= \overline{R}_{XY} + (\overline{\nabla}_X h)(Y, .) - (\overline{\nabla}_Y h)(X, .) + [h_X, h_Y] + h_{\overline{T}_X Y}$$

**QED.**

$2^0$. **DEFINITION 2.1**. We call a pair *(P(G),g)* a quasi homogeneous structure if $\overline{\nabla} h = 0$ on *M*.

From *Theorem 2.4* it follows that for a quasi homogeneous structure

(2.7) $R_{XY}Z = \overline{R}_{XY}Z + [h_X, h_Y]Z + h_{\overline{T}_X Y}Z$.

Examples of quasi homogeneous structures one can find in [41].
The following vector fields play an important role in our consideration.

**DEFINITION 2.2**. We call a vector field *X* a particular one, respectively a nearly particular one, if $h_X = 0$ on *M*, respectively $h_X^+ = 0$.

If $X_1, ..., X_n$ are orthonormal vector fields on a neighbourhood of some point of $M$, then (1,1) tensor fields $r^1, r^2$ are defined on this neighbourhood by

(2.8) $< r^1 X, Y > = \sum_{k=1}^{n} < h_X X_k, h_Y X_k >,$

(2.9) $\quad < r^2 X, Y >= \sum_{k=1}^{n} < h_X^+ X_k, h_Y^+ X_k > .$

Let $Y_1,...,Y_n$ be another vector fields, then for every point from our neighbourhood $Y_j = \sum_k a_j^k X_j$ , where $(a_j^i) \subset O(n)$, and

$$< r^1 X, Y >= \sum_{j=1}^{n} < h_X Y_j, h_Y Y_j > = \sum_{j,k=1}^{n} < a_j^k h_X X_k, a_j^k h_X X_k > = \sum_{k=1}^{n} < h_X X_k, h_Y X_k >.$$

So, we have correctly defined on $M$ tensor field $r^1$, for $r^2$ it is analogously. It is evident that for any $X,Y \in \mathbf{X}(M)$

(2.1O) $\qquad < r^i X, Y >=< X, r^i Y >, \quad i = 1,2 .$

**PROPOSITION 2.5.** A set of particular (nearly particular) vector fields on $M$ concides with Ker $r^1$ ( Ker $r^2$).

**Proof.** If $X \in Ker\, r^1$, then $< r^1 X, X >= \sum_{k=1}^{n} \| h_X X_k \|^2 = O$, hence $h_X X_k = 0$ for each $X_k$, $k = 1,...,n$, and from linear independence of $X_1,...,X_n$ it follows that $h_X Y = 0$ for any $Y \in \mathbf{X}(M)$. Conversely, if $h_X Y = O$, then $r^1 X = O$ from (2.8) and $X \in Ker\ r^1$. The proof for $r^2$ is the same.

**QED.**

**3⁰**. The following theorem describes a situation in a quasi homogeneous case.

**THEOREM 2.6**. Let $(P(G),g)$ be a quasi homogeneous structure i.e. $\overline{\nabla} h = O$. Then we have
I) $\overline{\nabla} r^1 = O$, $\overline{\nabla} r^2 = O$; II) there exist almost product structures
$T(M) = Ker\ r^1 \oplus Im\, r^1; T(M) = Ker\ r^2 \oplus Im\, r^2$
(direct sums) on $M$ and $Im\, r^i = (Ker\ r^i)^{\perp}, i = 1,2$.

**Proof**. I) Let $Z$ be a vector field and $\gamma$ a curve segment in $M$ defined by $Z$ or more precisely by a local 1-parameter group of transformations induced by $Z$. We denote by $X,Y,X_1,...,X_n$ the vector fields defined on some neighbourhood of $\gamma$ which are obtained by the parallel translation of $X_p, Y_p, X_{1p},...,X_{np} \in T_p(M)$ along $\gamma$, $p \in \gamma$, in the connection $\overline{\nabla}$. Obviously $X_1,...,X_n$ are remained orthonormal because the connection $\overline{\nabla}$ is metric. So we have

$(\overline{\nabla}_Z X)_p = (\overline{\nabla}_Z Y)_p = (\overline{\nabla}_Z X_k)_p = O$ for any $k = 1,...,n$ and

$((\overline{\nabla}_Z r^i)X)_p = (\overline{\nabla}_Z r^i X)_p$.

$$< \overline{\nabla}_Z r^1 X, Y >_p = Z < r^1 X, Y >_p = \sum_{k=1}^{n} (< \overline{\nabla}_Z h_X X_k, h_Y X_k >$$

$$+ < h_X X_k, \overline{\nabla}_Z h_Y X_k >)_p = O,$$

$(\overline{\nabla}_Z h^+)_X Y = \frac{1}{2} \overline{\nabla}_Z (h_X Y + h_Y X) = O$ and $\overline{\nabla}_Z r^2 = O$ too.

II) From I) it follows that $r^1, r^2$ are affinors, hence $Ker\ r^i$ define differentiable distributions on $M$.

From (2.1O) we have that $(Ker\ r^i)^{\perp} = Im\, r^i$, therefore

$T(M) = Ker\ r^i \oplus Im\, r^i, i = 1,2$, are almost product structures on $M$.

**QED.**

**DEFINITION 2.3**. A pair *(P(G),g)* is called a strict structure if $h_X \neq O$ for each $X \in X(M)$, $X \neq O$.

It is evident that *(P(G),g)* is strict if and only if $Ker\ r^1 = \{O\}$ for each point of $M$.

**4$^0$**. The following notions can be useful for the study of the structure *(P(G),g)* over the manifold $M$. We call

a) $r^1$ the induced Ricci mapping,

b) $< r^1 X, Y >$ the induced Ricci tensor,

c) tr $r^1$ the induced scalar curvature,

d) $< r^1 X, X > \|X\|^{-2}$ the inducud curvature in direction $X$ of the structure *(P(G),g)* over the manifold $M$.

The similar notions can be defined for the tensor field $r^2$.

## *§2. ISOMETRY GROUPS AND HOLONOMY FIBRE BUNDLES*

**1[0]**. We remind (see *Definition 1.5*) that a connection $\Gamma$ in *O(M)* is an assignment of a subspace $Q_u$ of $T_uO(M)$ to each $u \in P$ such that (a) $T_u(O(M)) = Q_u \oplus V_u$, where $V_u$ is the tangent space of the fibre, (b) $(R_a)_* Q_u = Q_{ua}$ and (c) $Q_u$ depends differentiably on *u*.

Then for every $\xi \in R^n$ and $u \in L(M)$ there exists unique vector $(B(\xi))_u$ in $Q_u$ such that $\pi((B(\xi))_u) = u(\xi)$. $B(\xi)$ is called the standard horizontal vector field corresponding to $\xi \in R^n$.

**PROPOSITION 2.7** [46]. Let $B(\xi)$ be a standard horizontal vector field and $\Theta$ the canonical form of *L(M)*. Then we have:
(1) $\Theta(B(\xi)) = \xi,\ \xi \in R^n$;
(2) $R_a(B(\xi)) = B(a^{-1}\xi),\ a \in GL(n,R),\ \xi \in R^n$;
(3) $B(\xi) \neq O$ on *M* for $\xi \neq O$.

We can define the Riemannian connection $\Gamma$ in *O(M)* by $u \to Q_u$ or by ***O***-valued form or by covarint derivative $\nabla$.

Let *A** be the fundamental vector field corresponding to $A \in$ ***O***. The Riemannian metric $g_\nabla$ is defined on *O(M)* by

$$g_\nabla(B(\xi),B(\eta)) = <\xi,\eta>,\quad \xi,\eta \in R^n;$$
$$g_\nabla(A_1^*,A_2^*) = -tr(A_1A_2),\ A_1,A_2 \in \mathbf{O};$$
$$g_\nabla(B(\xi),A^*) = O,\quad \xi \in R^n, A \in \mathbf{O}.$$

The action of the isometry group *I(M)* on *M* induces the action of *I(M)* as the isometry group on $(O(M),g_\nabla)$ (see [67]) by the formula:

(2.11) $\varphi u = \varphi(p;X_1,\ldots,X_n) = (\varphi p;\varphi_* X_1,\ldots,\varphi_* X_n)$,

where $u = (p;X_1,\ldots,X_n) \in O(M)$ and $\varphi \in I(M)$.

For our structure *(P(H),g), H=G∩O(n),* we can consider the canonical connection $\overline{\Gamma}$ defined by $\overline{\nabla}$, $\overline{\omega}$, $u \to \overline{Q}_u$ and the standard horizontal vector fields $\overline{B}(\xi)$ too. We introduce distributions $V^{\mathbf{h}}$, $V^{\mathbf{m}}$ by $V^{\underline{\mathbf{h}}} = \{ A^* : A \in \underline{h} \},\ V^{\underline{\mathbf{m}}} = \{ A^* : A \in \underline{m} \}$

It is obvious that $V^{\underline{\mathbf{m}}} \perp V^{\underline{\mathbf{h}}}$ in $V$ and $T_u P(H) \cap V_u = V_u^{\underline{h}}$ for each $u \in P(H)$.

**LEMMA 2.8**. The horizontal distribution $\overline{Q}_u$ of the canonical connection of the structure *(P(H),g)* is defined by the following formula

(2.11) $$u \to \overline{Q}_u = T_u(P(H)) \cap V_u^{\underline{h}\perp} = T_u(P(H)) \cap (V_u^{\underline{m}} \oplus Q_u),$$

where $u \in P(H)$.

**Proof**. As $\omega(\overline{Q}_u) \subset \underline{\mathbf{m}}$, therefore $\overline{\omega}(\overline{Q}_u) = O$. Since $T_u(P(H)) = V_u^{\underline{h}} \oplus \overline{Q}_u$, hence $dim\ \overline{Q}_u = n$ and $\overline{Q}_u = Ker\ \overline{\omega}$.

**QED.**

**2⁰**. Let $\overline{\boldsymbol{\omega}}$ be the canonical connection in *P(H)*, and $P(u_0)$ the holonomy bundle through $u_0$, where $u_0 \in P(H)$, i.e., the set of all $u \in P(H)$ which can be joined to $u_0$ by a ( piece-wise differentiable) horizontal curve. Further, let $\Phi(u_0)$ denote the holonomy group with reference frame $u_0$. Then the famous Reduction theorem, [46], sets that

(I) $P(u_0)$ is a differentiable subbundle of *P(H)* with the structure group $\Phi(u_0)$.

(II) The connection $\overline{\boldsymbol{\omega}}$ is reducible to a connection in $P(u_0)$.

**THEOREM 2.9**. The linear canonical connection of the structure $(P(u_0),g)$ coincides with the canonical connection $\overline{\omega}$ of the pair *(P(H),g)*.

**Proof**. From *Lemma 2.8* we see that the horizontal distribution of the canonical connection in $P(u_0)$ is the orthogonal complement of $V_u^{\underline{\bar{h}}}$ to $T_u(P(u_0))$, where $V_u^{\underline{\bar{h}}}$ is the vertical subspace corresponding to the holonomy algebra $\underline{\boldsymbol{h}}$. Since $\overline{Q}_u \subset T_u P(u_0)$ and $\overline{Q}_u \perp V^{\underline{h}} \supset V^{\underline{\bar{h}}}$, hence it is clear from *Lemma 2.8* that above - mentioned distribution coincides with $\overline{Q}_u$ for every $u \in P(u_0)$ and, therefore, for any $u \in L(M)$ from the right invariance.

**QED.**

The theorem shows that sometimes we can identify *P(H)* with the holonomy fibre bundle $P(u_0)$ in our consideration.

Now we give a geometric characterization of *Definition 2.3*.

**THEOREM 2.1O**. The structure *(P(H),g)* is strict if and only if $Q_u \cap \overline{Q}_u = \emptyset$ for each $u \in O(M)$.

**Proof**. If $\tau = \omega - \varpi$ and $u \in O(M)$, then, according to [7], we have

$$h_X Y = u\tau(\overline{X})u^{-1}Y$$

for $X, Y \in T_{\pi(u)}(M)$, $\overline{X} \in T_u(O(M))$, $\pi_*(\overline{X}) = X$.

If $\overline{X} \in \overline{Q}_u$, then the condition $h_X \neq O$ is equivalent to one that $\tau(\overline{X}) = \omega(\overline{X}) \neq O$ and $\overline{X} \notin Q$. Having taken all such $\overline{X} \notin Q_u$ we have got that $Q_u \cap \overline{Q}_u = \emptyset$. The converse it is evident.

**QED.**

From the proof of *Theorem 2.1O* it is obvious that for each $u \in O(n)$ $\pi_*(Q_u \cap \overline{Q}_u)$ conincides with $(Ker\ r^1)_{\pi(u)}$.

**3⁰**. Let $\overline{\Gamma}$ be the canonical connection of *(P(H),g)* and $\varphi$ a transformation of *M*. It is well known, [46], that the following conditions are equivalent :

1) $\varphi$ is an affine transformation of *M* with respect to $\overline{\nabla}$,
2) $\varphi^* \varpi = \varpi$
3) every standard horizontal vector field $\overline{B}(\xi)$ is invariant under $\varphi$,
4) $\varphi_* \overline{\nabla}_X Y = \overline{\nabla}_{\varphi_* X} \varphi_* Y$ for all $X, Y \in \mathfrak{X}(M)$.

If $\varphi$ is an isometry, then $\varphi$ is an affine transformation with respect to $\Gamma$ but we do not say the same about $\overline{\Gamma}$.

**DEFINITION 2.4** [48]. The group of all the affine transformations of *M* with respect to $\overline{\nabla}$ preserving each holonomy bundle *P(u)* , $u \in L(M)$, is called the transvection group and it is denoted by $Tr(M, \overline{\nabla})$ or simply by $Tr(\overline{\nabla})$.

It is evident that if preserves a fixed holonomy bundle $P(u_0)$, then it also preserves the holonomy bundle *P(u)* for each $u \in L(M)$. More geometrically, an affine transformation $\varphi$ of $(M,\overline{\nabla})$ belongs to the group $Tr(\overline{\nabla})$ if and only if the following holds : for every point $p\in M$ there is a piece-wise differentiable curve $\gamma$ joining $p$ to $\varphi(p)$ such that the tangent map $\varphi_{*p} : T_p(M) \to T_{\varphi(p)}(M)$ coincides with the parallel translation along $\gamma$.

The following theorem describes the situation for canonical connection $\overline{\nabla}$ of the structure $(P(H),g)$.

**THEOREM 2.11**. If $P(H)$ is invariant under an isometry $\varphi$ i.e. $\varphi(P(H))\subset P(H)$, then $\varphi$ is an affine transformation with respect to $\overline{\nabla}$. Conversely if $(P(H),g)$ is strict and for some frame $u\in P(u_0)\subset P(H)$ $\varphi(u)\in P(u_0)$, where an isometry $\varphi$ is an affine transformation with respect to $\overline{\nabla}$, then $\varphi(P(H))\subset P(H)$ and $\varphi\in Tr(\overline{\nabla})\subset I(M)$.

**Proof**. Since $\varphi$ is an isometry, then $\varphi(O(M))\subset O(M)$. From (2.11) it follows that $\varphi_*(V)=V$. We know, [46], that $\varphi$ preserves fibres, canonical form $\Theta$ is invariant under $\varphi$ and $\varphi$ is an isometry for $g_{\nabla}$ on $O(M)$, therefore $\varphi_*(Q)=\varphi_*(V^{\perp})=V^{\perp}=Q$.

$\varphi_*(T_u(P(H))=T_{\varphi(u)}(P(H))$ and $\varphi_*(V^{\underline{\mathbf{h}}})= V^{\underline{\mathbf{h}}}$ , $V^{\underline{\mathbf{h}}} =T_u(P(H))\cap V$, hence $\varphi_*(V^{\underline{\mathbf{h}}\perp}) = V^{\underline{\mathbf{h}}\perp}$ From *Lemma 2.8* we have got $\varphi_*(\overline{Q})=\overline{Q}$. The form $\omega$ of the Riemannian connection is invariant under $\varphi$ and it is obvious that the form $\overline{\omega}$ is also invariant under $\varphi$. So, $\varphi$ is an affine transformation with respect to $\overline{\nabla}$.

Conversely, if an isometry $\varphi$ is an affine transformation, then $\varphi_*(\overline{Q})=\overline{Q}$ and $\varphi_*(Q)=Q$. Using *Theorem 2.9*, let $P(H)=P(u_0)$ at first. $\varphi$ preserves

$$Q^{\perp} = V^{\underline{\mathbf{h}}} \oplus V^{\underline{\mathbf{m}}}, \quad \overline{Q}^{\perp} = V^{\underline{\mathbf{h}}} \oplus \overline{V}.$$

From *Theorem 2.1O* we have that $Q^{\perp}\cap\overline{Q}^{\perp}=V^{\underline{\mathbf{h}}}$ and $V^{\underline{\mathbf{h}}}$ is invariant under $\varphi$. Let $\tilde{\Lambda}$ be the following distribution

$$\tilde{\Lambda} : \Lambda_u = T_u(P(\tilde{H})), \quad u\in O(M),$$

where $P(\tilde{H})$ is one from the structures conjugate to $P(H)$, $u\in P(\tilde{H})$. $\tilde{\Lambda}$ is invariant under $\varphi$, therefore $\varphi$ transforms the fibres of the foliation $\tilde{\Lambda}$ onto the fibres. Hence, if $\varphi(u) \in P(u_0)$ for some $u \in P(u_0)$, then $\varphi(P(u_0)) = P(u_0)$ and $\varphi \in Tr(\overline{\nabla})$. If $\varphi$ preserves $P(u_0)$, then $\varphi$ preserves an extension of $P(u_0)$, that is, $P(H)$.

**QED.**

**4⁰**. We follow [46], [73]. Let $\overline{\nabla}$ be an affine connection on *M*. A vector field *X* on *M* is called an infinitesimal affine transformation if, for each $x \in M$, a local 1-parameter group of local transformations $\varphi_t$ of a neighbourhood *U* of *x* into *M* preserves the connection $\overline{\nabla}$.

A vector field *X* on *M* is called an infinitesimal isometry (or, a Killing vector field) if the local 1-parameter group of local transformations generated by *X* in a neighbourhood of each point of *M* consists of local isometries. An infinitesimal isometry is necessarily an infinitesimal affine transformation with respect to $\nabla$. *X* is an infinitesimal isometry if and ohly if $\mathcal{L}_X g = O$ ($\mathcal{L}_X$ is the Lie differentiation with respect to *X*).

Let $\overline{\nabla}$ be a complete affine connection on *M*. Then every infinitesimal affine transformation *X* of *M* is complete, that is, *X* generates a global 1-parameter group of transformations of *M*.

A vector field *X* is an infinitesimal affine transformation of *M* if and only if

(2.12) $\mathcal{L}_X \cdot \overline{\nabla}_Y - \overline{\nabla}_Y \cdot \mathcal{L}_X = \overline{\nabla}_{[X,Y]}$ for each $Y \in \mathfrak{X}(M)$.

Let $\overline{\nabla}$ be as usually the canonical connection of a structure *(P(H),g)*.

**THEOREM 2.12** [67]. A metric connection (in particular $\overline{\nabla}$) is complete on the complete Riemannian manifold *M*.

So we have got that any infinitesimal affine transformation of *M* with respect to $\overline{\nabla}$ is complete if *M* is a complete Riemannian manifold.

**THEOREM 2.13**. A vector field *X* is an infinitesimal isometry and an affine transformation with respect to $\overline{\nabla}$ if and only if $\mathcal{L}_X g = O$ and $\mathcal{L}_X h = O$, where $h = \nabla - \overline{\nabla}$.

Proof. For any $Y,Z \in \mathfrak{X}(M)$ we obtain

$$(\mathcal{L}_X h)(Y,Z) = [X, h_Y Z] - h_{[X,Y]}Z - h_Y[X,Z]$$

$$= ([X, \nabla_Y Z] - \nabla_{[X,Y]}Z - \nabla_Y[X,Y])$$

$$- ([X, \overline{\nabla}_Y Z] - \overline{\nabla}_{[X,Y]}Z - \overline{\nabla}_Y[X,Y])$$

$$= (\mathcal{L}_X \cdot \nabla_Y - \nabla_Y \mathcal{L}_X - \nabla_{[X,Y]})Z$$

$$- (\mathcal{L}_X \cdot \overline{\nabla}_Y - \overline{\nabla}_Y \mathcal{L}_X - \overline{\nabla}_{[X,Y]})Z.$$

From the above-mentioned results the theorem follows.

**QED.**

## *§3. HOMOGENEOUS STRUCTURES*

**1⁰**. The Riemannian manifold *(M,g)* is called a homogeneous one if the full isometry group *I(M)* is a transitive Lie group of transformations of *M*. The following theorem of W.Ambrose, I.M.Singer [3], [67], [48] plays an important role in the study of such manifolds.

**THEOREM (AS)**. A connected, complete, simply connected Riemannian manifold *(M,g)* is a homogeneous one if and only if there exists a tensor field $h$ of type (1,2) on $M$ such that for any $X,Y,Z \in \mathfrak{X}(M)$

(AS1) $\langle h_X Y, Z \rangle + \langle Y, h_X Z \rangle = O$,

(AS2) $(\nabla_X R)_{YZ} = [h_X, R_{YZ}] - R_{h_X YZ} - R_{Y h_X Z}$,

(AS3) $(\nabla_X h)_Y = [h_X, h_Y] - h_{h_X Y}$,

where $\nabla$ and $R$ denote the Riemannian connection and the Riemannian curvature tensor field respectively.

Every such a tensor field $h$ satisfying the conditions *(AS)* on $M$ is called a homogeneous Riemannian structure on $M$. We can define a canonical connection on $M$ by the formula

$\tilde{\nabla} = \nabla - h$

Then the conditions *(AS1), (AS2), (AS3)* are equivalent to the following

(2.13) $\tilde{\nabla} g = O,\ \tilde{\nabla} h = O,\ \tilde{\nabla}\tilde{R} = O$,

where $\tilde{R}$ is the Riemannian curvature tensor field of $\tilde{\nabla}$.

Let $(M, \tilde{\nabla})$ be a connected manifold with an affine connection $\tilde{\nabla}$. Then the following two conditions are equivalent :

(I) The transvection group $Tr(\tilde{\nabla})$ acts transitively on each holonomy bundle $P(u) \subset L(M)$.

(II) $M$ can be expressed as a reductive homogeneous space $K/K_0$ with respect to a decomposition $\underline{\mathbf{k}} = \underline{\mathbf{k}}_0 + \underline{\mathbf{m}}$, where $K$ is effective on $M$ and $\tilde{\nabla}$ is the canonical connection of $K/K_0$, $\underline{\mathbf{k}}$, $\underline{\mathbf{k}}_0$, are the Lie algebras of $K$, $K_0$ respectively.

If (I) is satisfied, then $Tr(\tilde{\nabla})$ is a connected Lie group and $M$ can be expressed in the form (II) with $K = Tr(\tilde{\nabla})$. For every expression of $M$ in the form (II), $Tr(\tilde{\nabla})$ is a normal Lie subgroup of $K$ and its Lie algebra is isomorphic to the ideal $\underline{\mathbf{m}} = [\underline{\mathbf{m}}, \underline{\mathbf{m}}]$ of $\underline{\mathbf{k}}$.

So, the space $(M, \tilde{\nabla})$ will be called an affine reductive space if the transvection group $Tr(M, \tilde{\nabla})$ acts transitively on each holonomy bundle.

Using results in [67], we can give the following characterisation of the homogeneous Riemannian structures :

Let $(M,g)$ be a connected, complete and simply connected Riemannian manifold, and $\nabla$ its Riemannian connection. A tensor field $h$ of type (1,2) on $(M,g)$ is a homogeneous structure on $(M,g)$ if and only if the new affine connection $\tilde{\nabla} = \nabla - h$ determines an affine reductive space $(M, \tilde{\nabla})$ and g is parallel with respect to $\tilde{\nabla}$. Each homogeneous Riemannian space $(M, g)$ admits at least one homogeneous structure ( but it can admit more than one).

Now, let $h$ be as usually the second fundamental tensor field of the pair $(P(H),g)$ on $M$.

**THEOREM 2.14**. Let $M$ be a simply connected, complete Riemannian manifold and let $\overline{\nabla}$ be the canonical connection of the structure $(P(H),g)$ and $\overline{\nabla} h = O$, $\overline{\nabla} R = O$. Then $h$ is a homogeneous Riemannian structure and $M \cong K/K_0$ is the Riemannian homogeneous space, where $K = Tr(\overline{\nabla})$ and $K_0$ is the isotropy subgroup of a fixed point $o \in M$. The structure $P(H)$ is invariant under $K$.

**Proof**. From *Theorem 2.12* we see that $\overline{\nabla}$ is complete on the complete Riemannian manifold $M$ and we have $\overline{\nabla} g=O$, $\overline{\nabla} h=O$, $\overline{\nabla} R=O$. From (2.3) it follows that $\overline{\nabla}\overline{T}=O$ and from (2.7) $\overline{\nabla} R=O$ implies $\overline{\nabla}\overline{R}=O$. So, there exists the minimal group of affine transformations $K=Tr(\overline{\nabla})$, which acts transitively on the holonomy fibre bundle $P(u_0)\subset P(H)$ (see [47]) . Since $\overline{\nabla}$ preserves $P(H)$ we see that $P(u)\subset P(H)$ for each $u\in P(H)$. As $P(u_0)$ is invariant under the action of $K$, therefore $P(H)$ is invariant too. Since $P(H)\subset O(M)$ it is evident that $K\subset I(M)$.

**QED.**

We describe now relations between canonical connections $\overline{\nabla}$ and $\tilde{\nabla}$.

**THEOREM 2.15**. Let $\tilde{\nabla}$ be the canonical connection of a Riemannian homogeneous space $M\cong K/K_0$, where $K=Tr(\tilde{\nabla})$, and let $P(u_0)$ be the holonomy fibre bundle of $\tilde{\nabla}$ containing $u_0\in O(M)$. Then for canonical connection $\overline{\nabla}$ of $G$-structure $(P(u_0),g)$ on $M$ we have $\overline{\nabla}=\tilde{\nabla}$.

**Proof**. Since $M\cong K/K_0$, where $K=Tr(\tilde{\nabla})$, then $K$ acts transitively on $P(u_0)$ and preserves $P(u_0)$, [48] . We consider a mapping

$$J_u : K\to O(M): \varphi\mapsto \varphi(u),\ u\in P(u_0).$$

Then $J_u(K)=P(u_0)$ and $K_u=(J_u)_*\underline{\boldsymbol{k}}$, where $\underline{\boldsymbol{k}}$ is the Lie algebra of $K$. Let $\underline{\boldsymbol{h}}$ be the Lie algebra of the holonomy group with the fixed frame $u_0$.

Using [67], we have from *Lemma 2.8* that the horizontal distribution of $\tilde{\nabla}$

$$\tilde{Q}_u=K_u\cap(K_u\cap V_u)^{\perp}=T_u(P(u_0))\cap(T_u(P(u_0))\cap V_u)^{\perp}$$

$$=T_u(P(u_0))\cap V_u^{\overline{h}\perp}=\overline{Q}_u$$

From the right invariance it is evident that $\tilde{Q}_u=\overline{Q}_u$ for each $u\in L(M)$, hence the connections $\tilde{\nabla}$ and $\overline{\nabla}$ coincide.

**QED.**

**REMARKS.**

1) Let $\tilde{\nabla}$ be the canonical connection of a Riemannian homogeneous space $M\cong K/K_0$, where $K=Tr(\tilde{\nabla})$, $P(H)$ be an invariant structure under the action of $K$, $\overline{\nabla}$ be the canonical connection of $(P(H),g)$. In general, $\tilde{\nabla}$ and $\overline{\nabla}$ are not the same. For example if $P(H)=P(O(n))=O(M)$, then $\overline{\nabla} = \nabla$ but $\tilde{\nabla} \neq \nabla$ if *M* is not a symmetric Riemannian space.

2) From *Theorems 2.14, 2.15* we see that the notion of the second fundamental tensor field *h* of *G*-structure generalises, in some sense, that of the homogeneous Riemannian structure on *M*.
3) The holonomy fibre bundle $P(u_0)$ gives a trivial example of a quasi homogeneous structure. More interetings examples one can find in [41].

**PROBLEM**. Let $M \cong K/K_0$ be a Riemannian homogeneous space with canonical connection $\tilde{\nabla}$, $K=Tr(\tilde{\nabla})$ and let *P(H)* be an invariant structure under the action of *K*. It seems interesting to obtain conditions on *M* and *P(H)*, when the connections $\tilde{\nabla}$ and $\overline{\nabla}$ necessarily coincide.

**2⁰**. We continue our consideration.

**THEOREM 2.16**. Let *P(G)* be an invariant structure over the homogeneous Riemannian space *(M,g)*, where $M \cong K/K_0$ .Then metric *g* is the associated one for *P(G)*.

**Proof**. For fixed point $o \in M$ in the class of *G*-structures conjugated with *P(G)* there exists a structure, let it be *P(G)* itself, such that $(P(G) \cap O(M))_o = P(H)_o$. Since *P(G)* and *O(M)* are invariant under the action of *K*, therefore for any $p \in M$

$$P(H)_{|p} = (P(G) \cap O(M))_{|p} = \varphi(P(H)_{|p}),$$

where $\varphi(o) = p,\ \varphi \in K$ . It is evident that *P(H)* is a reduction of *P(G)* and *H* is the maximal compact subgroup of *G*, $H \in O(n)$.

**QED.**

Let *M* be a simply connected homogeneous Riemannian space with the canonical connection $\tilde{\nabla}$. Then $M \cong K/K_0$ , where $K=Tr(\tilde{\nabla})$ and $K_0$ is the isotropy subgroup of the fixed point $o \in M$, [47] . We remark that $\nabla$ and $\tilde{\nabla}$ are complete. Let $M = M_1 \times M_2 \times ... \times M_r$ be the de Rham decomposition of *M*, then there naturally exists an almost-product structure *P(AP)*: $T_p(M) = T_p^1 \oplus ... \oplus T_p^r$.

**LEMMA 2.17**. The structure *P(AP)* is invariant under the action of *K* and $\tilde{\nabla}$.

**Proof**. Since $K = Tr(\tilde{\nabla}) \subset I^0(M)$, where $I^0(M)$ is the maximal connected isometry group, therefore $\varphi(P(AP)) = P(AP)$ for every $\varphi \in K$, [46] , and *P(AP)* is invariant under the action of *K*. We consider projectors $\pi_i$, $i = 1,...,r$ on $T^i$, $\pi_i^2 = \pi_i$, $\pi_i \pi_j = O$ for $i \neq j$. Since each $\pi_i$, $i = 1,...,r$, is invariant under the action of *K* it is parallel under $\tilde{\nabla}$, [47], and *P(AP)* is invariant under $\tilde{\nabla}$.

**QED.**

**THEOREM 2.18**. Let $M = M_1 \times ... \times M_r$ be the de Rham decomposition of simply connected homogeneous Riemannian space $M \cong K/K_0$ and let *P(G)* be invariant under the action of *K* on *M*. Then, in the class of *G*-structures conjugated with *P(G)* there exists a structure, let it be *P(G)* itself, which has a common subbundle with *P(AP)* i.e. *P(G)* induces the structure $P(G_i)$ on $M_i$, $i = 1,...,r$.

**Proof**. We can choose *P(G)* from the class of conjugated structures in such a way that there exists a frame $u_0 \in (P(G) \cap P(AP))|_0$. Since $K = Tr(\tilde{\nabla})$, then the set $K(u_0) = \{\varphi(u_0) : \varphi \in K\}$ coincides with holonomy fibre bundle $P(u_0)$. From the invariance of *P(G)* and *B(AP)* under the action of *K* it follows that $P(u_0) \subset P(G)$ and $P(u_0) \subset P(AP)$ that is $P(u_0)$ is a common subbundle of *P(G)* and *B(AP)*, hence *P(G)* is reduced to subgroup $\tilde{G} = G_1 \times ... \times G_r$. Therefore we have obtained the structures $P(G_i)$ on $M_i$, $i = 1,...,r$, where $M_i$ is identified with corresponding manifold passing through the point *o*.

**QED.**

It is well known that an invariant almost Hermitian structure on a symmetric Riemannian space M is a Kaehlerian one.

This result is generalised for *G*-structures.

**THEOREM 2.19**. Let $M \cong K/K_0$ be a symmetric Riemannian space, where $K = Tr(\nabla)$ and let *P(G)* be invariant under the action of *K* on *M* . Then the structure *(P(G),g)* is a particular structure.

**Proof**. For fixed frame $u_0 \in P(H) = P(G) \cap O(M)$ we can consider a set $K(u_0) = \{\varphi(u_0) : \varphi \in K\}$, which coincides with holonomy fibre bundle $P(u_0)$. Since *P(H)* is invariant under the action of *K,* then $P(u_0) \subset P(H)$. From *Theorem 2.15* $\nabla$ is the canonical connection of the pair *(P(u),g)*. For $u \in P(u_0)$ $Q_u \subset T_u(P(u_0) \subset T_u((P(H))$ and $\overline{Q}_u = Q_u$ from (2.11). From the right invariance $\overline{Q}_u = Q_u$ for any $u \in L(M)$, that is, $\nabla = \overline{\nabla}$ and *(P(G),g)* is a particular structure.

**QED.**

## *§4. SUBMANIFOLDS AND FOLIATIONS*

**$1^0$**. Let *M'* be a *k*-dimensional manifold isometrically immersed in our *n*-dimensional Riemannian manifold *M*. Since the discussion is local, we may assume, if we want, that *M'* is imbedded in *M*. The submanifold *M'* is also a Riemannian manifold with respect to the restriction of *g* on *M'* and for any $p \in M'$ $T_p(M) = T_p(M') \oplus T_p(M')^{\perp}$. For any $X, Y \in \mathcal{X}(M')$ and $Z \in T(M')^{\perp}$ we have, [47],[73], that

$$\text{(2.14)} \qquad \nabla_X Y = \nabla'_X Y + \alpha(X,Y), \quad \nabla_X Z = -A_Z X + D_X Z ,$$

where $\nabla'_X Y$ and $-A_Z X$ are the tangential components of $\nabla_X Y$ and $\nabla_X Z$, respectively, $\alpha(X,Y)$ and $D_X Z$ are the normal components of $\nabla_X Y$ and $\nabla_X Z$.

The first formula of (2.14) is called the Gauss formula and the second formula is called the Weingarten formula, $\alpha(X,Y)$ is called the second fundamental form of *M'* (or of the immersion). It is well known that $\nabla'$ is the Riemannian connection of *(M',g)*, *D* is a metric connection in the normal bundle $T(M)^{\perp}$, $\alpha$ is symmetric, $A_Z X$ is bilinear in *Z* and *X* and

$$< \alpha(X,Y), Z > = < A_Z X, Y > .$$

Let *P(H)*, $H = G \cap O(n)$ be a structure over *M* with the canonical connection $\overline{\nabla}$. It is easy to see, for example [73], that the above-mentioned results are independent on the torsion of $\nabla$, excluding the symmetry of $\alpha$, therefore we can rewrite them for the metric connection $\overline{\nabla}$

$$\text{(2.15)} \qquad \overline{\nabla}_X Y = \overline{\nabla}'_X Y + \overline{\alpha}(X,Y), \quad \overline{\nabla}_X Z = -A_Z X + \overline{D}_X Z,$$

$$< \overline{\alpha}(X,Y), Z > = < \overline{A}_Z X, Y >$$

We remark that equations of Gauss, Codazzi and Ricci can be also rewritten for the connection $\overline{\nabla}'$. So, $\overline{\nabla}'$ is a metric connection on the submanifold *M'* and we can consider the holonomy fibre bundle $\overline{P}'(u_0)$, of $\overline{\nabla}'$ over *M'*, where $u_0 \in O(M')$.

We shall say that the structure $\overline{P}'(u_0)$ over *M'* is induced by the structure *P(H)* over *M*. Of course, relations between $\overline{P}'(u_0)$ and *P(H)* are very weak or vanish in general case and it is necessary to have addititional restrictions on *M'* and *P(H)* to obtain more close connection.

Let *M'* be an autoparallel submanifold of *M* with respect to $\overline{\nabla}$, [47], that is, for any $X,Y \in X(M)$ $\overline{\nabla}_X Y \in T(M')$. It is equivalent to the condition $\overline{\alpha}(X,Y)=0$.

**THEOREM 2.2O**. Let *M'* be an autoparallel submanifold of *M* and for some frame $u_0 \in P(H)$ $u_0=\{p_0; X_1,...,X_k, X_{k+1},...,X_n\}$, where $p_0 \in M'$, $X_1,...,X_k \in T_{p_0}(M')$ and $X_{k+1},...,X_n \in T_{p_0}(M')^{\perp}$. Then the structure *P(H)* induces by the natural way structures *P(H')* over *M'* and $P(H'^{\perp})$ on the vector bundle $T(M')^{\perp}$.

**Proof**. We consider the holonomy fibre bundle $P'(u_0)$ of $\overline{\nabla}$ over *M'* in *O(M)*. Projections of all the horizontal curves belong to *M'* for $P'(u_0)$. Then $P'(u_0) \subset P(u_0)|_{M'} \subset P(H)|_{M'}$. Since $\overline{\nabla}_X Y \in T(M')$ for any $X,Y \in X(M')$ it is clear that $P'(u_0)=\{p; Y_1,...,Y_k, Y_{k+1},...,Y_n\}$, where $p \in M'$, $Y_1,...,Y_k \in T_p(M')$ and $Y_{k+1},...,Y_n \in T_p(M')^{\perp}$. We introduce $\overline{P}'(u_o)=\{p; Y_1,...,Y_k\}$ and $\overline{P}'^{\perp}(u_0)=\{p; Y_{k+1},...,Y_n\}$. It is evident that $\overline{P}'(u_0)$ is the holonomy fibre bundle of $\overline{\nabla}'=\overline{\nabla}$ over *M'* with a structure group *H'* and $\overline{P}'^{\perp}(u_0)$ is a principal fibre bundle with a structure group $H'^{\perp}$.

**QED.**

**2⁰**. Let *(P(G),g)* be a quasi homogeneous structure, then from *Theorem 2.6* $\overline{\nabla} r^1=0$ and $T(M)=Ker\, r^1 \oplus Im\, r^1$, where $r^1$ is defined by (2.8). We denote $Ker\, r^1$ and $Im\, r^1$ by $T^1$ and $T^2$ respectively.

**LEMMA 2.21**. The distribution $T^1$ is integrable and its maximal integral manifolds are totally geodesic submanifolds with respect to $\nabla$.

**Proof.** From *Theorem 2.6* ,I) it follows that the almost product structure $T(M) = T^1 \oplus T^2$ is invariant with respect to $\overline{\nabla}$ and for any $X, Y \in T^1$ $\overline{\nabla}_X Y = \nabla_X Y \in T^1$.

Therefore $[X, Y] = \nabla_X Y - \nabla_Y X \in T^1$ and $T^1$ is autoparallel under ∇, that is, its maximal integral manifolds are totally geodesic.

**QED.**

We can choose *P(G)* from the class of conjugated structures in such a way that there exists a frame

$$u_0 = \{ p_0; X_1, ..., X_k, X_{k+1}, ..., X_n \},$$

where $X_1, ..., X_k \in T^1$ and $X_{k+1}, ..., X_n \in T^2$.

**THEOREM 2.22**. The structure *P(G)* induces a particular structure $P(G_1)$ and a strict structure $P(G_2)$ on the vector bundle $T^1$ and $T^2$ respectively.

**Proof.** Since $T^1 \oplus T^2$ is invariant under $\overline{\nabla}$, then the holonomy fibre bundle of $\overline{\nabla}$ $P(u_0) = \{ p; Y_1, ..., Y_k, Y_{k+1}, ..., Y_n \}$, where $p \in M$, $Y_1, ..., Y_k \in T^1$ and $Y_{k+1}, ..., Y_n \in T^2$. We know that $P(u_0)$ is a subbundle of *P(G)*. So, we can consider the structures $P(G_1) = \{ p; Y_1, ..., Y_k \}$ on $T^1$ and $P(G_2) = \{ p; Y_{k+1}, ..., Y_n \}$ on $T^2$. In the class of conjugated with $P(G_1)$ structures there exists a structure, let it be $P(G_1)$ itself, that $P(G_1) \cap P(O(k)) = P(H_1)$. As $\overline{\nabla} = \nabla$ on $T^1$, then ∇ is a connection in $P(H_1)$ and the second fundamental tensor field of $P(H_1)$ vanishes. It is obvious that $P(G_2)$ is the strict structure on $T^2$ because each particular vector field belongs to $T^1$.

**QED.**

**REMARK**. It is evident from the proof of *Theorem 2.22* that if we have a general almost product structure $T(M) = T^1 \oplus T^2$, which is invariant with respect to $\overline{\nabla}$, then *P(G)* induces structures $P(G_1)$, $P(G_2)$ on $T^1$, $T^2$ respectively.

**3⁰**. Let $\widetilde{\Lambda}$ be a foliation on *M* defined by an integrable distribution $T^1$ and let *U* be a foliated chart of $\widetilde{\Lambda}$, [66]. We can consider the quotient manifold $U^* \cong U / \widetilde{\Lambda}$ defined by a differentiable projection

$$\pi : U \to U^* : x \mapsto U \cap \Lambda_x ,$$

where $\Lambda_x$ is the fibre passing through $x \in U$. If $T^2 = T^{1\perp}$, then it is evident that $\pi_*$ is an isomorphism between $T_x^2$ and $T_{\pi(x)}(U^*)$. Let $P(G)$ be a structure defined by $O$-deformable tensor fields $K_1,...,K_m$, that is, the frames of $P(G)$ are the same for which every $K_i$, $i = 1,...,m$, has the constant components on $M$. It is well known that the structure $P(G)$ defines a structure $\pi(P(G))$ over $U^*$ by the natural way if we have for each $X \in T^1$

(2.16) $\mathsf{L}_X K_1 = 0,..., \mathsf{L}_X K_m = O$.

It is similarly for the case when we can define the quotient manifold $M^* = M/\widetilde{\Lambda}$.

# CHAPTER 3

# PROBLEMS OF CLASSIFICATION OF G - STRUCTURES

In this chapter we give some results and theorems concerning a classification of *G*-structures on Riemannian manifolds.

In §1 , we consider a classification of Riemannian *G* - structures with respect to the orthogonal group *O(n)*. The classification consists of eight classes . Riemannian *G* - sructures of class $\boldsymbol{T_1}$, are discussed in §2 and nearly particular structures (class $\boldsymbol{T_3}$) in §3 . In §3, we also consider some algebraic constructions connected with a nearly particular quasi homogeneous structure and the main decomposition theorem for a manifold with such a structure.

We refer to [23] , [24] , [41] , [67] .

## *§1. CLASSIFICATION WITH RESPECT TO O(n)*

**$1^0$** . We continue the study of the fixed pair *(P(H),g)*, where $H=G\cap O(n)$, and its second fundamental tensor field $h=\nabla-\overline{\nabla}$, where $\nabla$ is the Riemannian connection, $\overline{\nabla}$ is the canonical that of the structure. We follow [67].

We consider a point $p\in M$ , $T = T_p(M)$ and a vector subspace $\boldsymbol{T}(T)$ of $\overset{3}{\otimes} T^*$

$$\boldsymbol{T}\,(T\,)=\left\{h\in\overset{3}{\otimes}T^*:h_{XYZ}=-h_{XZY}\,,\;X,Y,Z\in T\right\}.$$

If $E_1,...,E_n$ is an arbitrary orthonormal basis of *T*, then $\boldsymbol{T}(T)$ is an Euclidean vector space under the inner product

$$(3.1)\quad <h,h'>=\sum_{i,j,k}h_{E_iE_jE_k}\cdot h'_{E_iE_jE_k}\ .$$

The natural action of the orthogonal group *O(T)* on *T* induces the action on $\boldsymbol{T}(T)$ by the formula

$$(3.2)\ (ah)_{XYZ}=h_{a^{-1}Xa^{-1}Ya^{-1}Z}\,,$$

where $X,Y,Z\in T$ and $a\in O(T)$.

We define

$$\boldsymbol{T}_1(T) = \{h \in \boldsymbol{T} : h_{XYZ} = <X,Y> \beta(Z) - <X,Y> \beta(Y),\ \ \beta \in T^*\},$$

$$\boldsymbol{T}_2(T) = \{h \in \boldsymbol{T}(T) : \boldsymbol{\sigma}\ h_{XYZ} = 0,\ \ c_{12}(h) = 0\},$$

$$\boldsymbol{T}_3(T) = \{h \in \boldsymbol{T}(T) : h_{XYZ} = -h_{YXZ}\},$$

where $X,Y,Z \in T$, $c_{12}(h)Z = \sum_i h_{E_i E_i Z}$ and $\sigma$ denotes the cyclic sum over $X,Y,Z$.

**THEOREM 3.1** [67]. If $dim\ T \geq 3$, then $\boldsymbol{T}(T)$ is the direct sum of the subspaces $\boldsymbol{T}_1$, $\boldsymbol{T}_2$, $\boldsymbol{T}_3$ and $\boldsymbol{T}_i$ each , $i=1,2,3$ is invariant and irreducible under the action of $O(T)$.

It is obvious from this theorem that we can construct 8 classes which are invariant with respect to $O(T)$

TABLE 3.1

| Class | Defining condition |
|---|---|
| particular class | $h = 0$ |
| $\boldsymbol{T}_1$ | $h_{XYZ} = <X,Y> \beta(Z) - <X,Z> \beta(Y),\ \ \beta \in (T)^*$ |
| $\boldsymbol{T}_2$ | $\boldsymbol{\sigma}\ h_{XYZ} = 0,\ c_{12}(h) = 0$ |
| $\boldsymbol{T}_3$ | $h_{XYZ} = -h_{YXZ}$ |
| $\boldsymbol{T}_1 \oplus \boldsymbol{T}_2$ | $\boldsymbol{\sigma}\ h_{XYZ} = 0$ |
| $\boldsymbol{T}_1 \oplus \boldsymbol{T}_3$ | $h_{XYZ} + h_{YXZ} = 2<X,Y>\beta(Z) - <X,Y>\beta(Y) - <Y,Z>\beta(X),\ \beta \in T^*$ |
| $\boldsymbol{T}_2 \oplus \boldsymbol{T}_3$ | $c_{12}(h) = 0$ |
| $\boldsymbol{T}$ | $h_{XZY} = -h_{XZY}$ |

**$2^0$** . From (2.2) for every $p \in M$ the second fundamental tensor $h_p = (\nabla - \overline{\nabla})_p \in \boldsymbol{T}(T)$ and this table is very useful for our study.

**DEFINITION 3.1**. We shall say that the structure *(P(H),g)* has a type $\boldsymbol{T}_{\alpha}$, or belongs to the class $\boldsymbol{T}_{\alpha}$, on *M* if $h_p \in \boldsymbol{T}_{\alpha}(T_pM))$ for each $p \in M$. Here $\boldsymbol{T}_{\alpha}$ is one from the classes from the *table 3.1*.

**THEOREM 3.2**. A pair *(P(H),g)* is a nearly particular structure if and only if its second fundamental tensor field h belongs to $\boldsymbol{T}_3$.

**Proof** evidently follows from the formula

$$h^{+}_{XYZ} = 1/2(h_{XYZ} + h_{YXZ}),$$

where $X,Y,Z \in X(M)$.

**QED.**

Thus, nearly particular structures give examples of the structures having the type $\boldsymbol{T}_3$.

**THEOREM 3.3**. Let *(P(H),g)* be a quasi homogeneous structure having a type $\boldsymbol{T}_\alpha$ for some point $p \in M$. Then this structure belongs to the class $\boldsymbol{T}_\alpha$ on *M*.

**Proof**. We consider a curve $\gamma(t)$, $t \in [O;1]$, $\gamma(O) = p$, $\gamma(1) = q$ and the parallel translation $\bar{\tau}_t$ of $\overline{\nabla}$ along $\gamma(t)$. We denote by $X$, $Y$, $Z$, $E_i$ the vector fields defined on some neighborhood of $\gamma$ which are parallel along $\gamma(t)$. Then from the condition $\overline{\nabla} h = 0$ it follows that

$$(h_{XYZ})_q = h_{\bar{\tau}_q X \bar{\tau}_q Y \bar{\tau}_q Z} = (h_{XYZ})_p$$

Therefore

$$(c_{12}(h)Z)_q = (\sum_i h_{E_i E_i Z})_q = (c_{12}(h)Z)_p,$$

$$(\sigma\, h_{XYZ})_q = (\sigma\, h_{XYZ})_p.$$

If $h\in \boldsymbol{T}_1\oplus \boldsymbol{T}_2$ , then $\beta_q$ is correctly defined form as the parallel translation of $\beta_p$ along $\gamma(t)$. The form $\beta_q$ is independent on the choice of a curve because $h$ and $g$ are parallel along $\gamma(t)$ and unique in the point $q$.

Using defining conditions from the *table 3.1* it is easily to see that all the classes are invariant under the parallel translation.

**QED.**

**REMARK**. Let $h$ be a homogeneous Riemannian structure on a homogeneous Riemannian space $M$ and let $P(u_0)$ be the holonomy fibre bundle of $\tilde{\nabla} = \nabla - h$ containing the frame $u_0 \in O(M)$, (see *theorem (AS)* in §3 of chapter 2) . The examples homogeneous Riemannian structures of the all 8 types are considered in [67]. From *theorem 2.15* it follows that the second fundamental tensor field of $(P(u_0),g)$ coincides with $h$. So, if $h\in \boldsymbol{T}_\alpha$, then the holonomy fibre bundle of $\tilde{\nabla} = \overline{\nabla}$ gives an example of a $G$ - structure of the corresponding class.

## *§2. RIEMANNIAN G-STRUCTURES HAVING TYPE* $\boldsymbol{T}_1$ *AND* $\boldsymbol{T}_2$

**1⁰**. We say that a structure $(P(H),g)$ belongs to the class $\boldsymbol{T}_1$ on $M$ if $h_p\in \boldsymbol{T}_1(T_p(M))$ for each $p \in M$, that is,

$$(3.3)\quad h_{XYZ} = \langle X,Y\rangle\beta(Z) - \langle X,Z\rangle\beta(Y),$$

where $\beta$ is a nonzero 1-form on $M$ and $X,Y,Z \in X(M)$.

We can define the nonzero vector field $\xi$ on $M$ by the formula

$$\langle \xi,X\rangle = \beta(X)$$

and in this case from (3.3) it folows that

$$(3.4)\quad h_XY= \langle X,Y\rangle\xi - \langle \xi,Y\rangle X.$$

It follows from (2.3) that the torsion $\overline{T}$ of the canonical connection has the following form

$$(3.5)\quad \overline{T}_XY = -2h_{\bar{X}}Y = h_YX - h_XY = \langle Y,X\rangle\xi - \langle\xi,X\rangle Y - \langle X,Y\rangle\xi + \langle\xi,Y\rangle X$$

$$= \langle\xi,Y\rangle X - \langle\xi,X\rangle Y$$

Let $L = [\xi]$ be the one-dimensional distribution on $M$ defined by the vector field $\xi$ and let $V = L^{\perp}$ be the orthogonal complement of $L$. We have the almost-product structure

(3.6) $T(M) = L \oplus V.$

**PROPOSITION 3.4**. The structure (3.6) is invariant with respect to the canonical connection $\overline{\nabla}$.

**Proof**. We take $p \in M$, $\zeta = \|\xi\|^{-1}\xi$, $u = \{p;\ \zeta_p = E_0, E_1, \ldots, E_{n-1}\} \in O(M)$ and such $\overline{E}_o, \overline{E}_1, \ldots, \overline{E}_{n-1} \in \overline{Q}_u = (Ker\,\overline{\omega})_u$ that $\pi_*(\overline{E}_i) = E_i$, $i = O, \ldots, n-1$, where $\pi$ is the canonical projection and $\overline{\omega}$ is the $\underline{\boldsymbol{O}}$-valued form of the canonical connection in $P(H)$. According to [7] we have that

$$h_{E_i} E_j = u\tau(\overline{E}_i)u^{-1}E_j,$$

where $\tau = \omega - \overline{\omega}$ and ω is $\underline{\boldsymbol{O}}$ - valued form of the Riemannian connection. From (3.4) it follows that

$$\tau(\overline{E}_i) = \left[\begin{array}{c|c} 0 & 0 \ \cdots \ \underbrace{\|\xi\|}_{i} \ \cdots \ 0 \\ \hline 0 & \\ \vdots & \\ -\|\xi\| & 0 \\ \vdots & \\ 0 & \end{array}\right] \in \underline{\boldsymbol{m}} \subset \underline{\boldsymbol{O}},\ i = 1, \ldots, n-1;\ \tau(\overline{E}_0) = [O]$$

Let $\underline{\boldsymbol{m}}_0$ be the linear span of the $\tau(\overline{E}_1), \ldots, \tau(\overline{E}_{n-1})$, that is,

$$\underline{\boldsymbol{m}}_0 = \left[\begin{array}{c|c} 0 & B \\ \hline -B & 0 \end{array}\right] \subset \underline{\boldsymbol{m}} \subset \underline{\boldsymbol{O}},\quad \underline{\boldsymbol{h}}_0 = \underline{\boldsymbol{m}}_0^{\perp} = \left[\begin{array}{c|c} 0 & 0 \ \cdots \ 0 \\ \hline 0 & \\ \vdots & A \\ 0 & \end{array}\right] \subset \underline{\boldsymbol{O}}.$$

Since $\underline{\boldsymbol{m}}_0 \subset \underline{\boldsymbol{m}}$ then $\underline{\boldsymbol{h}} \subset \underline{\boldsymbol{h}}_0$ and $H$ is a Lie subgroup of the Lie group $K = 1 \times O(n-1)$, that is, $H = 1 \times H'$, where $H' \subset O(n-1)$. The structure $P(H)$ is invariant under the parallel translation of the connection $\overline{\nabla}$, hence its extension $P(K)$ is invariant too.

**QED.**

**PROPOSITION 3.5**. The distribution $p \to V_p$ is integrable on *M*.

**Proof**. We consider $X,Y \in V$, $X \perp Y$. It folows from (4.4) that $h_X Y \in V$ and $h_Y X \in V$, proposition 3.4 implies that $\overline{\nabla}_X Y \in V$ and $\overline{\nabla}_Y X \in V$, therefore

$$\nabla_X Y = \overline{\nabla}_X Y + h_X Y \in V, \quad \nabla_Y X = \overline{\nabla}_Y X + h_Y X \in V$$

and $[X,Y] = \nabla_X Y - \nabla_Y X \in V$.

**QED.**

**2⁰**. We have obtained the foliation $\tilde{\Lambda}$ on *M* defined by the integrable distribution *V*. Let *M'* be a fibre of this foliation containing a point $p \in M$ and let $\nabla'$ be the Riemannian connection on the submanifold *M'*.

**LEMMA 3.6**. $\nabla' = \overline{\nabla}$ on M.

**Proof**. If $X,Y \in \mathcal{X}(M')$ then, from (3,4), we have

$$\nabla'_X Y = \nabla_X Y - <\nabla_X Y, \zeta> \zeta = \nabla_X Y - <h_X Y, \zeta> \zeta$$
$$= \nabla_X Y - <X,Y> \|\xi\| \zeta = \nabla_X Y - <X,Y> \xi = \overline{\nabla}_X Y.$$

**QED.**

This lemma implies that the second fundamental from of submanifold *M′* is defined by the formula

(3.7) $\alpha(X,Y) = <X,Y>\xi$

for any vector fields *X,Y* tangent to *M′*.

**THEOREM 3.7**. The structure *P(H)* induces the particular structure *P(H′)* over *M′*.

**Proof**. From the proof of *proposition 3.4* we see that $H = 1 \times H'$ and we can consider the structure *P(H′)* over *M′*

$P(H') = \{\{p; E_1,...,E_{n-1}\}\}$,

where $p \in M'$ and $\{p; \zeta_p = E_0, E_1, ..., E_{n-1}\} \in P(H)$.

Since *P(H)* is invariant under $\overline{\nabla}$ over *M*, then *P(H′)* is invariant under $\overline{\nabla}_{|M'} = \nabla'$. So, from the construction of the second fundamental tensor field *h*, it is evident that $\overline{\omega}' = \omega'$ and $h' = \nabla' - \overline{\nabla}' = 0$ on $M'$.

**QED.**

For any vector fields *X,Y,Z,W* tangent to *M′* we have , [73] , equation of Gauss

(3.8)
$$<R_{XY}Z,W> = <R'_{XY}Z,W> - <\alpha(X,W),\alpha(Y,Z)> + <\alpha(Y,W),\alpha(X,Z)>,$$

where *R*, *R′* are the curvature tensors of $\nabla,\nabla'$ respectively and, using (3.7), we have obtained

(3.9) $$<R_{XY}Z,W> = <R'_{XY}Z,W> + \|\xi\|^2(<Y,W><X,Z> - <X,W><Y,Z>).$$

**3⁰**. Let *M* be a space of constant curvature *k*, [73] , then for any $X,Y,Z \in \mathfrak{X}(M)$ the curvature tensor of $\nabla$ is given by

(3.10) $R_{XY}Z = k(<Y,Z>X - <X,Z>Y)$

**TEOREM 3.8**. Every fibre *M′* of the foliation $\tilde{\Lambda}$ is a space of constant curvature $k + \|\xi\|^2$ and $\|\xi\|$ is constant on $M'$( $\|\xi\| = c$).

**Proof**. It follows from (4.9) and (4.1O) that

$$<R'_{XY}Z,W> = k(<Y,Z><X,W> - <X,Z><Y,W>) - \|\xi\|^2(<Y,W><X,Z> - <X,W><Y,Z>)$$
$$= (k + \|\xi\|^2)(<Y,Z><X,W> - <X,Z><Y,W>),$$

where $X,Y,Z,W \in \mathfrak{X}(M')$. For each plane $[X \wedge Y]$ in the tangent space $T_p(M')$, where *X,Y* is an orthonormal basis for $[X \wedge Y]$ ,the sectional curvature $k([X \wedge Y]) = <R'_{XY}Y,X> = k + \|\xi\|^2$ and it depends only on the point $p \in M'$. From the Schur\`s theorem, [73] , $(k + \|\xi\|^2)$ is constant on *M′* and $\|\xi\| = c$ on *M′*.

**QED.**

**REMARK**. Perhaps $\|\xi\|$ is not constant on *M*.

**$4^0$**. Now, we shall generalise a theorem of F. Tricerri and L. Vanhecke from [67] considered for h.R.s. on *M* on the case of quasi-homogeneous structures having the type $\boldsymbol{T}_1$.

Let $\|\xi\|$ be constant on *M* , then we shall obtain that *P(H)* is a quasi homogeneous structure.

**PROPOSITION 3.9**. $\|\xi\| = c$ on *M* if and only if $\overline{\nabla} h = O$.

**Proof**. It is obvious that $\|\xi\| = c$ if and only if $\overline{\nabla}_X \xi = c\overline{\nabla}_X \zeta = 0$ for any $X \in \boldsymbol{X}(M)$. For an integral curve $\gamma(t)$ of the vector field *X* we can consider vector fields *Y, Z* which are parallel along $\gamma(t)$, that is, $\overline{\nabla}_X Y = \overline{\nabla}_X Z = 0$. So, we have got

$$(\overline{\nabla}_X h)(Y,Z) = \overline{\nabla}_X h(Y,Z) = \overline{\nabla}_X (<Y,Z>\xi - <\xi,Z>Y) = <Y,Z>\overline{\nabla}_X \xi + <\overline{\nabla}_X \xi, Z>Y$$

From this identity the proposition follows.

**QED.**

**THEOREM 3.1O**. Let *(P(H),g)* be a quasi homogeneous structure over *M* belonging to the class $\boldsymbol{T}_1$. Then *(M,g)* is locally isometric to $R\times M'$ with the Riemannian metric $ds^2 = c^2 dt^2 + e^{-2c^2 t} g'$, where $g'$ is the induced Riemannian metric on $M'$ and $dt(\xi)=1$.

**Proof** of this theorem is the same to one considered in [67] for the case when h is a homogeneous Riemannian structure (h.R.s.) on *M* ( see *chapter* 2, §3).

**REMARK**. *Proposition 3.9* and *theorem 3.1O* allow us to construct examples of quasi homogeneous structures having the type $\boldsymbol{T}_1$.

**$5^0$**. Let *h* be a h.R.s. on *M* , that is, $\overline{\nabla} h = 0$, $\overline{\nabla}\overline{R} = 0$, where $\overline{R}$ is the curvature tensor field of $\overline{\nabla}$. The following results were considered in [67] .

1) All the nontrivial h.R.s. on surfaces have the type $\boldsymbol{T}_1$.
2) If *h* is a nonzero h.R.s. on a surface *M*, then *(M, g)* has a constant negative curvature.
3) Let *M* be a complete, simply connected surface, then *M* admits a nontrivial h.R.s. if and only if *(M, g)* is isometric to the hyperbolic surface $H^2$.
4) If a manifold *M* admits a h.R.s. of the class $\boldsymbol{T}_1$, then *(M, g)* has a constant

negative curvature.

5) Let $M$ be a complete, simply connected Riemannian manifold ,then there exists a nonzero h.R.s. $h \in \boldsymbol{T}_1$ if and only if $(M, g)$ is isometric to the hyperbolic space $H^n$.
6) For any fibre $M'$ of the foliation $\tilde{\Lambda}$ we have $R' = O$.

**6⁰**. O.Kowalski and F.Tricerri in [49] otained the following classification results:

I) Each connected, complete and simply connected Riemannian manifold $(M, g)$ of dimension $n = 3$ admitting a non-trivial homogeneous structure of class $\boldsymbol{T}_2$ is isometric to one of the following homogeneous Riemannian spaces:

(a) The ordinary sphere $S^3(R)$ with the sectional curvature $k = 1/R^2$.

(b) The group $E(1,1)$ with a left invariant metric $g$ such that $dim\, I(M, g) = 3$. The admissible metrics form a one-parameter family.

(c) The universal covering group $\widetilde{SL(2,R)}$ of $SL(2;R)$ with a left-invariant metric $g$ such that dim $I(M,g) = 3$. The admissible metrics form a two-parameter family.

(d) The group $\widetilde{SL(2,R)}$ with a left-invariant metric g such that dim $I(M,g) = 4$. The admissible metrics form a two-parameter family.

(e) The group $SU(2) \cong S^3$ with a left-invariant metric g such that dim $I(M,g)=4$. The admissible metrics form a two-parameter family.

(f) The Heisenberg group $H_3$ with any left-invariant metric. These metrics form a one-parameter family.

Moreover, the spaces sub(b) are 4-symmetric, but not naturally reductive; the spaces sub(c) are also not naturally reductive. The spaces sub(d)-(f) are non-symmetric naturally reductive spaces.

In dimension 3, a non-symmetric space admitting a non-trivial homogeneous structure of class $\boldsymbol{T}_3$ also admits a non-trivial homogeneous structure of class $\boldsymbol{T}_2$.

II) Each connected, complete and simply connected Riemannian manifold $(M, g)$ of dimension $n = 4$ admitting a non-trivial homogeneous structure of type $\boldsymbol{T}_2$ is isometric to one of the following Riemannian homogeneous spaces:

(a) $(M, g)$ is the Cartesian space $\boldsymbol{R}^4[x, y, z, t]$ provided with a Riemannian metric of the form

$$g = e^{-2\alpha t} dx^2 + e^{-2\beta t} dy^2 + e^{2(\alpha+\beta)t} dz^2 + dt^2$$

where $\alpha \neq 0,\ \beta \neq 0,\ \alpha \neq -\beta$ are real parameter.

(b) *(M, g)* is the Cartesian space $\mathbf{R}^4[x, y, u, v]$ provided with a Riemannian metric of the form

$$g = \{-x + \sqrt{x^2 + y^2 + 1}\}du^2 + \{x + \sqrt{x^2 + y^2 + 1}\}dv^2$$

$$-2y\, du\, dv + \lambda^2 (1 + x^2 + y^2)^{-1} \{(1 + y^2)dx^2 + (1 + x^2)dy^2\}$$

$$-2xy\, dx\, dy,$$

where $\lambda \succ 0$ is a real parameter.

(c) *(M, g)* is a Riemannian product $(M_3, g')\times\mathbf{R}$ , where $(M_3, g')$ is one of the spaces given in *I)*.

The space (a) and (b) are always irreducible as Riemannian manifolds.

The space (a) can be described as the matrix group *G* whose elements are

$$\begin{pmatrix} e^{\alpha t} & 0 & 0 & x \\ 0 & e^{\beta t} & 0 & y \\ 0 & 0 & e^{(\alpha+\beta)t} & z \\ 0 & 0 & 0 & t \end{pmatrix}$$

equipped with a special left-invariant metric.

The space (b) is the homogeneous space *G/H* , where *G* is the group of all positive equiaffine transformations of the plane and *H* is the subgroup of all rotations around the origin, equipped with a left-invariant metric from a special one-parameter family. All spaces of type (b) are 3-symmetric Riemannian spaces.

Let *h* be a h.R.s. of the type $\underline{T}_2$ on *M* , where *M* is one from the manifolds considered in I), II) and $P(u_0)$ the holonomy bundle of the canonical connection $\tilde{\nabla} = \nabla - h$. According to *Theorem 2.15,* the second fundamental tensor field of $P(u_0)$ coincides with *h* and $P(u_0)$ gives an example of a *G*-structure having the type $\boldsymbol{T}_2$.

## §3. NEARLY PARTICULAR STRUCTURES

**1⁰**. It follows from *Theorem 3.2* that the nearly particular strucures coincide with those having a type $\boldsymbol{T}_3$. In this case we see that for any $X,Y,Z \in \mathfrak{X}(M)$

(3.11) $$h_{XYZ} = -h_{XZY} = h_{ZXY}$$

If $dim\ M = 2$, then (3.11) implies that $h = O$ , thus, every nearly particular structure is a particular one.

Let $M \cong K/K_o$ be a naturally reductive Riemannian homogeneous space with $K$–invariant Riemannian metric $g = <\ ,\ >$ and with $ad(K_0)$ - invariant decomposition

$$\underline{\mathbf{k}} = \underline{\mathbf{k}} \oplus \underline{\mathbf{m}}$$

such that for any $X,Y,Z \in \underline{\mathbf{m}}$

$$<[X,Y]_{\underline{\mathbf{m}}}, Z> + <[X,Z]_{\underline{\mathbf{m}}}, Y> = O$$

where $K = Tr(\overline{\nabla}\ )$ and $<\ ,\ >$ is the inner product in $\underline{\mathbf{m}}$ induced by $g$, (see [47]). Then the canonical connection $\tilde{\nabla}$ satisfies (2.13) and defines such a h.R.s. $h = \nabla - \tilde{\nabla}$ on M that $h_X Y = \frac{1}{2}[X,Y]_{\underline{\mathbf{m}}}$. It is evident that $h_X Y = -h_Y X$. Let $P(u_0)$ be the holonomy fibre bundle of $\tilde{\nabla}$ containing a frame $u_0 \in O(M)$. Then from *Theorem 2.15* it follows that $\nabla - \overline{\nabla} = \nabla - \tilde{\nabla}$ and the second fundamenal tensor field of $P(u_0)$ coincides with the homogeneous Riemannian structure $h$. Thus $P(u_0)$ gives an example of a nearly particular structure on the space $M$.

If we consider a nearly particular structure on a complete, simply connected manifold $M$ for which $\overline{\nabla} h = 0$, $\overline{\nabla} R = 0$ then, from *Theorem 2.4*, $(M, g)$ is a naturally reductive Riemannian homogeneous space with h.R.s. $h$.

It is well known that all the irredusible symmetric spaces are naturally reductive. Every isotropy irreducible homogeneous space belongs to this class , [71]. Each nearly Kaehlerian 3-symmetric space is naturally reductive, [33] .

The case of compact Lie groups was discussed in [4].

See the survey about the recent research on naturally reductive Riemannian homogeneous space in [5O].

Let $M$ be a 3-dimensional manifold. Some results about this case one can find in [67].

A) Let $M$ be a connected, simply connected, complete manifold, $dim\ M = 3$, and there exists a h.R.s. $h$ on $M$ having a type $\boldsymbol{T}_3$, $h \neq 0$. Then $M$ is isometric to $R^3$, $S^3$, $H^3$ or to one of the following Lie groups with a left-invariant Riemannian metric: 1) $SU(2)$ , 2) $\widetilde{SL(2,R)}$(universal covering of the group $SL(2,R)$), 3) the Heisenberg group.

B) If there exists a h.R.s. $h$ on $M$, $dim\ M = 3$, having a type $\boldsymbol{T}_1 \oplus \boldsymbol{T}_3$ then $h$ belongs to $\boldsymbol{T}_1$ or $\boldsymbol{T}_3$ separately.

**REMARK**. As we shall see later on, a nearly particular structure generalizes a notion of the nearly Kaehlerian structure [32], [34] .

**$2^0$**. We consider now a nearly particular quasi homogeneous structure *P(H)* on a Riemannian manifold *(M, g)*, that is, $h = h^-$ and $\overline{\nabla} h = 0$.

Let *T* be an Euclidian vector space under an inner product < , > and there is defined a bilinear operation

$$* : V \times V \to V$$

such that for any $X,Y,Z \in T$

1) $X * X = O$,  2) $<X*Y, Z> = <X, Y*Z>$.

Constructed algebra is called a *QR*–algebra. The theory of *QR*–algebras was developed in [41] .

A *QR*–algebra is called simple if it have no ideals except *{O}* and *T*, semisimple if it have no Abelian ideals except *{O}*.

**DEFINITION 3.2**. We say that a *QR*–algebra *T* is a direct product of ideals $I_0, I_1, ..., I_r$ if $T = I_0 \oplus I_1 \oplus ... \oplus I_r$, where $I_i$ is orthogonal to $I_j$ for every $i \neq j$. In this case we denote $T = I_0 \overset{\bullet}{\times} I_1 \overset{\bullet}{\times} ... \overset{\bullet}{\times} I_r$.

**THEOREM 3.11** [41]. Each *QR* - algebra *T* is a direct product of its Abelian ideal $I_0$ and non-Abelian simple ideals.
This decomposition is unique up to an order of factors.

Let *T* be a *QR* - algedra and $T = I_0 \overset{\bullet}{\times} I_1 \overset{\bullet}{\times} ... \overset{\bullet}{\times} I_r$ is the direct product from *Theorem 3.11*. Let *dim T = n* and $E_1, ..., E_n$ be an arbitrary orthonormal basis of *T*.

**DEFINITION 3.3**. We shall call a linear mapping $A : T \to T$ a fundamental operator, if it is defined by

$$(3.11) \quad < AX, Y > = \sum_{k=1}^{n} < X * E_k, Y * E_k >$$

It is easily checked an independence of *A* from a choice of an orthonormal basis (see(2.9)).

**THEOREM 3.12**. Let $T$ be a $QR$ - algedra and $A$ be its fundamental operator. Then we have

1) $<AX,Y> = <X,AY>$, 2) $Ker\ A = I_0$, 3) $A(I_i)=I_i\ ,i=1,..,r$.

**Proof**. 1) evidently follows from (3.11).
2). We have from (3.11) that

$$< AX,X >= \sum_{k=1}^{n} \left\| X * E_k \right\|^2$$

and, if $X \in Ker\ A$, then $X*E_k = O$, $k = 1,...,n$; therefore $X*Y = O$ for any $Y \in T$ and $X \in I_0$. Conversely, if $X \in I_0$, then $X*E_k = O$, $k = 1,...,n$ and it follows from (3.11) that $<AX,Y> = O$ for every $Y \in T$ , that is $X \in Ker\ A$. Thus $I_0 = Ker\ A$.

3) Let $E_1^1,...,E_{S_1}^1,...,E_1^r,...,E_{S_r}^r$ be such an orthonormal basis of $T$ that $E_1^i,...,E_{S_i}^i$ is a basis of $I_i$. If $X \in I_i$ and $Y \in I_j$, where $i \neq j$, then

$$< AX,Y >= \sum_{k=1}^{S_i} < X * E_k^i , Y * E_k^i > = 0\,,$$

therefore $AI_i \subset I_i$. Since A is nonsingular on $I_0^{\perp}$, hence $AI_i = I_i$, $i = 1,...,r$.

**QED.**

**COROLLARY**. There exists an orthonormal basis of $T$ composed from eigenvectors of the operator $A$ that the matrix of $A$ has following form in this basis

$$\begin{bmatrix}
\begin{matrix} 0 & & 0 \\ & \ddots & \\ 0 & & 0 \end{matrix} & & & 0 \\
& \begin{matrix} \lambda_1 & & 0 \\ & \ddots & \\ 0 & & \lambda_1 \end{matrix} & & \\
& & \ddots & \\
0 & & & \begin{matrix} \lambda_e & & 0 \\ & \ddots & \\ 0 & & \lambda_e \end{matrix}
\end{bmatrix}$$

where $\lambda_j > 0, \quad j = 1,...,l$.

Thus $T$ is a direct sum of the orthogonal each other proper subspaces of $A$, $T = I_0 \oplus \Lambda_1 \oplus ... \oplus \Lambda_e$

**DEFINITION 3.4**. A *QR* - algebra $T$ is called a *QRA* – algebra if $l = r$ and $\Lambda_i = I_i, \ i = 1,...,r$

**EXAMPLE**. If $A$ has only one eigenvalue or two (one of which is equal to zero), then a *QR* - algebra $T$ is a *QRA* - algebra.

**PROPOSITION 3.13**. If $T$ is a *QRA* - algebra , then

$A(X*Y) = AX*Y = X*AY, \ X,Y \in T.$

**Proof**. The operation $*$ is bilinear and the operator $A$ is linear, therefore it is sufficienly to prove this formula for vectors of some basis of $T$. If $T$ is a *QRA* – algebra, then above-mentioned basis $E_1^1,...,E_{S_1}^1,...,E_1^r,...,E_{S_r}^r$ consists of the eigenvectors of the operator $A$. Thus we have

$$E_k^i * E_e^j = 0, \ i \neq j; \ AE_j^i * E_k^i = E_j^i * AE_k^i = A(E_j^i * E_k^i) = \lambda_i (E_j^i * E_k^i).$$

**QED.**

**$3^0$**. We apply this algebraic construction to our situaion.

**THEOREM 3.14**. Let $P(H)$ be a nearly particular quasi homogeneous structure on a manifold $M$ . There exists unique up to an isomorphism *QR* - algebra $T$ associated to $P(H)$.

**Proof**. We consider $p \in M$ , $T = T_p(M)$ and define an operation $*$ by $X * Y = h_X Y$ , where $X,Y \in T$ and $h_X Y$ is the second fundamental tensor of structure $P(H)$ on $T$. It follows from (3.11) that $(T,*)$ is a *QR* - algebra. Let $q$ be another point of $M$ and $T' = T_q(M)$. We have to show that the *QR* - algebra $(T',*)$ is isomorphic to $(T,*)$. Let $\gamma$ be a curve segment in $M$ connecting $p$ and $q$ , and let $\bar{\tau}$ be the parallel translation along $\gamma$ in the connection $\overline{\nabla}$ . $\bar{\tau}$ is an isomorphism of the vector spaces $T_p(M)$ and $T_q(M)$. Since $\overline{\nabla} g=0$ then $<\overline{\tau} X, \overline{\tau} Y> = <X,Y>$ for $X,Y \in T_p(M)$. From the identity $\overline{\nabla} h=0$ it follows that $(\overline{\tau} X)*(\overline{\tau} Y)=\overline{\tau}(X,Y)$.

Thus, we have obtained that $\overline{\tau}$ is an isomorphism of $QR$ - algebras $(T, *)$ and $(T', *)$.

**QED.**

It follows from (2.6) that in our case for any $X,Y,Z \in \mathcal{X}(M)$

(3.12) $$R_{XY}Z = \overline{R}_{XY}Z + [h_X, h_Y]Z + h_{\overline{T}_X Y} Z .$$

We can define so-called induced Ricci tensor by the formula

(3.13) $$ri(X,Y) = \sum_{k=1}^{n} <(\overline{R} - R)_{XE_k} Y, E_k>,$$

where $E_1,...,E_n$ are orthonormal vector fields on some neighbourhood.

**LEMMA 3.15**. $ri(X,Y) = <r^1 X, Y> = \sum_{k=1}^{n} <h_X E_k, h_Y E_k>$, where $r^1$ is defined by (2.8) .

**Proof**. Using (2.3), (3.11), (3.12) we have

$$ri(X,Y) = \sum_{k=1}^{n} <(\overline{R} - R)_{XE_k} Y, E_k> = \sum_{k=1}^{n} (h_{E_k h_X YE_k} - h_{X h_{E_k} YE_k} + 2h_{h_X E_k YE_k})$$

$$= \sum_{k=1}^{n} (2h_{YE_k h_X E_k} - h_{XE_k h_Y E_k} - h_{E_k E_k h_X Y})$$

$$= \sum_{k=1}^{n} (2<h_Y E_k, h_X E_k> - <h_X E_k, h_Y E_k>)$$

$$= \sum_{k=1}^{n} <h_X E_k, h_Y E_k> = <r^1 X, Y>$$

**QED.**

For every $p \in M$ and any $X,Y \in T_p(M)$ we can define $X*Y = h_X Y$ and obtain a $QR$–algebra $(T, *)$. It follows from (3.11) that the fundamental operator $A$ of this algebra coincides with $r^1$ and $\overline{\nabla} r^1 = 0$ from *Theorem 2.6*. The proper subspaces of $A = r^1$ define the differentiable distributions $p \to \Lambda_{ip}$, $i = 1,...,l$, which are orthogonal each other and $Ker\, r^1 = I_0$ is the differentiable ideal on *M*. If $Y \in \Lambda_i$, then $A\overline{\nabla}_X Y = \overline{\nabla}_X AY = \lambda_i \overline{\nabla}_X Y$ and $\Lambda_i$, $i = 1,...,l$ are invariant with respect to $\overline{\nabla}$. Every ideal $I_j$, $j = 1,...,r$ is a direct sum of some proper subspaces from *Theorem 3.12*, 3), therefore it defines a differentiable distribution $p \to I_{jp}$ on *M* and we have

(3.14) $T_p(M) = I_{0p} \oplus I_{1p} \oplus ... \oplus I_{rp}$

**PROPOSITION 3.16**. The distributions $I_j, j = O,...,r$ are invariant with respect to $\nabla$, in particular, they are integrable and the maximal integral manifold of $I_j$ is a totally geodesic submanifold of *M*.

**Proof**. For any $Y \in I_j$, $j = O,...,r$ we can write $\nabla_X Y = \overline{\nabla}_X Y + h_X Y$. $I_j$ is invariant with respect to $\overline{\nabla}$, hence $\overline{\nabla}_X Y \in I_j$ and $h_X Y = X * Y \in I_j$ because $I_j$ is the ideal. Thus $\nabla_X Y \in I$ and $I_j$ is invariant with respect to $\nabla$.

**QED.**

The structure (3.14) and *P(H)* are invariant with respect to $\overline{\nabla}$, therefore, in the class of *G* - structures conjugated with *P(H)* there exists a sructure $P(\overline{H})$, which has a common subbundle with (3.14) , that is $\overline{H} = H_0 \times ... \times H_r \subset O(n)$ and $P(\overline{H})$ induces the structure $P(H_j)$ on $I_j$. $P(H_0)$ is invariant with respect to $\overline{\nabla} = \nabla$ on $I_0$, therefore it is particular. For $P(H_j)$ , $j = 1,...,r$ $h_X \neq 0$ for any $X \in X(M)$, $X \neq O$, hence $P(H_j)$ is strict.

We can unite now all the above-mentioned facts in the main decomposition theorem for a manifold with a nearly particular quasi homogeneous structure *P(H)*.

**THEOREM 3.17**. Let *M* be a complete, simply connected manifold and *P(H)* be a nearly particular quasi homogeneous structure over *M*. Then we have

1) *M* is isometric to direct product $M_0 \times M_1 \times ... \times M_r$, $i = 0,...,r$, where each $M_j$ is a totally geodesic submanifold of *M* ($M_j$ is the maximal integral submanifold of $I_j$, see *Proposition 3.16*, passing through a fixed point of *M*);

2) there exists a structure $P(\overline{H})$ conjugated with *P(H)* in *O(M),* which induces a particular structure $P(H_0)$ over $M_0$ and a strict, nearly particular, quasi homogeneous one $P(H_j)$, $j = 1,...,r$ over $M_j$;
3) if associated to *P(H)* the *QR* - algebra *(T, ∗)* is the *QRA* - algebra, then $ri(X,Y) = \lambda_j < X,Y >$ on $M_j$, where $\lambda_j$ is an eigenvalue of *A* (see (3.13)).

**4^0^**. Let *P(H)* be a nearly particular structure, maybe $\overline{\nabla} h \neq 0$ and *R*, $\overline{R}$ denote the curvature tensor fields of $\nabla, \overline{\nabla}$ respectively.

**PROPOSITION 3.18**. $<(R-\overline{R})_{XY}Y,X> = \|h_X Y\|^2$, $X,Y \in \mathcal{X}(M)$.

**Proof**. The condition $h_X X = O$ implies that $\nabla_X X = \overline{\nabla}_X X$ for any $X \in \mathcal{X}(M)$. From (2.3) $\overline{\nabla}_X Y - \overline{\nabla}_Y X - [X,Y] = -2h_X Y$ and $h_{XYZ} = h_{ZXY}$ from (3.11).

$$<(R-\overline{R})_{XY}Y,X> = <\nabla_X \nabla_Y Y,X> - <\nabla_Y \nabla_X Y,X> - <\nabla_{[X,Y]}Y,X>$$

$$- <\overline{\nabla}_X \overline{\nabla}_Y Y,X> + <\overline{\nabla}_Y \overline{\nabla}_X Y,X> + <\overline{\nabla}_{[X,Y]}Y,X>$$
$$= X<\nabla_Y Y,X> - <\nabla_Y Y,\nabla_X X> - Y<\nabla_X Y,X> + <\nabla_X Y,\nabla_Y X>$$
$$- X<\overline{\nabla}_Y Y,X> + <\overline{\nabla}_Y Y,\overline{\nabla}_X X> + Y<\overline{\nabla}_X Y,X> - <\overline{\nabla}_X Y,\overline{\nabla}_Y X>$$
$$- h_{[X,Y]YX}$$
$$= -Y[X<Y,X> - <Y,\nabla_X X> - X<Y,X> + <Y,\overline{\nabla}_X X>]$$
$$+ <\nabla_X Y,\nabla_Y X> - <\overline{\nabla}_X Y,\overline{\nabla}_Y X> - h_{[X,Y]YX}$$
$$= <\overline{\nabla}_X Y + h_X Y,\overline{\nabla}_Y X + h_Y X> - <\overline{\nabla}_X Y,\overline{\nabla}_Y X> - h_{[X,Y]YX}$$
$$= <\overline{\nabla}_X Y,h_Y X> + <h_X Y,\overline{\nabla}_Y X> - \|h_X Y\|^2 - h_{[X,Y]YX}$$
$$= <\overline{\nabla}_Y X + [X,Y] - 2h_X Y,h_Y X> + <\overline{\nabla}_Y X,h_X Y>$$
$$- \|h_X Y\|^2 - h_{[X,Y]YX}$$
$$= \|h_X Y\|^2 + h_{YX[X,Y]} - h_{[X,Y]YX} = \|h_X Y\|^2$$

**QED.**

This proposition can be useful when we study concrete structures on spaces of constant curvature *k*.

# CHAPTER 4
# RIEMANNIAN REGULAR σ-MANIFOLDS

Symmetric spaces and their generalizations play an important role in modern differential geometry and its applications, [48],[58]. In this chapter we introduce and study the so-called Riemannian regular σ-manifolds, which generalize on the one hand the spaces with reflexions [56] and on the other hand the Riemannian regular s-manifolds [48]. We want to point out that the term "subsymmetry" was first used in [62] .

In §1, we give the axioms of Riemannian (locally) regular σ-manifold and determine the canonical connection and the foliation of mirrors. The structure of Riemannian regular σ-manifold is considered in §2. Every such a manifold can be described as a fibre bundle over a regular *s*-manifold. §3 is devoted to the Lie algebra of infinitesimal automorphisms of Riemannian regular σ-manifolds. In §4, we consider orbits under the action of the structural group *G* and the main examples of the Riemannian regular σ-manifolds. Riemannian locally regular σ-manifolds are discussed in §5.

The works [26] , [27] , [48] , [56] , [57] are close to this chapter.

## *§1. BASIC NOTIONS*

**1⁰. DEFINITION 4.1**. We call a conneced Riemannian manifold *(M, g)* with a family of local isometries $\{s_x : x\in M\}$ a Riemannian locally regular σ-manifold (R.l.r. σ-m.), if

1) $s_x(x)=x$, 2) the tensor field $S : S_x=(s_x)_{*x}$ is smooth and invariant under any subsymmetry $s_x$, 3) there exists a connection $\tilde{\nabla}$ on *M* invariant under any $s_x$, such that

$$\tilde{\nabla}S=\tilde{\nabla}g=0.$$

As $S_x=(s_x)_{*x}$, it is evident that

$$g(SX,SY) = g(X,Y)\ ,\ X,Y\in \mathfrak{X}(M)\ .$$

Let *M* be a R.l.r. σ-m. and suppose all the subsymmetries $s_x$ are determined globally. Then, we call $M$ a Riemannian regular $\sigma$-manifold (R.r. $\sigma$-m.).

Let the closure $G = CL(\{s_x\})$ of the group generated by the set $\{s_x : x \in M\}$ in the full isometry group $I(M)$ of a R.r. $\sigma$-m. $M$ be a transitive Lie group of transformations. Then $M$ is a Riemannian homogeneous space with the canonical connection $\tilde{\nabla}$ ($M$ is reductive, [47]).

$S$ is $G$-invariant, $S$ is invariant under every $s_x$, and it follows that $\tilde{\nabla} S = \tilde{\nabla} g = 0$.

The following example shows that the axiom 3) of *Definition 1.1* is significant.

**EXAMPLE**. Let $M = R^2 \times R$ be three - dimensional Euclidian space with the standard flat Riemannian metric $g$. For every point $x\ (x_1, x_2, t)$ we define $s_x$ as the rotation around the axis $x_1 = c_1,\ x_2 = c_2$ on the angle $t$ of any plane $R^2 : t = const$, hence $s_x(x) = x$. It is clear that $s_x$ is an isometry for each $x \in M$ and $S_x = (s_x)_{*x}$ has the following matrix on $T_x(M) = R^3$

$$\begin{bmatrix} cost & -sint & 0 \\ sint & cost & 0 \\ 0 & 0 & 1 \end{bmatrix}$$

The affinor $S$ is smooth on $M$. One can verify with a help of compasses that

$$s_x \cdot s_y = s_w \cdot s_x, \quad w = s_x(y), \quad x,y,w \in M$$

Differentiating this equality at the point $y$ we obtain

$$(s_x \cdot s_y)_{*y} = (s_x)_{*y} \cdot (s_y)_{*y} = (s_x)_{*y} \cdot S_y,$$

$$(s_w \cdot s_x)_{*y} = S_w \cdot (s_x)_{*y},$$

therefore $S$ is invariant under any $s_x$. It is clear that $S$ is not $O$-deformable and a canonical connection $\tilde{\nabla}$ does not exist.

**$2^0$**. The condition $\tilde{\nabla} S = O$ on a R.l.r. $\sigma$-m. $M$ implies that $S$ has on $M$ a constant Jordan normal form. An almost product structure can be defined on $M$: $T(M) = T^1(M) \oplus T^2(M)$, where $T^1$ is a distribution corresponding to the eigenvalue 1 , $T^2 = T^{1\perp}$

In the case when $T^1 = \{O\}$, $M$ is a Riemannian locally regular $s$-manifold, [48]. Further on, we assume $T^1 \neq \{O\}$.

**THEOREM 4.1** . Let $M$ be a R.l.r. $\sigma$-m. Then the distribution $T^1$ is integrable and its maximal integral manifolds are totally geodesic submanifolds with respect to the Riemannian connection $\nabla$.

**Proof**. From the fact that connections $\nabla, \tilde{\nabla}$ are invariant it follows that the tensor field $h = \nabla - \tilde{\nabla}$ is also invariant under every $s_x$. Since h is invariant and $S_x = (s_x)_{*x}$, it follows that

$$h_{SX}SY = Sh_X Y,\ X,Y \in \mathcal{X}(M)$$

Let $X,Y \in T^1$, then $Sh_X Y = h_{SX}SY = h_X Y$ and $h_X Y = \nabla_X Y - \tilde{\nabla}_X Y \in T^1$. Since $\tilde{\nabla}S = 0$, $T^1$ is invariant under $\tilde{\nabla}$ and we get

$$\tilde{\nabla}_X Y \in T^1,\ \nabla_X Y = \tilde{\nabla}_X Y + h_X Y \in T^1,\ [X,Y] = \nabla_X Y - \nabla_Y X \in T^1.$$

$T^1$ is autoparallel under $\nabla$ and it follows that its maximal integral submanifolds are totally geodesic.

**QED.**

The distribution $T^1$ defines the foliation $\tilde{\Lambda} = \{\Lambda_x : x \in M\}$. The fibres of $\tilde{\Lambda}$ will be called the mirrors.

The canonical connection is unique for any Riemannian locally regular $s$-manifold, [48]. For R.l.r. $\sigma$-m. we have

**PROPOSITION 4.2**. Let $\tilde{\nabla}, \tilde{\nabla}'$ be canonical connections from *Definition 1.1* and $X \in T^1$ . Then $\tilde{\nabla}_X = \tilde{\nabla}'_X$ on $M$.

**Proof**. $S$ has no fixed vectors exept the null vector in $T^2$, hence $(I - S)$ is an isomorphism on $T^2$ and for $X \in T^2$, $X \neq O$, $(I - S)X \neq O$. Let $X \in T^2$ , $Y \in \mathcal{X}(M)$ and let $\tilde{\nabla}, \tilde{\nabla}'$ be canonical connections from *Definition 1.1*, $E = \tilde{\nabla} - \tilde{\nabla}'$. Then for $X = (I - S)X_1$, $Y = SY_1$ we have

$$E_X Y = E_{(I-S)X_1}SY_1 = E_{X_1}SY_1 - E_{SX_1}SY_1 - SE_{X_1}Y_1 - SE_{X_1}Y_1 = 0,$$

therefore $\tilde{\nabla}_X = \tilde{\nabla}'_X$ $(SE_{X_1}Y_1 = E_{X_1}SY_1$ because $\tilde{\nabla}(S) = \tilde{\nabla}'(S) = 0$,

$SE_{X_1}Y_1 = E_{SX_1}SY_1$ because $E$ is invariant under every $s_x$ ).

**QED.**

It is desirable to get an explicit form of a canonical connection.

**THEOREM 4.3**. Let $M$ be a R.l.r. σ-m., $\nabla$ the Riemannian connection, and $\tilde{\nabla}'$ canonical connection, $\pi_1, \pi_2$ the projections on $T^1$, $T^2$ respectively. Then some a new canonical connection

$$(4.1)\ \tilde{\nabla}_X Y = \nabla_X Y - (\nabla_{(I-S)^{-1}\pi_2 X} S)(S^{-1}Y),\ X,Y \in \mathcal{X}(M).$$

is determined on $M$.

**Proof**. As $(I - S)$ is a nonsingular on $T^2$, it is obvious that $\tilde{\nabla}$ is a connection. Let $h_X Y = (\nabla_{(I-S)^{-1}\pi_2 X} S)S^{-1}Y$. Since $\nabla$ and $S$ are invariant under every $s_x$, it follows from (4.1) that $\tilde{\nabla}$ is also invariant under every $s_x$. So, both the tensor fields $h = \nabla - \tilde{\nabla}$ and $h' = \nabla - \tilde{\nabla}'$ are invariant under every $s_x$, hence $Sh_X Y = h_{SX} SY$ and $Sh'_X Y = h'_{SX} SY$ for $X,Y \in \mathcal{X}(M)$. Using that $S \cdot \pi_1 = \pi_1$ we obtain

$$\begin{aligned}(\tilde{\nabla}_X S)Y &= (\nabla_X S)Y - (h_X SY - Sh_X Y) = (\nabla_X S)Y - (h_X SY - h_{SX} SY) = (\nabla_X S)Y - h_{(I-S)X} SY \\ &= (\nabla_X S)Y - (\nabla_{\pi_2 X} S)Y = (\nabla_{\pi_1 X} S)Y = \nabla_{\pi_1 X} SY - S\nabla_{\pi_1 X} Y - (\tilde{\nabla}'_{\pi_1 X} SY - S\tilde{\nabla}'_{\pi_1 X} Y) \\ &= h'_{\pi_1 X} SY - Sh'_{\pi_1 X} Y = h'_{\pi_1 X} SY - h'_{S\pi_1 X} SY = 0, \quad X,Y \in \mathcal{X}(M)\end{aligned}$$

Further on, we have

$$h_X Y = (\nabla_{(I-S)^{-1}\pi_2 X} S)(S^{-1}Y) = h'_{(I-S)^{-1}\pi_2 X} S \cdot S^{-1}Y - Sh'_{(I-S)^{-1}\pi_2 X} S^{-1}Y$$

$$= h'_{(I-S)^{-1}\pi_2 X} Y - h'_{S(I-S)^{-1}\pi_2 X} Y = h'_{(I-S)(I-S)^{-1}\pi_2 X} Y = h'_{\pi_2 X} Y$$

For $X \in T^1$ $h_X Y = O$ and $\tilde{\nabla}_X = \nabla_X$, for $X \in T^2$ $h_X Y = h'_X Y$ and $\tilde{\nabla}_X = \tilde{\nabla}'_X$. The identity $\nabla g = \tilde{\nabla}' g = 0$ implies $\tilde{\nabla} g = 0$.

**QED.**

**REMARK.** From *Theorem 4.3* it follows that the connection $\tilde{\nabla}$ defined by (4.1) is necessary canonical, i.e., it has to be realized

(4.2) $\tilde{\nabla} S = 0, \quad Sh_X Y = h_{SX} SY, \quad X, Y \in \mathfrak{X}(M)$

on R.l.r. σ-m. Therefore, if the conditions (4.2) are not fulfilled on a Riemannian manifold *(M, g)* with an affinor *S*, then the affinor *S* is not defined by some structure of a R.l.r. σ-m.

*Definition 4.1.* can be rewritten in the following form

**DEFINITION 4.2.** We call a connected Riemannian manifold *(M, g)* with a family of local isometries $\{s_x : x \in M\}$ a Riemannian locally regular σ-manifold, if 1) $s_x(x) = x$, 2) the tensor field $S : S_x = (s_x)_{*x}$ is smooth, *O*-deformable and invariant under any subsymmetry $s_x$, 3) $\tilde{\nabla} S = \tilde{\nabla} g = 0$, where $\tilde{\nabla}$ is the canonical connection defined by (4.1).

**3⁰. DEFINITION 4.3.** We call a connected Riemannian manifold *(M, g)* with a family of local isometries $\{s_x : x \in M\}$ a Riemannian locally regular σ-manifold of order *k* (R.l.r. σ-m.o.k), if 1) $s_x(x) = x$, 2) the tensor field *S* determined by the formula $S_x = (s_x)_{*x}$ is smooth, invariant under any $s_x$ and satisfies the condition $S^k = I$.

Let *M* be a R.l.r. σ-m.o.k and suppose all the symmetries are determined globally. Then we call *M* a Riemannian regular σ-manifold of order *k* (R.r. σ-m.o.k).

The following theorem shows that any R.l.r. σ-m.o.k is a R.l.r. σ-m.

**THEOREM 4.4.** Let *M* be a R.l.r. σ-m.o.k, $S^k = I$, and $\nabla$ the Riemannian connection of *g*. Then the connection

(4.3) $$\tilde{\nabla}_X Y = \nabla_X Y - \frac{1}{k}\sum_{j=1}^{k-1} \nabla_X (S^j) S^{k-j} Y = \frac{1}{k}\sum_{j=0}^{k-1} S^j \nabla_X S^{k-j} Y, \quad X, Y \in \mathfrak{X}(M)$$

is determined on *M*, $\tilde{\nabla} S = \tilde{\nabla} g = 0$, and $\tilde{\nabla}$ is invariant under every $s_x$.

**Proof.** $\tilde{\nabla}$ is obviously a connection. Further, we have

$$(\tilde{\nabla}_X S)Y = \frac{1}{k}\sum_{j=0}^{k-1}(S^j \nabla_X S^{k-j+1} Y - S^{j+1} \nabla_X S^{k-j} Y) = \frac{1}{k}(\nabla_X S^{k+1} Y - S^k \nabla_X SY) = 0$$

$$g(\tilde{\nabla}_X Y, Z) + g(Y, \tilde{\nabla}_X Z) = \frac{1}{k}\sum_{j=0}^{k-1}[\,g(S^j \nabla_X S^{k-j} Y, Z) + g(Y, S^j \nabla_X S^{k-j})Z\,]$$

$$= \frac{1}{k}\sum_{j=0}^{k-1}[\,g(\nabla_X S^{k-j} Y, S^{k-j} Z) + g(S^{k-j} Y, \nabla_X S^{k-j} Z)]$$

$$= \frac{1}{k}\sum_{j=0}^{k-1} Xg(S^{k-j} Y, S^{k-j} Z) = Xg(Y, Z),\ X, Y \in \mathfrak{X}(M)$$

that is, $\tilde{\nabla} g = 0$. As $\nabla$ and $S$ are invariant under every $s_x$, it follows from (4.3) that $\tilde{\nabla}$ is also invariant under every $s_x$.

**QED.**

**PROPOSITION 4.5.** Let $M$ be a R.l.r. σ-m.o.k. Then the canonical connections defined by (4.1), (4.3) coincide.

**Proof.** For $X \in T^2$ the coincidence follows from *Proposition 4.2*. If $X \in T^1$ and $\tilde{\nabla}$ is defined by (4.1), then we see from (4.1) that $\tilde{\nabla}_X = \nabla_X$ and $\nabla_X S = \tilde{\nabla}_X S = 0$. Finally, from (4.3), *Proposition 4.5* follows.

**QED.**

## *§2. RIEMANNIAN REGULAR σ - MANIFOLDS*

**1**$^0$. In this paragraph, we consider a R.r. σ-m. and its foliation of mirrors.

**LEMMA 4.6** [46]. Let $\varphi$ and $\psi$ be isometries on *(M, g)*, $\varphi(x) = \psi(x)$, $\varphi_{*x} = \psi_{*x}$ for some $x \in M$. Then $\varphi = \psi$ on $M$.

**LEMMA 4.7.** All the subsymmetries $s_x$ are affine transformations with respect to $\tilde{\nabla}$.

**Proof** obviously follows from *Definition 4.1*.

**PROPOSITION 4.8.** Let $M$ be a R.r. σ-m. and $s_x$ a subsymmetry on $M$. Then we have $s_{x|\Lambda_x} = id_{|\Lambda_x}$ and if $x_1 \in \Lambda_x$, then $s_x = s_{x_1}$ on $M$.

**Proof.** Since $s_x$ and $S$ commute, $T^1$ and $\tilde{\Lambda}$ are invariant under $s_x$ and it follows that $s_x(\Lambda_x) = \Lambda_x$. For the restriction $s_{x|\Lambda_x}$ we have $s_x(x) = x$, $s_{x*x} = I$. According to *Lemma 4.6*, $s_x = id$ on $\Lambda_x$. Let $x_1 \in \Lambda_x$, then $s_{x_1|\Lambda_x} = id$ and $s_{x_1}(x) = s_x(x) = x$. Consider $v \in T_x(M)$ and a curve $\tau_t$ connecting $x$ and $x_1$. Denote the parallel translation with respect to the connection $\tilde{\nabla}$ by $\tilde{\tau}_t$. According to *Lemma 4.7*, all the subsymmetries commute with the parallel translation; the parallel translation commutes with *S*, because $\tilde{\nabla} S = 0$. Thus,

$$\tilde{\tau}_t((s_{x_1})_{*x}(v)) = (s_{x_1})_{*x_1}(\tilde{\tau}_t(v)) = S\tilde{\tau}_t(v) = \tilde{\tau}_t(Sv)$$

and we get $(s_{x_1})_{*x} = (s_x)_{*x} = S$. According to *Lemma 4.6*, $s_{x_1} = s_x$ on *M*.

**QED.**

**THEOREM 4.9**. Let $M$ be a R.r. σ-m., $N = \{\Lambda_x : x \in M\}$, $\pi : M \to N : x \mapsto \Lambda_x$. Then *N* is a smooth manifold and $\pi$ is a differentiable submersion.

**Proof.** According to [61], it is sufficient to show that the foliation is regular. Let *U(x)* be a convex neighbourhood of *x* in which there exists a foliated chart of the foliation $\tilde{\Lambda}$, [66] , and let $x_1 \in U(x)$. Suppose that $\bar{\Lambda}_{x_1}$, $\bar{\Lambda}_{x_2}$ are connected components of $\Lambda_{x_1} \cap U(x)$ which do not coincide ($x_2 \in U(x)$). Then there exists a unique minimizing geodesic $\gamma(t)$ in *U(x)*, where $t \in [t_1, t_2]$, $\gamma(t_1) = x_1$, $\gamma(t_2) = x_2$. The isometry $s_x$ transforms $\gamma$ into a geodesic $\gamma' \subset U(x)$ and $\gamma'$ is a minimizing geodesic too, [46]. *Proposition 4.8* yields that $s_{x_1}(\Lambda_{x_1}) = \Lambda_{x_1}$ and $s_{x_1}(x_1) = x_1$, $s_{x_1}(x_2) = x_2$. Since the minimizing geodesic which connects $x_1$ and $x_2$ is unique we have $\gamma' = \gamma$. Thus $s_{x_1}(\gamma) = \gamma$ and $(s_{x_1})_{*x_1}(\dot{\gamma}) = S_{x_1}(\dot{\gamma}) = \dot{\gamma}$, hence $\dot{\gamma}_{x_1} \in T^1_{x_1}$. According to *Theorem 4.1*, $\Lambda_{x_1}$ is a totally geodesic submanifold of *M*, so $\gamma \subset \Lambda_{x_1}$. Because $\bar{\Lambda}_{x_1}, \bar{\Lambda}_{x_2}$ are arcwise connected in *U(x)*, they coincide. The contradiction obtained proves the theorem.

**QED.**

**2⁰**. We want to show that any R.r. σ-m. can be considered as a fibre bundle. Let *I(M)* be the full isometry group of a R.r. σ-m. equipped with the compact open topology and let $G = CL(\{s_x\})$ be the closure in *I(M)* of the group generated by the set $\{s_x : x \in M\}$. Then *G* is a Lie group of transformations.

**LEMMA 4.1O.** The foliation $\tilde{\Lambda}$ is invariant under all the transformations of the group *G*, i.e., *G* transforms mirrors into mirrors.

**Proof.** Consider a sequance $\{ a_n \} \to a \in G$, where $a_n \in G$. As *S* is invariant under subsymmetries, *S* is also invariant under each $a_n$. But then $a_* \cdot S = S \cdot a_*$. Since the tensor field *S* is invariant under the group *G*, $T^1$ is also invariant under *G*. It follows that *G* transforms mirrors of the foliation $\tilde{\Lambda}$ into mirrors.

**QED.**

**LEMMA 4.11** [48]. If $G \subset I(M)$ is a closed subgroup then all *G*-orbits are closed in *M*.

Let us define the action of the group *G* on the manifold *N*: $G \times N \to N : (a,y) \mapsto \pi(a \cdot x)$, where $y=\pi(x)$. From *Lemma 4.1O* we see that this definition is correct. The action is obviously differentiable.

**THEOREM 4.12.** Let *M* be a R.r. σ-m. and *N* the corresponding manifold of mirrors. Then the group *G* is a transitive Lie group of transformations of the manifold *N*.

**Proof.** Let $x_0 \in M$ and $U(x_0)$ be a convex neighbourhood of $x_0$ with respect to ∇, which is a foliated chart of the foliation $\tilde{\Lambda}$. Suppose that *x* is an arbitrary point in $U(x_0)$, $x \notin \Lambda_{x_0}$, *r* is a distance from $x_0$ to the *G*-orbit $G(x)$ of the point *x*: $r = \inf_{a \in G} d(x_0, a(x))$. Since $G(x)$ is closed, one can find $z \in G(x)$ such that $r = d(x_0, z)$. Let us suppose that $z \notin \Lambda_{x_0}$. Then there exists a geodesic segment of the length *r* joining $x_0$ and *z*. Let *w* be a point of this segment between $x_0$ and *z*. Then $\dot{\gamma}_w \notin T^1$ because otherwise, according to *Theorem 4.1*, the whole segment would lie in $\Lambda_w$ and $z \in \Lambda_w = \Lambda_{x_0}$. Thus, $s_w(z) \neq z$, $s_w(z) \in G(x)$. Hence, all the points $x_0$, *z*, *w*, $s_w(z)$ lie in $U(x)$. Using the triangle inequality we get

$$d(x_0, s_w(z)) < d(x_0, w) + d(w, s_w(z)) = d(x_0, w) + d(s_w(w), s_w(z)) = d(x_0, w) + d(w,z)$$

$$= d(x_0, z) = r.$$

The contradiction obtained shows that $z \in \Lambda_{x_0}$. Thus, for any mirror $y = \Lambda_x$,

$y \in \pi\,(U(x_0))$, one can find an element of the group $G$ transforming $y$ into $y_0 = \Lambda_{x_0}$ , and for any $y_1$, $y_2 \in \pi\,(U(x_0))$ there exists such a transformation $a \in G$ that $y_2 = a(y_1)$.

Covering a segment of the curve between two arbitrary points of $N$ by a finite number of neighbourhoods like $\pi\,(U(x_0))$ we conclude that the group $G$ is a transitive Lie group of transformations of $N$.

**QED.**

**COROLLARY 4.13.** All the fibres of the foliation $\tilde{\Lambda}$ are diffeomorphic to the standard fibre $\Lambda = \Lambda_p$, where $p \in M$ is a fixed point.

It is well known that the component of identity of a Lie group acting transitively on a manifold $N$ is also transitive on $N$, so, later on, we shall assume the group $G$ to be connected.

**COROLLARY 4.14.** Let $p \in M$ and let $H$ be the isotropy subgroup of $\Lambda_p \in N$. The mapping $G/H \to N : aH \mapsto \Lambda_{a(p)}$ is a diffeomorphism of the manifold $G/H$ and $N$.

Let $G(G/H, H)$ be the principal fibre bundle with the base $G/H$ and the structure group $H$. Since $H$ acts on the manifold $\Lambda = \Lambda_p$ to the left, it is possible to consider $G \times_H \Lambda$, which is the fibre bundle over the base space $G/H$ with the standard fibre $\Lambda$ and the structure group $H$ associated with the principal fibre bundle.

Let $a \otimes x$ be the equivalence class containing $(a, x)$, where $(ab, x) \sim (a, bx)$, $b \in H$.

**THEOREM 4.15.** Let $M$ be a R.r. σ-m. The mappings $\Phi : G \times_H \Lambda \to M : a \otimes x \mapsto a(x)$ and $G/H \to N : aH \mapsto \Lambda_{a(p)}$ are diffeomorphisms. The following diagram is commutative :

$$\begin{array}{ccc} G \times_H \Lambda & \to & M \\ \downarrow & & \downarrow \\ G/H & \to & N \end{array} \qquad (4.4)$$

**Proof.** $\Phi$ is obviously a correctly defined differentiable mapping, $\Phi$ is surjective because $G$ is transitive on $N$. Let us check the injectivity of $\Phi$. Let

$a_1(x_1)=a_2(x_2)$, then

$$a_1^{-1}a_2 = b \in H \text{ and } a_1 \otimes x_1 = a_1 b \otimes b^{-1}x_1 = a_2 \otimes x_2 .$$

The mapping $G \times \Lambda \to M : (a,x) \mapsto a(x)$ is a submersion and the following diagram is commutative :

$$\begin{array}{ccc} G\times\Lambda & \longrightarrow & M \\ & \searrow \quad \swarrow & \\ & G\times_H \Lambda & \end{array}$$

Thus $\Phi$ is a diffeomorphism and *diagram (4.4)* is evidently commutative.

**QED.**

**3⁰**. In this paragraph we consider a manifold of mirrors as a regular *s*-manifold. Let $p\in M$ be again a fixed point, $y_p=\Lambda_p\in N$. According to *Proposition 4.8* every subsymmetry $s_x$ defines a diffeomorphism $s_y$ of the manifold $N$, where $y=\pi(x)$. It is clear that $s_y(y)= y$ and $s_{y*y} = \overline{S}$, where the Jordan normal form $\overline{S}$ coincides with the normal form of the tensor field $S$ restricted to $T^2$. It is also evident that $\overline{S}$ is invariant under the group $G$ acting transitively on $N$.

**LEMMA 4.16.** Let $a(\Lambda_p)=\Lambda_x$, where $x= a(p)$. Then $s_x= a\cdot s_p\cdot a^{-1}$ on $M$, $a\in G$.

**Proof.** $s_x(x)=x$ and $(a\cdot s_p\cdot a^{-1})(x)= x$. Then $(s_x)_{*x} = S_x$ and $(a\cdot s_p \cdot a^{-1})_{*x} = a_{*p}\cdot(s_p)_{*p}\cdot a_{*x}^{-1} = a_{*p}\cdot S_p \cdot a_{*x}^{-1} = S_x$, because $S$ is $G$-invariant. According to *Lemma 4.6*, $s_x$ coincides with $a\cdot s_p\cdot a^{-1}$ on $M$.

**QED.**

**PROPOSITION 4.17.** Let $M$ be a R.r. σ-m. and let $N$ be a manifold of its mirrors. Then $\mu : M\times N \to N : (y_1,y_2) \mapsto s_{y_1}(y_2)$ is a real analytic mapping.

**Proof.** $N\cong G/H$ has a structure of a real analytic manifold such that the action of $G$ on $N$ and the projection $\overline{p} : G \to G/H$ are analytic, [46]. One can find a neighbourhood $W\subset N$ of a point $y_0$ for which there exists an analytic section $v: W \to G$ of the fibre bundle $\overline{p}: G \to G/H$. According to *Lemma 4.16* $s_y = \pi(s_x) = \pi(a \cdot s_p \cdot a^{-1}) = a \cdot s_{y_p} \cdot a^{-1}$. Therefore, for any $y \in W$, $s_y = v(y) \cdot s_{y_p} \cdot (v(y))^{-1}$, where $s_{y_p} \in G$ is analytic. Thus, the mapping $(y_1, y_2) \mapsto s_{y_1}(y_2)$ is analytic on $W \times N$ and, in fact, on $N \times N$.

**QED.**

**DEFINITION 4.4** [48]. A regular *s*-manifold is a manifold *N* with a multiplication $\mu: N \times N \to N$ such that the mappings $s_y: N \to N$, $y \in N$ given by $s_y(z) = \mu(y, z)$ satisfy the following axioms:

1) $s_y(y) = y$,

2) each $s_y$ is a diffeomorphism,

3) $s_y \cdot s_z = s_w \cdot s_y$, where $w = s_y(z)$,

4) for each $y \in N$, $s_{y*y}: T_y(N) \to T_y(N)$ has no fixed vectors exept the null vector.

**THEOREM 4.18.** Let *M* be a R.r. σ-m. and *N* its manifold of mirrors. Then *N* is a regular *s*-manifold.

**Proof.** According to *Proposition 4.17* $\mu$ is differentiable, the axioms 1) and 2) are evident, 4) folllows from the fact that $S_{|T^2}$ has no fixed vectors except the null one. We consider the axiom 3). Let $x, u, v \in M$, $\pi(x) = y$, $\pi(u) = z$, $\pi(v) = w$, $v = s_x(u)$. Let us prove that $s_x \cdot s_u = s_v \cdot s_x$. We have

$$(s_x \cdot s_u)(u) = (s_v \cdot s_x)(u) = v,$$

$$(s_x \cdot s_u)_{*u} = (s_x)_{*u} \cdot (s_u)_{*u} = (s_x)_{*u} \cdot S_u = S_v \cdot (s_x)_{*u} = (s_v)_{*v} \cdot (s_x)_{*u} = (s_v \cdot s_x)_{*u}.$$

According to *Lemma 4.6* we get $s_x \cdot s_u = s_v \cdot s_x$. Projecting this equality onto *N* we obtain that $s_y \cdot s_z = s_w \cdot s_y$, where $w = s_y(z)$.

**QED.**

**THEOREM 4.19.** Let a R.r. σ-m. *M* be compact. Then its manifold of mirrors *N* is a Riemannian regular *s*-manifold.

**Proof.** Since the group *I(M)* of all the isometries of *M* is compact, the group

*G* is also compact. Assume < , >* is an arbitrary Riemannian metric on *N*, $X,Y\in T_y(N)$. The elements of the group *G* are isometries with respect to the following metric < , > on *N*:

$$<X,Y>=\int_{a\in G}<a_*X,a_*Y>^*$$

The rest follows from *Theorem 4.18*.

**QED.**

**REMARK**. If *H* is not compact then *G/H* can not be a Riemannian regular *s*-manifold because according to [47], the isotropy subgroup of a homogeneous Riemannian space must be compact.

## *§3. LIE ALGEBRA OF INFINITESIMAL AUTOMORPHISMS OF RIEMANNIAN REGULAR σ-MANIFOLDS*

**1⁰**. Let $(M,\mu)$ be a manifolld with a multiplication (no other properties exept the differentiability of the mapping $\mu : M \times M \to M$ are required). Let $x, y \in M$, $X \in T_x(M)$, $Y \in T_y(M)$. The point $\mu(x, y)$ is denoted by $x \cdot y$ as usual. The products $x \cdot Y$ and $X \cdot y$ in $T_{x \cdot y}(M)$ are defined by the following formulas

(4.5) $x \cdot Y = d/dt|_0(x \cdot \alpha(t)), \quad X \cdot y = d/dt|_0(\beta(t) \cdot y),$

where $\alpha(t)$, $\beta(t)$ are parametrized curves in $M$ such that $Y = d\alpha(0)/dt$, $X = d\beta(0)/dt$.

**LEMMA 4.2O.** Let $\alpha, \beta : (-l, l) \to M$ be two parametrized curves in $M$, and let $\alpha \cdot \beta$ denote the curve given by $(\alpha \cdot \beta)(t) = \alpha(t) \cdot \beta(t)$ $t \in (-l, l)$. Then

$d(\alpha \cdot \beta)(o)/dt = d\alpha(o)/dt \cdot \beta(o) + \alpha(o) \cdot d\beta(o)/dt$

**Proof** evidently follows from the "Leibniz formula", [46], p.1O.

**DEFINITION 4.5.** An automorphism of $(M, \mu)$ is defined as a diffeomorphism $\varphi : M \to M$ such that $\varphi(x \cdot y) = \varphi(x) \cdot \varphi(y)$ for every $x,y \in M$.

We shall often denote the tangent mapping of $\varphi$ by the same symbol $\varphi$.

**LEMMA 4.21.** Let $\varphi$ be an automorphism of $(M, \mu)$ and $x \in M$, $X \in T(M)$. Then

(4.6) $\varphi(x \cdot X) = \varphi(x) \cdot \varphi(X), \quad \varphi(X \cdot x) = \varphi(X) \cdot \varphi(x).$

**Proof** is obvious.

We consider the elements of $\mathfrak{X}(M)$ as the cross-sections $M \to T(M)$.

**DEFINITION 4.6** [58], p.51. Let $(M, \mu)$ be a manifold with multipliplication. A derivation of $(M, \mu)$ is a vector field $X \in \mathfrak{X}(M)$ such that

$X(p \cdot q) = X(p) \cdot q + p \cdot X(q)$ for all $p, q \in M$.

The set of all the derivations of $(M, \mu)$ will be denoted by $Der(M, \mu)$.

**PROPOSITION 4.22** [48], p.49. a) The derivations of $(M, \mu)$ form a Lie subalgebra of the Lie algebra $\boldsymbol{X}(M)$.

b) A one-parameter group of transformations of $M$ is a group of automorphisms if and only if the corresponding vector field is a derivation.

**$2^0$**. Let $(M, \{s_x\})$ be a R.r. σ-m. One can define a multiplication on $M$ by the formula

$$\mu : M\times M \to M:(x,y)\mapsto x\cdot y=\mu (x,y)=s_x(y).$$

**LEMMA 4.23.** Subsymmetries $s_x$ and $s_x^{-1}$ are automorphisms of $(M, \mu)$.

**Proof.** From the regularity condition (see proof of *Theorem 4.18*) we obtain

$$s_x(y\cdot z)=s_x(s_y(z))=s_{s_x(y)}(s_x(z))=s_x(y)\cdot s_x(z),$$

$$s_x^{-1}(y\cdot z)=s_x^{-1}(s_x(y')\cdot s_x(z'))=s_x^{-1}(s_x(y'\cdot z'))=y'\cdot z'=s_x^{-1}(y)\cdot s_x^{-1}(z).$$

**QED.**

We shall write simply $x^{-1}\cdot y$ instead of $s_x^{-1}(y)$ for $x, y\in M$.

**LEMMA 4.24.** Let $p\in M$, $X\in T_p(M)$, $x,y\in M$; then

(4.7) $X = p\cdot X + X\cdot p,$

(4.8) $X\cdot (x\cdot y) = (X\cdot x)\cdot (p\cdot y) + (p\cdot x)\cdot (X\cdot y).$

The proof is analogous to that considered in [48], p.49.

It is well-known that a diffeomorphism $\varphi$ of $M$ induces an automorphism $\tilde{\varphi}$ of the algebra $\boldsymbol{T}(M)$ of all the global tensor fields on $M$, [46], p.28.

**LEMMA 4.25.** If $\varphi$ is an isometry and automorphism of $(M, \mu)$, then $\tilde{\varphi}$ preserves the tensor fields $S$, $S^{-1}$, $I - S$, $(I-S)^{-1}\circ\pi_2$, where $\pi_1,\pi_2$ are projections on distributions $T^1$, $T^2$.

**Proof.** Let $p \in M$, $X \in T_p(M)$. Then $\varphi(S_p X)=\varphi(p \cdot X)=\varphi(p) \cdot \varphi(X)=S_{\varphi(p)} \cdot \varphi(X)$, i.e., $\varphi \cdot S_p = S_{\varphi(p)} \cdot \varphi$ on $T_p(M)$ and $\tilde{\varphi}$ preserves the tensor field $S$, therefore $S^{-1}$, $I - S$ are also invariant under $\tilde{\varphi}$. Further, $\varphi(\pi_1 X)= \varphi(S_p \pi_1 X)=S_{\varphi(p)} \varphi(\pi_1 X)$ and $\varphi(\pi_1 X) \in T^1$ i.e. $\varphi(T^1)=T^1$. Since $\varphi$ is an isometry and $T^2=T^{1\perp}$, then $\varphi(T^2)=T^2$ and $\varphi \cdot \pi_2 = \pi_2 \cdot \varphi$. The rest is obvious.

**QED.**

**PROPOSITION 4.26.** Let $X$ be an arbitrary tangent vector from $T_p(M)$ and a mapping $L(X):M \to T(M)$ is defined by the formula

(4.9) $L(X)(x)=(I_p - S_p)^{-1} \pi_2 X \cdot (p^{-1} \cdot x), \quad x \in M.$

Then $L(X)$ is a derivation of $(M,\mu)$ and $L(X)(p)=\pi_2 X$. The mapping $L:T(M) \to Der(M,\mu)$ is linear and injective on each vector space $T_p^2(M)$.

**Proof.** We can see easily that $L(X)(x) \in T_x(M)$ for each $x$, and thus, $L(X) \in X(M)$. Using formula (4.8) we get

$$\begin{aligned} L(X)(x \cdot y) &= (I_p - S_p)^{-1} \pi_2 X \cdot [p^{-1} \cdot (x \cdot y)] = (I_p - S_p)^{-1} \pi_2 X \cdot [(p^{-1} \cdot x) \cdot (p^{-1} \cdot y)] \\ &= [(I_p - S_p)^{-1} \pi_2 X \cdot (p^{-1} \cdot x)] \cdot [p \cdot (p^{-1} \cdot y)] + \\ &\quad [p \cdot (p^{-1} \cdot x)] \cdot [(I_p - S_p)^{-1} \pi_2 X \cdot (p^{-1} \cdot y)] = L(X)(x) \cdot y + x \cdot L(X)(y) \end{aligned}$$

and the first statement follows. Further, using formula (4.7) we obtain

$$\begin{aligned} L(X)(p) &= (I_p - S_p)^{-1} \pi_2 X \cdot (p^{-1} \cdot p) = (I_p - S_p)^{-1} \pi_2 X \cdot p \\ &= (I_p - S_p)^{-1} \pi_2 X - p \cdot (I_p - S_p)^{-1} \pi_2 X \\ &= I_p (I_p - S_p)^{-1} \pi_2 X - S_p (I - S_p)^{-1} \pi_2 X = \pi_2 X. \end{aligned}$$

The rest is evident.

**QED.**

Let $\underline{\mathbf{m}} = \{ L(X) : X \in T_p^2(M) \}$. From *Proposition 4.26* it follows that the mapping $L : T_p^2(M) \to \underline{\mathbf{m}}$ is an isomorphism of vector spaces.

**PROPOSITION 4.27.** For every isometric automorphism $\varphi$ of $(M,\mu)$ we have $\varphi_* \cdot L = L \cdot \varphi_*$ on $T(M)$, with the values in $Der(M,\mu)$.

**Proof.** Let $p,x \in M$ and $X \in T_p(M)$. Using formula (4.6) and *Lemma 4.2* we get

$$L(\varphi_* X)(x)=(I_{\varphi(p)} - S_{\varphi(p)})^{-1}(\varphi_* \pi_2 X) \cdot ([\varphi(p)]^{-1} \cdot x) = \varphi_* ((I_p - S_p)^{-1} \pi_2 X) \cdot ([\varphi(p)]^{-1} \cdot x)$$

$$=\varphi_*((I_p - S_p)^{-1}\pi_2 X)\cdot\varphi[p^{-1}\cdot\varphi^{-1}(x)]=\varphi_*\,(L(X)(\varphi^{-1}(x)))=(\varphi_* L(X))(x).$$

**QED.**

**3⁰**. We consider now the canonical connection of a R.r. σ-m. *M.*

**PROPOSITION 4.28.** The formula

(4.1O)

$$\tilde{\nabla}_X Y=\nabla_{\pi_1 X} Y+[L(X),Y](p),\quad p\in M,\ X\in T_p(M),\ Y\in\mathfrak{X}(M)$$

defines a connection $\tilde{\nabla}$ on *M*. Each isometric automorphism of $(M,\mu)$ is an affine transformation of $(M,\tilde{\nabla})$, and every infinitesimal isometry $Z\in Der(M,\mu)$ is an infinitesimal affine transformation of $(M,\tilde{\nabla})$.

**Proof.** It is clear that the introduced operation $\tilde{\nabla}$ is linear with respect to *X* and *Y*. We have

$$\tilde{\nabla}_X(f\,Y)=\nabla_{\pi_1 X}(f\,Y)+[L(X),f\,Y](p)=(\pi_1 X)(f)Y(p)+f(p)\nabla_{\pi_1 X}Y$$
$$+\,(L(X)(p)(f)Y(p)+f(p)[L(X),Y](p)$$
$$=f(p)(\nabla_{\pi_1 X}Y+[L(X),Y](p))+(\pi_1 X)(f)Y(p)+(\pi_2 X)(f)Y(p)$$
$$=f(p)\tilde{\nabla}_X Y+X(f)Y(p).$$

We have got that $\tilde{\nabla}$ is a connection on *M*.
Further, if an isometry $\varphi\in Aut(M,\mu)$, then using *Proposition 4.27* we obtain

$$\tilde{\nabla}_{\varphi_*(X)}(Y)=\nabla_{\pi_1\varphi_*(X)}\varphi_*(Y)+[L(\varphi_*(X)),\varphi_*Y](\varphi(p))$$

$$=\varphi_*\nabla_{\pi_1 X}Y+(\varphi_*[L(X),Y])(\varphi(p))=\varphi_*(\tilde{\nabla}_X Y).$$

Finally, if an infinitesimal isometry $Z\in Der(M,\mu)$, then *Z* generates a local group of local automorphisms of $(M,\mu)$ (an easy modification of *Proposition 4.22*) and thus, a local group of local affine transformations of $(M,\tilde{\nabla})$. So, *X* is an infinitesimal affine transformation of $(M,\tilde{\nabla})$.

**QED.**

**COROLLARY 4.29.** The connection $\tilde{\nabla}$ is invariant under all $s_x$, $x\in M$.

**Proof** immediately follows from *Lemma 4.23.*

**PROPOSITION 4.3O.** $\tilde{\nabla}S=0$.

**Proof.** 1) If $X\in T^1$, then $\tilde{\nabla}_X S=\nabla_X S=0$ (see proof of *Proposition 4.3*).

2) For $X\in T_p^2(M)$ and $Y\in \mathfrak{X}(M)$ we have $\tilde{\nabla}_X Y=(\mathcal{L}_{L(X)}Y(p))$ and $(\tilde{\nabla}_X S)Y=\tilde{\nabla}_X SY-S\tilde{\nabla}_X Y=[\mathcal{L}_{L(X)}SY-S(\mathcal{L}_{L(X)}Y](p)=[(\mathcal{L}_{L(X)}S)(Y)]\cdot(p)$. On the other hand, $L(X)$ is an infinitesimal automorphism and, according to the proof of *Lemma 4.25*, $S$ is invariant with respect to the local automorphisms of $(M,\mu)$. Therefore, $\mathcal{L}_{L(X)}S=0$ and the Proposition follows.

**QED.**

**PROPOSITION 4.31**. The connections defined by the formulas (4.10), (4.1) coincide, i.e., for $X\in T_p(M)$ and $Y\in \mathfrak{X}(M)$ we have

(4.11) $\tilde{\nabla}_X Y=\nabla_{\pi_1 X}Y+[L(X),Y](p)=\nabla_X Y-(\nabla_{(I-S)^{-1}\pi_2 X}S)(S^{-1}Y)$

**Proof.** We denote by $\tilde{\nabla}$ the connection (4.10) and by $\tilde{\nabla}'$ the second one. Both the connections are invariant under all $s_x$ and $\tilde{\nabla}S=\tilde{\nabla}'S=0$. Applying the proof of *Proposition 4.2* we see that $\tilde{\nabla}_X=\tilde{\nabla}'_X$ for $X\in T^2$ (metric properties are not used in the proof).

For $X\in T^1$ $\tilde{\nabla}_X=\tilde{\nabla}'_X=\nabla_X$ and the Proposition follows.

**QED.**

**PROPOSITION 4.32.** For every $X\in T_p(M)$ the derivation $L(X)$ is an infinitesimal isometry.

**Proof.** Using (4.11) and the equality $\nabla g=\tilde{\nabla}g=0$ we have

$$X<Y,Z>=<\tilde{\nabla}_X Y,Z>+<Y,\tilde{\nabla}_X Z>=<\nabla_{\pi_1 X}Y,Z>+<Y,\nabla_{\pi_1 X}Z>$$
$$+<[L(X),Y](p),Z>+<Y,[L(X),Z](p)>,\quad Y,Z\in\mathfrak{X}(M),$$

therefore

$$\pi_2 X<Y,Z>=<[L(X),Y](p),Z>+<Y,[L(X),Z](p)>$$

and

$$(L_{L(X)(p)}\, g)(Y,Z)=(L_{\pi_2 X}\, g)(Y,Z)=\pi_2 X<Y,Z>-<[L(X),Y](p),Z>$$
$$-<Y,[L(X),Z](p)>=0.$$

The rest follows from [46], p.237.

**QED.**

**4⁰**. Let $N$ be a manifold of mirrors of a R.r σ-m. $M$. It was shown in *Theorem 4.18* that $N$ is a regular $s$-manifold (maybe non-Riemannian). Using [48] we shall consider some information about affine regular $s$-manifolds. Let $q\in N$, then for $X'\in T_q(N)$ a mapping $L(X):N\to T(N)$ is defined by the formula

$$(4.12)\qquad L(X')(y)=(I_q-S_q)^{-1}X'\cdot(q^{-1}\cdot y),\quad y\in N,$$

$L(X')$ is a derivation of $(N,\mu)$ and $L(X')(q)=X'$.
The formula

$$(4.13)\quad \tilde{\nabla}'_{X'}Y'=[L(X'),Y'](q),\ q\in N,\ X'\in T_q(N),\quad Y'\in \mathfrak{X}(N),$$

defines an affine connection $\tilde{\nabla}'$ on $N$. Each automorphism of $(N,\mu)$ is an affine transformation of $(N,\tilde{\nabla}')$ and each derivation of $(N,\mu)$ is an infinitesimal affine transformation of $(N,\tilde{\nabla}')$.

The connection $\tilde{\nabla}'$ is invariant under all $s_y$, $y\in N$ and $\tilde{\nabla}'S=0$. For each regular $s$-manifold $(N,\mu)$ there exists a unique connection $\tilde{\nabla}'$ which is invariant with respect to all $s_y$ and such that $\tilde{\nabla}'S=0$. This connection coincides with one given by (4.13) and

$$(4.14)\quad \tilde{\nabla}'_{X'}Y'=\nabla'_{X'}Y'-(\nabla'_{(I-S)^{-1}X'}S)(S^{-1}Y'),\quad X',Y'\in \mathfrak{X}(N),$$

where $\nabla'$ is any $s_y$-invariant connection on $(N,\mu)$.

The canonical connection $\tilde{\nabla}'$ is complete and has parallel curvature and parallel torsion, i.e., $\tilde{\nabla}'\tilde{R}'=0$ and $\tilde{\nabla}'\tilde{T}'=0$ on $N$. $(N,\tilde{\nabla}')$ is an affine reductive space. Every derivation on $(N,\mu)$ is complete and thus, it determines a one-parameter group of automorphisms of $(N,\mu)$. Consequently, $Der(N,\mu)$ is the Lie algebra of the group $Aut(N,\mu)$ of all the automorphisms of $(N,\mu)$.

The elementary transvections of a regular $s$-manifold $(N,\{s_y\})$ are defined as the automorphisms of the form $s_y\circ s_z^{-1}$, $y,z\in N$. Let $G$ be the closure in the full isometry group $I(N)$ of the group generated by all the elementary transvections. $G$ is called the group of transvections of $(N,\{s_y\})$ and is denoted by $Tr(N,\{s_y\})$. The transvection group $Tr(N,\{s_y\})$ coincides with the transvection group $Tr(N)$ of the affine reductive space $(N,\tilde{\nabla}')$. Further, $Tr(N)$ is a normal subgroup of the group $Aut(N,\mu)$ and $Tr(N)$ is transitive on $N$.

Let $q\in N$ and $\underline{\mathbf{m}} =\{L(X'):X'\in T_q(N)\}$ , where $L(X')$ is defined by (4.12). Then the group $G = Tr(N)$ coincides with the group generated by the ideal

$$\underline{\mathbf{g}} = \underline{\mathbf{m}}+[\underline{\mathbf{m}}, \underline{\mathbf{m}}]$$

We have seen that the canonical projection $\pi : M \to N$ is a homomorphism of manifolds with multiplications $(M,\mu)$, $(N,\mu)$ and $\pi$ is an isomorphism of vector spaces $T_p^2(M)$ and $T_q(N)$.

**LEMMA 4.33.** Let $X \in T_p^2(M)$ and $X' = \pi_*(X) \in T_q(N)$. Then $\pi_*(L(X)) = L(X')$.

**Proof** is the same to that of *Proposition 4.27*, where we have to replace an automorphism $\varphi$ and the homomorphism $\pi$.

We can identify the action of the group $G \subset I(M)$ generated by the set $\{s_x: x\in M\}$ with the action of $G$ on $N$ by means of $\pi$, see $1^0$, $2^0$ §2.

So, the Lie algebra $\underline{\mathbf{g}}$ of the closed Lie group $G$ may be written in the following form

(4.15) $$\underline{\mathbf{g}} = \underline{\mathbf{m}}+[\underline{\mathbf{m}}, \underline{\mathbf{m}}],$$

where $\underline{\mathbf{m}} = \{ L(X) : X \in T_p^2(M), p \in M \}$. Thus, if $\dim T^2 = n_2$, then $\dim \underline{\mathbf{g}} \le n_2 + n_2^2$.

## *§4. ORBITS ON RIEMANNIAN REGULAR σ-MANIFOLDS AND SOME EXAMPLES*

**$1^0$**. We want to consider orbits under the action of the group $G$ on $M$. We consider a point $p\in M$ and its isotropy subgroup $H(p) \subset G$. According to *Lemma 4.10* the fibre $\Lambda_p$ is invariant with respect to $H(p)$, hence $H(p) \subset H$. If $M(p)$ is an orbit of a point $p$ under the action of a closed subgroup $G \subset I(M)$, then the homogeneous space $M(p) \cong G/H(p)$ is a closed submanifold of the manifold $M$, see *Lemma 4.11*.

Since $H(p)$ is a compact subgroup, then there exists a positive definite inner product $B$ on the Lie algebra $\underline{\mathbf{g}}$ of the isometry group $G$ which is invariant under the action of the group $ad_G(H)$. Let $\underline{\bar{\mathbf{m}}}$ be the orthogonal complement to $\underline{\bar{\mathbf{h}}}$ in $\underline{\mathbf{g}}$ with respect to $B$, where $\underline{\bar{\mathbf{h}}}$ is the Lie algebra of the group $H(p)$. Then the homogeneous space $G/H(p)$ is the reductive space with respect to decomposition $\underline{\mathbf{g}} = \underline{\bar{\mathbf{h}}} \oplus \underline{\bar{\mathbf{m}}}$, [47].

**PROPOSITION 4.34.** The orbit $M(p)$ is a R.r. $\sigma$-m. itself for any $p \in M$.

**Proof.** Since $s_x \in G$, then $s_x(M(p)) = M(p)$ for any $x \in M$ and $T(M(p))$ is invariant with respect to $S$. Hence the conditions 1),2) of *Definition 4.1* are realized. Let $\tilde{\nabla}$ be the canonical connection of the Riemannian homogeneous space $M(p) \cong G/H(p)$. Then $\tilde{\nabla} g = 0$, $\tilde{\nabla} S = 0$ because $S$ is invariant under the action of $G$ and $\tilde{\nabla}$ is also invariant with respect to $G$. So, the third axiom of *Definition 4.1* is fulfilled and $\tilde{\nabla}$ is a canonical connection of the R.r. $\sigma$-m. $M(p)$.

**QED.**

**LEMMA 4.35.** For each $x \in M(p)$, $T_x^2(M) \subset T_x(M(p))$.

**Proof.** Since $G$ is a transitive Lie group on the manifold of mirrors $N$, then $\pi(M(p)) = N$ and $\pi_*(T_p(M(p))) = T_q(N)$, where $q = \pi(p)$. Further, for each $X' \in T_q(N)$ there exists such a vector $X \in T_p(M(p))$ that $\pi_*(X) = X'$. Let $X = X_1 + X_2$, where $X_1 \in T_p^1$, $X_2 \in T_p^2$ then $\pi_*(X) = \pi_*(X_2)$. $T_p(M(p))$ and $T_p^2$ are invariant with respect to $S$, hence they are also invariant under $I - S$. So, we have

$$(I - S)X = X_1 + X_2 - X_1 - SX_2 = (I - S)X_2 \in T_p^2 \cap T_p(M(p)).$$

If we look over all the vectors $X'$ from $T_q(N)$, then we obtain all the vectors $(I - S)X_2 \in T_p^2$, i.e. $T_p^2$ itself, because $(I - S)$ is nonsingular on $T_p^2$. Thus, $T_p^2 \subset T_p(M(p))$. From the invariance of $T^2$ and $T(M(p))$ under the action of the group $G$ it follows that $T_p^2(M) \subset T_x(M(p))$ for every $x \in M(p)$.

**QED.**

**LEMMA 4.36.** The vector field $L(X)$ can be restricted to the manifold $M(p)$, i.e., $L(X)(x) \in T_x(M(p))$ for any $x \in M(p)$, $X \in T_p^2$.

**Proof.** Using (4.9) we have

$L(X)(x)=(I_p - S_p)^{-1}X\cdot(p^{-1}\cdot x)=X'\cdot y,$

where $X'=(I_p-S_p)^{-1}X\in T_p^2,\ y\in M(p)$. Let $\beta(t)$ be such a parametrized curve in $M(p)$ that $X'=d\beta(o)/dt$. Then $\beta(t)\cdot y=s_{\beta(t)}(y)\in M(p)$ and $X'\cdot y=d/dt_{|0}(\beta(t)\cdot y)\in T_x(M(p))$.

**QED.**

**PROPOSITION 4.37.** The canonical connection $\tilde{\nabla}$ of the Riemannian homogeneous space $M(p)\cong G/H(p)$ coincides with the connection $\tilde{\nabla}'$ given by (4.11).

**Proof.** From the formula (4.10) it follows that for $X\in T_p^2$, $Y\in \mathcal{X}(M(p))$ $\tilde{\nabla}'_X Y=[L(X),Y](p)\in T_p(M(p))$. According to *Proposition 4.34* $\tilde{\nabla}$ is a canonical connection of the R.r. σ-m. $M(p)$, therefore for $X\in T_p^2$ $\tilde{\nabla}_X Y=[L(X),Y](p)$ too, see *Proposition 4.2,* 4.31, and $\tilde{\nabla}_X=\tilde{\nabla}'_X$.

Let $\tilde{R},\tilde{R}'$ be corresponding qurvature tensors of $\tilde{\nabla},\tilde{\nabla}'$, and let $X$, $Y$ be linearly independent vectors in $T_p^2$. There exist such local coordinates $(x_1,x_2,\ldots,x_n)$ of $p$ that $X=\frac{\partial}{\partial x_1}_{|p}$, $Y=\frac{\partial}{\partial x_2}_{|p}$. Then for $Z\in\mathcal{X}(M(p))$ we get

$$(\tilde{R}_{XY}Z)_p=(\tilde{\nabla}_X\tilde{\nabla}_Y Z-\tilde{\nabla}_Y\tilde{\nabla}_X Z-\tilde{\nabla}_{[X,Y]}Z)_p=(\tilde{\nabla}_X\tilde{\nabla}_Y Z-\tilde{\nabla}_Y\tilde{\nabla}_X Z)_p$$

$$=(\tilde{\nabla}'_X\tilde{\nabla}'_Y Z-\tilde{\nabla}'_Y\nabla'_X Z)=(\tilde{R}'_{XY}Z)_p.$$

From the invariance under $G$ it is evident that

$\tilde{R}_{XY}Z=\tilde{R}'_{XY}Z$ for $X,Y\in T^2(M(p))$, $Z\in\mathcal{X}(M(p))$.

Furter,

$dim\ G = dim\ \underline{\mathbf{g}} = dim([\underline{\mathbf{m}},\underline{\mathbf{m}}]+\underline{\mathbf{m}}) \geq dim\ T_p(M((p))= dim\ (G/H(p))$

hence $T_p(M(p))=T_p^2(M(p))+[\underline{\mathbf{m}},\underline{\mathbf{m}}]$, where $\underline{\mathbf{m}}=\{L(X):X\in T_p^2(M)\}$. Using that for $X\in T_p^2$ $L(X)(p)=X$ we obtain

$$(\tilde{\nabla}_{[L(X),L(Y)]})_p = (\tilde{\nabla}_{L(X)}\tilde{\nabla}_{L(Y)} - \tilde{\nabla}_{L(Y)}\tilde{\nabla}_{L(X)} - \tilde{R}_{L(X)L(Y)})_p$$

$$= (\tilde{\nabla}'_{L(X)}\tilde{\nabla}'_{L(Y)} - \tilde{\nabla}'_{L(Y)}\tilde{\nabla}'_{L(X)} - \tilde{R}'_{L(X)L(Y)})_p = (\tilde{\nabla}'_{[L(X),L(Y)]})_p$$

and $\tilde{\nabla}_X = \tilde{\nabla}'_X$ for any $X \in T_p(M(p))$. From the invariance $\tilde{\nabla}, \tilde{\nabla}'$ under $G$ the Proposition follows.

**QED.**

**COROLLARY 4.38** [47]. 1) The canonical connection $\tilde{\nabla}$ defined by (4.11) is a complete connection on the Riemannian homogeneous space $M(p)$;
2) $\tilde{\nabla}g = \tilde{\nabla}S = \tilde{\nabla}h = \tilde{\nabla}\tilde{T} = \tilde{\nabla}\tilde{R} = 0$ on $M(p)$;
3) every vector field $L(X)$, $X \in T_p(M(p))$ is complete on $M(p)$, see [46], p.234.

**REMARK.** In general case a type of orbit $M(p)$ depends on $p \in \Lambda$ and can be different in various points of $M$.

**2⁰. DEFINITION 4.7.** We call $M$ a Riemannian regular $\sigma$-manifold of maximal torsion if for some $p \in M$ $M(p)=M$, i.e., $M$ is a Riemannian homogeneous space, and of minimal torsion if $T^2$ is an integrable distribution on $M$.

We recall here a definition of covering space. Given a connected locally arcwise connected topological space $N$, a connected space $M$ is called a covering space over $N$ with projection $\pi : M \to N$ if every point $x$ of $N$ has a connected open neighbourhood $U$ such that each connected component of $\pi^{-1}(U)$ is open in $M$ and is mapped homeomorphically onto $U$ by $\pi$.

**PROPOSITION 4.39.** Let $M$ be a R.r. $\sigma$-m. of minimal torsion and $N$ be a manifold of its mirrors. Then $M(p)$ is a covering space over $N$ with canonical projection $\pi : M(p) \to N$.

**Proof.** Since the distribution $T^2$ is integrable we can consider its maximal integral manifold $M^2(p)$ containing $p$. At first, we verify that $M^2(p)$ is invariant under the multiplication $\mu(x,y)=s_x(y)$, $x,y \in M^2(p)$. Let $\gamma(t)$, $t \in [0,1]$ be such a parametrized curve in $M^2(p)$ that $\gamma(0)=x$, $\gamma(1)=y$. Since for any $t \in [0,1]$ $\gamma(t) \in M^2(p)$, then $\dot{\gamma}_t \in T_{\gamma(t)}(M^2(p)) = T^2_{\gamma(t)}(M)$. The distribution $T^2$ is invariant

under $s_x$, therefore

$$X(t)=(s_x)_{*\gamma(t)}(\dot\gamma_t)\in T^2_{s_x(\gamma(t))}(M)=T_{s_x(\gamma(t))}(M^2(p))$$

where $X(t)$ is a tangent vector at the point $s_x(\gamma(t))$ to the curve $s_x(\gamma(t))$. So, the curve $s_x(\gamma(t))$ lies in $M^2(p)$ and $s_x(y)=s_x(\gamma(1))\in M^2(p)$. We have got that $M^2(p)$ is invariant under multiplication.

Further, for any $x\in M$ $s_x=s_{x_1}$, where a point $x_1$ belongs to intersection of the mirror $\Lambda_x$ and $M^2(p)$, see *Proposition 4.8*, hence $M^2(p)$ is invariant under any subsymmetry $s_x$. As the group $G$ is the closure in $I(M)$ of the group generated by the set $\{s_x: x\in M\}$, then $M^2(p)$ is also invariant with respect to the action of $G$. So, the orbit $M(p)\subset M^2(p)$ but, according to *Lemma 4.35*, for any $x\in M(p)$ $T^2_x(M)=T_x(M^2(p))\subset T_x(M(P))$, hence $T_x(M(p))=T_x(M^2(p))$ and $M(p)=M^2(p)$ because they are connected.

Since $\pi_*$ is an isomorphism of $T^2_x$ and $T_{\pi(x)}(N)$ the rest is evident.

**QED.**

**COROLLARY 4.40.** If $M$ is a R.r. σ-m. of minimal torsion, then the action of $H$ is discrete on $\Lambda$ and $L(X)(x)\in T^2_x$ for any $x\in M$.

**Proof.** The first statement is clear because $M(p)$ intersects $\Lambda_p$ in a discrete set of points of $M$. The second one follows from *Lemma 4.36*.

**QED.**

**PROPOSITION 4.41.** Let $M$ be a R.r. σ-m. of maximal torsion. Then the Lie group $H$ acts transitively on $\Lambda$ and $\Lambda\cong H/H(p)$ is a Riemannian symmetric space with canonical connection $\nabla$.

**Proof.** For any $x_1, x_2\in\Lambda_p$ there exists an element $a\in G$ that $x_2=a(x_1)$. Since the foliation $\tilde\Lambda$ is invariant under the action of the group $G$, then $a(\Lambda)=\Lambda$ and $a\in H$. The distribution $T^1$ is autoparallel with respect to $\tilde\nabla$, therefore, according to *Proposition 4.37*, $\tilde\nabla_{|\Lambda}=\nabla$ is the canonical connection of the Riemannian homogeneous space $\Lambda\cong H/H(p)$, i.e., $\Lambda$ is a Riemannian symmetric space.

**QED.**

**PROPOSITION 4.42.** Let $M$ be a R.r. σ-m. and let its mirrors be

one-dimensional. Then $M$ is a R.r. $\sigma$-m. of minimal or maximal torsion.

**Proof.** In this case $dim\, M = n_2+1$, where $n_2 = dim\, T^2$. According to *Lemma 4.35*, for each point $p\in M$ $dim\, M(p)\geq n_2$.

If for some point $p_0\in M$ $dim\, M(p_0)= n_2+1$, then $T_{p_0}(M(p_0))=T_{p_0}(M)$ and $M(p_0)=M$, therefore $M$ is a R.r. $\sigma$-m. of maximal torsion.

If $dim\, M(p)= n_2$ for every $p\in M$, then $T_x(M(p))=T_x^2(M)$, $x\in M(p)$ and the distribution $T^2$ is integrable.

**QED.**

**3⁰**. We consider the main example of a Riemannian regular $\sigma$-manifold of order $k$. Let $(N,g^2)$ be a Riemannian regular homogeneous $s$-manifold of order $k$, [48], then $N\cong G/H$, where $G_0^\sigma \subset H \subset G^\sigma$, $G^\sigma=\{a\in G : \sigma(a)= a\}$, $G_0^\sigma$ is the component of the identity of $G^\sigma$, $\sigma$ is the automorphism of the group $G$ $(\sigma^k= id)$. (Here $G$ is a connected group of isometries which acts transitively on $N$). Let $G(G/H, H)$ be a principal fibre bundle with the base $G/H$ and the structure group $H$. Let $(\Lambda, g^1)$ be a Riemannian manifold and let $H$ act on $\Lambda$ to the left. We consider the fibre bundle $G\times_H \Lambda$ which is associated with $G(G/H, H)$ and again denote by $a\otimes x$ the equivalence class containing $(a, x)$, where $(ah, x) \sim (a, hx)$, $h\in H$.

Now we state the main theorem of this section.

**THEOREM 4.43.** $M= G\times_H\Lambda$ is a R.r. $\sigma$-m.o.k.

**Proof** will be given step by step in the next paragraphs.

**LEMMA 4.44** [58] . The formulas

$$aH\cdot bH=a(a^\sigma)^{-1}b^\sigma H,\ a^\sigma=\sigma(a),\ b^\sigma=\sigma(b),\quad a,b\in G$$

define a regular multiplication on $N$.

**LEMMA 4.45.** The formula

$$(a\otimes u)\cdot(b\otimes v)=a(a^\sigma)^{-1}b^\sigma\otimes v$$

defines a regular multiplication on $M\cong G\times_H\Lambda$. The projection $\pi : G\times_H \Lambda\to G/H$ is a homomorphism of spaces with multiplications.

**Proof** is analogous to that considered in [56] for the case $\sigma^2 = id$.

We have a family of symmetries $\{s_y : y \in N\}$ on $N$, $s_y(z) = y \cdot z$, and a tensor field $\bar{S}_y = (s_y)_{*y}$ which is invariant under all $s_y$. It is clear that $\bar{S}^k = I$. The family of subsymmetries $\{s_x : x \in M\}$, $s_x(z) = x \cdot z$, and the tensor field $S_x = (s_x)_{*x}$ are defined on $M$. $S$ is invariant under all $s_x$ from the regularity condition. Since $\pi$ is a homomorphism of spaces with multiplications, we have

(4.16) $\qquad \pi \cdot s_x = s_{\pi(x)}, \quad \pi_* \cdot S = \bar{S}$

**LEMMA 4.46.** Let $\Lambda_x$ be the fibre which contains $x \in M$. Then $s_x = id$ on $\Lambda_x$ and if $x_1 \in \Lambda_x$, then $s_x = s_{x_1}$.

**Proof.** Let $x = a \otimes u$, $z = b \otimes v \in \Lambda_x$, $a,b \in G$, then $a = bh$, $h \in H$, because $\pi(x) = \pi(z)$ and we obtain

$$x \cdot z = (a \otimes u) \cdot (b \otimes v) = (b \otimes hu) \cdot (b \otimes v) = b(b^\sigma)^{-1} b^\sigma \otimes v = b \otimes v = z.$$

If $x_1 = a_1 \otimes u_1 \in \Lambda_x$, then $a_1 = ah$ because $\pi(x) = \pi(x_1)$ and $x_1 = a_1 \otimes u_1 = a \otimes hu_1$, $h \in H$, and for any $\bar{z} \in M$

$$x \cdot \bar{z} = (a \otimes hu_1) \cdot (\bar{b} \otimes \bar{v}) = a(a^\sigma)^{-1} \bar{b} \otimes \bar{v} = x \cdot \bar{z}.$$

**QED.**

The foliation $\tilde{\Lambda} = \{\Lambda_x : x \in M\}$ defines the distribution $T^1$ on $M$. According to *Lemma 4.46* $S_{|T^1} = I$ and , since $\bar{S}$ has no fixed vectors exept the null vector, the eigenspace of S corresponding to the eigenvalue $1$ coincides with $T^1_x$. Let $T^2_x$ be the direct sum of all the eigenspaces of $S_x$ except $T^1_x$. From (4.16) we get $S^k = I$ and $\pi_* : T^2_x \to T_{\pi(x)}(N)$ is an isomorphism. The structure of the almost product $T(M) = T^1 \oplus T^2$ is defined on $M$. The action of the group $G$ on the homogeneous space $N \cong G/H$ induces the action of $G$ on $M \cong G \times_H \Lambda : (a, b \otimes u) \mapsto ab \otimes u$ and we have $\pi(a \cdot x) = a \cdot \pi(x)$, $a,b \in G$, $x \in M$.

**LEMMA 4.47.** The tensor field $S$ is invariant under all the elements of $G$ on $M$.

**Proof.** We shall show that $(b\cdot s_x)(z)=(s_{b(x)}\cdot b)(z)$, $b\in G$, $x,z\in M$. Indeed, $b\cdot(x\cdot z)=ba(a^\sigma)^{-1}c^\sigma\otimes v$, $(ba\otimes u)\cdot(bc\otimes v)=(ba)(b^\sigma a^\sigma)^{-1}b^\sigma c^\sigma\otimes v=ba(a^\sigma)^{-1}c^\sigma\otimes v$, where $x=a\otimes u$, $z=c\otimes v$. Considering the tangent mappings we get $b_*\cdot S_x=S_{b(x)}\cdot b_{*x}$.

**QED.**

According to *Lemma 4.47* the distributions $T^1$, $T^2$ are invariant under $G$, hence the foliation $\tilde{\Lambda}$ is also $G$ - invariant.

We define the following Riemannian metric on the distribution $T^2$:

$$g_x^2(X,Y)=g_{\pi(x)}^2(\pi_*X,\pi_*Y),\quad X,Y\in T_x^2.$$

Then $g^2(a_*X,a_*Y)=g^2(\pi_*(a_*X),\pi_*(a_*Y))=g^2(a_*(\pi_*X),a_*(\pi_*Y))$

$$=g^2(\pi_*X,\pi_*Y)=g^2(X,Y), X,Y\in T^2,\ a\in G.$$

Thus, the elements of the group $G$ are isometries on $T^2$. Let $p\in M$ be a fixed point and $\Lambda_p=\Lambda$. We define a Riemannian metric on the distribution $T^1$ as follows:

$$g_x^1(X,Y)=g^1(a_*X,a_*Y),\ a\in G,\ a(x)\in\Lambda,\ X,Y\in T^1.$$

The element $a$ exists because $G$ is a transitive Lie group of transformations of $N$. Let $b\in G$, $b(x)\in\Lambda$, then $\Lambda$ is invariant under $h=ab^{-1}$ and $h\in H$. Since $H$ acts on $\Lambda$ as an isometry group, we get

$$g^1(b_*X,b_*Y)=g^1(h_*(b_*X),h_*(b_*Y))=g^1(a_*X,a_*Y),\quad X,Y\in T^1.$$

It follows that the metric $g^1$ is well-defined on $T^1$. It is clear that the elements of the group $G$ are isometries on $T^1$.

We define a Riemannian metric on $M$ as follows: $g_{|T^1}=g^1$, $g_{|T^2}=g^2$ and $T^1$, $T^2$ are orthogonal in the metric $g$. From the above we see that $G$ is an isometry group with respect to $g$. A transformation $s_x$ is identified with an element of $G$ and $s_x$ is an isometry , too.

Hence *Theorem 4.43* follows.

**QED.**

**REMARK.** If the action of $H$ on $\Lambda$ is trivial, then $M$ is simply a direct Riemannian product of the Riemannian manifold $\Lambda$ and the Riemannian regular homogeneous $s$-manifold of order $k$ $N$ , $M\cong G/H\times\Lambda$.

**4⁰**. We consider another example of a Riemannian regular σ-manifold.

**DEFINITION 4.8** [57]. A manifold $M$ with a differentiable multiplication is called a reflexion space, if the following axioms are satisfied:
(i) $x\cdot x= x$; (ii) $x\cdot (x\cdot y)= y$; (iii) $x\cdot (y\cdot z)=(x\cdot y)\cdot (x\cdot z)$;
where $x,y,z\in M$.

**PROPOSITION 4.48.** Let $M$ be a reflexion space and let a mapping $s_x : M \to M : y \mapsto s_x(y) = x\cdot y$ be an isometry for any $x\in M$. Then $M$ is a Riemannian regular σ-manifold of order 2.

**Proof.** From *Definition 4.8* we see that $s_x(x)=x\cdot x$, $s_x \cdot s_y = s_{s_x(y)} \cdot s_x$ and the isometry $S_x = (s_x)_{*x}$ is smooth and invariant under every $s_x$. From (ii) it follows that $s_x^2(y) = y$, hence $S^2= I$ on $M$. All the axioms of *Definition 4.3* are realized.

**QED.**

**DEFINITION 4.9** [48]. A connected Riemannian manifold $(M, g)$ with a familly of local isometries $\{s_x : x\in M\}$ is called a locally $k$-symmetric Riemannian space (k-s.l.R.s.) if the following axioms are fullfilled:
a) $s_x(x)=x$ and $x$ is the isolated fixed point of the local symmetry $s_x$;
b) the tensor field $S : S_x = (s_x)_{*x}$ is smooth and invariant under any local isometry s;
c) $S^k=I$ and $k$ is the least of such positive integers.
If all the symmetries are determined globally, then $(M, g)$ is called $k$-symmetric Riemannian space (k-s.R.s.).

Comparing with *Definition 4.3* we conclude that every k-s.l.R.s. $M$ is a R.l.r. σ-m.o.k . If $M$ is a $k$-s.l.R.s., then the canonical connection $\tilde{\nabla}$ can be defined by (4.14) or (4.3). $\tilde{\nabla}$ is unique and $\tilde{\nabla}S = \tilde{\nabla}g = 0$, see [48].

**5⁰**. Let $M$ be a 2k-s.R.s. We define the family of isometries $\{\sigma_x : x\in M\}$, where $\sigma_x =(s_x)^k$, $P_x = (\sigma_x)_{*x} = S_x^k$, $\pi_1=1/2(I+P)$, $\pi_2=1/2(I-P)$, $T^1=\pi_1(T(M))$, $T^2=\pi_2(T(M))$.

**THEOREM 4.49.** Let $M$ be a 2k-s.R.s. Then $M$ is a reflexion space and the natural multiplication is analytic.

**Proof.** The closure $\overline{G} = GL(\{s_x\})$ of the group generated by the set $\{s_x : x \in M\}$ in the full isometry group $I(M)$ is a transitive Lie group of transformations of $M$, [48] , and there exists such analytic structures on $\overline{G}$ and $M \cong \overline{G}/H_p$, where $H_p$ is the isotropy subgroup of $p \in M$, that the action of $\overline{G}$ on $M$ is analytic. There exists a neighbourhood $U \ni p$ and such an analytic section $\nu : U \to \overline{G}$ of the fibre bundle $\pi : G \to M$ that for $x \in U$ $s_x = \nu(x) \cdot s_p \cdot (\nu(x))^{-1}$. It is evident that $\sigma_x = \nu(x) \cdot \sigma_p \cdot (\nu(x))^{-1}$. $\sigma_p \in \overline{G}$ is analytic, therefore the mapping $(x,y) \mapsto \sigma_x(y)$ is also analytic on $U \times M$ and, as a result, on $M \times M$. So, we have obtained an analytic mapping

$$\mu : M \times M \to M : (x,y) \mapsto x \cdot y = \sigma_x(y).$$

Since $\sigma_x(x) = x$ and $\sigma_x^2 = id$, then $x \cdot x = x$ and $x \cdot (x \cdot y) = y$. As $s_{x*}(S) = S$, then $\sigma_{x*}(P) = P$ and it follows that $(\sigma_x)_{*y} \cdot P = P_z \cdot (\sigma_x)_{*y}$, where $z = \sigma_x(y)$. An isometry is uniquely defined by its tangent mapping, therefore $\sigma_x \cdot \sigma_y = \sigma_z \cdot \sigma_x$, $z = \sigma_x(y)$, and this identity is equivalent to the last axiom (iii) of a reflexion space.

**QED.**

Using *Proposition 4.48* and *Theorems 4.12 - 4.15* we have the following

**THEOREM 4.50.** Let $(M,g)$ be a 2k-s.R.s. Then the distribution $T^1$ is involutory and a set of its integral submanifolds $N = \{\Lambda_x : x \in M\}$ is a smooth manifold. If $G = CL(\{\sigma_x\})$ is the closure in $I(M)$ of a group generated by the set $\{\sigma_x : x \in M\}$ and $H$ is the isotropy subgroup of the fixed mirror $\Lambda_p \in N$, then the mappings

$$\Phi : G \times_H \Lambda \to M : a \otimes x \mapsto a(x), \; G/H \to N : aH \mapsto \Lambda_{a(p)}$$

are isomorphisms and the following diagram is commutative:

$$\begin{array}{ccc} G \times_H \Lambda & \to & M \\ \downarrow & & \downarrow \\ G/H & \to & N \end{array}.$$

The manifold $N$ is a symmetric space, maybe non - Riemannian.

If the distribution $T^2$ is involutary and $N$ is a simply connected manifold, then, [57], $M$ is isomorphic to $N\times\Lambda$ as reflexion spaces.

It is evident that $P=S^k$ is invariant under $\overline{G}=CL(\{s_x\})$, therefore the foliaion $\tilde{\Lambda}=\{\Lambda_x : x\in M\}$ is also invariant under $\overline{G}$, see *Lemma 4.10*.

It is defined the natural action of $\overline{G}$ on $N$ and this action is transitive because $G\subset\overline{G}$. Let $H_p$ be the isotropy subgroup of $p\in M$, then $M\cong\overline{G}/H_p$ and it is clear that $H_p\subset\overline{H}$, where $\overline{H}$ is a subgroup of $\overline{G}$ consisting of such elements that $\Lambda$ is invariant with respect to them. It is evident that $\overline{H}=CL(\{s_x\})$, where $x\in\Lambda$ and $\overline{H}$ is a transitive group of transformations of the submanifold $\Lambda=\Lambda_p$. Since $S^k=I$ on $\Lambda$, then $s_x^k=id$ on $\Lambda$, $x\in\Lambda$, and all the usual axioms are satisfied, therefore $\Lambda$ with the family $\{s_x : x\in\Lambda\}$ is a $k$-symmetric space. So, we obtain the following

**THEOREM 4.51.** Let $(M, g)$ be a 2k-s.R.s. Then such a sequence of closed Lie groups $\overline{G}\supset\overline{H}\supset H_p$ is defined that $M\cong\overline{G}/H_p$; $N\cong\overline{G}/\overline{H}$ is a symmetric space; $\Lambda\cong\overline{H}/H_p$ is a k-symmetric space.

## *§5. RIEMANNIAN LOCALLY REGULAR σ-MANIFOLDS*

**$1^0$.** We consider now a Riemannian locally regular σ-manifold $M$ with the canonical connection $\tilde{\nabla}$ defined by

$$(4.17)\quad \tilde{\nabla}_X Y=\nabla_X Y-(\nabla_{(I-S)^{-1}\pi_2 X}S)(S^{-1}Y),\ X,Y\in\mathfrak{X}(M)$$

**PROPOSITION 4.52.** Let $\{s_x : x\in M\}$ be the family of local subsymmetries of the manifold $(M,g)$. If a tensor field $P$ is invariant under any $s_x$, then $\tilde{\nabla}_X P=0$ for each $X\in T^2$. In particular $\tilde{\nabla}_X h=\tilde{\nabla}_X\tilde{T}=\tilde{\nabla}_X\tilde{R}=0$ for every $X\in T^2$, where $h=\nabla-\tilde{\nabla}$ and $\tilde{T}$, $\tilde{R}$ are the torsion and curvature tensor fields of $\tilde{\nabla}$.

**Proof.** We consider the integral curve $\gamma(t)$ of a vector field $X\in T^2$, $\gamma(0)=x_0$, and denote by $X_1, ..., X_l$, $\omega^1, ..., \omega^m$ vector and covector fields defined on some neighbourhood of $\gamma(t)$ which are parallel along $\gamma(t)$ with respect to $\tilde{\nabla}$. As $S$ is parallel along $\gamma(t)$ under $\tilde{\nabla}$, then $SX_i$, $S^*\omega^j$, $i=1, ..., l$; $j=1, ..., m$, are also parallel.

Since $P$ is invariant under any $s_x$ we have

$$P(S^*\omega^1,...,S^*\omega^m,X_1,...,X_l)=P(\omega^1,...,\omega^m,SX_1,...,SX_l),$$
$$\tilde{\nabla}_X P(S^*\omega^1,...,S^*\omega^m,X_1,...,X_l)=(\tilde{\nabla}_{SX}P)(\omega^1,...,\omega^m,SX_1,...,SX_l).$$

Taking the covariant derivative of the first equality in direction of $SX$ at the point $x_0$ we obtain

$$\tilde{\nabla}_{SX}P(S^*\omega^1,...,S^*\omega^m,X_1,...,X_l)=(\tilde{\nabla}_{SX}P)(\omega^1,...,\omega^m,SX_1,...,SX_l)$$

and subtracting two last identities we get

$$(\tilde{\nabla}_{(I-S)X}P)(S^*\omega^1,...,S^*\omega^m,X_1,...,X_l)=0.$$

Since $(I-S)$ is a nonsingular on $T^2$ this equality implies $\tilde{\nabla}_X P=0$ for any $X\in T^2$.

As any $s_x$ is a locally affine transformation of $\tilde{\nabla}$, then $\tilde{T}$, $\tilde{R}$ are invariant under $s_x$, hence $\tilde{\nabla}_X\tilde{T}=\tilde{\nabla}_X\tilde{R}=0$ for $X\in T^2$, $\tilde{\nabla}_X h=0$ because $\tilde{\nabla}$ and $\nabla$ are also invariant with respect to $s_x$.

**QED.**

**COROLLARY 4.53.** $S(\tilde{T})=\tilde{T}$, $S(\tilde{R})=\tilde{R}$ on $M$.

**Proof** is evident because $S_x=(s_x)_{*x}$.

An arrangement of a local subsymmetry $s_x$, $x\in M$ is described below.

**PROPOSITION 4.54.** Let $M$ be a R.l.r. σ-m., $x\in M$ and let $B_x$ be such an open geodesic ball with a center $x$ that $s_x$ and $e\tilde{x}p_x^{-1}$ are defined on $B_x$. Then

(4.18) $s_x=e\tilde{x}p_x\cdot S_x\cdot e\tilde{x}p_x^{-1}$,

where $e\tilde{x}p$ is the exponential mapping of $\tilde{\nabla}$ at the point $x$.

**Proof.** Let $y$ be a point in $B_x$, $y=e\tilde{x}p_x v$. Since $s_x$ is a local isometry and affine transformation with respect to $\tilde{\nabla}$, then it transforms the geodesic $e\tilde{x}p_x(tv)$ of $\tilde{\nabla}$ onto the geodesic $e\tilde{x}p_x(tSv)$ preserving the length of a segment of a curve. As $S_x=(s_x)_{*x}$, then

$$(s_x \cdot e\tilde{x}p_x)(v) = (e\tilde{x}p_x \cdot S_x)(v)$$

and the proposition follows.

**QED.**

By similar arguments we obtain

(4.19) $$s_x = exp_x \cdot S_x \cdot exp_x^{-1},$$

where *exp* is the exponential mapping of $\nabla$ at the point $x$, $s_x(B_x)=B_x$, where $s_x$ and $exp_x^{-1}$ are defined on an open geodesic ball $B_x$.

**2⁰**. Let $(M,\{s_x\})$ and $(M',\{s_y\})$ be two R.l.r. σ-m. with the canonical connections $\tilde{\nabla}$, $\tilde{\nabla}'$ defined by (4.17). A local isometry $\varphi$ of $M$ into $M'$ is called a local isomorphism if $\varphi_*(S) = S'$ in a domain of definition of $\varphi$.

**LEMMA 4.55.** A local isomorphism $\varphi$ is necessarily a local affine mapping of $(M,\tilde{\nabla})$ into $(M',\tilde{\nabla}')$.

**Proof.** Since $\varphi$ is a local isometry, then it is an affine mapping with respect to $\nabla$. From the invariance of $S$ under $\varphi$ and (4.17) our lemma follows.

**QED.**

A vector field $X \in \mathfrak{X}(M)$ is called an infinitesimal automorphism of R.l.r. σ-m. $M$ if, for each $x \in M$, a local 1-parameter group of local transformations $\varphi_t$ of a neighbourhood $U$ of $x$ into $M$ generated by $X$ consists of local automorphisms of $(M,\{s_x\})$.

**THEOREM 4.56.** (1) A vector field $X \in \mathfrak{X}(M)$ is an infinitesimal automorphism of a R.l.r. σ-m. if and only if $L_X g=0$, $L_X S=0$.

(2) The set of all the infinitesimal automorphisms of a R.l.r. σ-m. $M$ is a subalgebra $\sigma(M)$ in the algebra $i(M)$ of all the infinitesimal isometries of $M$.

(3) An infinitesimal automorphism $X$ on $M$ is an infinitesimal affine transformation of $(M,\tilde{\nabla})$.

**Proof.** (1) Let $\varphi_t$ be the local 1-parameter group of local transformations generated by $X$. A tensor field $P$ is invariant under $\varphi_t$ for every $t$ if and only if $L_X P=0$, [46].

(2) For $X,Y\in\sigma(M)$ we have

$$L_X g=L_X S=L_Y g=L_Y S=0,$$

therefore

$$L_{[X,Y]}g=L_X\cdot L_Y g-L_Y\cdot L_X g=0,$$

$$L_{[X,Y]}S=L_X\cdot L_Y S-L_X\cdot L_X S=0.$$

So, $\sigma(M)$ is a subalgebra in $i(M)$.

(3) For every $t$ $\varphi_t$ is a local isometry and $\varphi_{t*}(S)=S$, hence $\varphi_{t*}(S^{-1})=S^{-1}$. If $X\in T^1$, then $SX=X$ and $S\varphi_{t*}(X)=\varphi_{t*}(SX)=\varphi_{t*}(X)$, therefore $\varphi_{t*}(X)\in T^1$ and $T^1$ is invariant under $\varphi_t$. For $Y\in T^2$ we have

$$<X,\varphi_{t*}(Y)>=<\varphi_{t*}(X'),\varphi_{t*}(Y)>=<X',Y>=0,$$

because $X'=\varphi_{t*}^{-1}(X)\in T^1$. So, $\varphi_{t*}(Y)\in T^2$ and $\pi_1$, $\pi_2$ are invariant under $\varphi_t$. If $X,Y\in X(M)$, then from (4.17) we obtain

$$\begin{aligned}\varphi_{t*}(\tilde\nabla_X Y)&=\varphi_{t*}(\nabla_X Y-\nabla_{(I-S)^{-1}\pi_2 X}Y+S\nabla_{(I-S)^{-1}\pi_2 X}S^{-1}Y)\\&=\nabla_{\varphi_{t*}(X)}\varphi_{t*}(Y)-\nabla_{\varphi_{t*}((I-S)^{-1}\pi_2 X)}\varphi_{t*}(Y)\\&+S\nabla_{\varphi_{t*}((I-S)^{-1}\pi_2 X)}S^{-1}\varphi_{t*}(Y)=\tilde\nabla_{\varphi_{t*}(X)}\varphi_{t*}(Y)\end{aligned}$$

**QED.**

Since the torsion and curvature tensor fields $\tilde T$ and $\tilde R$ are invariant with respect to $\varphi_t$, then

(4.20) $$L_X\tilde T=L_X\tilde R=0,$$

where $X$ is any infinitesimal automorphism of R.l.r. $\sigma$-m. $M$.

**3⁰. LEMMA 4.57** [48],[58]. Let $(M,\tilde\nabla)$ be a manifold with a connection $\tilde\nabla$ and $\pi:T(M)\to M$ the canonical projection. Then there exists such a neighbourhood $U_0$ of the null section $O_M$ in $T(M)$ that the mapping

$$\pi \times e\tilde{x}p : v \mapsto (\pi(v), e\tilde{x}p_{\pi(v)} v)$$

is a diffeomorphism of $U_0$ onto a neighbourhood $U$ of the diagonal in $M\times M$.

It is evident that for every point $p$ we can choose such an open geodesic ball $B_p$ (a closed ball $\overline{B}_p$) with the center $p$ under the Riemannian connection $\nabla$ that $B_p \times B_p \subset U \quad (\overline{B}_p \times \overline{B}_p \subset U)$.

**LEMMA 4.58.** Let $M$ be a R.l.r. $\sigma$-m. and $p\in M$. Then there exists such an open ball $B_p$ and a differentiable mapping

$$\mu : B_p \times B_p \to M : \mu(x;y) - s_x(y) = (e\tilde{x}p_x \cdot S_x \cdot e\tilde{x}p_x^{-1})(y),\ x,y \in B_p$$

that every local subsymmetry $s_x$, $x\in B_p$, is an isometry of $B_p$ on $s_x(B_p)$.

**Proof.** Let $U_0^{'}$ be a neighbourhood of the null section , given in *Lemma 4.57*. We see that a set $S^{-1}U_0^{'} = \{ S_{\pi(v)}^{-1} : v \in U_0^{'} \}$ is also a neighbourhood of the null section. Let $U_0 = U_0^{'} \cap S^{-1}U_0^{'}$ and let $U$ be a corresponding neighbourhood of the diagonal in $M\times M$. We can take such an open geodesic ball $B_p$ (a closed ball $\overline{B}_p$) that $B_p \times B_p \subset U \ (\overline{B}_p \times \overline{B}_p \subset U)$. Then $\mu(x:y) = (e\tilde{x}p_x \cdot S_x \cdot e\tilde{x}p_x^{-1})(y)$ is defined on $B_p\times B_p$. The rest follows from *Proposition 4.54.*

**QED.**

It is evident that a set $\sigma(M)$ of all the local automorphisms of a R.l.r. $\sigma$-m. $(M,\{s_x\})$ is a pseudogroup of transformations on $M$, [46], p.1.

Let $\overline{B}_p(R)$ be a closed ball of the radius $R$, $\overline{\Lambda}_x$ be the connected component of $\Lambda_x \cap \overline{B}_p(R)$ containing $x\in B_p(R)$ and $B_p(R,\Lambda)=\{ \overline{\Lambda}_x : x \in B_p(R)\}$.

**THEOREM 4.59.** If $M$ is a R.l.r. $\sigma$-m., then for every $p\in M$ there exists such a geodesic ball $B_p(R)$ that $\sigma(M)$ acts locally transitive on $B_p(R,\Lambda)$.

**Proof.** For every point $p\in M$ we can choose such a closed ball $\overline{B}_p(2R)$ of the radius $2R$ that $\overline{B}_p(2R)\times\overline{B}_p(2R)\subset U$, where $U$ was considered in *Lemma 4.57, 4.58*, therefore $s_x(y)= \mu(x, y)$ is defined for any $x,y \in \overline{B}_p(2R)$. If $d(x_0,\overline{\Lambda}_x)$ is a distance between $x_0 = p$ and $\overline{\Lambda}_x$ then, since $\overline{\Lambda}_x$ is closed, there exists such a point $y_0 \in \overline{\Lambda}_x$ that $d(x_0,\overline{\Lambda}_x)=d(x_0,y_0)$. Let $\gamma_0(t)$ $t\in[0;t_0]$ be the unique geodesic segment with respect to $\nabla$ in $\overline{B}_p(R)$ that $\gamma_0(0)=x_0$, $\gamma_0(t_0)=y_0$ and $x_1=\gamma_0(t/2)$ the midpoint of $\gamma_0$. If $y_1=s_{x_1}(x_0)$, then

$$R>d(x_0,y_0)=d(x_0,x_1)+d(x_1,y_0)=d(s_{x_1}(x_0),s_{x_1}(x_1))+d(x_1,y_0)$$

$$=d(y_1,x_1)+d(x_1,y_0)\geq d(y_0,y_1)$$

and $y_1\in\overline{B}_p(2R)$.

Let $\gamma_1(t)$, $t\in[0;t_1]$ be the geodesic segment joining $y_1$ with $y_0$ and $x_2$ its midpoint. If $y_2=s_{x_2}(y_1)$, then

$$R>d(y_0,y_1)=d(y_0,x_2)+d(x_2,y_1)=d(s_{x_2}(x_2),s_{x_2}(y_1))+d(y_0,x_2)$$

$$=d(y_0,x_2)+d(x_2,y_2)\geq d(y_0,y_2)$$

Containing this process we get two sequences $\{y_n\}$, $\{x_n\}$ in $\overline{B}_p(2R)$ that

$$d(y_0,y_1)\geq d(y_0,y_2)\geq\ldots\geq d(y_0,y_n)\geq\ldots$$

and

$$d(y_0,x_2)\geq d(y_0,x_3)\geq\ldots\geq d(y_0,x_{n+1})\geq\ldots$$

because $x_{i+1}$ is the midpoint of segment $\gamma_i$ joining $y_0$ and $y_i$.

Let $\bar{x}$, $\bar{y}$ be the limit points of $\{x_n\}$, $\{y_n\}$ correspondingly. Since $y_n=s_{x_n}(y_{n-1})=\mu(x_n,y_{n-1})$ and $\mu$ is continuous in $B_p(2R)$ then, taking the passage to the limit, we obtain $\bar{y}=s_{\bar{x}}(\bar{y})$, where $\bar{x}$ is the midpoint of segment $\bar{\gamma}$ connecting $y_0$ and $\bar{y}$. As $s_{\bar{x}}(\bar{x})=\bar{x}$ the uniqueness of the geodesic segment joining two points in $\overline{B}_p(2R)$ implies that $s_{\bar{x}}(\bar{\gamma}_1)=\bar{\gamma}_1$, where $\bar{\gamma}_1$ is a part of $\bar{\gamma}$ between $\bar{x}$ and $\bar{y}$, $(s_{\bar{x}})_{*\bar{x}}(\dot{\bar{\gamma}})=S_{\bar{x}}(\dot{\bar{\gamma}})=\dot{\bar{\gamma}}$ and $\dot{\bar{\gamma}}\in T^1$. From *Theorem 4.1* it follows that $\bar{\gamma}\in\overline{\Lambda}_{y_0}$ and $\bar{y}\in\overline{\Lambda}_{y_0}$.

If we denote $\varphi_n=s_{x_n}\cdot\ldots\cdot s_{x_1}$, then every $\varphi_n$ is a local automorphism of $(M,\{s_x\})$ defined on some neighbourhood of $x_0=p$ and the sequence $\{\varphi_n\}$ has a limit point $\varphi$. It is clear that $\varphi$ is also a local automorphism of $(M,\{s_x\})$ i.e. $\varphi\in\sigma(M)$.

Similarly, for any other $\overline{\Lambda}_z\in B_p(R,\Lambda)$ there exists a local automorphism $\psi$

of $(M,\{s_x\})$ that $\psi(x_0)\in\overline{\Lambda}_z$. In this case $\psi\cdot\varphi^{-1}\in\sigma(M)$ is an automorphism locally transforming $\overline{\Lambda}_x$ in $\overline{\Lambda}_z$.

**QED.**

**PROBLEM.** To construct an example of a R.l.r. σ-m. which has such a foliation of mirrors that $\sigma(M)$ does not act locally transitive on the set $N=\{\Lambda_x : x\in M\}$.

$4^0$. **DEFINITION 4.10.** Let $(M,\{s_x\})$ be a R.l.r. σ-m. A submanifold $M'\subset M$ is said to be invariant if the following holds:
for every two points $p,q\in M'$ and every local automorphism $\varphi$ of $M$ such that $\varphi(p)=q$ we have $\varphi(M'\cap\varphi^{-1}(M))\subset M'$.

**THEOREM 4.60.** Let $M'\subset M$ be an invariant submanifold of $(M,\{s_x\})$. Then $(M',\{s_{x|M'}\})$ is naturally a R.l.r. σ-m.

**Proof.** From *Definition 4.10* it follows that for every $x\in M'$, $s_{x|M'}$ is a local isometry of $(M',g)$. If we consider an integral curve $\gamma(t)$ in M′ of a vector $X\in T_x(M')$, $x=\gamma(t)$, then $s_x(\gamma(t))\subset M'$ and $S_x(X)=(s_x)_{*x}(X)\in T_x(M')$.

So, $T_x(M)$ is invariant with resect to $S$ and the conditions 1), 2) of *Definition 4.1* are fulfilled.

Let $\pi$ be the projection to $T_x(M')$ of $T_x(M)=T_x(M')\oplus T_x(M')^{\perp}$. It is well-known that $\nabla'_X Y=\pi\nabla_X Y,\quad X,Y\in\mathfrak{X}(M)$, is the Riemannian connection of $(M',g)$ and it is clear that $\tilde{\nabla}'_X Y=\pi\tilde{\nabla}_X Y$ is also a metric connection, i.e., $\tilde{\nabla}' g=0$ on $M'$. As $T_x(M')$ is invariant under $S$, then $\pi\cdot S=S\cdot\pi$. Using (4.1) for $X,Y\in\mathfrak{X}(M')$ we have

$$(\tilde{\nabla}'_X S)Y=\tilde{\nabla}'_X SY-S\tilde{\nabla}'_X Y=\pi\tilde{\nabla}_X SY-S\pi\tilde{\nabla}_X Y=\pi\tilde{\nabla}_X SY-\pi S\tilde{\nabla}_X Y=\pi(\tilde{\nabla}_X S)Y=0$$

.

So, $\tilde{\nabla}'$ is a canonical connection of $(M',\{s_x\})$ and the axiom 3) of *Definition 4.1* is realised.

**QED.**

If $M$ is a R.l.r.$s$-m., then every invariant submanifold $M'$ of $M$ is autoparallel with respect to the canonical connection $\tilde{\nabla}$ of $(M,\{s_x\})$, see [48]. It is not so for a R.l.r. σ-m.

**EXAMPLE.** Let $M$ be a R.l.r.$s$-m. and $\Lambda' \subset \Lambda$ be a submanifold of a Riemannian manifold $\Lambda$ and also $\Lambda'$ is not an autoparallel submanifold with respect to the Riemannian connection $\nabla$. By a natural way $M\times\Lambda$ is a R.l.r. σ-m., where

$$s_{(x,u)}(y,v) = (s_x(y),v), \quad x,y \in M, \quad u,v \in \Lambda.$$

It is clear that $M\times\Lambda'$ is an invariant submanifold of $M\times\Lambda$ which is not autoparallel with respect to the canonical connection $\tilde{\nabla}$ defined by (4.1) because $\tilde{\nabla}_{|(x,\Lambda')} = \nabla$.

# CHAPTER 5
# CLASSICAL STRUCTURES

In this chapter, we have our methods illustrated for concrete structures which are called classical.

§1 is devoted to the study of canonical connection $\overline{\nabla}$ and the second fundamental tensor field *h* of a Riemannian almost product structure *(P,g)* and their relations with geometry of a manifold *M*. A behaviour of *P* along geodesics and ties with another structures as, for example, a structure of reflection space are discussed here too.

In §2, we conduct the same policy for almost Hermitian manifolds. Also, the classification of A.Gray and L.M.Hervella has been rewritten in terms of the field *h*, that is, the canonical connection $\overline{\nabla}$ for every from 16 classes is constructed.

Some examples of almost Hermitian manifolds are discussed in §3. We consider a projection of the classes on submanifolds, conformal changes of a metric *g*, a quasi homogeneous structure *(J, g)* of the class $\boldsymbol{U_4}$, a behaviour of *J* along a curve, almost Hermitian structures on Riemannian locally regular *s*-manifolds.

In the last §4, a structure defined by an affinor *F* satisfying $F^3+F=0$ and an associated to *F* Riemannian metric *g* is considered. For *(F, g)* we have obtained the canonical connection $\overline{\nabla}$ and the tensor field *h*. Some relations with an affinor $S_x = (s_x)_{*x}$ of a locally 4 - symmetric Riemannian space are also looked through.

We refer to [19] , [22] , [26] , [34] , [35] , [37].

## *§1. ALMOST PRODUCT MANIFOLDS*

**1⁰**. A tensor field P of type *(1,1)* such that $P^2=I$ is called an almost product structure on a manifold *M*. We put

$\pi_1=1/2(I+P),\ \pi_2=1/2(I-P).$

Then

(5.1) $\pi_{1_2} + \pi_2 = I,\ \pi_1^2 = \pi_1,\ \pi_2^2 = \pi_2,\ \pi_1\pi_2 = \pi_2\pi_1 = 0;\ P = \pi_1 - \pi_2 .$

For any Riemannian metric $\tilde{g}$ on *M* a new Riemannian metric *g* is defined by the formula

$g(X,Y)= \tilde{g}(X,Y)+ \tilde{g}(PX,PY),$

where $X,Y\in\mathcal{X}(M)$. Every such a metric $g= < , >$ satisfies the formula

(5.2) $<PX, PY> = <X, Y>$

Let $\nabla$ be the Riemannian connection of such a fixed metric $g=< , >$ on $M$. We define a connection $\overline{\nabla}$ on M by

(5.3) $$\overline{\nabla}_X Y = \nabla_X Y - \frac{1}{2}\nabla_X (P)PY = \pi_1\nabla_X \pi_1 Y + \pi_2\nabla_X \pi_2 Y,$$

where $X,Y\in\mathcal{X}(M)$.

Then

$$\overline{\nabla}_X Y = \frac{1}{2}(\nabla_X Y + P\nabla_X PY)$$

and

$$2\overline{\nabla}_X (P)Y = 2(\overline{\nabla}_X PY - P\overline{\nabla}_X Y) = \overline{\nabla}_X PY + P\nabla_X Y - P\nabla_X Y - P\nabla_X Y = 0$$

Using (5.2) we obtain

$$\begin{aligned} 2(<\overline{\nabla}_X Y,Z> + <Y,\overline{\nabla}_X Z>) &= <\nabla_X Y,Z> + <Y,\nabla_X Z> + <P\nabla_X PY,Z> + <Y,P\nabla_X PZ> \\ &= X<Y,Z> + <\nabla_X PY,PZ> + <PY,\nabla_X PZ> \\ &= X<Y,Z> + X<PY,PZ> = 2X<Y,Z> \end{aligned}$$

for $X,Y,Z\in\mathcal{X}(M)$, i.e., $\overline{\nabla}g = 0$.

According to *Theorem 1.2* a tensor field $P$, $P^2 = I$, is always $0$-deformable, that is, $P$ is an affinor, and every Riemmannian metric $g$ satisfying (5.2) is an associated one to the corresponding $G$-structure.

The space of all such associated metrics is infinite dimensional. Later on, we shall consider only associated metrics satisfying (5.2).

We have $\overline{\nabla}\pi_1 = \overline{\nabla}\pi_2 = 0$, therefore $\pi_1$ and $\pi_2$ define two complementary distributions $T^1=\pi_1(T(M))$ and $T^2=\pi_2(T(M))$, where $T^1$ is a distribution corresponding to the eigenvalue $1$ of $P$ and $T^2$ is a distribution corresponding to the eigenvalue -1 of P, hence

$$T(M)= T^1(M)\oplus T^2(M)$$

We can consider a set *P(G)* of linear frames over *M* such that for every $u \in P(G)$ the affinor *P* has the following matrix

$$(P) = \begin{bmatrix} E_{n_1} & 0 \\ 0 & -E_{n_2} \end{bmatrix}$$

where $n_1$ and $n_2$ are dimensions of $T^1$ and $T^2$ respectively. The set *P(G)* is a *G*-structure with the structure group *G* consisting of matrixes of the following form

$$\begin{bmatrix} A_1 & 0 \\ 0 & A_2 \end{bmatrix},$$

where $A_k \in GL(n_k, R)$, $k = 1,2$.

From *Theorem 1.2* it follows that there exists a reduction of *G* to its maximal compact subgroup *H* consisting of matrixes of the form above, where $A_k \in O(n_k)$, $k = 1,2$. A structure *P(H)* defines a Riemannian metric $g = \langle \, , \, \rangle$ on *M*, which is evidently an associated metric to *P(G)* and $P(H) = P(G) \cap O(M)$. Every such defined a metric *g* satisfies (5.2) and a pair *(P, g)* is called an almost-product Riemannian structure (a.p.R.s.).

Further. we consider the canonical connection $\overline{\nabla}$ and the second fundamental tensor field of the structure *(P(H), g)*.

**THEOREM 5.1.** The canonical connection $\overline{\nabla}$ of the *G*-structure *(P(H),g)* corresponding to an a.p.R.s. *(P,g)* is defined by (5.3). The second fundamental tensor field *h* of *(P(H),g)* is determined by

$$(5.4) \quad h_X Y = \pi_1 \nabla_X \pi_2 Y + \pi_2 \nabla_X \pi_1 Y = \frac{1}{2} \nabla_X (P) PY,$$

where $X, Y \in \mathfrak{X}(M)$ and $\nabla$ is the Riemannian connection of the associated Riemannian metric *g*.

**Proof.** It is clear that the Lie algebra $\underline{\boldsymbol{h}}$ of the structure group *H* of *P(H)* has the following form

$$\underline{\boldsymbol{h}} = \{ \tilde{x} = \begin{bmatrix} B_1 & 0 \\ 0 & B_2 \end{bmatrix} : \tilde{x} \in \underline{\boldsymbol{o}} \}.$$

It is evident that $\underline{\boldsymbol{o}} = \underline{\boldsymbol{h}} \oplus \underline{\boldsymbol{m}}$, where

$$\underline{\mathbf{m}} = \{ \tilde{x} = \begin{bmatrix} 0 & C_1 \\ C_2 & 0 \end{bmatrix} : \tilde{x} \in \underline{\mathbf{o}} \}$$

and

$$tr( \begin{bmatrix} B_1 & 0 \\ 0 & B_2 \end{bmatrix} \cdot \begin{bmatrix} 0 & C_1 \\ C_2 & 0 \end{bmatrix}^T ) = 0 ,$$

therefore $\underline{\mathbf{m}} = \underline{\mathbf{h}}^{\perp}$ with respect to Killing form . For each $\omega \in \underline{\mathbf{o}}$ we can define the natural decomposition $\omega = \omega_{|\underline{\mathbf{h}}} + \omega_{|\underline{\mathbf{m}}}$ by the formula

$$\omega = ( p_1 \omega p_1 + p_2 \omega p_2 ) + ( \omega - p_1 \omega p_1 - p_2 \omega p_2 ),$$

where

$$p_1 = \begin{bmatrix} E_{n_1} & 0 \\ 0 & 0 \end{bmatrix}, \quad p_2 = \begin{bmatrix} 0 & 0 \\ 0 & E_{n_2} \end{bmatrix}$$

Let $\varphi$ be a cross section of $P(H)$ over some neighbourhood $U$ which assigns to each $x \in U$ the linear frame $((X_1)_x, \ldots, (X_n)_x)$ and $X, Y = \sum_k f^k X_k$ be vector fields on $M$. Then from (1.6) and *Definition 1.9* it follows that

$$\overline{\nabla}_{X_x} Y = \varphi(x) \overline{\omega}(\varphi_* X_x) \varphi(x)^{-1} Y_x + \sum_k (X f^k)(x)(X_k)_x$$

$$= \varphi(x)[ p_1 \omega(\varphi_* X_x) p_1 + p_2 \omega(\varphi_* X_x) p_2 ] \varphi(x)^{-1} Y_x + \sum_k (X f^k)(x)(X_k)_x$$

where $\varphi(x)$ is considered as a mapping of $R^n$ onto $T_x(M)$. It is obvious that

$$\varphi(x) \cdot p_1 = \pi_1 \cdot \varphi(x), \quad \varphi(x) \cdot p_2 = \pi_2 \cdot \varphi(x)$$

and

$$(\pi_1 + \pi_2)[(X f^k)(x)(X_k)_x] = (X f^k)(x)(\pi_1 X_k)_x + (X f^k)(x)(\pi_2 X_k)_x$$

We obtain that

$$\overline{\nabla}_{X_x} Y = [\pi_1 \varphi(x)\omega(\varphi_* X_x)\varphi(x)^{-1}\pi_1 Y_x + \sum_k (Xf^k)(x)(\pi_1 X_k)_x]$$
$$+ [\pi_2 \varphi(x)\omega(\varphi_* X_x)\varphi(x)^{-1}\pi_2 Y_x + \sum_k (Xf^k)(x)(\pi_2 X_k)_x]$$
$$= \pi_1 \nabla_{X_x} \pi_1 Y + \pi_2 \nabla_{X_x} \pi_2 Y$$
$$\overline{\nabla}_X Y = \frac{1}{4}(I+P)\nabla_X (I+P)Y + \frac{1}{4}(I-P)\nabla_X (I-P)Y$$
$$= \frac{1}{4}(\nabla_X Y + P\nabla_X Y + \nabla_X PY + P\nabla_X PY + \nabla_X Y - P\nabla_X Y - \nabla_X PV + P\nabla_X PY)$$
$$= \frac{1}{2}(\nabla_X Y + P\nabla_X PY) = \nabla_X Y - \frac{1}{2}\nabla_X(P)PY$$
$$h_X Y = \frac{1}{2}(\nabla_X Y - P\nabla_X PY) = \frac{1}{2}\nabla_X(P)PY = \pi_1 \nabla_X \pi_2 Y + \pi_2 \nabla_X \pi_1 Y$$

**QED.**

Later on, we shall call $\overline{\nabla}$ and $h$ the canonical connection and the second fundamental tensor field of an a.p.R.s. *(P, g)* respectively.

**2⁰**. We consider now the torsion tensor field $\overline{T}$ of the canonical connection $\overline{\nabla}$ of an a.p.R.s. *(P, g)*. It follows from (2.3) and from *Proposition 2.3* that $\overline{T} = -2h^-$ and $\overline{T} = 0$ if and only if the pair *(P, g)* is a particular structure, that is, $\nabla P = 0$. We have, see [46], [73], that the condition $\nabla P = 0$ implies a manifold $M$ to be a locally decomposable Riemannian manifold or, if $M$ is a complete, simply connected Riemannian manifold, to be a globally decomposable Riemannian manifold $M = M_1 \times M_2$, where $M_1$, $M_2$ are the maximal integral manifolds of the distributions $T^1$, $T^2$ passing through a fixed point of $M$.

Thus, a particular structure *(P, g)* is a structure of local Riemannian product.

From *Theorem 1.5* it follows that an almost product structure is integrable if and only if the Nijenhuis tensor field

$$N(P)(X,Y) = [PX, PY] - P[X, PY] - P[PX, Y] + [X, Y]$$

vanishes on $M$, where $X,Y \in \mathcal{X}(M)$.

**PROPOSITION 5.2.** For $X,Y \in \mathcal{X}(M)$

$$N(P)(X,Y) = -2((\overline{T}_{PX} PY + \overline{T}_X Y).$$

**Proof.** We remark that from (5.3)

$$\overline{T}_X Y=\overline{\nabla}_X Y-\overline{\nabla}_Y X-[X,Y]=\frac{1}{2}(\nabla_X Y-\nabla_Y X-[X,Y]+P\nabla_X PY-P\nabla_Y PX-[X,Y])$$

$$=1/2(P\nabla_X PY-P\nabla_Y PX-[X,Y]).$$

It follows from this identity that

$$-N(P)(X,Y)=P[PX,Y]+P[X,PY]-[PX,PY]-[X,Y]$$

$$=P\nabla_{PX}Y-P\nabla_Y PX+P\nabla_X PY-P\nabla_{PY}X-[PX,PY]-[X,Y]$$

$$=(P\nabla_{PX}Y-P\nabla_{PY}X-[PX,PY])+(P\nabla_X PY-P\nabla_Y PX-[X,Y])$$

$$=2[\overline{T}_{PX}PY+\overline{T}_X Y]$$

**QED.**

Thus, an almost product structure $P$ is integrable if and only if

$$\overline{T}_{PX}PY=-\overline{T}_X Y\,.$$

**PROPOSITION 5.3.** The distribution $T^1(T^2)$ is integrable if and only if $\overline{T}_X Y=0$ for any $X,Y\in T^1(T^2)$.

**Proof.** Using (5.3) we obtain

$$\overline{T}_{\pi_1 X}\pi_1 Y=\pi_1\nabla_{\pi_1 X}\pi_1 Y-\pi_1\nabla_{\pi_1 Y}\pi_1 X-[\pi_1 X,\pi_1 Y]=\pi_1(\nabla_{\pi_1 X}\pi_1 Y-\nabla_{\pi_1 Y}\pi_1 X$$
$$-[\pi_1 X,\pi_1 Y])-\pi_2[\pi_1 X,\pi_1 Y]=-\pi_2[\pi_1 X,\pi_1 Y].$$

The rest follows from the Frobenius theorem, [64].

**QED.**

Let $\overline{R}$ and $R$ be the curvature tensor fields of the connection $\overline{\nabla}$ and $\nabla$ respectively, then for $X,Y,Z\in\mathcal{X}(M)$ we have

$$\overline{R}_{XY}Z=\frac{1}{2}R_{XY}Z+\frac{1}{2}PR_{XY}PZ+h_Y h_X Z-h_X h_Y Z\,.$$

Really, it is easily to check that

$$4h_Y h_X Z = \nabla_Y \nabla_X Z - \nabla_Y P \nabla_X PZ - P \nabla_Y P \nabla_X Z + P \nabla_Y \nabla_X PZ$$
$$4h_X h_Y Z = \nabla_X \nabla_Y Z - \nabla_X P \nabla_Y PZ - P \nabla_X P \nabla_Y Z + P \nabla_X \nabla_Y PZ,$$

and

$$4\overline{R}_{XY} Z = 4(\overline{\nabla}_X \overline{\nabla}_Y Z - \overline{\nabla}_Y \overline{\nabla}_X Z - \overline{\nabla}_{[X,Y]} Z)$$
$$= (\nabla_X \nabla_Y Z + \nabla_X P \nabla_Y PZ + P \nabla_X P \nabla_Y Z + P \nabla_X \nabla_Y PZ)$$
$$-(\nabla_Y \nabla_X Z + \nabla_Y P \nabla_X PZ + P \nabla_Y P \nabla_X Z + P \nabla_Y \nabla_X PZ) - 2(\nabla_{[X,Y]} Z + P \nabla_{[X,Y]} PZ)$$

$$= 2R_{XY} Z + 2PR_{XY} PZ + 4h_Y h_X Z - 4h_X h_Y Z.$$

**$3^0$**. Let $T^1$ be integrable and $M_1$ be the maximal integral manifold of $T^1$ passing through some point of $M$. The manifold $M_1$ is a $n_1$- dimensional submanifold of $M$. For any $X,Y \in X(M_1)$ we know that

$$\nabla_X Y = \nabla^1_X Y + \alpha(X,Y) = \pi_1 \nabla_X Y + \pi_2 \nabla_X Y = \overline{\nabla}_X Y + h_X Y,$$

where $\nabla^1$ is the Riemannian connection and $\alpha(X,Y)$ is the second fundamental form of $M_1$. Equating the tangential and normal components of this equality we have got that

$$\nabla^1_X Y = \overline{\nabla}_X Y, \ \alpha(X,Y) = h_X Y.$$

Thus, the notion of the second fundamental tensor field $h$ of an almost product Riemannian structure (Riemannian $H$-structure) is a generalization of that of the second fundamental form of a submanifold and it is the same for the notion of the canonical connection $\overline{\nabla}$.

Since $\overline{T} = -2h^-$, then $h^- = 0$ on $M_1$ and $h_X Y = h^+_X Y$ on $M_1$. It is well known, [47], that the following conditions are equivalent

1) $M_1$ is an autoparallel submanifold with respect to $\nabla$,
2) $M_1$ is totally geodesic with respect to $\nabla$,
3) $\alpha(X,Y) = h_X Y = 0$ on $M_1$.

Now it is evident the following

**PROPOSITION 5.4.** The distribution $T^1(T^2)$ is integrable and defines the foliation of totally geodesic with respect to $\nabla$ maximal integral submanifolds if and only if $h_X Y = 0$ for any $X,Y \in T^1(T^2)$.

The following proposition describes the similar situation for both the distributions.

**PROPOSITION 5.5.** The distributions $T^1$ and $T^2$ are integrable and define the foliations of totally geodesic with respect to $\nabla$ maximal integral submanifolds if and only if $h=0$ on $M$, that is, $(P, g)$ is the local Riemannian product structure .

**Proof.** It follows from *Proposition 5.4* that $h_X Y=0$, when $X,Y\in T^1$ or $X,Y\in T^2$. Let a vector field $X$ be from $T^1$ and a vector field $Y$ be from $T^2$. Then $h_{XYZ} = <\pi_1 \nabla_X Y, Z> = -h_{XZY}$ from (2.2). Thus, if $Z\in T^2$ then $<\pi_1 \nabla_X Y, Z>=0$, if $Z\in T^1$ then $h_X Z=0$ because $X,Z\in T^1$ , and we have obtained that $h_{XYZ}=0$ for any $Z\in \mathfrak{X}(M)$. The case when $X\in T^2$, $Y\in T^1$ is the same.

The converse is evident.

**QED.**

Let $(P, g)$ be a nearly particular structure, that is, $h_X Y = -h_Y X$ for any $X,Y\in \mathfrak{X}(M)$. For an a.p.R.s. this condition is equivalent to one that $\nabla_X(P)X=0$. Really, we have

$$2h_X X=\nabla_X X-P\nabla_X PX=P(P\nabla_X X-\nabla_X PX)=-P\nabla_X(P)X.$$

Since $P$ is a nonsingular affinor hence $h_X X=0$ if and only if $\nabla_X(P)X=0$.

**THEOREM 5.6.** An a.p.R.s. $(P, g)$ is a nearly particular structure if and only if it is a particular structure (a local Riemannian product structure).

**Proof.** 1) If $X\in T^1$ and $Y\in T^2$, then from (5.4) we have

$$h_X Y=\pi_1\nabla_X Y= -h_Y X = -\pi_2\nabla_Y X=0$$

because $h_X Y\in T^1\cap T^2=\{0\}$.

2) If $X,Y\in T^2$, then we obtain

$$h_X Y+h_Y X=\pi_1(\nabla_X Y+\nabla_Y X)=0, \quad \pi_1(\nabla_X Y-\nabla_Y X)= \pi_1([X,Y]).$$

From these equalities it follows that

$$h_X Y=1/2\pi_1([X,Y]).$$

and

$$h_{XYZ} = \langle h_X Y,Z\rangle = 1/2\langle \pi_1([X,Y]),Z\rangle = -h_{XZY} = -\langle h_X Z,Y\rangle.$$

If $Z\in T^1$, then from 1) it follows that $h_X Z=0$ because $X\in T^2$.
If $Z\in T^2$, then $\langle \pi_1([X,Y]),Z\rangle=0$.
Thus, we see that $h_X Y=0$ in this case.
3) If $X,Y\in T^1$, then $h_X Y=1/2\pi_2([X,Y])$ and the rest is similar to 2).
The converse is evident.

**QED.**

It is obvious that $g'=e^{2\rho} g$ is an associated to $P$ Riemannian metric satisfying (5.2) too.

**PROPOSITION 5.7.** The tensor field h of a structure $(P, g)$ is invariant under conformal changes of the metric $g$.

**Proof.** If $g'=e^{2\rho} g$, then

(5.5) $\nabla'_X Y = \nabla_X Y + X(\rho)Y + Y(\rho)X - \langle X,Y\rangle \, grad\, \rho$,

where $\langle grad\, \rho, X\rangle = X\rho$. Really, we have

$$\nabla'_X Y - \nabla'_Y X = \nabla_X Y - \nabla_Y X = [X,Y],$$
$$\langle\nabla'_X Y,Z\rangle' + \langle Y,\nabla'_X Z\rangle' = e^{2\rho}[\langle\nabla_X Y,Z\rangle + \langle Y,\nabla_X Z\rangle + X(\rho)\langle Y,Z\rangle + Y(\rho)\langle X,Z\rangle$$
$$-\langle X,Y\rangle Z(\rho) + X(\rho)\langle Z,Y\rangle + Z(\rho)\langle X,Y\rangle - \langle X,Z\rangle Y(\rho)]$$
$$= e^{2\rho}[X\langle Y,Z\rangle + 2X(\rho)\langle Y,Z\rangle] = X[e^{2\rho}\langle Y,Z\rangle] = X\langle Y,Z\rangle'$$

From (5.4), (5.5) it follows that

$h'_X Y = \pi_1\nabla'_X \pi_2 Y + \pi_2\nabla'_X \pi_1 Y = h_X Y$, where $X,Y,Z\in\mathfrak{X}(M)$.

**QED.**

**4⁰**. We consider now a behavior of an almost product structure $P$ along geodesics.

**LEMMA 5.8.** Let $\gamma(t)$ be a geodesic with respect to $\overline{\nabla}$ $(\nabla)$ and $X=\gamma'(t)$.

Then $\nabla_X \overline{\nabla}_X Y = \overline{\nabla}_X \nabla_X Y$ if and only if $(\overline{\nabla}_X h)(X,Y)=0$ $((\nabla_X h)(X,Y)=0)$, where a vector field $Y$ is defined on some neighbourhood of $\gamma(t)$.

**Proof.** We have $\overline{\nabla}_X X = 0$ $(\nabla_X X = 0)$ and

$$(\overline{\nabla}_X h)(X,Y) = \overline{\nabla}_X h_X Y - h_{\overline{\nabla}_X X} Y - h_X \overline{\nabla}_X Y$$

$$= \overline{\nabla}_X \nabla_X Y - \overline{\nabla}_X \overline{\nabla}_X Y - \nabla_X \overline{\nabla}_X Y + \overline{\nabla}_X \overline{\nabla}_X Y = \overline{\nabla}_X \nabla_X Y - \nabla_X \overline{\nabla}_X Y$$

.

It is the same for $\nabla$.

**QED.**

If $(P, g)$ is a quasi homogeneous structure, then $\overline{\nabla} h = 0$ and from *Lemma 5.8* it follows that

(5.6) $\nabla_X \overline{\nabla}_X Y = \overline{\nabla}_X \nabla_X Y$.

Let $M$ and $P$ be analytic, $\gamma(t)$ be such a curve on $M$ that $\nabla_X \overline{\nabla}_X Y = \overline{\nabla}_X \nabla_X Y$, where $X = \gamma'(t)$ and $Y$ is any vector field on some neighbourhood of $\gamma(t)$. We denote $P^{(k)}(t) = (\nabla^k_{X\ldots X} P)$.

**LEMMA 5.9.** $\nabla_X \overline{\nabla}_X Y = \overline{\nabla}_X \nabla_X Y$ if and only if $\overline{\nabla}_X h_X Y = h_X \overline{\nabla}_X Y$ or $\nabla_X h_X Y = h_X \nabla_X Y$.

**Proof.** We can compare the identities $\overline{\nabla}_X h_X Y = \overline{\nabla}_X \nabla_X Y - \overline{\nabla}^2_{XX} Y$ and $h_X \overline{\nabla}_X Y = \nabla_X \overline{\nabla}_X Y - \overline{\nabla}^2_{XX} Y$. The proof of the second equivalence is the same.

**QED.**

**PROPOSITION 5.10.** If $\nabla_X \overline{\nabla}_X Y = \overline{\nabla}_X \nabla_X Y$, then we have

(5.7) $P\nabla_X(P)Y + \nabla_X(P)PY = 0,$

(5.8) $-P\nabla^2_{XX}(P)Y = \nabla_X(P)^2 Y$,

(5.9) $P\nabla_X(P)^k Y = (-1)^k \nabla_X(P)^k PY,$

(5.10) $\nabla_X(\nabla_X(P)^2)=0.$

**Proof.** (5.7). Differentiating the equality $P^2=I$ we obtain $\nabla_X(P)P+P\nabla_X(P)=0$ and it follows from (5.4) that $\nabla_X(P)Y=2h_X PY=-2Ph_XY.$

(5.8)

$$-P\nabla^2_{XX}(P)Y=-2P[\nabla_X h_X PY-h_X P\nabla_X Y]=--2P[h_X\nabla_X PY-h_X P\nabla_X Y]$$

$$=-2Ph_X\nabla_X(P)Y = \nabla_X(P)^2Y.$$

(5.9). Using the method of mathematical induction we get that for $k=1$, $P\nabla_X(P)Y=-\nabla_X(P)PY$, i.e., (5.7).

Let (5.9) be true for $k-1$, then

$$P\nabla_X(P)^kY=(-1)^{k-1}\nabla_X(P)^{k-1}\cdot P\nabla_X(P)Y=(-1)^k\nabla_X(P)^kPY.$$

(5.10) $$P\nabla_X(\nabla_X(P)^2)=P\nabla^2_{XX}(P)\nabla_X(P)+P\nabla_X(P)\nabla^2_{XX}(P)$$
$$=-\nabla_X(P)^2\nabla_X(P)-\nabla_X(P)P\nabla^2_{XX}(P)$$
$$=-\nabla_X(P)^3+\nabla_X(P)^3=0.$$

**QED.**

(5.7) - (5.10) can be rewrritten in the following form

(5.11) $P\cdot P'+P'\cdot P = 0,$

(5.12) $P''= -P(P')^2=-(P')^2P,$

(5.13) $P\cdot(P')^k=(-1)^k\cdot(P')^k\cdot P,$

(5.14) $[(P')^2]'=0.$

We describe now the $k$-th derivative of $P$.

**PROPOSITION 5.11.** If $\nabla_X\overline{\nabla}_XY=\overline{\nabla}_X\nabla_XY$, then we have

(5.15) $P^{(2m)}=(-1)^mP(P')^{2m},$

(5.16) $P^{(2m+1)}= (-1)^m(P')^{2m+1}.$

**Proof.** Using the method of mathematical induction from (5.12), (5.14) we obtain (5.15). Really, if $m=1$, then $P''=P(P')^2$. Let (5.15) be true for $m$, then

$$P^{(2m+2)}=[P^{(2m)}]''=(-1)^m\{P[(P')^2]^m\}''=(-1)^m P''\cdot(P')^{2m}=(-1)^{m+1}P(P')^{2m+2}.$$

To get (5.16) we use (5.15) and obtain

$$P^{(2m+1)}=[P^{(2m)}]'=(-1)^m\{P[(P')^2]^m\}'=(-1)^m(P')^{2m+1}.$$

**QED.**

We can define *cos P* and *sin P* by the formulas

$$\cos P=\sum_{m=0}^{\infty}(-1)^m P^{2m}/(2m)!,\quad \sin P=\sum_{m=0}^{\infty}(-1)^m P^{2m+1}/(2m+1)!.$$

Let $\tau_t$ be the parallel translation from $\gamma(0)$ to $\gamma(t)$ along a curve $\gamma(t)$ with respect to $\nabla$.

**PROPOSITION 5.12.** Let $P(t)$ be an almost product structure in $\gamma(t)$ and $\nabla_X\overline{\nabla}_X Y=\overline{\nabla}_X\nabla_X Y$ along $\gamma(t)$. Then

$$P(t)=\tau_t\cdot[\cos(tP'(o))P(o)+\sin(tP'(o))]\cdot\tau_t^{-1}.$$

**Proof.** Using (5.15), (5.16) we expand $\tau_t^{-1}\cdot P(t)\cdot\tau_t$ in a power series in the case, when $t=o$.

$$\tau_t^{-1}\cdot P(t)\cdot\tau_t=\sum_{k=0}^{\infty}t^k/k!P^{(k)}(o)=\sum_{m=0}^{\infty}[t^{2m}/(2m)!P^{(2m)}(o)+t^{2m+1}/(2m+1)!P^{(2m+1)}(o)]$$

$$=\sum_{m=0}^{\infty}[(-1)^m t^{2m}/(2m)!P'(o)^{2m}P(o)+(-1)^m t^{2m+1}/(2m+1)!P'(o)^{2m+1}]$$

$$=\cos(tP'(o))P(o)+\sin(tP'(o)).$$

**QED.**

**THEOREM 5.13.** Let a) $-\lambda^2$, b) $\lambda^2(\lambda\neq o)$ be an eigenvalue of $P'(o)^2$ and let $X(o)$ be a corresponding eigenvector, $X(t)=\tau_t(X(o))$. Then

a) $P(t)X(t)=\cosh(\lambda t)\tau_t(P(o)X(o))+1/\lambda\,\sinh(\lambda t)\tau_t(P'(o)X(o)),$

b) $P(t)X(t) = \cos(\lambda t)\tau_t(P(o)X(o))+1/\lambda \sin(\lambda t)\tau_t(P'(o)X(o))$.

**Proof.** Using *Proposition 5.12* we obtain

a) $\tau_t^{-1}(P(t)X(t))=(\tau_t^{-1}\cdot P(t)\cdot\tau_t)X(o)$

$$=\sum_{m=0}^{\infty}[(-1)^m\frac{t^{2m}}{(2m)!}P(o)P'(o)^{2m}X(o)+(-1)^m\frac{t^{2m+1}}{(2m+1)!}P'(o)^{2m+1}X(o)]$$

$$=\sum_{m=0}^{\infty}[(-1)^m\frac{t^{2m}}{(2m)!}(-1)^m\lambda^{2m}P(o)X(o)+(-1)^m\frac{t^{2m+1}}{(2m+1)!}(-1)^m\lambda^{2m}P'(o)X(o)]$$

$$=\sum_{m=0}^{\infty}[\frac{(t\lambda)^{2m}}{(2m)!}P(o)X(o)+\frac{1}{\lambda}\frac{(t\lambda)^{2m+1}}{(2m+1)!}P'(o)X(o)]$$

$$=\cosh(\lambda t)P(o)X(o)+\frac{1}{\lambda}\sinh(\lambda t)P'(o)X(o).$$

b) $\tau_t^{-1}(P(t)X(t))=\sum_{m=0}^{\infty}[(-1)^m\frac{(t\lambda)^{2m}}{(2m)!}P(o)X(o)+\frac{1}{\lambda}(-1)^m\frac{(t\lambda)^{2m+1}}{(2m+1)!}P'(o)X(o)]$

$$=\cos(\lambda t)P(o)X(o)+\frac{1}{\lambda}\sin(\lambda t)P'(o)X(o).$$

**QED.**

**COROLLARY.** Let $\gamma(t)$ be a periodic geodesic and $\Lambda(t)$ be a proper space corresponding to an eigenvalue $\lambda^2$ $(\lambda\neq o)$ of $P'(t)^2$. Then the restriction of $P(t)$ on $\Lambda(t)$ is parallel to the restriction of $P(t+2\pi k/\lambda)$ on $\Lambda(t+2\pi k/\lambda)$, $k\in Z$.

**Proof.** Let $X(t)$ and $X(t+2\pi k/\lambda)$ be parallel each other eigenvectors of $P'^2$. From *Theorem 5.13* b) it follows that $P(t+2\pi k/\lambda)X(t+2\pi k/\lambda)$ is parallel to $P(t)X(t)$.

**QED.**

Similar results were considered in [34] for an almost complex structures of nearly Kaehlerian manifolds.

**$5^0$.** A classification of almost product Riemannian manifolds was obtained in [59]. Using tensor fields $h$ of such structures, (5.4), an analogous classification can be got. In the next section we discuss a similar construction for almost Hermitian manifolds.

Let $\Phi$ be a pseudo-Riemannian metric on $M$ defined by a $G$-structure $P(G)$. It follows from *Theorem 1.1* that there exists a reduction of $G$ to the maximal compact subgroup $H=O(n, n_1)$, which defines an associated Riemannian metric $g= \langle\, ,\, \rangle$. The space of all such associated metrics is infinite dimensional. For a fixed $g= \langle\, ,\, \rangle$ one can consider the a.p.R.s. $P$ defined by condition

(5.17) $\Phi(X, Y) = \langle PX, Y\rangle$

It is evident that $P(H)$ is the same for both the structures $(\Phi, g)$ and $(P, g)$, hence to construct the canonical connection $\overline{\nabla}$ (tensor field $h$) of the pair $(\Phi, g)$ we can take the canonical connection (see (5.3)) of the pair $(P, g)$, where $P, \Phi, g$ are connected by (5.17).

**$6^0$.** As we have seen in *Chapter 4* the Riemannian almost product structure $(P, g){:}T(M)=T^1(M)\oplus T^2(M)$ is defined on a R.l.r. $\sigma$-m. $M$, where $T^1$ denotes the distribution of mirrors. Of course, on an arbitrary R.l.r. $\sigma$-m. the connection $\tilde{\nabla}$ defined by (4.1) and $\overline{\nabla}$ defined by (5.3) are not the same.

If $M$ is a R.l.r. $\sigma$-m.o.2 (reflexion space), see *Definitions 4.3, 4.8,* then $S^2=I$ on $M$ and $S=P$.

Using (4.3) we obtain

$$\tilde{\nabla}_X Y = \frac{1}{2}(\nabla_X S^2 Y + S\nabla_X SY) = \frac{1}{2}(\nabla_X Y + P\nabla_X PY)$$

$$= \nabla_X Y - \frac{1}{2}\nabla_X(P)PY = \overline{\nabla}_X Y,\ \ X,Y\in \mathfrak{X}(M).$$

So, the canonical connection of a Riemannian reflection space coincides with that of the Riemannian almost product strructure defined by (5.3). It follows from *Proposition 4.5, 4.52* and *Corollary 4.53* that

$$h_X = 0,\ X \in T^1;\ \overline{\nabla}_X h = \overline{\nabla}_X \overline{T} = \overline{\nabla}_X \overline{R} = 0,\ X \in T^2;\ P(\overline{T}) = \overline{T},\ P(\overline{R}) = \overline{R}$$

on a R.l.r. $\sigma$-m.o.2.

Thus, this conditions are necessary for a Riemannian almost product structure $(P, g)$ on a manifold $M$ to be the structure defined by that of a Riemannian reflection space.

## §2. ALMOST HERMITIAN MANIFOLDS

**1⁰**. A tensor field $J$ of type $(1,1)$ such that $J^2 = -I$ is called an almost complex structure on a manifold $M$. A manifold $M$ with a fixed almost complex structure $J$ is called an almost complex manifold.

Every almost complex manifold is of even dimensional.We refer to [73], [47], [34] for detailed information.

For any Riemannian metric $\tilde{g}$ on $M$ a new Riemannian metric $g$ is defined by the formula

$$g(X,Y) = \tilde{g}(X,Y) + \tilde{g}(JX,JY),$$

where $X,Y \in \mathcal{X}(M)$. For every such a metric $g = <\,,\,>$ we have

$$<JX, JY> = <X,Y>,$$

for any vector fields $X$ and $Y$ on $M$.

Let $\nabla$ be the Riemannian connection of such a fixed metric $g = <\,,\,>$, then one can define a connection $\overline{\nabla}$ on $M$ by

$$(5.18)\quad \overline{\nabla}_X Y = \frac{1}{2}(\nabla_X Y - J\nabla_X JY) = \nabla_X Y + \frac{1}{2}\nabla_X(J)JY$$

where $X,Y \in \mathcal{X}(M)$.

Further, we obtain

$$2\overline{\nabla}_X(J)Y = 2(\overline{\nabla}_X JY - J\overline{\nabla}_X Y) = \nabla_X JY + J\nabla_X Y - J\nabla_X Y - \nabla_X JY = 0$$

and

$$\begin{aligned} 2(<\overline{\nabla}_X Y,Z> + <Y,\overline{\nabla}_X Z>) = {} & <\nabla_X Y,Z> + <Y,\nabla_X Z> - <J\nabla_X JY,Z> \\ & - <Y,J\nabla_X JZ> = X<Y,Z> + <\nabla_X JY,JZ> + <JY,\nabla_X JZ> \\ & = X<Y,Z> + X<JY,JZ> = 2X<Y,Z> \end{aligned}$$

for any $X,Y,Z \in \mathcal{X}(M)$, i.e., $\overline{\nabla} g = 0$.

According to *Theorem 1.2*, the tensor field $J$ is $O$-deformable and the Riemannian metric $g$ is an associated one to the corresponding $G$-structure. Later on, we shall consider only associated metrics of this form.

For a fixed affinor $J$ a set of all the associated metrics is infinite dimensional.

A fixed pair *(J, g)* is called an almost Hermitian structure (a.H.s.) and *M* is called an almost Hermitian manifold.

We can consider a set *P(H)* of all the orthonormal frames over *M* such that for every $u \in P(H)$ the tensor field *J* has the following matrix

$$(J) = \begin{bmatrix} 0 & E_n \\ -E_n & 0 \end{bmatrix},$$

where *dim M= 2n*. The set *P(H)* is a *G*-structure with the structure group

$$H = U(n) = \left\{ \begin{bmatrix} A & B \\ -B & A \end{bmatrix} : A, B \in O(n) \right\}$$

We consider now the canonical connection $\overline{\nabla}$ and the second fundamental tensor field h of the structure *(P(H),g)*.

**THEOREM 5.14.** The canonical connection $\overline{\nabla}$ of the *G*-structure *(P(H),g)* corresponding to an a.H.s. *(J, g)* is defined by (5.18). For $X, Y \in \mathfrak{X}(M)$ we have

(5.19) $$h_X Y = -\frac{1}{2} \nabla_X (J) JY = \frac{1}{2} (\nabla_X Y + J \nabla_X JY)$$

where $\nabla$ is the Riemannian connection of the Riemannian metric *g*.

**Proof.** It is clear that the Lie algebra $\underline{\mathbf{h}}$ of the structure group *H* of *P(H)* has the following form

$$\underline{\mathbf{h}} = \{ \tilde{x} = \begin{bmatrix} A & B \\ -B & A \end{bmatrix} : \tilde{x} \in \underline{\mathbf{o}} \}.$$

It is evident that $\underline{\mathbf{o}} = \underline{\mathbf{h}} \oplus \underline{\mathbf{m}}$, where

$$\underline{\mathbf{m}} = \{ \tilde{x} = \begin{bmatrix} C & D \\ D & -C \end{bmatrix} : \tilde{x} \in \underline{\mathbf{o}} \},$$

and

$$tr(\begin{bmatrix} A & B \\ -B & A \end{bmatrix} \cdot \begin{bmatrix} C & D \\ D & -C \end{bmatrix}^T ) = 0$$

that is $\underline{\boldsymbol{m}} = \underline{\boldsymbol{h}}^{\perp}$,with respect to the Killing form . For any $\omega \in \underline{\boldsymbol{o}}$ we can obtain the natural decomposition $\omega = \omega_{|\underline{h}} + \omega_{|\underline{m}}$ by the formula

$$\omega = 1/2(\omega - j\omega j) + 1/2(\omega + j\omega j),$$

where

$$j = \begin{bmatrix} 0 & E_n \\ -E_n & 0 \end{bmatrix}$$

Let $\varphi$ be a cross section of $P(H)$ over some neighbourhood $U$ which assigns to each $x \in U$ the linear frame $((X_1)_x, ..., (X_{2n})_x)$ and $X, Y = \sum_k f^k X_k$ be vectors fields on $M$. Then from (1.6) and *Definition 1.9* it follows that

$$\begin{aligned}
\overline{\nabla}_{X_x} Y &= \varphi(x)\overline{\omega}(\varphi_* X_x)\varphi(x)^{-1} Y_x + \sum_k (Xf^k)(x)(X_k)_x \\
&= \frac{1}{2}\varphi(x)[\omega(\varphi_* X_x) - j\omega(\varphi_* X_x) j]\varphi(x)^{-1} Y + \sum_k (Xf^k)(x)(X_k)_x \\
&= \frac{1}{2}\nabla_{X_x} Y - \frac{1}{2} J\varphi(x)\omega(\varphi_* X)\varphi(x)^{-1} JY - \frac{1}{2}\sum_k J(Xf^k)(x)J(X_k)_x \\
&= \frac{1}{2}(\nabla_{X_x} Y - J\nabla_{X_x} JY) = \nabla_{X_x} Y + \frac{1}{2}\nabla_{X_x}(J)JY,
\end{aligned}$$

where $\varphi(x)$ is considered as the mapping of $R^n$ onto $T_x(M)$. It is obvious that $\varphi(x) \cdot j = J \cdot \varphi(x)$.

**QED.**

**$2^0$**. If $\overline{T}$ is the torsion tensor field of $\overline{\nabla}$, then for $X, Y \in X(M)$ we obtain

$$\begin{aligned}
\overline{T}_X Y &= \overline{\nabla}_X Y - \overline{\nabla}_Y X - [X,Y] = \frac{1}{2}(\nabla_X Y - \nabla_Y X - [X,Y] - J\nabla_X JY + J\nabla_Y JX - [X,Y]) \\
&= \frac{1}{2}(J\nabla_Y JX - J\nabla_X JY - [X,Y])
\end{aligned}$$

From *Theorem 1.5* it follows that an almost complex structure is integrable if and only if the Nijenhuis tensor field

$$N(J)(X,Y)=[JX, JY]-J[X,JY]-J[JX,Y]-[X, Y]$$

vanishes on $M$, where $X,Y\in \mathfrak{X}(M)$.

We have

$$\begin{aligned}-N(J)(X,Y) &= J\nabla_X JY - J\nabla_{JY}X + J\nabla_{JX}Y - J\nabla_Y JX - [JX,JY] + [X,Y]\\ &= (J\nabla_{JX}Y - J\nabla_{JY}X - [JX,JY]) - (J\nabla_Y JX - J\nabla_X JY - [X,Y])\\ &= 2(\bar{T}_{JX}JY - \bar{T}_X Y)\end{aligned}$$

Using (2.3) we have obtained that for $X,Y\in \mathfrak{X}(M)$

$$N(J)(X,Y) = 4(h_{JX}^{-}JY - h_X^{-}Y).$$

Let $\bar{R}$ and $R$ be the curvature tensor fields of the connections $\bar{\nabla}$, $\nabla$ resectively, then for $X,Y,Z\in \mathfrak{X}(M)$ we have

$$\bar{R}_{XY}Z = \frac{1}{2}R_{XY}Z - \frac{1}{2}JR_{XY}JZ + h_Y h_X Z - h_X h_Y Z$$

Really, it is easily to verify that

$$\begin{aligned}4h_Y h_X Z &= \nabla_Y\nabla_X Z + \nabla_Y J\nabla_X JZ + J\nabla_Y J\nabla_X Z + J\nabla_Y J^2\nabla_X JZ,\\ 4h_X h_Y Z &= \nabla_X\nabla_Y Z + \nabla_X J\nabla_Y JZ + J\nabla_X J\nabla_Y Z + J\nabla_X J^2\nabla_Y JZ\end{aligned}$$

and

$$\begin{aligned}4\bar{R}_{XY}Z &= (\nabla_X\nabla_Y Z - \nabla_X J\nabla_Y JZ - J\nabla_X J\nabla_Y Z + J\nabla_X J^2\nabla_Y JZ)\\ &\quad - (\nabla_Y\nabla_X Z - \nabla_Y J\nabla_X JZ - J\nabla_Y J\nabla_X Z + J\nabla_Y J^2\nabla_X JZ)\\ &\quad - 2(\nabla_{[X,Y]}Z - J\nabla_{[X,Y]}JZ)\\ &= 2R_{XY}Z - 2JR_{XY}JZ + 4h_Y h_X Z - 4h_X h_Y Z.\end{aligned}$$

$3^0$. **DEFINITION 5.1.** [73]. An almost Hermitian manifold $M$ with almost complex structure $J$ is called a Kaehlerian manifold if $\nabla J=0$ and a nearly Kaehlerian manifold if $(\nabla_X J)X=0$, where $\nabla$ is the Riemannian connection of $g$, $X \in \mathcal{X}(M)$.

**PROPOSITION 5.15.** Let $(J,g)$ be an almost Hermitian structure on a manifold $M$. $(J,g)$ is a Kaehlerian (nearly Kaehlerian) structure if and only if it is a particular (nearly particular) one.

**Proof.** From (5.19) we have that

$$\nabla_X(J)Y=2h_X JY, \quad X,Y\in\mathcal{X}(M),$$

hence $\nabla J=0$ if and only if $h=0$.

$$h_X X = 1/2(\nabla_X X + J\nabla_X JX)$$

and

$$J(\nabla_X (J)X)=J\nabla_X JX+\nabla_X X=2h_X X,$$

therefore $(\nabla_X J)X=0$ if and only if $h_X X=0$ for any $X\in\mathcal{X}(M)$.

**QED.**

This proposition implies that the notion of a nearly particular (particular) structure is a generalization of that of a nearly Kaehlerian (Kaehlerian) structure.

**PROPOSITION 5.16.** If $(J,g)$ is a nearly Kachlerian structure then $\overline{\nabla}_X\nabla_X Y=\nabla_X\overline{\nabla}_X Y$ for any $X,Y\in\mathcal{X}(M)$.

**Proof.** We can find the following identity in [34]

$$<\nabla^2_{XX}(J)Y,JZ>=<\nabla^2_{YZ}(J)X,JX>=-<\nabla_X(J)Y,\nabla_X(J)Z>,$$

$X,Y,Z\in\mathcal{X}(M)$

This equality implies that

$$<\nabla^2_{XX}(J)Y,JZ>=-4<h_X JY,h_X JZ>=4<h_X h_X JY,JZ>=2<h_X\nabla_X(J)Y,JZ>$$
$$=-<\nabla_X(J)J\nabla_X(J)Y,JZ>=<J\nabla^2_X(J)Y,JZ>.$$

We obtain that

$$\nabla^2_{XX}(J)Y = J\nabla^2_X(J)Y$$

and

$$\nabla^2_{XX}(J)Y = \nabla_X(\nabla_X JY - J\nabla_X Y) - (\nabla_X J\nabla_X Y - J\nabla_X \nabla_X Y)$$
$$= \nabla^2_{XX} JY - \nabla_X J\nabla_X Y - \nabla_X J\nabla_X Y + J\nabla^2_{XX} Y,$$

$$J\nabla^2_X(J)Y = J(\nabla_X J\nabla_X JY - J\nabla^2_{XX} JY - \nabla_X J^2\nabla_X Y + J\nabla_X J\nabla_X Y)$$
$$= J\nabla_X J\nabla_X JY + \nabla^2_{XX} JY + J\nabla^2_{XX} Y - \nabla_X J\nabla_X Y.$$

Hence

$$-\nabla_X J\nabla_X Y = J\nabla_X J\nabla_X JY,$$

or

$$J\nabla_X J\nabla_X Y = \nabla_X J\nabla_X JY.$$

Thus

$$2\nabla_X \overline{\nabla}_X Y = \nabla_X \nabla_X Y - \nabla_X J\nabla_X JY, \quad 2\overline{\nabla}_X \nabla_X Y = \nabla_X \nabla_X Y - J\nabla_X J\nabla_X Y$$

and we have got that $\nabla_X \overline{\nabla}_X Y = \overline{\nabla}_X \nabla_X Y$.

**QED.**

One can find a more detailed information about the nearly Kaehlerian manifolds in [73], [32] , [34]. There are many results about Kaehlerian manifolds. A bibliography about such structures can be found in [73], [45], [47].

**4$^0$**. Using the results in [35] we consider now a classification of almost Hermitian structures with respect to the tensor field *h*.

**LEMMA 5.17.** If $h$ is the second fundamental tensor field of an almost Hermitian structure $(J,g)$ then for any $X,Y,Z \in \mathcal{X}(M)$ we have

(5.20) $$h_{XYZ} = -h_{XZY} = -h_{XJYJZ}.$$

**Proof.** The first equality follows from (2.2). To prove the second one we use (5.19)

$$2h_{XJYJZ} = <\nabla_X JY - J\nabla_X Y, JZ> = -<J\,\nabla_X (J)Y, Z>$$

Differentiating the identity $J^2=-I$ we obtain

$$\nabla_X (J)JY + J\nabla_X (J)Y = 0$$

and

$$-<J\nabla_X (J)Y, Z> \;=\; <\nabla_X (J)JY, Z> = -2h_{XYZ}.$$

**QED.**

**LEMMA 5.18.** If $\Phi(X,Y)=<JX,Y>$ for any $X,Y\in \mathcal{X}(M)$, then

(5.21) $(\nabla_X \Phi)(Y,Z)=2h_{XYJZ}$ .

**Proof.** Since $\nabla$ is a metric connection, then

$$X<Y,Z> = <\nabla_X Y, Z> + <Y, \nabla_X Z>$$

and

$$(\nabla_X\Phi)(Y,Z) = X<JY,Z> - <J\nabla_X Y, Z> - <JY, \nabla_X Z> = <\nabla_X JY, Z> - <J\nabla_X Y, Z>$$
$$= <\nabla_X (J)Y, Z> = 2h_{XJYZ} = 2h_{XYJZ}.$$

from (5.19) and (5.2O) .

**QED.**

We can also consider a notion of the codifferential $\delta T$ of a tensor field $T$ (see [73]). For example, let $T$ be a tensor field of type $(0,3)$. Then T is the tensor field of type $(0,2)$ defined by

$$\delta T = -g^{tj}\nabla_t T_{jik} = -\nabla^j T_{jik},$$

where we have put $\nabla^j = g^{tj}\nabla_t$ and $T_{jik}$ are local components of $T$. If $\{E_i, JE_i\}$, $i=1,...,n$ is a local orthonormal basis, defined on an open neighbourhood, then codifferential of $\Phi$ is computed by the formula

(5.22) $\delta\Phi(X)=-\sum_{i=1}^{n}[(\nabla_{E_i}\Phi)(E_i,X)+(\nabla_{JE_i}\Phi)(JE_i,X)]$, $X\in\mathfrak{X}(M)$.

Also, for $h$ let $\beta=c_{12}(h)$ be defined by

(5.23) $\beta(X)=c_{12}(h)X=\sum_{i=1}^{n}[h_{E_iE_iX}+h_{JE_iJE_iX}]$, $X\in\mathfrak{X}(M)$.

From (5.21), (5.22), (5.23) it follows that

(5.24) $\delta\Phi(X)=-2\sum_{i=1}^{n}[h_{E_iE_iJX}+h_{JE_iJE_iJX}]=-2\beta(JX)$,

$\delta\Phi(JX)=2\beta(X)$.

Let $p$ be a fixed point and $T=T_p(M)$. We consider a vector subspace $\overline{\boldsymbol{T}}(T)$ of $\overset{3}{\otimes}T^*$

(5.24) $\overline{T}(T)=\{h\in\overset{3}{\otimes}T^*:h_{XYZ}=-h_{XZY}=-h_{XJYJZ},\ X,Y,Z\in T\}$

and define four subspaces of $\overline{\boldsymbol{T}}(T)$ as follows:

(5.25) $\overline{\boldsymbol{T}}_1=\{h\in\overline{\boldsymbol{T}}:h_{XXZ}=0,\ X,Z\in T\}$,

$\overline{\boldsymbol{T}}_2=\{h\in\overline{\boldsymbol{T}}:\boldsymbol{\sigma}h_{XYZ}=0,\ X,Y,Z\in T\}$,

$\overline{\boldsymbol{T}}_3=\{h\in\overline{\boldsymbol{T}}:h_{XYZ}-h_{JXJYJZ}=\beta(Z)=0,\ X,Y,Z\in T\}$,

$\overline{T}_4=\{h\in\overline{T}:h_{XYZ}=\frac{1}{2(n-1)}[<X,Y>\beta(Z)-<X,Z>\beta(Y)]$
$-<X,JY>\beta(JZ)+<X,JZ>\beta(JY),\ X,Y,Z\in T\}$.

We consider now a decomposition of the space $\overline{\boldsymbol{T}}(T)$ on irreducible components under the action of $U(n)$, where the action is defined by (3.2) and the inner product by (3.1)

**THEOREM 5.19** [35]. We have $\overline{\boldsymbol{T}}=\overline{\boldsymbol{T}}_1\oplus\overline{\boldsymbol{T}}_2\oplus\overline{\boldsymbol{T}}_3\oplus\overline{\boldsymbol{T}}_4$. This direct sum is orthogonal, and it is preserved under the induced representation of $U(n)$ on $\overline{\boldsymbol{T}}$ . The induced representation of $U(n)$ on $\overline{\boldsymbol{T}}$ is irreducible.
For $n=1$, $\overline{\boldsymbol{T}}=\{0\}$; for $n=2$, $\overline{\boldsymbol{T}}_1=\overline{\boldsymbol{T}}_3=\{0\}$ and $\overline{\boldsymbol{T}}=\overline{\boldsymbol{T}}_2\oplus\overline{\boldsymbol{T}}_4$. For $n=2$, $\overline{\boldsymbol{T}}_2$ and $\overline{\boldsymbol{T}}_4$ are nontrivial, and for $n\geq 3$ all of the $\overline{\boldsymbol{T}}_i$ are nontrivial.

**REMARK.** If $\dim T=2n$, then $\dim\overline{\boldsymbol{T}} = \dim Hom(R^{2n}, \underline{\mathbf{m}})=2n(n^2-n)$ and $\dim\overline{\boldsymbol{T}}_1=1/3n(n-1)(n-2)$, $\dim\overline{\boldsymbol{T}}_2=2/3n(n-1)(n+1)$, $\dim\overline{\boldsymbol{T}}_3=n(n+1)(n-2)$ (for $n\geq 2$), and $\dim\overline{\boldsymbol{T}}_4=2n$.

Let $\overline{\boldsymbol{T}}_\alpha$ denote a direct sum of all the classes $\overline{\boldsymbol{T}}_i$ such that $i\in\alpha\subset\{1,2,3,4\}$. We can form $2^4=16$ invariant subspaces of $\overline{\boldsymbol{T}}$ (including $\{0\}$ and $\overline{\boldsymbol{T}}$ ).

**DEFINITION 5.2.** We call that the structure $(J,g)$ has a type $\overline{\boldsymbol{T}}_\alpha$ or $(M, J, g)$ belongs to the class $\overline{\boldsymbol{T}}_\alpha$, if $h_p\in\overline{\boldsymbol{T}}_\alpha(T_p(M))$, for every $p\in M$, where $h$ is the second fundamental tensor field of $(J,g)$ on $M$.

**REMARK.** One can replace $\overline{\boldsymbol{T}}_\alpha$ and $\boldsymbol{U}_\alpha$, $h$ and $\nabla\Phi$, $\beta$ and $\delta\Phi$ in (5.24), (5.25), *Theorem 5.19*, *Definition 5.2* and have got a classification of A.Gray and L.M.Hervella, [35], given by

**Table 5.1**

| Class | Defining condition |
|---|---|
| K | $\nabla\Phi=0$ |
| $\boldsymbol{U}_1=\boldsymbol{NK}$ | $\nabla_X(\Phi)(X,Y)=0$ (or $3\nabla\Phi = d\Phi$) |
| $\boldsymbol{U}_2=\boldsymbol{AK}$ | $d\Phi=0$ |
| $\boldsymbol{U}_3=\boldsymbol{SK}\cap\boldsymbol{H}$ | $\delta\Phi=N(J)=0$ (or $\nabla_X(\Phi)(Y,Z)-\nabla_{JX}(\Phi)(JY,Z)=\delta\Phi=0$ ) |
| $\boldsymbol{U}_4$ | $\nabla_X(\Phi)(Y,Z)=\frac{-1}{2(n-1)}\{<X,Y>\delta\Phi(Z)-<X,Z>\delta\Phi(Y)$ $-<X,JY>\delta\Phi(JZ)+<X,JZ>\delta\Phi(JY)\}$ |
| $\boldsymbol{U}_1\oplus\boldsymbol{U}_2=\boldsymbol{QK}$ | $\nabla_X(\Phi)(Y,Z)+\nabla_{JX}(\Phi)(JY,Z)=0$ |
| $\boldsymbol{U}_3\oplus\boldsymbol{U}_4=\boldsymbol{H}$ | $N(J)=0$ (or $\nabla_X(\Phi)(Y,Z)-\nabla_{JX}(\Phi)(JY,Z)=0$ ) |
| $\boldsymbol{U}_1\oplus\boldsymbol{U}_3$ | $\nabla_X(\Phi)(X,Y)-\nabla_{JX}(\Phi)(JX,Y)=\delta\Phi=0$ |
| $\boldsymbol{U}_2\oplus\boldsymbol{U}_4$ | $d\Phi=\Phi\wedge\Theta$ (or $\sigma\{\nabla_X(\Phi)(Y,Z)-\frac{1}{n-1}\Phi(X,Y)\delta\Phi(JZ)\}=0$) |
| $\boldsymbol{U}_1\oplus\boldsymbol{U}_4$ | $\nabla_X(\Phi)(X,Y)=\frac{-1}{2(n-1)}\{\lVert X\rVert^2\delta\Phi(Y)-<X,Y>\delta\Phi(X) - <JX,Z>\delta\Phi(JX)\}$ |
| $\boldsymbol{U}_2\oplus\boldsymbol{U}_3$ | $\sigma\{\nabla_X(\Phi)(Y,Z)-\nabla_{JX}(\Phi)(JY,Z)\}=\delta\Phi=0$ |
| $\boldsymbol{U}_1\oplus\boldsymbol{U}_2\oplus\boldsymbol{U}_2=\boldsymbol{SK}$ | $\delta\Phi=0$ |
| $\boldsymbol{U}_1\oplus\boldsymbol{U}_2\oplus\boldsymbol{U}_4$ | $\nabla_X(\Phi)(Y,Z)+\nabla_{JX}(\Phi)(JY,Z)=\frac{-1}{h-1}\{<X,Y>\delta\Phi(Z) - <X,Z>\delta\Phi(Y)-<X,JY>\delta\Phi(JZ)+<X,JZ>\delta\Phi(JY)\}$ |
| $\boldsymbol{U}_1\oplus\boldsymbol{U}_2\oplus\boldsymbol{U}_4=G_1$ | $\nabla_X(\Phi)(X,Y)-\nabla_{JX}(\Phi)(JX,Y)=0$ (or $<N(J)(X,Y),X>=0$ ) |
| $\boldsymbol{U}_2\oplus\boldsymbol{U}_3\oplus\boldsymbol{U}_4=G_2$ | $\sigma\{\nabla_X(\Phi)(Y,Z)-\nabla_{JX}(\Phi)(JY,Z)\}=0$ (or $\sigma<N(J)(X,Y),JZ>=0$) |
| U | No condition |

**5⁰**. Our next aim is to show that the both classifications are the same. We can define a mapping $C{:}\boldsymbol{U}\to\overline{\boldsymbol{T}}:(\nabla\Phi)_p\to h_p$.

If follows from (5.21) that

(5.26) $(\nabla_X\Phi)(Y,Z)=2h_{XYJZ}$ and $h_{XYZ}=-1/2(\nabla_X\Phi)(Y,JZ)$,

therefore $C$ is a one –to-one correspondence. From (5.26) we see that $C(\nabla\Phi_1+\nabla\Phi_2)=C(\nabla\Phi_1)+C(\nabla\Phi_2)$, $C(\alpha\nabla\Phi)=\alpha C(\nabla\Phi)$, hence $C$ is an isomorphism.

**THEOREM 5.2O.** We have $C(\boldsymbol{U}_i)=\overline{\boldsymbol{T}}_i$, where $i=1,2,3,4$.

**Proof.** 1) $C(\boldsymbol{U}_1)=\overline{\boldsymbol{T}}_1$ because $\nabla_X(\Phi)(X,Y)=2h_{XJXY}=2h_{XXJY}$.

2) For $\boldsymbol{U}_2\subset\boldsymbol{U}_1\oplus\boldsymbol{U}_2$ it follows, [35], that $(\nabla_X\Phi)(Y,Z)=-(\nabla_{JX}\Phi)(JY,Z)$ or $2h_{XYJZ}=-2h_{JXJYJZ}=2h_{JXYZ}$, hence $h_{XYJZ}=h_{XJYZ}=h_{JXYZ}$.

$\boldsymbol{U}_2: d\Phi=O$, or $(\nabla_X\Phi)(Y,Z)+(\nabla_Y\Phi)(Z,X)+(\nabla_Z\Phi)(X,Z)$

$$=2(h_{XYJZ}+h_{YZJX}+h_{ZXJY})=2\sigma h_{XYJZ}=0$$

Thus $\sigma h_{XYZ}=0$ and $C(\boldsymbol{U}_2)\subset\overline{T}_2$. But $\boldsymbol{U}_2$ and $\overline{\boldsymbol{T}}_2$ are equal dimensional, therefore $C(\boldsymbol{U}_2)=\overline{\boldsymbol{T}}_2$.

3) We have $(\nabla_X\Phi)(Y,Z)-(\nabla_{JX}\Phi)(JY,Z)=2(h_{XYJZ}-h_{JXJYZ})$ and $\delta\Phi(Z)=-2\beta(JZ)$, hence $C(\boldsymbol{U}_3)=\overline{\boldsymbol{T}}_3$.

4) Using the defining conditions for $\boldsymbol{U}_4$ and $\overline{\boldsymbol{T}}_4$ we obtain

$$\begin{aligned}(\nabla_X\Phi)(Y,Z)=-1/2(n-1)[&<X,Y>\delta\Phi(Z)-<X,Z>\delta\Phi(Y)-<X,JY>\delta\Phi(JZ)\\ &+<X,JZ>\delta\Phi(JY)]=2/2(n-1)[<X,Y>\beta(JZ)-<X,Z>\beta(JY)\\ &+<X,JY>\beta(Z)-<X,JZ>\beta(Y)]=2h_{XYJZ}\end{aligned}$$

and $C(U_4)=\overline{\boldsymbol{T}}_4$.

**QED.**

Thus, the classifications coincide and the classification in *Table 5.1* from [35] can be rewritten in terms of the tensor field *h*. Let $\dim M\geq 6$, then we have got

Table 5.2

| Class | Defining condition |
|---|---|
| K | $h=0$ |
| $\overline{\boldsymbol{T}}_1\equiv\boldsymbol{U}_1=\boldsymbol{NK}$ | $h_XX=0$ |
| $\overline{\boldsymbol{T}}_2\equiv\boldsymbol{U}_2=\boldsymbol{AK}$ | $\sigma h_{XYZ}=0$ |
| $\overline{T}_3 \equiv U_3 = SK \cap H$ | $h_{XYZ} - h_{JXJYZ} = \beta(Z) = 0$ |
| $\overline{T}_4 \equiv U_4$ | $h_{XYZ} = 1/2(n-1)[< X,Y > \beta(Z) - < X,Z > \beta(Y) - < X,JY > \beta(JZ) + < X,JZ > \beta(JY)]$ |
| $U_1 \oplus U_2 = QK$ | $h_{XYJZ} = h_{JXYZ}$ |
| $U_3 \oplus U_4 = H$ | $N(J)=0$ or $h_{XYJZ} = -h_{JXYZ}$ |
| $U_1 \oplus U_3$ | $h_{XXY} - h_{JXJXY} = \beta(Z) = 0$ |
| $U_2 \oplus U_4$ | $\sigma[h_{XYJZ} - 1/(n-1) < JX,Y > \beta(Z)] = 0$ |
| $U_1 \oplus U_4$ | $h_{XXY} = -1/2(n-1)[< X,Y > \beta(X) - \|X\|^2 \beta(Y) - < X,JY > \beta(JX)]$ |
| $U_2 \oplus U_3$ | $\sigma[h_{XYJZ} + h_{JXYZ}] = \beta(Z) = 0$ |
| $U_1 \oplus U_2 \oplus U_2 = SK$ | $\beta = 0$ |
| $U_1 \oplus U_2 \oplus U_4$ | $h_{XYJZ} - h_{JXYZ} = 1/(n-1)[< X,Y > \beta(JZ) - < X,Z > \beta(JY) + < X,JY > \beta(Z) - < X,JZ > \beta(Y)]$ |
| $U_1 \oplus U_3 \oplus U_4$ | $h_{XJXY} + h_{JXXY} = 0$ |
| $U_2 \oplus U_3 \oplus U_4$ | $\sigma[h_{XYJZ} + h_{JXYZ}] = 0$ |
| $U$ | No condition |

**REMARK.** In fact, for every class of the classification of A.Gray and L.M. Hervella the canonical connection $\overline{\nabla} = \nabla - h$ has been constructed and to study these classes of almost Hermitian manifolds the torsion tensor field $\overline{T}$ ,the curvature tensor field $\overline{R}$ and other usual characteristics of $\overline{\nabla}$ can be applied.

## *§3. SOME EXAMPLES OF ALMOST HERMITIAN MANIFOLDS*

**$1^0$**. Let $M'$ be a submanifold of an almost Hermitian manifold $(M,J,g)$ which is invariant with respect to $J$, that is, for every $X\in T(M')$, $JX\in T(M')$ too.

We call $M'$ strongly invariant if $h_X Y\in T(M')$ for any $X,Y\in \mathcal{X}(M')$. It is easy to see from (2.14),(2.15) that $M'$ is strongly invariant if and only if $\alpha(X,Y)=\overline{\alpha}(X,Y)$, where $X,Y\in \mathcal{X}(M')$.

**PROPOSITION 5.21.** $\alpha(X,Y)=\overline{\alpha}(X,Y)$ if and only if $J\alpha(X,Y)=-\alpha(X,JY)$ for $X,Y\in \mathcal{X}(M')$.

**Proof.** If $p\in M'$ then $T_p(M)=T_p(M')\oplus T_p(M')^{\perp}$ and $T_p(M')^{\perp}$ is invariant under $J$ too.

Let $\pi$ be the projection on $T_p(M')^{\perp}$, then we have

$$\nabla_X Y=\nabla'_X Y+\alpha(X,Y)=\nabla'_X Y+\pi\nabla_X Y,\ \overline{\nabla}_X Y=\overline{\nabla}'_X Y+\overline{\alpha}(X,Y)=\overline{\nabla}'_X Y+\pi\overline{\nabla}_X Y,$$

hence

$$\alpha(X,Y)-\overline{\alpha}(X,Y)=\pi h_X Y=1/2(\pi\nabla_X Y+J\pi\nabla_X JY)=1/2(\alpha(X,Y)+J\alpha(X,JY)).$$

**QED.**

**PROPOSITION 5.22.** 1) Every totally geodesic with respect to $\nabla$ invariant submanifold $M'$ of $(M,J,g)$ is strongly invariant. 2) Each autoparallel with respect to $\overline{\nabla}$ strongly invariant submanifold $M'$ is totally geodesic with respect to $\nabla$.

**Proof.** 1) $M'$ is an autoparallel submanifold of $M$, that is, $\nabla_X Y\in T(M')$ for $X,Y\in \mathcal{X}(M')$, [47], and $JY\in T(M')$. It follows from (5.19) that $h_X Y=1/2(\nabla_X Y-J\nabla_X JY)\in T(M')$.

2) It is evident from the formula $\nabla_X Y=\overline{\nabla}_X Y+h_X Y$.

**QED.**

**THEOREM 5.23.** Let $(M,J,g)$ belong to a class $\boldsymbol{U}_\alpha$ from the *table 5.1*, where $\alpha=1,2,4,(1,2),(3,4)(2,4),(1,4),(1,2,4),(1,3,4),(2,3,4)$ and $M'$ be a strongly invariant submanifold of $M$. Then $(M',J,g)$ belongs to a subclass of $\boldsymbol{U}_\alpha$.

**Proof.** Let $p\in M'$, $T=T_p(M)$, $T'=T_p(M')$ and let

$$f:\overline{\boldsymbol{T}}(T)\to\overline{\boldsymbol{T}}(T'):h\mapsto h_{|T'}$$

be the linear mapping of restriction. It is correct because $h_X Y\in T'$ for $X,Y\in T'$.

1) For $\alpha=1,2,(1,2),(3,4),(1,3,4),(2,3,4)$ a proof follows from the *table 5.1*.

2) $\alpha=4$. We take an orthonormal basis

$$E_1,\dots,E_k,\ E_{k+1},\dots E_n,\ JE_1,\dots,JE_k,\ JE_{k+1},\dots,JE_n$$

of $T$ in such a vay that $E_1,\dots,E_k,\ JE_1,\dots,JE_k$ is a basis of $T'$. If $h\in\overline{\boldsymbol{T}}_4(T)$, then

$$h_{JE_iJE_iE_l}=h_{E_iE_iE_l}=\begin{cases}0 & ,\ l=i\\ 1/2(n-1)\beta(E_l), & l\neq i\end{cases},\ i,l=1,\dots,n.$$

It is the same for $h_{E_iE_iJE_l}$. From (5.23) we have

$$\beta'(E_l)=\sum_{i=1}^{k}(h_{E_iE_iE_l}+h_{JE_iJE_iE_l})=2(k-1)/2(n-1)\beta(E_l)$$

$$\beta'(JE_l)=(k-1)/(n-1)\beta(JE_l),\quad l=1,\dots,k$$

and $\beta(X)=(n-1)/(k-1)\beta'(X)$ for $X\in T'$, therefore

$$\begin{aligned}h_{XYZ}&=1/2(n-1)[<X,Y>\beta(Z)-<X,Z>\beta(Y)-<X,JY>\beta(JZ)+<X,JZ>\beta(JY)]\\&=1/2(k-1)[<X,Y>\beta'(Z)-<X,Y>\beta'(Y)-<X,JY>\beta'(JZ)+<X,JZ>\beta'(JY)],\end{aligned}$$

where $X,Y,Z\in T'$.

3) $\alpha=(2,4),(1,4),(1,2,4)$. It is clear that

$$f(\overline{\boldsymbol{T}}_i(T)\oplus\overline{\boldsymbol{T}}_j(T))\subset f(\overline{\boldsymbol{T}}_i(T))\oplus f(\overline{\boldsymbol{T}}_j(T))\subset\overline{\boldsymbol{T}}_i(T')\oplus\overline{\boldsymbol{T}}_j(T')$$

**QED.**

**REMARK.** If we have such a basis as in 2) of *Theorem 5.23*, then

$$\beta(X)=\sum_{i=1}^{k}(h_{E_iE_iX}+h_{JE_iJE_iX})+\sum_{i=k+1}^{n}(h_{E_iE_iX}+h_{JE_iJE_iX})=\beta'(X)+\beta'^{\perp}(X),$$

where $X\in T'$ and the situation is not evident for $\alpha=3,(1,3),(2,3),(1,2,3)$. Really, if $\beta=0$ on $T$ it does not imply that $\beta'=0$ on $T$.

**$2^0$**. Let $(J,g)$ and $(J,g')$ be locally conformally related, that is, $g'=e^{2\rho}g$, then it

follows from (5.5) that

$$\nabla'_X Y = \nabla_X Y + X(\rho)Y + Y(\rho)X - <X,Y> \; grad \; \rho$$

and we have from (5.19) that

$$h'_X Y = 1/2(\nabla'_X Y + J\nabla'_X JY) = h_X Y + 1/2(X(\rho)Y + Y(\rho)X - <X,Y> grad\,\rho)$$
$$+ JX(\rho)JY + J(JY)(\rho)X - <X,JY> J\; grad\,\rho),$$

(5.27) $$h'_{XYZ} = e^{2\rho}[h_{XYZ} + 1/2(Y(\rho) <X,Y> - <X,Y> Z(\rho) + (JY)(\rho) <JX,Z>$$
$$+ <X,JY> (JZ)(\rho))].$$

If $(J,g')$ is a Kaehlerian structure and $\beta=(n-1)d\rho$, then

(5.28) $h_{XYZ}=1/2(n-1)[<X,Y>\beta(Z)-<X,Z>\beta(Y)-<X,JY>\beta(JZ)+<X,JZ>\beta(JY)]$,

therefore $h \in \overline{\boldsymbol{T}}_4$.

**3⁰**. Using [35] we illustrate the structures of the different classes.

1) The class of nearly Kachlerian manifolds $\boldsymbol{U}_1$=***NK***, see [73], [32], [34]. The most well-known example in this class is the sphere $S^6$, [33]. $S^6$ belongs to the class $\boldsymbol{U}_1$ but not to the ***K***.
2) The class of almost Kaehlerian manifolds $\boldsymbol{U}_2$=***AK***. The tangent bundle *T(M)* of a Riemannian manifold always has a naturally defined complex structure and a metric such that $T(M)\in\boldsymbol{U}_2$. If M is not flat, then $T(M)\notin\boldsymbol{K}$. An example of a compact 4-dimensional manifold in $\boldsymbol{U}_2$ was given in [65], which has no Kaehlerian metric.
3) The class of Hermitian semi-Kaehlerian manifolds $\boldsymbol{U}_3=\boldsymbol{H}\cap\boldsymbol{SK}$. Any complex parallelizable manifold is in $\boldsymbol{U}_3$.
4) As we have got above $\boldsymbol{U}_4$ is a class which contains locally conformal Kaehlerian manifolds.

One can find more detaited information about the following classes in literature.

$\boldsymbol{U}_1\oplus\boldsymbol{U}_2=\boldsymbol{QK}$ is the class of quasi-Kaehlerian manifolds.

$\boldsymbol{U}_3\oplus\boldsymbol{U}_4=\boldsymbol{H}$ is the class of Hermitian manifolds.

$\boldsymbol{U}_2\oplus\boldsymbol{U}_4$ is a class which contains locally conformal almost Kaehlerian manifolds.

$\boldsymbol{U}_1\oplus\boldsymbol{U}_2\oplus\boldsymbol{U}_3=\boldsymbol{SK}$ is the class of semi-Kaehlerian manifolds.

$\boldsymbol{U}_1\oplus\boldsymbol{U}_3\oplus\boldsymbol{U}_4$ and $\boldsymbol{U}_2\oplus\boldsymbol{U}_3\oplus\boldsymbol{U}_4$ are classes studied in [37].

**4$^0$**. Let *(J,g)* be a quasi homogeneous structure, i. e., $\overline{\nabla}h=0$.

**THEOREM 5.24.** Let *(J,g)* be a quasi homogeneous structure having a type $\overline{T}_\alpha$ for some point $p\in M$. Then this structure belongs to the class $\overline{T}_\alpha$ on *M*.

Proof is similar to that of *Theorem 3.3*.

We consider now a quasi homogeneous structure *(J,g)* having the type $\overline{T}_4\equiv U_4$. Let $\xi$ be such a vector field on *M* that $\beta(X)=<\xi,X>$ for any $X\in X(M)$ and let $L=[\xi\wedge J\xi]$ be the 2-dimensional distribution defined by $\xi$, $J\xi$. We take $V=L^\perp$ and obtain that $T(M)=L\oplus V$, $\beta(J\xi)=0$, $\beta(Y)=0$ for each $Y\in V$.

**PROPOSITION 5.25.** $\overline{\nabla}h=0$ if and only if $\overline{\nabla}_W\xi=0$ for any $W\in X(M)$.

**Proof.** For an integral curve $\gamma(t)$ of the vector field *W* we can consider vector fields *X,Y,Z* which are parallel along $\gamma(t)$, i.e., $\overline{\nabla}_W X=\overline{\nabla}_W Y=\overline{\nabla}_W Z=0$.

Then it follows that

$$(\overline{\nabla}_W h)(X,Y,Z)=\overline{\nabla}_W h_{XYZ}=1/2(n-1)[<X,Y><\overline{\nabla}_W\xi,Z>-<X,Z><\overline{\nabla}_W\xi,Y>$$
$$-<X,JY><\overline{\nabla}_W\xi,JZ>+<X,JZ><\overline{\nabla}_W\xi,JY>]=0 \text{ if } \overline{\nabla}_W\xi=0$$

Conversely, let $\overline{\nabla}h=0$. $\overline{\nabla}_W JZ=J\overline{\nabla}_W Z=0$ and $[Z\wedge JZ]$, $[Z\wedge JZ]^\perp$ are invariant with respect to $\overline{\nabla}$, therefore if $X\in[Z\wedge JZ]_p^\perp$ for $p=\gamma(o)$, then $\overline{\tau}_t(X)\in[Z\wedge JZ]_{\gamma(t)}^\perp$, where $\overline{\tau}_t$ is the parallel translation with respect to $\overline{\nabla}$ along $\gamma(t)$. Thus,

$$(\overline{\nabla}_W h)(X,X,Z)=\overline{\nabla}_W h_{XXZ}=1/2(n-1)\|X\|^2<\overline{\nabla}_W\xi,Z>=0 \text{ for every } Z,$$

hence $\overline{\nabla}_W\xi=0$ and $\overline{\nabla}_W J\xi=J\overline{\nabla}_W\xi=0$ too.

**QED.**

If $\xi=\|\xi\|\zeta$, where $\|\zeta\|=1$, then $\overline{\nabla}_W\xi=\|\xi\|\overline{\nabla}_W\zeta+(W(\|\xi\|))\zeta=0$, therefore $\overline{\nabla}_W\zeta=0$ for any $W\in X(M)$ and $\|\xi\|=c$ on *M*.

The formula (5.28) can be rewritten in the following form

$$(5.29)\quad h_{XYZ}=1/2(n-1)[<X,Y><\xi,Z>-<X,Z><\xi,Y>-<X,JY><JZ,\xi>$$
$$+<X,JZ><JY,\xi>].$$

Using (5.29) it is easily to check that for any $X,Y,Z\in V$

(5.30) $h_{\xi\xi Z}=h_{\xi J\xi Z}=h_{\xi\xi J\xi}=h_{J\xi\xi\xi}=h_{J\xi\xi J\xi}=h_{\xi J\xi J\xi}=h_{\xi YZ}=h_{J\xi YZ}=h_{XYZ}=0,$

$$h_{XY\xi}=c^2/2(n-1)<X,Y>,\ \ h_{XYJ\xi}=c^2/2(n-1)<X,JY>.$$

From *Proposition 5.25* it follows that $L$ is invariant with respect to $\overline{\nabla}$, hence $V$ is invariant under too.

**PROPOSITION 5.26.** The distributions $L$ and $V$ are integrable, $[\xi, J\xi]=0$. Every maximal integral manifold of L or V is a Kaehlerian one.

**Proof.** For any $X,Y\in V$ we have

$$<[X,Y],\xi>=<\nabla_X Y-\nabla_Y X,\ \xi>=h_{XY\xi}-h_{YX\xi}=0,$$

$$<[\xi,J\xi],Z>=<\nabla_\xi J\xi-\nabla_{J\xi}\xi,Z>=h_{\xi J\xi Z}-h_{J\xi\xi Z}=0,$$

hence $V$ and $L$ are integrable,

$$<[\xi, J\xi],\xi>=h_{\xi J\xi\xi}-h_{J\xi\xi\xi}=0,$$

$$<[\xi, J\xi], J\xi>=h_{\xi J\xi J\xi}-h_{J\xi\xi J\xi}=0,$$

therefore $[\xi,J\xi]=0$.

Since $h_{XYZ}=0$, then from (5.19) it is clear that every maximal integral manifold of $V$ is Kaehlerian. For L this follows from (5.30).

**QED.**

**THEOREM 5.27.** For a quasi homogeneous structure $(J,g)$ having the type $\overline{T}_4$ the metric $g$ can be locally conformally changed into a Kaehlerian one.

**Proof.** The affinor $J$ is integrable because $\overline{T}_4\subset \boldsymbol{H}=\overline{T}_3\oplus\overline{T}_4$ and for every $p\in M$ there exists such a coordinate neighbourhood $U(x_1,\dots,x_{2n-2},t,v)$ of $p$ that

$$X_1=\frac{\partial}{\partial x_1},\dots,X_{n-1}=\frac{\partial}{\partial x_{n-1}},\ JX_1=\frac{\partial}{\partial x_n},\dots,JX_{n-1}=\frac{\partial}{\partial x_{2n-2}}$$

belong to $V$ and $\xi=\dfrac{\partial}{\partial t}$, $J\xi=\dfrac{\partial}{\partial v}\in L$. Thus, we obtain $dt(\xi)=1$, $dt(J\xi)=dt(X_i)=dt(JX_i)=0$, where $i=1,\dots,n-1$. If $\rho=c^2/n-1\cdot t$ on $U$, then

$\beta(Z)=<\xi ,Z>=(n-1)\ d\rho(Z)$.

Let $g'=e^{2\rho}g$ on $U$, then from (5.27), (5.28) it follows that $h'_{XYZ}=0$ for any vector fields on $U$ and $g'$ is a Kaehlerian metric on $U$.

**QED.**

**REMARK.** Using these propositions examples of quasi homogeneous structures having the type $\bar{T}_4$ can be constructed.

**$5^0$. DEFINITION 5.3.** [35]. Let $(M,J,g)$ be an almost Hermitian manifold.Then $\mu$ is the tensor field of type (2.1) defined by

(5.31)

$$<\mu(X,Y),Z>=\nabla_X(\Phi)(Y,Z)+1/2(n-1)[<X,Y>\delta\Phi(Z)-<X,Z>\delta\Phi(Y)$$
$$-<X,JY>\delta\Phi(JZ)+<X,JZ>\delta\Phi(JY)],$$

for $X,Y,Z\in \boldsymbol{X}(M)$.

The tensor field $\mu$ measures the failure of an almost Hermitian manifold to be conformally Kaehlerian.

**THEOREM.5.28.** [35]. Let $(M,J,g)$ and $(M,J,g')$ be locally conformally related almost Hermitian manifolds. Then the corresponding tensor field $\mu$ and $\mu'$ satisfy $\mu'=\mu$ and for any class $\boldsymbol{U}_\alpha$ given in *table 5.1* we have $U'_\alpha\subset U_4\oplus U_\alpha$. Thus $\boldsymbol{U}_\alpha=\boldsymbol{U}'_\alpha$ if and only if $\boldsymbol{U}_4\subset\boldsymbol{U}_\alpha$. Using *table 5.1* the defining relation for each of the conformally invariant classes can be rewritten in terms of $\mu$ and we have

$M\in\boldsymbol{U}_4$ if and only if $\mu=0$,

$M\in\boldsymbol{U}_1\oplus\boldsymbol{U}_4$ if and only if $\mu(X,X)=0$ for all $X\in\boldsymbol{X}(M)$,

$M\in\boldsymbol{U}_2\oplus\boldsymbol{U}_4$ if and only if $\sigma<\mu(X,Y),Z>=0$ for all $X,Y,Z\in\boldsymbol{X}(M)$,

$M\in\boldsymbol{U}_3\oplus\boldsymbol{U}_4$ if and only if $\mu(X,Y)-\mu(JX,JY)=0$ for all $X,Y\in\boldsymbol{X}(M)$,

$M\in\boldsymbol{U}_1\oplus\boldsymbol{U}_2\oplus\boldsymbol{U}_4$ if and only if $\mu(X,Y)+\mu(JX,JY)=0$ for all $X,Y\in\boldsymbol{X}(M)$,

$M\in\boldsymbol{U}_1\oplus\boldsymbol{U}_3\oplus\boldsymbol{U}_4$ if and only if $\mu(X,X)-\mu(JX,JX)=0$ for all $X\in\boldsymbol{X}(M)$,

$M\in\boldsymbol{U}_2\oplus\boldsymbol{U}_3\oplus\boldsymbol{U}_4$ if and only if $\sigma<\mu(X,Y)-\mu(JX,JY),Z>=0$

for all $X,Y,Z\in\boldsymbol{X}(M)$.

Let $M$ and $J$ be analytic, $\gamma(t)$ be such a curve on $M$ that

$$(5.32)\ \nabla_X \overline{\nabla}_X Y = \overline{\nabla}_X \nabla_X Y$$

where $X=\gamma'(t)$ and $Y$ is any vector field defined on some neighbourhood of $\gamma(t)$.

If $(J,g)$ is a quasi homogeneous structure and $\gamma(t)$ is a geodesic with respect to $\overline{\nabla}$, then (5.32) is fulfilled from *Lemma 5.8.*

*Proposition 5.16* implies (5.32) for nearly Kaehlerian structures. We can define *cos J, sin J* by the formulas

$$\cos J = \sum_{m=0}^{\infty}(-1)^m J^{2m}/(2m)!\,,\ \sin J = \sum_{m=0}^{\infty}(-1)^m J^{2m+1}/(2m+1)!$$

**THEOREM 5.29.** If is the parallel translation from $\gamma(o)$ to $\gamma(t)$ with respect to $\nabla$ along such a curve $\gamma(t)$ that (5.32) is fulfilled and $J(t)$ is the almost Hermitian structure in $\gamma(t)$, then

$$J(t)=\tau_t \circ [\cos(tJ'(o))J(o)+\sin(tJ'(o))]\circ \tau_t^{-1},$$

where $J'=\nabla_X J$.

Let a) $-\lambda^2$, b) $\lambda^2$ $(\lambda\neq 0)$ be an eigenvalue of $(J'(o))^2$ and let $X(o)$ be a corresponding eigenvector, $X(t)=\tau_t(X(o))$. Then

a) $J(t)X(t)=\cosh(\lambda t)\tau_t(J(o)X(o))+1/\lambda\,\sinh(\lambda t)\tau_t(J'(o)X(o))$,

b) $J(t)X(t)=\cos(\lambda t)\tau_t(J(o)X(o)) + 1/\lambda\,\sin(\lambda t)\tau_t(J'(o)X(o))$.

Proofs are similar to ones considered in *Propositions 5.8 – 5.13*.

Let $(M,\Phi)$ be an almost symplectic manifold, that is, a manifold with a 2-form $\Phi$ which has maximal rank. It is well-known, see [64], that $\Phi$ is defined by a $G$-structure $P(G)$, $G=Sp(n,R)$, $\dim M=2n$. We see from *Theorem 1.1* that there exists a reduction of $G$ to the maximal compact subgroup $H=U(n)$, which defines an associated metric $g=< , >$. The space of all such associated metrics is infinite dimensional and for a fixed metric $g=< , >$ the almost Hermitian structure $J$ is determined by

(5.33) $\Phi(X,Y)=<JX,Y>$.

It is clear that $P(H)$ is the same for both the structures $(\Phi,g)$ and $(J,g)$, therefore the canonical connection $\overline{\nabla}$ of the pair $(\Phi,g)$ coincides with that of $(J,g)$, where $J$, $\Phi$, $g$ are related by (5.33).

**6⁰**. If $M$ is a 3-s.l.R.s., see *Definition 4.9*, then $S^3=I$ and $S$ has only two eigenvalues $-\frac{1}{2}\pm\frac{\sqrt{3}}{2}i$. Since $I-S^3=(I-S)(S^2+S+I)$ and $(I-S)$ is nonsingular on $M$, then we obtain $S^2+S+I=0$. An almost-complex structure $J$ on $M$ is defined by the formula

$$(5.34)\quad J=\frac{1}{\sqrt{3}}(2S+I),\quad S=-\frac{1}{2}I+\frac{\sqrt{3}}{2}J$$

We have $J^2=1/3(4S^2+4S+4I-3I)=-I$ and for $X,Y\in \mathfrak{X}(M)$

$$<JX,JY>=1/3(4<SX,SY>+2<SX,Y>+2<X,SY>+<X,Y>)$$
$$=1/3(3<X,Y>+2<X,S^2Y>+2<X,SY>+2<X,Y>)=<X,Y>.$$

So, $(J,g)$ is an almost Hermitian structure on $M$.

We know, [48], that for a R.l.r.s-m. $(M,\{s_x\})$ there exists the unique connection $\tilde{\nabla}$ which is invariant under every $s_x$, $x\in M$, and such that $\tilde{\nabla}S=\tilde{\nabla}g=0$. $\tilde{\nabla}$ is defined by (4.1), where $\pi_2=I$, or (4.3).

**THEOREM 5.30** [48]. Let $(M,\{s_x\})$ be a R.l.r.s-m. and $\tilde{\nabla}$ its canonical connection. Then

$$(5.35)\quad \tilde{\nabla}g=\tilde{\nabla}\tilde{R}=\tilde{\nabla}h=\tilde{\nabla}S=0,\quad S(\tilde{R})=\tilde{R},\ S(h)=h,\ S(g)=g,$$

where $h=\nabla-\tilde{\nabla}$ and $\tilde{R}$ is the curvature tensor field of $\tilde{\nabla}$. Conversely, if $\tilde{\nabla}$ is a metric connection on a Riemannian manifold $(M,g,S)$, where $S$ and $(I-S)$ are nonsingular affinors on $M$, and (5.35) are fulfilled for $\tilde{\nabla}$, $g$, $h$, $R$, $S$, then there exists a Riemannian locally regular structure $\{s_x\}$ on $M$, that is the canonical connection of $(M,\{s_x\})$ and $S_x=(S_x)_{*x}$.

Returning to our 3-s.l.R.s. we can consider the canonical connection $\overline{\nabla}$ of the a.H.s. $(J,g)$. It is well- known that $\overline{\nabla}g=0$ and from (5.34) it is clear that $\overline{\nabla}S=0$. From the invariance of $\nabla$ under any $s_x$ and from (5.18), (5.34) it follows that $\overline{\nabla}$ is also invariant under $s_x$. So, the uniqueness of the canonical connection of $(M,\{s_x\})$ implies that $\overline{\nabla}=\tilde{\nabla}$ on $M$. Thus, for $\overline{\nabla}$ we get

$$(5.36)\quad \overline{\nabla}h=\overline{\nabla}\overline{R}=0,\ S(h)=h,\ S(\overline{R})=\overline{R}.$$

Using *Theorem 5.30* we obtain

**PROPOSITION 5.31.** Let *(M,J,g)* be an a.H.m. with canonical connection $\overline{\nabla}, h = \nabla - \overline{\nabla}$. If $S = -\frac{1}{2} I + \frac{\sqrt{3}}{2} J$ and (5.36) are fulfilled, then there exists such a Riemannian locally regular *s*-structure $\{s_x\}$ on *M*, that $(M,\{s_x\})$ is a 3-s.l.R.s. and $S_x = (s_x)_{*x}$.

One can find detailed information about 3-s.R.s. in [33].

We consider now more general case.

Let $(M,\{s_x\})$ be a R.l.r.s-m. and $S_x = (s_x)_{*x}$ has only complex eigenvalues $a_1 \pm b_1 i, ..., a_r \pm b_r i$. We define distributions

$D_i$, $i=1,...,r$ by $D_i = ker\,(S^2 - 2a_i S + I)$

It is clear that every $X \in X(M)$ has the unique decomposition $X = X_1 + ... + X_r$, where $X_i \in D_i$, $i=1,...,r$.

An almost complex structure *J* on *M* is defined by

(5.37) $$JX = \sum_{i=1}^{r} 1/b_i(S - a_i I)X_i$$

Really, if $\pi_i$ is the projection on $D_i$, then $\pi_i \pi_j = 0,\ i \neq j,\ \pi_i^2 = \pi_i$ and

$$J^2 = (\sum_{i=1}^{r} 1/b_i(S - a_i I)\pi_i)^2 = \sum_{i=1}^{r} 1/b_i^2(S^2 - 2a_i S + a_i^2 I)\pi_i^2 = \sum_{i=1}^{r} 1/b_i^2(S^2 - 2a_i S + I - b_i^2 I)\pi_i$$
$$= -\sum_{i=1}^{r} \pi_i = -I,$$

It is clear that $D_i$ and $D_j$ are orthogonal each other for $i \neq j$. Further, using that $S^2 = 2a_i S - I$ on $D_i$ we obtain

$$< JX, JY > = \sum_{i=1}^{r} 1/b_i^2 < (S - a_i I)X,\ (S - a_i I)Y_j >$$
$$= \sum_{i=1}^{r} 1/b_i^2 (< X_i, Y_i > - a_i < SX_i, S^2 Y_i > - a_i < SX_i, Y_i > + a_i^2 < X_i, Y_i >)$$
$$= \sum_{i=1}^{r} 1/b_i^2 (< X_i, Y_i > - a_i^2 < X_i, Y_i >) = < X, Y >.$$

Thus, the almost Hermitian structure *(J,g)* is defined on *M*. Let $\overline{\nabla}$ be its canonical connection, see (5.18). Since $S = a_i I + b_i J$ on $D_i$ and $\overline{\nabla} J = 0$, then $\overline{\nabla}_Y SX_i = S\overline{\nabla}_Y X_i$, $i=1,...,r$ and $\overline{\nabla}S = \overline{\nabla}g = 0$ on *M*. It is obvious that $D_i$ and *J* are invariant under any $s_x$, hence it follows from (5.18) that $\overline{\nabla}$ is also invariant under $s_x$. The uniqueness of the canonical connection $\tilde{\nabla}$ of $(M,\{s_x\})$ implies that $\overline{\nabla} = \tilde{\nabla}$ on *M* and (5.36) are fulfilled. So, $h = \nabla - \overline{\nabla}$ defines a locally homogeneous Riemannian structure, see *§3, Chapter 2.*

**REMARK.** If *(M,J,g)* is an a.H.m. and $\overline{\nabla}h = \overline{\nabla}\overline{R} = 0$, where $\overline{\nabla}$ is the canonical connection of *(J,g)*, then a searching for a suitable affinor *S* is a sufficiently difficult problem. We must require $\overline{\nabla}S = 0$, *S(h)=h,* $S(\overline{R}) = \overline{R}$, *S(g)=g* and *S, (I–S)* have to be nonsingular. After that *Theorem 5.30* can be applied.

## *§4. STRUCTURE DEFINED BY AFFINOR F SATISFYING $F^3+F=0$*

**$1^0$**. A structure on an n-dimensional manifold *M* given by a non-null *(1,1)* tensor *F* field satisfying

$$F^3+F=0$$

in called an *f*-structure, [72]. Later on, we shall see that such a tensor field *F* is always *O*-deformable. From *Theorem 1.2* it follows that there exists an associated metric *g* for *P(F)=P(G)* defined by a structure *P(H),* where $H=G\cap O(n)$. The group *H* can be chosen, [72], as

$$H=O(n-2n_2)\times U(n_2),$$

where $2n_2$ is the rank of *F*. If *M* is orientiable and $n-2n_2=1$, then an *f*-structure gives an almost contact structure. We put

$$\pi_1=F^2+I\,,\quad \pi_2=-F^2$$

and consider $L=\pi_1(T(M))$, $V=\pi_2(T(M))$. It is easy to verify that

$$I=\pi_1+\pi_2,\ \pi_1^2=\pi_1,\ \pi_2^2=\pi_2,\ F\pi_1=\pi_1F=0,\ \pi_2F=F\pi_2=F\,,$$

hence *L* and *V* are complementary distributions on *M,* $T(M)=L\oplus V$, where $dim\,V=2n_2$, $dim\,L=n_1=n-2n_2$, and *F* determines the almost complex structure on the distribution $V$.

**PROPOSITION 5.32.** For every $f$-structure $F$ on a manifold $M$ a periodic affinor $S$, $S^4=I$, can be constructed. Conversely, any affinor $S$, $S^4=I$ defines the $f$-structure .

**Proof.** If we put $S_{|V}=F_{|V}$ and $S_{|L}=I_{|L}$ (or $S_{|L}=-I_{|L}$), then it is evident that $S^4=I$. Conversely, if $S^4=I$, then we can take $P=S^2$ and

$$\pi_2=1/2(I+P),\ \pi_2=1/2(I-P),\ L=\pi_1(T(M)),\ V=\pi_2(T(M)).$$

The $f$-structure is defined by the conditions $F=0$ on $L$ and $F=S$ on $V$ . It is clear that $F^3+F=0$.

**QED.**

If $g=<\ ,\ >$ is a Riemannian associated metric determined by a structure $P(H)$, then $V=L^{\perp}$ and $<FX,FY>=\ <X,Y>$ for any $X,Y\in V$.

For any Riemannian metric $\tilde{g}$ on $M$ an associated metric $g$ can be defined by the formula

$$g(X,Y)=\tilde{g}(X,Y)+\tilde{g}(SX,SY)+\tilde{g}(S^2X,S^2Y)+\tilde{g}(S^3X,S^3Y),\ X,Y\in \mathbf{X}(M),$$

where $S_{|V}=F_{|V}$ and $S_{|L}=I_{|L}$ (or $S_{|L}=-I_{|L}$), $S^4=I$. Since $g(SX,SY)=g(X,Y)$ and $g(PX,PY)=g(X,Y)$, where $P=S^2$, for any $X,Y\in \mathbf{X}(M)$, it is obvious that $g$ is associated . The space of all such associated metrics is infinite dimensional.

**$2^0$**. We consider now a fixed pair $(F,g)$ (or $(S,g)$) and calculate its canonical connection $\overline{\nabla}$ and the second fundamental tensor field $h$.

**THEOREM 5.33.** For a Riemannian $f$-structure $(F,g)$ we have

$$(5.38)\quad \overline{\nabla}_X Y=\pi_1\nabla_X\pi_1 Y+\frac{1}{2}\pi_2(\nabla_X\pi_2 Y-F\nabla_X F\pi_2 Y),$$

$$(5.39)\quad h_X Y=\pi_1\nabla_X\pi_2 Y+\pi_2\nabla_X\pi_1 Y+\frac{1}{2}\pi_2(\nabla_X\pi_2 Y+F\nabla_X F\pi_2 Y),$$

$$(5.40)\quad \overline{\nabla}_X Y=\frac{1}{4}(\nabla_X Y+S\nabla_X S^3 Y+S^2\nabla_X S^2 Y+S^3\nabla_X SY),$$

where $X,Y,Z\in \mathcal{X}(M)$.

**Proof.** It is clear that the Lie algebra $\underline{\mathbf{h}}$ of the structure group $H$ of $P(H)$ has the following form

$$\underline{\mathbf{h}}=\{\,\tilde{x}=\left[\begin{array}{c|cc} A' & 0 & \\ \hline 0 & A & B \\ & -B & A \end{array}\right]:\tilde{x}\in\underline{\mathbf{o}}\,\},$$

and $\underline{\mathbf{o}}=\underline{\mathbf{h}}\oplus\underline{\mathbf{m}}$, where

$$\underline{\mathbf{m}}=\{\,\tilde{x}=\left[\begin{array}{c|cc} 0 & -C^T & \\ \hline C & D & K \\ & K & -D \end{array}\right]:\tilde{x}\in\underline{\mathbf{o}}\,\}.$$

We see that $\underline{\mathbf{m}}=\underline{\mathbf{h}}^{\perp}$ with respect to Killing form. For every $\omega\in\underline{\mathbf{o}}$ the natural decomposition $\omega=\omega_{|\underline{\mathbf{h}}}+\omega_{|\underline{\mathbf{m}}}$ can be defined by the formulas

$$\omega_{|\underline{\mathbf{h}}}=p_1\omega p_1+1/2(p_2\omega p_2-jp_2\omega p_2j),\ \ \omega_{|\underline{\mathbf{m}}}=\omega-\omega_{|\underline{\mathbf{h}}}.$$

where

$$p_1=\begin{bmatrix} E_{n_1} & 0 \\ 0 & 0 \end{bmatrix},\ \ p_2=\begin{bmatrix} 0 & 0 \\ 0 & E_{2n_2} \end{bmatrix},\ \ j=\left[\begin{array}{c|cc} 0 & 0 & \\ \hline 0 & 0 & E_{n_2} \\ & -E_{n_2} & 0 \end{array}\right].$$

Let $\varphi$ be a cross section of $P(H)$ over some neighbourhood $U$, which assigns to each $x\in U$ the linear frame $(X_1)_x,\dots,(X_n)_x$ and $X$, $Y=\sum_k f^kX_k$ be vector fields on $M$. Then from (1.6) it follows that

$$\begin{aligned}\overline{\nabla}_{X_x}Y&=\varphi(x)\overline{\omega}(\varphi_*X_x)\varphi(x)^{-1}Y_x+\sum_k(Xf^k)(x)(X_k)_x\\&=\varphi(x)[\,p_1\omega(\varphi_*X_x)p_1+1/2(p_2\omega(\varphi_*X_x)p_2)\\&\quad-jp_2\omega(\varphi_*X_x)p_2j)]\varphi(x)^{-1}Y_x+\sum_k(Xf^k)(x)(X_k)_x\ ,\end{aligned}$$

where $\varphi(x)$ is considered as the mapping of $R^n$ onto $T_x(M)$. It is clear that $\varphi(x)\cdot p_1=\pi_1\cdot\varphi(x),\ \varphi(x)\cdot p_2=\pi_2\cdot\varphi(x),\ \varphi(x)\cdot j=F\cdot\varphi(x)$, hence we obtain

$$\overline{\nabla}_{X_x} Y = [\pi_1\varphi(x)\omega(\varphi_* X_x)\varphi(x)^{-1}\pi_1 Y_x + \sum_k (Xf^k)(x)(\pi_1 X_k)_x]$$
$$+1/2[\pi_2\varphi(x)\omega(\varphi_* X_x)\varphi(x)^{-1}\pi_2 Y_x + \sum_k (Xf^k)(x)(\pi_2 X_k)_x]$$
$$-1/2[F\pi_2\varphi(x)\omega(\varphi_* X_x)\varphi(x)^{-1}\pi_2 FY_x + \sum_k F(Xf^k)(x)F(\pi_2 X_k)_x]$$
$$= \pi_1\nabla_{X_x}\pi_1 Y + 1/2(\pi_2\nabla_{X_x}\pi_2 Y - F\pi_2\nabla_{X_x}F\pi_2 Y).$$

Using the fact that $\nabla_X Y=(\pi_1+\pi_2)\nabla_X(\pi_1+\pi_2)Y$ and obtained formula of $\overline{\nabla}$ we get (5.39).

Since $S_{|V} = F_{|V}$ and $\pi_1=1/2(I+S^2)$, $\pi_2=1/2(I-S^2)$ we have

$$\overline{\nabla}_X Y = 1/2\pi_1\nabla_X\pi_1 Y + 1/2 S\pi_1\nabla_X S\pi_1 Y + \frac{1}{2}\pi_2\nabla_X\pi_2 Y - \frac{1}{2}S\pi_2\nabla_X S\pi_2 Y$$
$$= 1/8(\nabla_X Y + S^2\nabla_X Y + \nabla_X S^2 Y + S^2\nabla_X S^2 Y + \nabla_X Y - S^2\nabla_X Y$$
$$-\nabla_X S^2 Y + S^2\nabla_X S^2 Y) + 1/8(S\nabla_X S^3 YS\nabla_X SY + S^3\nabla_X SY$$
$$+ S^3\nabla_X S^3 Y - S\nabla_X SY + S\nabla_X S^3 Y + S^3\nabla_X SY - S^3\nabla_X S^3 Y)$$
$$= 1/4(\nabla_X Y + S\nabla_X S^3 Y + S^2\nabla_X S^2 Y + S^3\nabla_X SY).$$

**QED.**

**REMARK.** It is easy to verify that the formula (5.38) can be rewritten in the following form

(5.41) $\overline{\nabla}_X Y = \nabla_X Y - 1/2 F\nabla_X FY + \nabla_X F^2 Y + F^2\nabla_X Y + 3/2 F^2\nabla_X F^2 Y$.

For a *f*-structure $F$, $F^3+F=0$ and for every connection $\nabla$ on $M$ (5.41) determines the connection and one can simply check that $\overline{\nabla}F = 0$, that is, the tensor field $F$ is $O$-deformable and always defines a $G$-structure , $\overline{\nabla}S = 0$ too.

**LEMMA 5.34.** For $X\in X(M)$ and $Y,Z\in V$ we have

(5.42) $h_{XYZ} = -h_{XFYFZ}$.

**Proof.** It follows from (5.39) that

$h_{XYZ} = \langle h_X Y, Z\rangle = \langle\pi_2 h_X Y, Z\rangle = 1/2\langle\nabla_X Y + F\nabla_X FY, Z\rangle$,

and

$h_{XFYFZ} = \langle h_X FY, FZ\rangle = \langle \pi_2 h_X FY, FZ\rangle = 1/2\langle \nabla_X FY - F\nabla_X Y, FZ\rangle = -1/2\langle F\nabla_X FY + \nabla_X Y, Z\rangle$.

**QED.**

**3⁰**. Let *(M,{$s_x$ })* be a 4-s.l.R.s. (or R.l.r. σ-m.o.4), see *Definitions 4.9, 4.3,* then for affinor $S : S_x = (s_x)_{*x}$ we have $S^4=I$.

It follows from *Theorem 4.4, Proposition 4.5* and ((5.4O) that the canonical connection $\tilde{\nabla}$ of *(M,{$s_x$ })* (formula (4.3)) coincides with $\overline{\nabla}$ defined by (5.4O). Therefore, if *(M,{$s_x$ })* is a 4-s.l.R.s., then the corresponding structure *(F,g)* is determined on *M* and for its $\overline{\nabla}$ and $h = \nabla - \overline{\nabla}$ (see (5.38),(5.39)) conditions (5.35) have to be fulfilled.

Conversely, using *Theorem 5.30* we obtain

**PROPOSITION 5.35.** Let *(F,g)* be a Riemannian *f*-structure on a manifold *M* with the canonical connection $\overline{\nabla}$, $h = \nabla - \overline{\nabla}$. If *S* is the affinor constructed in *Proposition 5.32* ($S^4=I,\ S_{|L} = -I_{|L}$) and

$$\overline{\nabla} h = \overline{\nabla}\,\overline{R} = 0,\ S(h) = h,\ S(\overline{R}) = \overline{R},$$

then there exists such a Riemannian locally regular *s*-structure *{$s_x$ }* on *M* that *(M,{$s_x$ })* is a 4-s.l.R.s. and $S_x = (s_x)_{*x}$.

**PROPOSITION 5.36.** Let *(M,{$s_x$ })* be a 4-s.l.R.s. and $h = \nabla - \overline{\nabla}$

$P=S^2,\ \pi_1 = 1/2(I+P),\ \pi_2 = 1/2(I-P),\ L=\pi_1(T(M)),\ V=\pi_2(T(M)),$

Then we have

1) $h_X Y=0$ for $X,Y\in L$;

2) $h_X Y=\pi_1\nabla_X Y$ for $X,Y\in V$;

3) $h_X Y=\pi_2\nabla_X Y$ for $X\in V$ and $Y\in L$;

4) $h_X Y=1/2(\nabla_X Y+S\nabla_X SY)\in V$ for $X\in L$ and $Y\in V$.

**Proof.** 1) It follows from (5.35) that $h_{SX}SY=Sh_X Y$ for $X,Y\in X(M)$, so that for $X,Y\in L$ $h_X Y = h_{S^2X} S^2 Y = S^2 h_X Y$ and $h_X Y\in L$, hence $\nabla_X Y\in L$. From another side, from (5.38) we obtain that $\overline{\nabla}_X Y = \pi_1 \nabla_X Y = \nabla_X Y$ and $h_X Y = \nabla_X Y - \overline{\nabla}_X Y = 0$.

2) $h_X Y = h_{S^2 X} S^2 Y = S^2 h_X Y$, that is, $h_X Y \in L$ for $X,Y \in V$. Using (5.38) we get $h_X Y = \pi_1 \nabla_X Y$.

3) $h_X Y = \nabla_X Y - \pi_1 \nabla_X Y = \pi_2 \nabla_X Y$, for $X \in V$ and $Y \in L$.

4) For $X \in L$ and $Y \in V$, we have $h_{S^2 X} S^2 Y = S^2 h_X Y = -h_X Y$, so that $h_X Y \in V$. Since $\nabla_X Y \in V$, then $\nabla_X Y \in V$.

**QED.**

From this proposition it follows that $L$ is integrable and its maximal integral manifolds are totally geodesic submanifolds with respect to the Riemannian connection $\nabla$.

**PROPOSITION 5.37**. Let $(M,\{s_x\})$ be a 4-s.l.R.s. and $h = \nabla - \overline{\nabla}$ Then $h \in \boldsymbol{T}_2 \oplus \boldsymbol{T}_3$ (see Table 3.1).

**Proof.** We can choose orthonormal vector fields by the folllowing way $E_1,\dots,E_k, SE_1,\dots,SE_k \in V$; $E_{k+1},\dots,E_n \in L$ and using *Proposition 5.36, Table 3.1* we obtain

$$c_{12}(h) = \sum_{i=1}^{k} (h_{E_i} E_i + h_{SE_i} SE_i) = \sum_{i=1}^{k} (h_{E_i} E_i + Sh_{E_i} Sh_{E_i})$$

$$= (I + S)\pi_1 \sum_{i=1}^{k} h_{E_i} E_i = 1/2(I + S + S^2 + S^3) \sum_{i=1}^{k} h_{E_i} E_i = 0$$

because $I - S \neq 0$ on $V$.

**QED.**

**REMARK.** If $(M,\{s_x\})$ is a R.l.r. σ-m.o.4 , then $Sh_X Y = h_{SX} SY$, where $X,Y \in \boldsymbol{X}(M)$ (see *Corollary 4.53*), so that *Propositions 5.36, 5.37* are also true in this case.

Let $(M,\{s_x\})$ be a R.l.r.s-m. and $S_x = (s_x)_{*x}$. We denote by $L=[-1]$ a distribution of eigenvector fields corresponding to the eigenvalue –1 and $V = L^{\perp}$. For $Y \in L$ we have $S\tilde{\nabla}_X Y = \tilde{\nabla}_X SY = -\tilde{\nabla}_X Y$, hence $\tilde{\nabla}_X Y \in L$ and $L,V$ are invariant under $\tilde{\nabla}$.

Let a $a_1 \pm b_1 i,\dots,a_r \pm b_r i$ be complex eigenvalues and $D_i = Ker\,(S^2 - 2a_i S + I)$, $i=1,\dots,r$. It is clear that $V = \sum_{i=1}^{r} D_i$ and every $X \in \boldsymbol{X}(M)$ has an unique decomposition

$X=X_0+X_1+...+X_r$, where $X_0 \in L$, $X_i \in D_i$. An affinor $F$ on $M$ is defined by

(5.43) $FX_0 = 0, \quad F(X - X_0) = \sum_{i=1}^{r} 1/b_i(S - a_i I)X_i$ .

By similar arguments as in **$6^0$**, §3 we get that *(F,g)* is a Riemannian *f*-structure on *M*, i.e., $F^3+F=0$.

Let $\overline{\nabla}$ be the canonical connection of *(F,g)* (see (5.38)) . Since $S=a_i I+b_i F$ on $D_i$, then $\overline{\nabla} S = \overline{\nabla} g = 0$ and $\overline{\nabla}$ is invariant under any $s_x$. The uniqueness of the canonical connection $\overline{\nabla}$ of $(M,\{s_x\})$ implies that $\overline{\nabla} = \tilde{\nabla}$ on $M$ and (5.36) are realized.

In particular, $\overline{\nabla} h = \overline{\nabla} R = 0$ and $h = \nabla - \overline{\nabla}$ defines a locally homogeneous Riemannian structure.

# CHAPTER 6

# A CLASSIFICATION OF ALMOST CONTACT METRIC STRUCTURES

In this chapter, using our interpretation of the classification of A.Gray, L.M.Hervella and the second fundamental tensor field $h$ of an almost contact metric structure, we get a classification of such structures in terms of $h$. There are $2^{12}$ classes of almost contact metric structures. A similar classification was considered by D.Chinea and C.Gonzalez (A.A.Alexiev and G.Ganchev) by a different method. Good relations have been found between both the classifications and this allows to describe the canonical connection $\overline{\nabla}$ for every class.

§1 is devoted to the obtaining of a new classification with help of the tensor field h.

In §2, we prove a theorem which states that both the classifications are the same up to an isomorphism. We follow especially closely to [14].

## *§1. ABOUT A CLASSIFICATION OF ALMOST CONTACT METRIC STRUCTURES*

**1⁰**. Let $M$ be a $(2n+1)$-dimensional manifold. An almost contact metric structure (a.c.m.s.) on $M$ is called $H$-structure $P(H)$, where $H=U(n)\times 1$. The extension of $P(H)$ to $O(2n+1)$ defines the Riemannian metric $g$ on $M$. $P(H)$ also determines the vector field $\xi(\|\xi\|=1)$; the almost product structure $T(M)=V\oplus L$, where $L=[\xi]$ $(dim\ L=1)$, $V=L^{\perp}$; the tensor field $F$ of type (1.1), $F^3+F=0$, $F(\xi)=0$, $F^2=-I$ on $V$, $<FX,FY>=<X,Y>$ for $X,Y\in V$; the $1$-form $\eta$, where $\eta(X)=<X,\xi>$ for $X\in \mathfrak{X}(M)$. We denote the projections to $L$ and $V$ by $\pi_1$ and $\pi_2$ correspondingly. ($\pi_1+\pi_2=I$, $\pi_i^2=\pi_i$, $i=1,2$).

**THEOREM 6.1.** For a.c.m.s. on $M$ and $X,Y\in \mathfrak{X}(M)$ we have

$$(6.1)\quad \overline{\nabla}_X Y=\pi_1\nabla_X\pi_1 Y+\frac{1}{2}\pi_2(\nabla_X\pi_2 Y-F\nabla_X F\pi_2 Y),$$

$$(6.2)\quad h_X Y=\pi_1\nabla_X\pi_2 Y+\pi_2\nabla_X\pi_1 Y+\frac{1}{2}\pi_2(\nabla_X\pi_2 Y+F\nabla_X F\pi_2 Y).$$

Proof is the same to that of Theorem 5.33.

If $h_{XYZ}=<h_X Y,Z>$, $X,Y,Z\in \mathbf{X}(M)$, then from (2.2) we have

(6.3) $h_{XYZ}=-h_{XZY}$.

It follows from Lemma 5.34 that for $X\in \mathbf{X}(M)$ and $Y,Z\in V$

(6.4) $h_{XYZ}=-h_{XFYFZ}$.

Since $\overline{\nabla}$ is the canonical connection of a.c.m.s. $P(H)$, then for $X,Y\in \mathbf{X}(M)$

$\overline{\nabla}F=\overline{\nabla}g=0$;

$\overline{\nabla}_X\xi=\pi_1\nabla_X\xi=0$ because $\|\xi\|=1$;

$\overline{\nabla}\eta=0$ because $\overline{\nabla}g=0$, $\overline{\nabla}\xi=0$;

$\overline{\nabla}\pi_i=0,\ i=1,2$, because $\pi_1=I+F^2$, $\pi_2=-F^2$.

**2⁰**. Let $p\in M$, $T=T_p(M)$, $V_p=V$, $L_p=L$. $T$ is an Euclidian vector space over $R$ with respect to the inner product < , > induced by $g$. We consider the subspace $\overline{\boldsymbol{T}}(T)$ formed by all the tensors of the type (0.3) which satisfy the identities (6.3), (6.4), i.e.,

(6.5)

$$\overline{T}(T)=\{h\in\overset{3}{\otimes}T^*:h_{XYZ}=-h_{XZY},\ X,Y,Z\in T;\ h_{XYZ}=-h_{XFYFZ}, X\in T,\ Y,Z\in V\}$$

$\overline{\boldsymbol{T}}(T)$ is an Euclidian vector space under the inner product defined by (3.1). The action of the group $H=U(n)\times 1$ is defined by (3.2). Later on, we shall choose an orthonormal basis of $T$ in the following way

(6.6) $E_1,\dots,E_n, FE_1,\dots, FE_n\in V$; $E_{2n+1}=\xi\in L$.

The action of $H$ on $V$ coincides with the ordinary action of $U(n)$. We consider now the induced action on $V^*\otimes V^*$ defined by

(6.7) $(ah)_{XY}=h_{a^{-1}Xa^{-1}Y}$,

where $a\in U(n)$, $X,Y\in V$, $h\in V^*\otimes V^*$.

The inner product is also determined on $V^*\otimes V$ by

(6.8) $<h^1,h^2>=\sum_{i,k}h^1_{E_iE_k}h^2_{E_iE_k}$,

where $h^1, h^2 \in V^* \otimes V^*$ and $E_1,...,E_{2n}$ form an orthonormal basis of $V$. Using the basis of type (6.6) we define

$$tr(h) = \sum_{i=1}^{n}(h_{E_iE_i} + h_{FE_iFE_i}), \quad \overline{tr}(h) = \sum_{i=1}^{n}(h_{E_iFE_i} - h_{FE_iE_i}).$$

**LEMMA 6.2.** If $n \geq 2$, then $V^* \otimes V^* = \bigoplus_{i=1}^{6} \overline{V_i}$, where $\overline{V_i}$ are invariant and irreducible under the action of $U(n)$. One can find the spaces $\overline{V_i}$, $i = 1,...,6$ in *Table 6.1*.

**Proof.** It is well known that $V^* \otimes V^*$ has the following decomposition to irreducible components with respect to the action of $O(2n)$.

$$V^* \otimes V^* = \Lambda^2 V \oplus S_0^2 V \oplus A_V,$$

where $\Lambda^2 V$ is the subspace of antisymmetric tensors, $S_0^2 V$ is the subspace of symmetric tensors of a zero trace. If we consider the elements of $V^* \otimes V^*$ to be the bilinear forms $\mu : V \times V \to R$, then $\mu \in A_V$ if and only if $\mu(X,Y) = 1/2n <X,Y> tr\mu$. All the components are invariant under the action of $U(n)$. Having fixed a basis of type (6.6) in $V$ we can identify $V^* \otimes V^*$ with the space of matrices. Thus, we have

$$\Lambda^2 V = \overline{V_1} \oplus \overline{V_2} \oplus \overline{V_3},$$

where

$$\overline{V_1} = \{ \begin{bmatrix} A & B_0 \\ -B_0 & A \end{bmatrix} : A^T = -A, B_0^T = B_0,\ trB_0 = 0 \}\ (dim\overline{V_1} = n^2 - 1),$$

$$\overline{V_2} = \{ \begin{bmatrix} C & D \\ D & -C \end{bmatrix} : C^T = -C,\ D^T = -D \} \quad (dim\overline{V_2} = n^2 - n),$$

$$\overline{V_3} = \{ \begin{bmatrix} 0 & E \\ -E & 0 \end{bmatrix} trB : B = B_0 + (trB)E \}\ (dim\overline{V_3} = 1).$$

The subspace of symmetric tensors of a zero trace is decomposed as

$$S_0^2 V = \overline{V_4} \oplus \overline{V_5},$$

where

$$\overline{V}_4 = \{ \begin{bmatrix} \overline{A}_0 & \overline{B} \\ -\overline{B} & \overline{A}_0 \end{bmatrix} : \overline{A}_0^T = \overline{A},\ \overline{B}^T = -\overline{B},\ tr\overline{A}_0 = 0 \} \quad (dim\overline{V}_4 = n^2 - 1),$$

$$\overline{V}_5 = \{ \begin{bmatrix} \overline{C} & \overline{D} \\ \overline{D} & -\overline{C} \end{bmatrix} : \overline{C}^T = \overline{C},\ \overline{D}^T = \overline{D} \} \quad (dim\overline{V}_5 = n^2 + n).$$

We denote

$$A_V = \overline{V}_6 \ (dim\overline{V}_6 = 1).$$

Using (6.8) it is easy to verify that all $\overline{V}_i$ are mutually orthogonal. Considering the matrix realization of *U(n)*

$$U(n) = \{ \begin{bmatrix} A & B \\ -B & A \end{bmatrix} \in O(2n) \}$$

we can check that all $\overline{V}_i$ are invariant and irreducible under the action of *U(n)*.

**QED.**

Let *X,Y* be arbitrary vectors in *V*. We have

Table 6.1

| Class | Dimension | Defining condition |
|---|---|---|
| $\overline{V}_1$ | $n^2-1$ | $h_{XY} = -h_{YX},\ h_{FXFY} = h_{XY},\ \overline{tr}(h) = 0$ |
| $\overline{V}_2$ | $n^2-n$ | $h_{XY} = -h_{YX},\ h_{FXFY} = -h_{XY}$ |
| $\overline{V}_3$ | $1$ | $h_{XY} = \frac{1}{2n} < X, FY > \overline{tr}(h)$ |
| $\overline{V}_4$ | $n^2-1$ | $h_{XY} = h_{YX},\ h_{FXFY} = h_{XY},\ tr(h) = 0$ |
| $\overline{V}_5$ | $n^2+n$ | $h_{XY} = h_{YX},\ h_{FXFY} = -h_{XY}$ |
| $\overline{V}_6$ | $1$ | $h_{XY} = 1/2n < X, Y > tr(h)$ |

**$3^0$**. We consider now the main theorem of the section about the decomposition of the space $\overline{\boldsymbol{T}}(T)$ to invariant and irreducible components with respect to the action of the group H=U(n)×1.

**THEOREM 6.3.** If $n \geq 3$, then there exists the orthogonal decomposition $\overline{\boldsymbol{T}}(T) = \bigoplus_{i=1}^{12} \overline{\boldsymbol{T}}_i$, where $\overline{\boldsymbol{T}}_i$ are invariant and irreducible under the action of *H*. The spaces $\overline{\boldsymbol{T}}_i$, $i=1,...,12$ are adduced in *Table 6.2*.

**Proof.** We have

$$T^* \otimes T^* \otimes T^* = T^* \otimes (V^* \oplus L^*) \otimes (V^* \oplus L^*) = T^* \otimes (V^* \otimes V^* \oplus V^* \otimes L^* \oplus L^* \otimes V^* \oplus L^* \otimes L^*)$$

From (6.3) we obtain

$$(T^* \otimes V^* \otimes L^*) \cap \overline{\boldsymbol{T}} = (T^* \otimes L^* \otimes V^*) \cap \overline{\boldsymbol{T}}$$

and as $h_{X\xi\xi} = 0$ $(T^* \otimes L^* \otimes L^*) \cap \overline{\boldsymbol{T}} = \emptyset$.

Using the invariance of $\overline{\boldsymbol{T}}$ and of the components under the action of the group *H* we get

$$\overline{\boldsymbol{T}} = (\overset{3}{\otimes} T^*) \cap \overline{\boldsymbol{T}} = (T^* \otimes V^* \otimes V^* \oplus T^* \otimes V^* \otimes L^*) \cap \overline{\boldsymbol{T}}$$
$$= (V^* \otimes V^* \otimes V^* \oplus V^* \otimes V^* \otimes L^* \oplus L^* \otimes V^* \otimes L^* \otimes L^* \otimes V^* \otimes V^*) \cap \overline{\boldsymbol{T}}$$

We denote $\overline{T}^1 = (V^* \otimes V^* \otimes V^*) \cap \overline{T}$, $\overline{T}^2 = (V^* \otimes V^* \otimes L^*) \cap \overline{T}$,
$\overline{T}_{11} = (L^* \otimes V^* \otimes L^*) \cap \overline{T}$, $\overline{T}_{12} = (L^* \otimes V^* \otimes V^*) \cap \overline{T}$.

The group *H* preserves the sum $T=V \oplus L$ invariant and the action of *H* on *V* coincides with that of *U(n)*, on *L* it is trivial. Thus *H* preserves the decomposition

$$\overline{\boldsymbol{T}} = \overline{\boldsymbol{T}}^1 \oplus \overline{\boldsymbol{T}}^2 \oplus \overline{\boldsymbol{T}}_{11} \oplus \overline{\boldsymbol{T}}_{12}$$

invariant. From (3.1) we can notice that the components of this decomposition are orthogonal . We define $\beta, \overline{\beta} \in T^*$ by

(6.9) $\beta(U) = \sum_{i=1}^{n} (h_{E_i E_i U} + h_{FE_i FE_i U})$

$$\overline{\beta}(U)=\sum_{i=1}^{n}(h_{E_iFE_iU}-h_{FE_iE_iU}),\ \ U\in T.$$

From (5.25) we have obtained the following decomposition of the space $\overline{\boldsymbol{T}}^1$ to orthogonal irreducible components under the action of *U(n)*

$$\overline{\boldsymbol{T}}^1=\bigoplus_{i=1}^{4}\overline{\boldsymbol{T}}_i,$$

where $\overline{\boldsymbol{T}}_1=\{h\in\overline{\boldsymbol{T}}^1: h_{XYZ}=-h_{YXZ},\ X,Y,Z\in V\},$

$\overline{\boldsymbol{T}}_2=\{h\in\overline{\boldsymbol{T}}^1: \boldsymbol{\sigma}\ h_{XYZ}=0,\ X,Y,Z\in V\},$

$\overline{T}_3=\{h\in\overline{T}^1: h_{XYZ}-h_{FXFXZ}=\beta(Z)=0,\ X,Y,Z\in V\},$

$$\overline{\boldsymbol{T}}_4=\{h\in\overline{\boldsymbol{T}}^1: h_{XYZ}=1/2(n-1)[<X,Y>\beta(Z)-<X,Z>\beta(Y)$$
$$-<X,FY>\beta(FZ)+<X,FZ>\beta(FY)],\ X,Y,Z\in V\},$$

The action of the group *H* on *L* is trivial. Using *Lemma 2.1* (*Table 6.1*) we can get the following decomposition of $\overline{\boldsymbol{T}}^2$ to orthogonal irreducible components under the action of *H*

$\overline{\boldsymbol{T}}^2=\bigoplus_{i=5}^{10}\overline{\boldsymbol{T}}_i$, where $\overline{\boldsymbol{T}}_i=\overline{V}_{i-4}\otimes L^*$ and *i=5,...,10.*

Since the action of *U(n)* on *V* is irreducible, therefore

$\overline{\boldsymbol{T}}_{11}=L^*\otimes V^*\otimes L^*$ is irreducible.

It follows from (6.3) , (6.4) that
$\overline{\boldsymbol{T}}_{12}=L^*\otimes\overline{V}_2$ is irreducible.

Thus, we have got $\overline{\boldsymbol{T}}=\bigoplus_{i=1}^{12}\overline{\boldsymbol{T}}_i$ and

*dim* $\overline{\boldsymbol{T}}$ *=(2n+1)n(n+1)=dim Hom (T,* ***m****),*

where ***m*** is the subspace of ***o****(2n+1)* defined in the proof of *Theorem 5.33*.

**QED.**

Let *X,Y,Z* be arbitrary vectors in *V*. Then, we have

Table 6.2

| Class | Dimension | Defining condition |
|---|---|---|
| $\overline{T}_0$ | 0 | $h=0$ |
| $\overline{T}_1$ | $1/3n(n-1)(n-2)$ | $h_{XYZ}=-h_{YXZ}$, $h_{XY\xi}=h_{\xi X\xi}=h_{\xi XY}=0.$ |
| $\overline{T}_2$ | $2/3n(n-1)(n+1)$ | $\sigma\, h_{XYZ}=0,\ h_{XY\xi}=h_{\xi X\xi}=h_{\xi XY}=0.$ |
| $\overline{T}_3$ | $n(n+1)(n-2)$ | $h_{XYZ}-h_{FXFYZ}=\beta(Z)=0,$ $h_{XY\xi}=h_{\xi X\xi}=h_{\xi XY}=0.$ |
| $\overline{T}_4$ | $2n$ | $h_{XYZ}=1/2(n-1)(<X,Y>\beta(Z)-<X,Z>\beta(Y)$ $-<X,FY>\beta(FZ)+<X,FZ>\beta(FY)),$ $h_{XY\xi}=h_{\xi X\xi}=h_{\xi XY}=0.$ |
| $\overline{T}_5$ | $n^2-1$ | $h_{XY\xi}=-h_{YX\xi},\ h_{FXFY\xi}=h_{XY\xi},\ \overline{\beta}(\xi)=0,$ $h_{XYZ}=h_{\xi X\xi}=h_{\xi XY}=0.$ |
| $\overline{T}_6$ | $n^2-n$ | $h_{XY\xi}=-h_{YX\xi},\ h_{FXFY\xi}=-h_{XY\xi},$ $h_{XYZ}=h_{\xi X\xi}=h_{\xi XY}=0.$ |
| $\overline{T}_7$ | 1 | $h_{XY\xi}=\dfrac{1}{2n}<X,FY>\overline{\beta}(\xi),\ h_{XYZ}=h_{\xi X\xi}=h_{\xi XY}=0.$ |
| $\overline{T}_8$ | $n^2-1$ | $h_{XY\xi}=h_{YX\xi},\ h_{FXFY\xi}=h_{XY\xi},\ \beta(\xi)=0,$ $h_{XYZ}=h_{\xi X\xi}=h_{\xi XY}=0.$ |
| $\overline{T}_9$ | $n^2+n$ | $h_{XY\xi}=h_{YX\xi},\ h_{FXFY\xi}=-h_{XY\xi},$ $h_{XYZ}=h_{\xi X\xi}=h_{\xi XY}=0.$ |
| $\overline{T}_{10}$ | 1 | $h_{XY\xi}=1/2n<X,Y>\beta(\xi)$ $h_{XYZ}=h_{\xi X\xi}=h_{\xi XY}=0.$ |
| $\overline{T}_{11}$ | $2n$ | $h_{XYZ}=h_{XY\xi}=h_{\xi XY}=0.$ |
| $\overline{\boldsymbol{T}}_{12}$ | $n^2 - n$ | $h_{\xi XY} = -h_{\xi YX} = -h_{\xi FXFY}$, $h_{XYZ} = h_{XY\xi} = h_{\xi X\xi} = 0.$ |
| $\overline{\boldsymbol{T}}$ | $(2n+1)n(n+1)$ | ------------ |

This classification was announced in [25] and given in [28]. In fact, for every class from *Table 6.2* the canonical connection $\overline{\nabla} = \nabla - h$ have been constructed and to study these classes the torsion tensor field $\overline{T}$ and various curvature characteristics of $\overline{\nabla}$ can be applied.

It is directly follows from the analysis of *Table 6.2* that

$$\overline{\boldsymbol{T}} = \overline{\boldsymbol{T}}_7 \oplus \overline{\boldsymbol{T}}_9 \oplus \overline{\boldsymbol{T}}_{10} \oplus \overline{\boldsymbol{T}}_{11}$$

when $n=1$, and

$$\overline{\boldsymbol{T}} = \overline{\boldsymbol{T}}_2 \oplus \overline{\boldsymbol{T}}_4 \bigoplus_{i=5}^{12} \overline{\boldsymbol{T}}_i$$

when $n=2$.

We denote the set of all the kinds of combinations of $\{1,...,12\}$ by $\boldsymbol{A}_{12}$. It is evident that every combination $\alpha = \{\alpha_1,...,\alpha_l\}$ determines the invariant subspace $\overline{\boldsymbol{T}}_\alpha = \overline{\boldsymbol{T}}_{\alpha_1} \oplus ... \oplus \overline{\boldsymbol{T}}_{\alpha_l}$ of $\overline{\boldsymbol{T}}$.

**DEFINITION 6.1.** We say that a.c.m.s. is of the class $\overline{\boldsymbol{T}}_\alpha$, $\alpha \in \boldsymbol{A}_{12}$, or has a type $\overline{\boldsymbol{T}}_\alpha$ on $M$ if $h_x \in \overline{\boldsymbol{T}}_\alpha(T_x(M))$, $\forall x \in M$, and it is strictly of the class $\overline{\boldsymbol{T}}_\alpha$ if $h_X \neq 0$, $\forall X \in \boldsymbol{X}(M)$, $X \neq 0$.

If $\overline{\boldsymbol{T}}(M)$ denotes the vector bundle over $M$ with the typical fibre $\overline{\boldsymbol{T}}(M)$, then $\overline{\boldsymbol{T}}_\alpha(M)$, $\alpha \in \boldsymbol{A}_{12}$ are the vector subbundles in $\overline{\boldsymbol{T}}(M)$ and a.c.m.s. is strictly of the class $\overline{\boldsymbol{T}}_\alpha$ if $x \to h_x$ is a nonzero cross - section of the corresponding subbundle. Since the power of $\boldsymbol{A}_{12}$ is equel to $2^{12}$, then there are 4096 classes of a.c.m.s. on $M$.

**THEOREM 6.4.** Let a.c.m.s. be a quasi homogeneous structure having a type $\overline{\boldsymbol{T}}_\alpha$ for some point $p \in M$. Then this structure is of the class $\overline{\boldsymbol{T}}_\alpha$ on $M$.

**Proof** is similar to that of *Theorem 3.3*.

## §2. ABOUT A CLASSIFICATION OF D.CHINEA AND C.GONZALEZ

**$1^0$**. Let $(F,\xi,\eta, g)$ be an almost contact metric structure (a.c.m.s.) on $M$, $\dim M=2n+1$ and let

$$\Phi(X,Y) = <X,FY>, \quad X,Y\in \mathfrak{X}(M).$$

Using a tensor field $\nabla\Phi$ A.A.Alexiev, G.Ganchev in [1] and D.Chinea, G.Gonzalez in [14] obtained a classification of a.c.m.s. on $M$ by a different method.

The main point of this section is to identify our classification of §1 with that in [14].

The following equalities one can find in [14].

(6.10) $(\nabla_X\Phi)(Y,Z)=<Y,\nabla_X(F)Z>,$

(6.11) $(\nabla_X\Phi)(Y,Z)+ (\nabla_X\Phi)(FY,FZ)=\eta(Z)(\nabla_X\eta)FY-\eta(Y)(\nabla_X\eta)FZ,$

(6.12) $(\nabla_X\eta)Y=<Y,\nabla_X\xi>=(\nabla_X\Phi)(\xi,FY),$

(6.13) $2d\eta(X,Y)=(\nabla_X\eta)Y-(\nabla_Y\eta)X,$

(6.14) $3d\Phi(X,Y,Z)=\sigma(\nabla_X\Phi)(Y,Z)$

where $d\eta$, $d\Phi$ are the exterior derivatives of $\eta$ and $\Phi$, $X,Y,Z\in\mathfrak{X}(M)$.

(6.15)

$$\delta\Phi(X)=-\sum_{i=1}^{n}[(\nabla_{E_i}\Phi)(E_i,X)+(\nabla_{FE_i}\Phi)(FE_i,X)]-(\nabla_\xi\Phi)(\xi,X),$$

(6.16) $$\delta\eta=-\sum_{i=1}^{n}[(\nabla_{E_i}\eta)E_i+(\nabla_{FE_i}\eta)FE_i],$$

where $\delta\Phi$ and $\delta\eta$ are the coderivatives of $\Phi$ and $\eta$, $\{E_i,FE_i,\xi\}$, $i=1,...,n$, is a a local orthonormal basis of the type (6.6) defined on an open subset of $M$.

Let $p$ be a fixed point of $M$ and $T=T_p(M)$, $X,Y,Z\in T$. We consider a vector subspace $\mathbf{C}(T)$ in $\overset{3}{\otimes}T^*$

(6.17) $$\mathbf{C}(T)=\{\alpha\in\overset{3}{\otimes}T^*:\alpha(X,Y,Z)=-\alpha(X,Z,Y)=-\alpha(X,FY,FZ)$$
$$+\eta(Y)\alpha(X,\xi,Z)+\eta(Z)\alpha(X,Y,\xi)\}$$

**THEOREM 6.5** [14]. $(\nabla\Phi)_p \in \boldsymbol{C}(T)$ and there exists an orthogonal decomposition $\boldsymbol{C}(T)=\bigoplus_{i=1}^{12} \boldsymbol{C}_i$, where $\boldsymbol{C}_i$ are invariant and irreducible under the action of the group $H=U(n)\times 1$. The spaces $\boldsymbol{C}_i$, $i=1,...,12$ are adduced in *Table 6.3*.

Let $X,Y,Z$ be arbitrary vectors in $T$. Then we have $\mathbf{2^0}$. We consider now relations between $h$ and the other tensor fields.

Table 6.3

| Class | Dimension | Defining condition |
|---|---|---|
| $\boldsymbol{C}_1$ | $1/3n(n-1)(n-2)$ | $(\nabla_X\Phi)(X,Y)=0,\ \nabla\eta=0$ |
| $\boldsymbol{C}_2$ | $2/3n(n-1)(n+1)$ | $d\Phi=\nabla\eta=0$ |
| $\boldsymbol{C}_3$ | $n(n+1)(n-2)$ | $(\nabla_X\Phi)(Y,Z)-(\nabla_{FX}\Phi)(FY,Z)=0,\ \delta\Phi=0$ |
| $\boldsymbol{C}_4$ | $2n$ | $(\nabla_X\Phi)(Y,Z)=-1/2(n-1)[<FX,FY>\delta\Phi(Z) - <FX,FZ>\delta\Phi(Y)-\Phi(X,Y)\delta\Phi(FZ) + \Phi(X,Z)\delta\Phi(FY)],\ \delta\Phi(\xi)=0$ |
| $\boldsymbol{C}_5$ | $1$ | $(\nabla_X\Phi)(Y,Z)=1/2n[\Phi(X,Z)\eta(Y)-\Phi(X,Y)\eta(Z)]\delta\eta$ |
| $\boldsymbol{C}_6$ | $1$ | $(\nabla_X\Phi)(Y,Z)=1/2n[<X,Z>\eta(Y)-<X,Y>\eta(Z)]\delta\Phi(\xi)$ |
| $\boldsymbol{C}_7$ | $n^2-1$ | $(\nabla_X\Phi)(Y,Z)=\eta(Z)(\nabla_Y\eta)FX+\eta(Y)(\nabla_{FX}\eta)Z,$ $\delta\Phi=0$ |
| $\boldsymbol{C}_8$ | $n^2-1$ | $(\nabla_X\Phi)(Y,Z)=-\eta(Z)(\nabla_Y\eta)FX+\eta(Y)(\nabla_{FX}\eta)Z,$ $\delta\eta=0$ |
| $\boldsymbol{C}_9$ | $n(n+1)$ | $(\nabla_X\Phi)(Y,Z)=\eta(Z)(\nabla_Y\eta)FX-\eta(Y)(\nabla_{FX}\eta)Z$ |
| $\boldsymbol{C}_{10}$ | $n(n-1)$ | $(\nabla_X\Phi)(Y,Z)=-\eta(Z)(\nabla_Y\eta)FX-\eta(Y)(\nabla_{FX}\eta)Z$ |
| $\boldsymbol{C}_{11}$ | $n(n-1)$ | $(\nabla_X\Phi)(Y,Z)=-\eta(X)(\nabla_\xi\Phi)(FY,FZ)$ |
| $\boldsymbol{C}_{12}$ | $2n$ | $(\nabla_X\Phi)(Y,Z)=\eta(X)\eta(Z)(\nabla_\xi\eta)FY-\eta(X)\eta(Y)(\nabla_\xi\eta)FZ$ |

**LEMMA 6.6.** For a.c.m.s. we have the following identities

(6.18) $h_{XYZ}=-1/2<\nabla_X(F)FY,Z>,\ X\in \mathcal{X}(M),\ Y,Z\in V$;

(6.18) $\nabla_X(F)Y=h_X FY-Fh_X Y,\ X,Y\in \mathcal{X}(M)$.

**Proof.** It follows from (6.2) that for $Y,Z\in V$

$h_{XYZ}=1/2<\nabla_X Y+F\nabla_X FY,Z>=-1/2<\nabla_X(F)FY,Z>$.

To obtain (6.19) we subtract the identity $\overline{\nabla}_X FY-F\overline{\nabla}_X Y=0$ from the identity $\nabla_X(F)Y=\nabla_X FY-F\nabla_X Y$.

**QED.**

**PROPOSITION 6.7.** We have the following identity

(6.20) $(\nabla_X\Phi)(Y,FZ)=2h_{XYZ}-<Y,\xi><\nabla_X\xi,Z>+2<Z,\xi><Y,\nabla_X\xi>$

where $X,Y,Z\in \mathcal{X}(M)$.

**Proof.** We consider the following cases
1) If $Y,Z\in V$ then, using (6.3), (6.10), (6.18), we get

$$(\nabla_X\Phi)(Y,FZ)=<Y,(\nabla_X F)FZ>=-2h_{XZY}=2h_{XYZ};$$

2) If $Y\in V$ and $Z=\xi$, then we obtain

$(\nabla_X\Phi)(Y,FZ)=(\nabla_X\Phi)(Y,F\xi)=0$,

$2h_{XYZ}=2<\nabla_X Y-\overline{\nabla}_X Y,\xi>=2<\nabla_X Y,\xi>=-2<\xi,\xi><Y,\nabla_X\xi>$;

3) For $Y=\xi$ and $Z\in V$ from (6.10), (6.19), (6.3) we have

$$(\nabla_X\Phi)(Y,FZ)=(\nabla_X\Phi)(\xi,FZ)=<\xi,(\nabla_X F)FZ>$$
$$=<\xi,h_X F^2 Z-Fh_X FZ>=-h_{XZ\xi}=h_{X\xi Z},$$

$$2h_{X\xi Z}-<\xi,\xi><\nabla_X\xi,Z>=2h_{X\xi Z}-<\nabla_X\xi-\overline{\nabla}_X\xi,Z>=2h_{X\xi Z}-h_{X\xi Z}=h_{X\xi Z}.$$

4) For $Z=Y=\xi$ we get $(\nabla_X\Phi)(Y,F\xi)=2h_{X\xi\xi}=0$.

**QED.**

Using (6.12) we obtain

$$(6.21)\ h_{XYZ} = \frac{1}{2}\eta(Y)(\nabla_X \eta)Z - \eta(Z)(\nabla_X \eta)Y + \frac{1}{2}(\nabla_X \Phi)(Y,FZ)$$
$$= \frac{1}{2}\eta(Y)(\nabla_X \Phi)(\xi,FZ) - \eta(Z)(\nabla_X \Phi)(\xi,FY)$$
$$+ \frac{1}{2}(\nabla_X \Phi)(Y,FZ), \quad X,Y,Z \in \boldsymbol{X(M)}$$

From (6.10) and (6.19) it follows

$$(6.22)\ (\nabla_X \Phi)(Y,Z) = \langle Y,(\nabla_X F)Z \rangle = \langle h_X FZ - Fh_X Z, Y \rangle,\ X,Y,Z \in \boldsymbol{X(M)}$$

Thus, for a fixed point $p$ of $M$ we can consider a mapping

$$B : \boldsymbol{C(T)} \rightarrow \overline{\boldsymbol{T}}(T) : (\nabla\Phi)_p \mapsto h_p$$

**THEOREM 6.8.** The mapping $B$ is a linear isomorphism between the vector spaces $\boldsymbol{C(T)}$ and $\overline{\boldsymbol{T}}(T)$.

**Proof.** From the second part of (6.21) it is evident that

$$B(\mu\nabla\Phi_1 + \nu\nabla\Phi_2) = \mu B(\nabla\Phi_1) + \nu B(\nabla\Phi_2)$$

and $B$ is a linear mapping. From (6.21) and (6.22) it follows that $B$ is a one-to-one correspondence.

**QED.**

From (6.12) we have that $(\nabla_X \eta)Y = \langle Y,\nabla_X \xi \rangle$, $(\nabla_X \eta)\xi = \langle \xi,\nabla_X \xi \rangle = 0$ and

$$(6.23)\ (\nabla_X \eta)Y = \langle \nabla_X \xi - \overline{\nabla}_X \xi, Y \rangle = h_{X\xi Y}.$$

Using (6.16), (6.23), (6.9) we obtain

$$(6.24)\ \delta\eta = -\sum_{i=1}^{n} [(\nabla_{E_i} \eta)E_i + (\nabla_{FE_i} \eta)FE_i] = -\sum_{i=1}^{n} [h_{E_i \xi E_i} + h_{FE_i \xi FE_i}]$$
$$= \sum_{i=1}^{n} [h_{E_i E_i \xi} + h_{FE_i FE_i \xi}] = \beta(\xi)$$

It follows from (6.15), (6.20), (6.10), (6.19), (6.9) that for $X \in V$

$$\delta\Phi(FX) = -\sum_{i=1}^{n}[(\nabla_{E_i}\Phi)(E_i,FX) + (\nabla_{FE_i}\Phi)(FE_i,FX)] - (\nabla_{\xi}\Phi)(\xi,FX)$$

$$= -\sum_{i=1}^{n}[2h_{E_iE_iX} - <E_i,\xi><\nabla_{E_i}\xi,X> + 2<X,\xi><E_i,\nabla_{E_i}\xi>$$

$$+ 2h_{FE_iFE_iX} - <FE_i,\xi><\nabla_{FE_i}\xi,X> + 2<X,\xi><FE_i,\nabla_{FE_i}\xi>]$$

$$- <\xi,(\nabla_{\xi}F)\xi> = -2\sum_{i=1}^{n}[h_{E_iE_iX} + h_{FE_iFE_iX}] - <\xi,\nabla_{\xi}FX>$$

$$= -2\beta(X) - h_{\xi\xi FX}.$$

From (6.10), (6.15), (6.19), (6.9) we obtain

$$\delta\Phi(\xi) = -\sum_{i=1}^{n}[<E_i,(\nabla_{E_i}F)\xi> + <FE_i,(\nabla_{FE_i}F)\xi>] - <\xi,(\nabla_{\xi}F)\xi>$$

$$= -\sum_{i=1}^{n}[-<E_i,Fh_{E_i}\xi> - <FE_i,Fh_{FE_i}\xi>] = -\sum_{i=1}^{n}[h_{E_i\xi FE_i} - h_{E_i\xi E_i}]$$

$$= \sum_{i=1}^{n}[h_{E_iFE_i\xi} - h_{FE_iE_i\xi}] = \overline{\beta}(\xi).$$

We have got

(6.25) $\delta\Phi(X) = 2\beta(FX) + h_{\xi\xi X}, \quad X \in V;$

**3⁰**. We are ready now to identify classes in both the classifications.

**LEMMA 6.9.** $B(\boldsymbol{C}_1) = \overline{\boldsymbol{T}}_1$

**Proof.** For $X,Y \in V$ we have from (6.20), (6.23) that

$$(\nabla_X\Phi)(X,Y) = -2h_{XXFY} = 0,\ (\nabla_X\eta)Y = -h_{XY\xi} = 0,\ (\nabla_{\xi}\eta)Y = -h_{\xi Y\xi} = 0$$

and from (6.12) that $(\nabla_X\eta)Y = <Y,\nabla_X\xi> = 0$, therefore

$$(\nabla_{\xi}\Phi)(X,Y) + (\nabla_X\Phi)(\xi,Y) = -2(h_{\xi XFY} + h_{X\xi FY}) = -2(h_{\xi XFY} - h_{XFY\xi}) = -2h_{\xi XFY} = 0.$$

We have got $h_{XXY}=0$, $h_{XY\xi}=h_{\xi Y\xi}=h_{\xi XY}=0$, hence $B(\boldsymbol{C}_1) \subset \overline{\boldsymbol{T}}_1$. Since $dim\,\boldsymbol{C}_1 = \frac{1}{3}n(n-1)(n-2) = dim\,\overline{\boldsymbol{T}}_1$, then $B(\boldsymbol{C}_1) = \overline{\boldsymbol{T}}_1$

**QED.**

**LEMMA 6.10.** $B(\boldsymbol{C}_2) = \overline{\boldsymbol{T}}_2$

**Proof.** It follows from (6.12), (6.20), (6.23) that for $X,Y,Z \in V$

$$h_{XY\xi} = h_{\xi X\xi} = 0, \quad (\nabla_X \Phi)(Y,Z) = -2h_{XYFZ},$$

$$3d\Phi(X,Y,Z) = \sigma(\nabla_X \Phi)(Y,Z) = -2(h_{XYFZ} + h_{YZFX} + h_{ZXFY}) = 0.$$

By analogy with 2) of *Theorem 5.20* we get $h_{XYFZ}=h_{XFYZ}=h_{FXYZ}$, hence

$$\sigma h_{XYZ}=0, \quad 3d\Phi(\xi,Y,Z)=-2h_{\xi YFZ}=0$$

and we have obtained that $B(\boldsymbol{C}_2) \subset \overline{\boldsymbol{T}}_2$.

Since $dim\,\boldsymbol{C}_2 = \frac{2}{3}n(n-1)(n+1) = dim\,\overline{\boldsymbol{T}}_2$, then $B(\boldsymbol{C}_2) = \overline{\boldsymbol{T}}_2$

**QED.**

**LEMMA 6.11.** $B(\boldsymbol{C}_3) = \overline{\boldsymbol{T}}_3$

**Proof.** Using (6.20) for $Y,Z \in V$ we get

$$(\nabla_X\Phi)(Y,Z)-(\nabla_{FX}\Phi)(FY,Z)=-2(h_{XYFZ}-h_{FXFYFZ})=0,$$

hence $h_{\xi XY}=0$; $h_{XYZ}=h_{FXFYZ}$, $X \in V$, and

$$(\nabla_X \Phi)(\xi,Y) = -2h_{X\xi FY} + <\nabla_X \xi, FY> = 2h_{XFY\xi} - h_{XFY\xi} = h_{XFY\xi} = 0$$

therefore $h_{\xi Y\xi} = 0;\ h_{XY\xi} = 0$, for $X \in V$. It follows from (6.25) that

$$\delta\Phi(Z) = 2\beta(FZ) - h_{\xi Z\xi} = 2\beta(FZ) = 0.$$

Thus, we have obtained that $B(\boldsymbol{C}_3) \subset \overline{\boldsymbol{T}}_3$.

Since $dim\,\boldsymbol{C}_3 = n(n+1)(n-2) = dim\,\overline{\boldsymbol{T}}_3$, then $B(\boldsymbol{C}_3) = \overline{\boldsymbol{T}}_3$

**QED.**

**LEMMA 6.12.** $B(\boldsymbol{C}_4) = \overline{\boldsymbol{T}}_4$

**Proof.** One can write from *Table 6.3* that for any $X,Y,Z\in \boldsymbol{X(M)}$

$$(6.26)\quad (\nabla_X \Phi)(Y,Z) = -\frac{1}{2(n-1)}[<FX,FY>\delta\Phi(Z) - <FX,FZ>\delta\Phi(Y) - <x,FY>\delta\Phi(FZ) + <X,FZ>\delta\Phi(FY)],\quad \delta\Phi(\xi)=0.$$

For $X,Y\in V$ it follows from (6.21), (6.26) that

$$h_{XY\xi} = -(\nabla_X\Phi)(\xi,FY)=0,\ h_{\xi XY}=\frac{1}{2}(\nabla_\xi\Phi)(X,FY)=0,\ h_{\xi Y\xi}=-(\nabla_\xi\Phi)(\xi,FY)=0;$$

from (6.25) $\delta\Phi(X)=2\beta(FX)$. Using (6.21) and (6.26) for $X,Y,Z\in V$ we have

$$h_{XYZ} = \frac{1}{2}(\nabla_X\Phi)(Y,FZ) = -\frac{1}{2(n-1)}[-<X,Y>\beta(Z) + <FX,Z>\beta(FY) + <X,FY>\beta(FZ) + <X,Z>\beta(Y)]$$

and $B(\boldsymbol{C}_4)\subset\overline{\boldsymbol{T}}_4$. Since $\dim \boldsymbol{C}_4=2n=\dim \overline{\boldsymbol{T}}_4$, then $B(\boldsymbol{C}_4)=\overline{\boldsymbol{T}}_4$.

**QED.**

**LEMMA 6.13.** $B(\boldsymbol{C}_5)=\overline{T}_{10}$

**Proof.** From Table 6.3 and (6.24) it follows that

$$(\nabla_X\Phi)(Y,Z) = \frac{1}{2n}[\Phi(X,Y)\eta(Y) - \Phi(X,Y)\eta(Z)]\delta\eta = \frac{1}{2n}[<X,FZ><Y,\xi> - <X,FY><Z,\xi>]\beta(\xi).$$

From this identity and from (6.21) for $X,Y,Z\in V$ we have

$$h_{XYZ} = \frac{1}{2}(\nabla_X\Phi)(Y,FZ)=0,\ h_{\xi YZ}=\frac{1}{2}(\nabla_\xi\Phi)(Y,FZ)=0,$$

$h_{\xi Y\xi} = -(\nabla_\xi\Phi)(\xi,FY)=0$ and

$$h_{XY\xi} = -(\nabla_X\Phi)(\xi,FY) = -\frac{1}{2n}<X,F^2,Y>\beta(\xi) = \frac{1}{2n}<X,Y>\beta(\xi).$$

Thus, $B(\boldsymbol{C}_5) \subset \overline{T}_{10}$ and since $dim\ \boldsymbol{C}_5 = 1 = dim\ \overline{T}_{10}$, then $B(\boldsymbol{C}_5) = \overline{T}_{10}$.

**QED.**

**LEMMA 6.14.** $B(\boldsymbol{C}_6) = \overline{\boldsymbol{T}}_7$

**Proof.** From Table 6.3 and (6.25) it follows that

$$(\nabla_X \Phi)(Y,Z) = \frac{1}{2n}[<X,Z><Y,\xi> - <X,Y><Z,\xi>]\overline{\beta}(\xi)$$

and from (6.21) for $X,Y,Z \in V$ we obtain

$$h_{XYZ} = \frac{1}{2}(\nabla_X \Phi)(Y,FZ) = 0, \quad h_{\xi YZ} = \frac{1}{2}(\nabla_\xi \Phi)(Y,FZ) = 0,$$

$$h_{\xi Y\xi} = -(\nabla_\xi \Phi)(\xi,FY) = 0, \quad h_{XY\xi} = -(\nabla_X \Phi)(\xi,FY) = \frac{1}{2n} <X,FY> \overline{\beta}(\xi).$$

Thus, $B(\boldsymbol{C}_6) \subset \overline{\boldsymbol{T}}_7$ and since $dim\ \boldsymbol{C}_6 = 1 = dim\ \overline{\boldsymbol{T}}_7$, then $B(\boldsymbol{C}_6) = \overline{\boldsymbol{T}}_7$.

**QED.**

**LEMMA 6.15.** $B(\boldsymbol{C}_7) = \overline{\boldsymbol{T}}_5$

**Proof.** From Table 6.3, (6.23) it follows that

$$(\nabla_X \Phi)(Y,Z) = \eta(Z)(\nabla_Y \eta)FX + \eta(Y)(\nabla_{FX}\eta)Z = <Z,\xi> h_{Y\xi FX} + <Y,\xi> h_{FX\xi Z},$$

therefore for $X,Y,Z \in V$ we get from (6.21) that

$$h_{XYZ} = \frac{1}{2}(\nabla_X \Phi)(Y,FZ) = 0, \ h_{\xi YZ} = \frac{1}{2}(\nabla_\xi \Phi)(Y,FZ) = 0,$$

$$h_{\xi Y\xi} = -(\nabla_\xi \Phi)(\xi,FY) = 0, \ h_{XY\xi} = -(\nabla_X \Phi)(\xi,FY) = -h_{FX\xi FY} = h_{FXFY\xi}$$

and from (6.25) $\overline{\beta}(\xi) = \delta\Phi(\xi) = 0$. Using (6.10), (6.19) one can verify that

$$(\nabla_X \Phi)(FY,\xi) = <FY,(\nabla_X F)\xi> = <FY, h_X F\xi - Fh_X \xi>$$
$$= - <Y, h_X \xi> = -h_{X\xi Y} = h_{XY\xi}.$$

On the other hand we see from Table 6.3 that

$$(\nabla_X \Phi)(FY,\xi) = <\xi,\xi> h_{FY\xi FX} = -h_{FYFX\xi} = -h_{YX\xi}$$

and for $X,Y\in V$ $h_{XY\xi}=-h_{YX\xi}$.

Since $dim\ \boldsymbol{C}_7=n^2-1=dim\ \overline{\boldsymbol{T}}_5$ and $B(\boldsymbol{C}_7)\subset\overline{\boldsymbol{T}}_5$, hence $B(\boldsymbol{C}_7)=\overline{\boldsymbol{T}}_5$.

**QED.**

**LEMMA 6.16.** $B(\boldsymbol{C}_8)=\overline{\boldsymbol{T}}_8$.

**Proof.** From Table 6.3 and (6.23) it follows that

$$(\nabla_X\Phi)(Y,Z)=-\eta(Z)(\nabla_Y\eta)FX+\eta(Y)(\nabla_{FX}\eta)Z=-<Z,\xi>h_{Y\xi FX}+<Y,\xi>h_{FX\xi Z},$$

therefore for $X,Y,Z\in V$ from (6.21) we obtain

$$h_{XYZ}=\frac{1}{2}(\nabla_X\Phi)(Y,FZ)=0,\ h_{\xi YZ}=\frac{1}{2}(\nabla_\xi\Phi)(Y,FZ)=0,$$

$$h_{\xi Y\xi}=-(\nabla_\xi\Phi)(\xi,FY)=0,\ h_{XY\xi}=-(\nabla_X\Phi)(\xi,FY)=-h_{FX\xi FY}=h_{FXFY\xi}$$

and from (6.24) we have $\delta\eta=\beta(\xi)=0$.

Using (6.1O),(6.19) one can check that

$$(\nabla_X\Phi)(FY,\xi)=<FY,(\nabla_XF)\xi>=<FY,h_XF\xi-Fh_X\xi>=-h_{X\xi Y}=h_{XY\xi}.$$

On the other hand we have seen from *Table 6.3* that

$$(\nabla_X\Phi)(FY,\xi)=-<\xi,\xi>h_{FY\xi FX}=h_{FYFX\xi}=h_{YX\xi},$$

hence for $X,Y\in V$ $h_{XY\xi}=h_{YX\xi}$.

Since $dim\ \boldsymbol{C}_8=n^2-1=dim\ \overline{\boldsymbol{T}}_8$ and $B(\boldsymbol{C}_8)\subset\overline{\boldsymbol{T}}_8$, therefore $B(\boldsymbol{C}_8)=\overline{\boldsymbol{T}}_8$

**QED.**

**LEMMA 6.17.** $B(\boldsymbol{C}_9)=\overline{\boldsymbol{T}}_9$.

**Proof.** One can write from Table 6.3 and (6.23) that

$$(\nabla_X\Phi)(Y,Z)=\eta(Z)(\nabla_Y\eta)FX-\eta(Y)(\nabla_{FX}\eta)Z=<Z,\xi>h_{Y\xi FX}-<Y,\xi>h_{FX\xi Z},$$

therefore from (6.21) for $X,Y,Z\in V$ we obtain

$$h_{XYZ}=\frac{1}{2}(\nabla_X\Phi)(Y,FZ)=0,\ h_{\xi YZ}=\frac{1}{2}(\nabla_\xi\Phi)(Y,FZ)=0,$$
$$h_{\xi Y\xi}=-(\nabla_\xi\Phi)(\xi,FY)=0,\quad h_{XY\xi}=-(\nabla_X\Phi)(\xi,FY)=-h_{FXFY\xi},$$
$$(\nabla_X\Phi)(FY,\xi)=<\xi,\xi>h_{FY\xi FX}=-h_{FYFX\xi}=h_{YX\xi}$$

and from (6.10), (6.19)

$$(\nabla_X\Phi)(FY,\xi)=<FY,(\nabla_X F)\xi>=<FY,h_X F\xi-Fh_X\xi>=h_{XY\xi}.$$

Thus, $B(\boldsymbol{C}_9)\subset\overline{\boldsymbol{T}}_9$ and, since $dim\ \boldsymbol{C}_9=n(n+1)=dim\ \overline{\boldsymbol{T}}_9$, then $B(\boldsymbol{C}_9)=\overline{\boldsymbol{T}}_9$

**QED.**

**LEMMA 6.18.** $B(\boldsymbol{C}_{10})=\overline{\boldsymbol{T}}_6$.

**Proof.** It follows from Table 6.3 and (6.23) that

$$(\nabla_X\Phi)(Y,Z)=-\eta(Z)(\nabla_Y\eta)FX-\eta(Y)(\nabla_{FX}\eta)Z=-<Z,\xi>h_{Y\xi FX}-<Y,\xi>h_{FX\xi Z}$$

The rest is similar to ones in Lemma 6.15 – 6.17.

**QED.**

**LEMMA 6.19.** $B(\boldsymbol{C}_{11})=\overline{\boldsymbol{T}}_{12}$

**Proof.** One can write from Table 6.3 that

$$(\nabla_X\Phi)(Y,Z)=-<X,\xi>(\nabla_\xi\Phi)(FY,FZ)$$

and from (6.21) for $X,Y,Z\in V$ we have

$$h_{XYZ}=\frac{1}{2}(\nabla_X\Phi)(Y,FZ)=0,\ h_{XY\xi}=-(\nabla_X\Phi)(\xi,FY)=0,$$
$$h_{\xi Y\xi}=-(\nabla_\xi\Phi)(\xi,FY)=0.$$

It follows from (6.5) that $h_{\xi XY}=-h_{\xi YX}=-h_{\xi FXFY}$ and $B(\boldsymbol{C}_{11})\subset\overline{\boldsymbol{T}}_{12}$. Since dim $\boldsymbol{C}_{11}$ =n(n–1)=dim $\overline{\boldsymbol{T}}_{12}$, then $B(\boldsymbol{C}_{11})=\overline{\boldsymbol{T}}_{12}$.

**QED.**

**LEMMA 6.20.** $B(\boldsymbol{C}_{12})=\overline{\boldsymbol{T}}_{11}$

**Proof.** It follows from *Table 6.3* and (6.23) that

$$(\nabla_X \Phi)(Y,Z) = <X,\xi><Z,\xi> h_{\xi\xi FY} - <X,\xi><Y,\xi> h_{\xi\xi FZ}$$

and from (6.21) for $X,Y,Z \in V$ we obtain

$$h_{XYZ} = \frac{1}{2}(\nabla_X \Phi)(Y,FZ) = 0,\ h_{XY\xi} = -(\nabla_X \Phi)(\xi,FY) = 0,$$

$h_{\xi YZ} = \frac{1}{2}(\nabla_\xi \Phi)(Y,FZ) = 0$ and $B(\boldsymbol{C}_{12}) \subset \overline{\boldsymbol{T}}_{11}$.

Since $dim\ \boldsymbol{C}_{12} = 2n = dim\ \overline{\boldsymbol{T}}_{11}$, then $B(\boldsymbol{C}_{12}) = \overline{\boldsymbol{T}}_{11}$.

**QED.**

So, we have got the final

**THEOREM 6.21.** The classifications given by *Tables 6.2, 6.3* are the same up to an isomorphism and the correspondence between the classes is described by

$$\boldsymbol{C}_1 \cong \overline{\boldsymbol{T}}_1,\ \boldsymbol{C}_2 \cong \overline{\boldsymbol{T}}_2,\ \boldsymbol{C}_3 \cong \overline{\boldsymbol{T}}_3,\ \boldsymbol{C}_4 \cong \overline{\boldsymbol{T}}_4,\ \boldsymbol{C}_5 \cong \overline{\boldsymbol{T}}_{10},\ \boldsymbol{C}_6 \cong \overline{\boldsymbol{T}}_7,\ \boldsymbol{C}_7 \cong \overline{\boldsymbol{T}}_5,$$
$$\boldsymbol{C}_8 \cong \overline{\boldsymbol{T}}_8,\ \boldsymbol{C}_9 \cong \overline{\boldsymbol{T}}_9,\ \boldsymbol{C}_{10} \cong \overline{\boldsymbol{T}}_6,\ \boldsymbol{C}_{11} \cong \overline{\boldsymbol{T}}_{12},\ \boldsymbol{C}_{12} \cong \overline{\boldsymbol{T}}_{11}.$$

# CHAPTER 7

# REMARKS ON GEOMETRY OF ALMOST CONTACT METRIC MANIFOLDS

In this chapter, we consider some of the almost contact metric structures.

In §1, examples of all the classes adduced in *Table 6.2* are given. In particular, α-Sasakian and α-Kenmotsu structures are identified and it is shown that, when $\alpha = const$, they are quasi homogeneous.

§2 is devoted to the conditions of integrability, normality and to the fundamental tensor fields $N^{(1)}$, $N^{(2)}$, $N^{(3)}$, $N^{(4)}$ of a.c.m.s. We identify some of the classes studied by various authors with those obtained from *Tables 6.3,6.4*.

Riemannian locally regular σ-manifolds with one-dimensional foliations of mirrors are discussed in §3. We consider necessary and sufficient conditions for M to be a R.l.r. σ-m. and the induced a.c.m.s. on *M*. In this case, the canonical connection $\tilde{\nabla}$ of R.l.r. σ-m. and that $\overline{\nabla}$ of induced a.c.m.s. are the same. R.l.r. σ-m. of order 3,4 are studied more explicitly.

We refer to [1], [2], [11], [14], [43], [48], [73].

## *§1. ABOUT BASIC CLASSES OF CLASSIFICATION*

We give examples of structures of basic classes given in *Table 6.2*.

**$1^0$. LEMMA 7.1.** Let $T(M) = L\oplus V$ be an almost product structure determined by an a.c.m.s. $(F, \xi, \eta, g)$ and $h$ be the second fundamental tensor field of the $G$-structure corresponding to a.c.m.s. Then, $L\oplus V$ is invariant with respect to $\nabla$ if and only if

$h_{ZY\xi} = 0$ for $Z\in \mathfrak{X}(M)$, $X\in V$;

$L\oplus V$ is integrable if and only if

$h_{XY\xi} = h_{YX\xi}$ for $X,Y\in V$.

**Proof.** It is evident from the following identities

$$h_{ZX\xi} = <\nabla_Z X, \xi> = -<\nabla_Z \xi, X>,$$

$$<[X,Y],\xi>=<\nabla_X Y,\xi>-<\nabla_Y X,\xi>=h_{XY\xi}-h_{YX\xi}$$

**QED.**

**PROPOSITION 7.2.** A manifold $M$ has an a.c.m.s. of class $\overline{\boldsymbol{T}}_\alpha$, where $\alpha\in\boldsymbol{A}_4$, if and only if $M$ is a local Riemannian product $M_1\times\tilde{M}$, where $dim\,M_1=1$ and $\tilde{M}$ is an almost Hermitian manifold of the corresponding class.

**Proof.** If $h\in\overline{\boldsymbol{T}}_\alpha$, $\alpha\in\boldsymbol{A}_4$, then according to *Lemma 7.1* we have from *Table 6.2* that the almost product structure $L\oplus V$ is invariant under $\nabla$, hence the manifold $M$ is the local Riemannian product of $M_1$ and $\tilde{M}$, where $M_1$, $\tilde{M}$ are the maximal integral manifolds of the distributions $L$, $V$ passing through a fixed point of $M$. The pair $(F,g)_{|V}$ determines an almost Hermitian structure on $\tilde{M}$ and from comparison of *Tables 5.1, 6.2* it follows that $(\tilde{M},F_{|\tilde{M}},g_{|\tilde{M}})$ belongs to the corresponding class $\overline{\boldsymbol{T}}_\alpha$ as an almost Hermitian manifold. For example, if $h\in\overline{\boldsymbol{T}}_1$, then $(\tilde{M},F_{|\tilde{M}},g_{|\tilde{M}})$ is a nearly Kaehlerian manifold, etc.

Conversely, if $M$ is the local Riemannian product of $M_1$ and $\tilde{M}$, where $dim\,M_1=1$ and $(\tilde{M},J,g)$ is an almost Hermitian manifold belonging to a class $\boldsymbol{U}_\alpha$, then one can define the natural a.c.m.s. on $M$, where $\xi$ is a tangent vector field to $M_1$, $\|\xi\|=1$, $F\xi=0$, $V=[\xi]^\perp$, $F_{|V}=J$. Using *Tables 5.1, 6.2* and *Lemma 7.1* it is easy to check that the constructed natural a.c.m.s. has a type $\overline{\boldsymbol{T}}_\alpha$.

**QED.**

This proposition makes possible to construct examples of a.c.m.s. of types $\overline{\boldsymbol{T}}_\alpha$, where $\alpha\in\boldsymbol{A}_4$

**DEFINITION 7.1** [42] . A.c.m.s. is called nearly - $K$ - cosymplectic if

$(\nabla_X F)Y+(\nabla_Y F)X=0$ and $\nabla_X\xi=0$ for $X,Y\in\mathfrak{X}(M)$.

**PROPOSITION 7.3.** A.c.m.s. is nearly - $K$ - cosymplectic if and only if it is of the class $\overline{\boldsymbol{T}}_1$.

**Proof.** If a.c.m.s. is nearly - $K$ - cosymplectic, then from (6.18) we have

$$h_{XYZ}=-\frac{1}{2}<(\nabla_X F)FY,Z>=-\frac{1}{2}<(\nabla_X F)Y,FZ>=\frac{1}{2}<(\nabla_Y F)X,FZ>$$

where $X\in\boldsymbol{X}(M)$, $Y,Z\in V$, therefore $h_{XYZ}=-h_{YXZ}$ for $X,Y,Z\in V$ and $h_{\xi YZ}=0$. If $X,Y\in\boldsymbol{X}(M)$, then

$$<\nabla_X\xi,Y>=<\nabla_X\xi-\overline{\nabla}_X\xi,Y>=h_{X\xi Y}=-h_{XY\xi}=0$$

So, $h_{\xi X\xi}=h_{XY\xi}=0$ for $X,Y\in V$ and $h\in\overline{\boldsymbol{T}}_1$.

The converse is easily verified.

**QED.**

$2^0$. We consider now structures of the class $\boldsymbol{C}_5(\overline{\boldsymbol{T}}_{10})$, $\boldsymbol{C}_6(\overline{\boldsymbol{T}}_7)$. Using the results in [43] one can state the folllowing

**DEFINITION 7.2.** A.c.m.s. $(F,\xi,\eta,g)$ is called :

1) $\alpha$–Sasakian if $(\nabla_X F)Y=\alpha\{<X,Y>\xi-\eta(Y)X\}$,

2) α–Kenmotsu if $(\nabla_X F)Y=\alpha\{<FX,Y>\xi-\eta(Y)FX\}$,

where $\alpha$ is a differentiable function on $M$ and $X,Y\in\boldsymbol{X}(M)$.

**PROPOSITION 7.4.** A.c.m.s. is $\alpha$–Sasakian if and only if it is a structure of the class $\boldsymbol{C}_6(\overline{\boldsymbol{T}}_7)$.

**Proof.** If $(\nabla_X F)Y=\alpha\{<X,Y>\xi-\eta(Y)X\}$, then from (6.10) we obtain

$$(\nabla_X\Phi)(Y,Z)=<Y,(\nabla_X F)Z>=\alpha\{<X,Z>\eta(Y)-\eta(Z)<X,Y>\}.$$

Using (6.15) we have

$$\delta\Phi(\xi)=-\sum_{i=1}^{n}\{(\nabla_{E_i}\Phi)(E_i,\xi)+(\nabla_{FE_i}\Phi)(FE_i,\xi)\}-(\nabla_\xi\Phi)(\xi,\xi)=2n\alpha.$$

One can see from *Table 6.3* that the structure is of the class $\boldsymbol{C}_6$.

Conversely, if a.c.m.s. is of the class $\boldsymbol{C}_6$ and $\alpha=\frac{1}{2n}\delta\Phi(\xi)$, then for any $Z\in\boldsymbol{X}(M)$

$$<(\nabla_X F)Y, Z> = \alpha\{<X,Y><\xi,Z> - \eta(Y)<X,Z>\},$$

therefore $(\nabla_X F)Y = \alpha\{<X,Y>\xi - \eta(Y)X\}$ and the structure is $\alpha$–Sasakian.

**QED.**

**PROPOSITION 7.5.** A.c.m.s. is $\alpha$–Kenmotsu if and only if it is a structure of the class $\boldsymbol{C}_5(\overline{\boldsymbol{T}}_{10})$.

**Proof.** If $(\nabla_X F)Y = \alpha\{<FX,Y>\xi - \eta(Y)FX\}$, then it follows from (6.10) that

$$(\nabla_X \Phi)(Y,Z) = <Y,(\nabla_X F)Z> = \alpha\{<FX,Z>\eta(Y) - \eta(Z)<FX,Y>\}$$
$$= -\alpha\{\Phi(X,Z)\eta(Y) - \Phi(X,Y)\eta(Z)\}.$$

Using (6.16) , (6.12) we obtain

$$\delta\eta = -\sum_{i=1}^{n}\{(\nabla_{E_i}\eta)E_i + (\nabla_{FE_i}\eta)FE_i\} = -\sum_{i=1}^{n}\{(\nabla_{E_i}\Phi)(\xi,FE_i) - (\nabla_{FE_i}\Phi)(\xi,E_i)\} = -2n\alpha$$

Thus, $-\alpha = 1/2n\delta\eta$ and it follows from *Table 6.3* that the structure has the type $\boldsymbol{C}_5$.

Conversely, if a.c.m.s. is of the class $\boldsymbol{C}_5$ and $\alpha = -1/2n\delta\eta$ then for any $Z \in \boldsymbol{X}(M)$

$$<(\nabla_X F)Y, Z> = \alpha\{<FX,Y><\xi,Z> - \eta(Y)<FX,Z>\},$$

therefore $(\nabla_X F)Y = \alpha\{<FX,Y>\xi - \eta(Y)FX\}$ and the structure is $\alpha$–Kenmotsu.

**QED.**

**PROPOSITION 7.6.** Let a.c.m.s. be $\alpha$–Sasakian or $\alpha$–Kenmotsu and $\alpha$ be a constant. Then the corresponding pair $(P(U(n)\times 1),\ g)$ is a quasi homogeneous structure.

**Proof.** From *Table 6.2*, (6.24) , (6.25) we have

$$h_{XYZ} = h_{\xi X\xi} = h_{\xi XY} = 0, \quad X,Y \in V;$$

$h_{XY\xi} = \frac{1}{2n}<X,FY>\overline{\beta}(\xi) = \frac{1}{2n}<X,FY>\delta\Phi(\xi)$ for $\boldsymbol{C}_6 \cong \overline{\boldsymbol{T}}_7$ and

$$h_{XY\xi} = \frac{1}{2n} < X,Y > \beta(\xi) = \frac{1}{2n} < X,Y > \delta\eta \text{ for } \boldsymbol{C}_5 \cong \overline{\boldsymbol{T}}_{10}$$

From the proofs of *Propositions 7.4, 7.5* we obtain that

$h_{XY\xi} = \alpha < X,FY >$ for $\boldsymbol{C}_6$ and $h_{XY\xi} = -\alpha < X,Y >$ for $\boldsymbol{C}_7$.

As $g$ and $F$ are invariant under $\overline{\nabla}$ and $\alpha$ is constant, therefore it is obvious that $\overline{\nabla} h = 0$ in these cases.

**QED.**

**$3^0$.** The following examples one can find in [14].

1) Let $H(p,1)$, $p \geq 1$, be the generalized Heisenberg group, i.e., the group of matrices of real numbers of the form

$$a = \begin{bmatrix} 1 & A & c \\ 0 & E_p & B^T \\ 0 & 0 & 1 \end{bmatrix}$$

where $E_p$ denotes the identity $p \times p$ matrix, $A=(a_1,...,a_p)$, $B=(b_1,...,b_p) \in R^p$, $c \in R$. $H(p,1)$ is a connected simply connected nilpotent Lie group of dimension $2p+1$ which is called a generalized Heisenberg group.

A global system of coordinates $(x_i, x_{p+i}, z)$, $1 \leq i \leq p$, on $H(p,1)$ is defined by

$$x_i(a)=a_i,\ x_{p+i}(a)=b_i,\ z(a)=c.$$

A basis for the left invariant $1$–forms on $H(p,1)$ is given by

$$\alpha_i = dx_i,\ \ \alpha_{p+i} = dx_{p+i},\ \gamma = dz - \sum_{j=1}^{p} x_j dx_{p+j},$$

and its dual basis of left invariant vector fields on $H(p,1)$ is obtained by

$$X_i = \frac{\partial}{\partial x_i},\ \ X_{p+i} = \frac{\partial}{\partial x_{p+i}},\ \ Z = \frac{\partial}{\partial z},\ \ i = 1,...,p$$

A left invariant metric on $H(p,1)$ is defined by

$$g = \sum_{k=1}^{2p} \alpha_k \otimes \alpha_k + \gamma \otimes \gamma$$

and the basis $\{X_k, Z\}$, $k=1,...,2p$, is orthonormal with respect to $g$.

Let $(F, \xi, \eta, g)$ be an a.c.m.s. on $H(p,1)$ and $F_n^m$ the components of $F$ with respect to basis $\{X_k, Z\}$. Using the Riemannian connection of $g$ it is obtained:

If $Z = \xi, F_n^m$ are constant and $F_j^{p+i} = -F_{p+j}^{i}$, $F_{p+j}^{p+i} = F_j^i$, $1 \le i,j \le p$, then $(F, \xi, \eta, g)$ is stricly of the class $\boldsymbol{C}_6 \oplus \boldsymbol{C}_7$. Moreover, it is of $\boldsymbol{C}_7(\overline{\boldsymbol{T}}_5)$ if and only if $\sum_{i=1}^{p} F_i^{p+i} = 0$, and it is of $\boldsymbol{C}_6(\overline{\boldsymbol{T}}_7)$ if and only if $F_i^{p+i} = F_j^{p+j} = \lambda$, where $\lambda$ is a nonzero constant, and the other components of $F$ are zero.

2) The generalized Heisenberg group $H(1,r)$, $r>1$, is the Lie group of real matrices of the form

$$a = \begin{bmatrix} E_r & A^T & c \\ 0 & 1 & B^T \\ 0 & 0 & 1 \end{bmatrix}$$

where $E_r$ denotes the identity $r \times r$ matrix, $A=(a_1,...,a_r)$, $B=(b_1,...,b_r) \in R^r$ and $c \in R$. $H(1,r)$ is a connected simply connected nilpotent Lie group of dimension $2r+1$ and the dimension of its center is $r>1$.

A global system of coordinates $(x_i, x_{r+i}, z)$, $1 \le i \le r$, on $H(1,r)$ is defined by

$$x_i(a) = a_i, \quad x_{r+i}(a) = b_i, \quad z(a) = c$$

A basis for the left invariant $1$–forms of $H(1,r)$ is given by

$$\alpha_i = dx_i, \quad \alpha_{r+i} = dx_{r+i}, \quad \gamma = dz$$

and its dual basis by

$$X_i = \frac{\partial}{\partial x_i}, \quad X_{r+i} = \frac{\partial}{\partial x_{r+i}}, \quad Z = \frac{\partial}{\partial z} + \sum_{j=1}^{r} x_j \frac{\partial}{\partial x_{r+j}}$$

This basis is orthonormal with respect to the left invariant metric defined by

$$g = \sum_{k=1}^{2r} \alpha_k \otimes \alpha_k + \gamma \otimes \gamma .$$

Let $(F, \xi, \eta, g)$ be an a.c.m.s. on $H(1,r)$ and $F_n^m$ the components of $F$ with respect to basis $\{X_k, Z\}$, $k=1,...,2r$. Using the Riemannian connection of the metric $g$ it is obtained :

if $Z = \xi$, $F_n^m$ are constant and $F_j^{r+i} = -F_{r+j}^i = 0$, $F_{r+j}^{r+i} = F_j^i$, then $(F, \xi, \eta, g)$ is of the class $\boldsymbol{C}_8(\overline{\boldsymbol{T}}_8)$;

if $Z = \xi$, $F_n^m$ are constant and $F_j^{r+i} = -F_{r+j}^i$, $F_{r+j}^{r+i} = F_j^i = 0$, then $(F, \xi, \eta, g)$ is of the class $\boldsymbol{C}_9(\overline{\boldsymbol{T}}_9)$.

3) Let $G$ be the Lie group of real matrices of the form

$$a = \begin{bmatrix} e^{-z} & 0 & x \\ 0 & e^{z} & y \\ 0 & 0 & 1 \end{bmatrix}$$

with the left invariant metric

$$g = e^{2z}dx^2 + e^{-2z}dy^2 + \lambda^2 dz^2, \quad \lambda > 0 .$$

$(G,g)$ is a $4$–symmetric space, which is isomorphic to the semi-direct product of $R$ and $R^2$, both with the additive group structure, and where the action of $R$ and $R^3$ is given by the matrix

$$\begin{bmatrix} e^{z} & 0 \\ 0 & e^{-z} \end{bmatrix}$$

i.e., the group $E(1,1)$ of rigid motions of the Minkowski $2$–space. With respect to the metric $g$, the basis of invariant vector fields $\{X_1, X_2, X_3\}$ given by

$$X_1 = e^{-z}\frac{\partial}{\partial x}, \quad X_2 = e^{z}\frac{\partial}{\partial y}, \quad X_3 = \frac{1}{\lambda}\frac{\partial}{\partial z}$$

is orthonormal.

It is verified that an a.c.m.s. $(F, \xi, \eta, g)$ on $G$ is of the class $\boldsymbol{C}_{12}(\overline{\boldsymbol{T}}_{11})$ if $\xi = X_1$ or $\xi = X_2$; and it is of the class $\boldsymbol{C}_9(\overline{\boldsymbol{T}}_9)$ if $\xi = X_3$.

4) Let $G$ be the complex matrix group $G$ of the form

$$a=\begin{bmatrix} e^{it} & 0 & z \\ 0 & e^{-it} & w \\ 0 & 0 & 1 \end{bmatrix}$$

Here $z$, $w$ denote complex variables and $t$ a real variable. This Lie group is diffeomorphic to $C^2(z,w)\times R(t)$. A left invariant metric on $G$ is defined by

$$g = dz\cdot d\bar{z} + dw\cdot d\overline{w} + dt^2$$

The vector fields $\{Z_1, \overline{Z}_1, Z_2, \overline{Z}_2, W\}$ given by

$$Z_1 = e^{it}\frac{\partial}{\partial z},\quad Z_2 = e^{-it}\frac{\partial}{\partial w},\quad W=\frac{\partial}{\partial t}$$

are invariant under the action of $G$ and they form an orthonormal basis of the Lie algebra of $G$. Put

$$X_1=\sqrt{2}\,Re(Z_1+Z_2),\qquad X_2=\sqrt{2}\,Im(Z_1+Z_2),$$
$$X_3=\sqrt{2}\,Im(Z_2-Z_1),\qquad X_4=\sqrt{2}\,Re(Z_1-Z_2).$$

Identifying $C^2\times R$ with $R^5$ with invariant Riemannian metric obtained from $g$, it follows that $\{X_1, X_2, X_3, X_4, W\}$ is an orthonormal basis on this space.

Let $(F, \xi, \eta, g)$ be an a.c.m.s. on $G$ and $F^i_j$ the components of $F$ with respect to $\{X_1, X_2, X_3, X_4, W\}$. It is obtained

a) If $\xi = W$, $F^i_j$ are constant and $F^3_2 = F^4_1$, $F^1_2 = F^3_4$, then $(F, \xi, \eta, g)$ is cosymplectic.

b) If $\xi = W$, $F^i_j$ are constant and $F^3_2 \neq F^4_1$, $F^1_2 \neq F^3_4$, then $(F, \xi, \eta, g)$ is of the class $\boldsymbol{C}_{11}(\overline{\boldsymbol{T}}_{12})$.

An example of a.c.m.s. satisfying the last condition, and so of the class $\boldsymbol{C}_{11}$, is the following:

$$FZ_1 = iZ_1,\quad FZ_2 = iZ_1,\quad F\overline{Z}_1 = -i\overline{Z}_2$$
$$F\overline{Z}_2 = -i\overline{Z}_1,\quad \xi = W,\qquad \eta = dt.$$

## *§2. ON SOME CLASSES OF ALMOST CONTACT METRIC MANIFOLDS*

We identify some classes given in *Tables 6.2, 6.3* with those studied in literature.

**1⁰**. Let *N(F)* be the Nijenhuis tensor field of *F*. From [42] and (6.19) it follows that for $X,Y\in \mathfrak{X}(M)$

(7.1) $N(F)(X,Y)=(\nabla_{FX}F)Y-F(\nabla_X F)Y-(\nabla_{FY}F)X+F(\nabla_Y F)X$

$$=h_{FX}FY-h_{FY}FX+F^2(h_XY-h_YX)+F(h_YFX-h_{FX}Y+h_{FY}X-h_XFY)$$

From *Theorem 1.5* we see that *F* is integrable if and only if *N(F)=0* on *M*.

**THEOREM 7.7.** An a.c.m.s. is integrable if and only if

$$h\in\overline{\boldsymbol{T}}_{In}=\overline{\boldsymbol{T}}_3\oplus\overline{\boldsymbol{T}}_4\oplus\overline{\boldsymbol{T}}_8\oplus\overline{\boldsymbol{T}}_{10}\oplus\overline{\boldsymbol{T}}_{11}.$$

**Proof.** If *F* is integrable, then the structure of almost product is also integrable and $h_{XY\xi}=h_{YX\xi}$, $X,Y\in V$, according to *Lemma 7.1*. Let $\widetilde{M}$ be some maximal integral manifold of the distribution *V*, then the restriction of *F* on $\widetilde{M}$ determines the almost complex structure $\widetilde{F}=F_{|\widetilde{M}}$ which is integrable too $(N(F)_{|V}=0)$, hence $(\widetilde{M},g,\widetilde{F})$ is Hermitian. In this case, one can see from *Table 5.1* that

(7.2) $h\in\overline{\boldsymbol{T}}_3\oplus\overline{\boldsymbol{T}}_4=\{h:h_{XYZ}=h_{FXFYZ},\ X,Y,Z\in V\}.$

Further, $N(F)(X,\xi)=F^2(h_X\xi-h_\xi X)+F(h_\xi FX-h_{FX}\xi)$,

$$<N(F)(X,\xi),Y>=-h_{X\xi Y}+h_{\xi XY}-h_{\xi FXFY}+h_{FX\xi FY}$$
$$=h_{XY\xi}-h_{FXFY\xi}+2h_{\xi XY}=0,\quad X,Y\in V.$$

By analogy, we have $h_{YX\xi}-h_{FYFX\xi}+2h_{\xi XY}=0$. Adding these equalities and taking into consideration that $h_{XY\xi}=h_{YX\xi}$ we obtain $h_{XY\xi}=h_{FXFY\xi}$ and $h_{\xi XY}=0$. Thus $h\in\overline{\boldsymbol{T}}_{In}$.

Conversely, if $h\in\overline{\boldsymbol{T}}_{In}$, then we have to check that *N(F)=0*. From (7.1),(7.2) it follows that

$$< N(F)(X,Y),\xi >= h_{FXFY\xi} - h_{FYFX\xi} = 0$$

$$< N(F)(X,Y),Z >= h_{FXFYZ} - h_{FYFXZ} - h_{XYZ} + h_{YXZ} - h_{YFXFZ} + h_{FXYFZ} - h_{FYXFZ}$$
$$+ h_{XFYFZ} = 2(h_{YXZ} - h_{FYFXZ} + h_{FXFYZ} - h_{XYZ}) = 0$$

Finally, $< N(F)(X,\xi),\xi >= 0$,

$$< N(F)(X,\xi),Y >= h_{XY\xi} - h_{FXFY\xi} + 2h_{\xi XY} = 0 .$$

**QED.**

So, there exist 32 classes of integrable a.c.m.s.

**2⁰**. The following tensor fields play an important role in geometry of a.c.m..s., see [8],[73].

$$N^{(1)}(X,Y)=N(F)(X,Y)+2d\eta(X,Y),\ \ N^{(2)}(X,Y)=(L_{FX}\eta)Y-(L_{FY}\eta)X,$$

$$N^{(3)}(X)=(L_{\xi}F)X,\quad N^{(4)}(X)=(L_{\xi}\eta)X,\quad X,Y\in \mathcal{X}(M)$$

**LEMMA 7.8.** $N^{(1)}(X,Y)=N(F)(X,Y)+(h_{YX\xi}-h_{XY\xi}),\ \ X,Y\in \mathcal{X}(M).$

**Proof.** Using (6.13) and (6.23) we have

$$2d\eta(X,Y) = (\nabla_X \eta)Y - (\nabla_Y \eta)X = h_{X\xi Y} - h_{Y\xi X} = h_{YX\xi} - h_{XY\xi} .$$

**QED.**

**LEMMA 7.9.** $N^{(2)}(X,Y) = h_{YFX\xi} - h_{FXY\xi} - h_{XFY\xi} + h_{FYX\xi}$, $X,Y\in \mathcal{X}(M)$

**Proof.**

$$(L_{FX}\,\eta)Y = FX<Y,\ \xi> - <[FX,Y],\ \xi>= <\nabla_{FX}Y - [FX,Y],\ \xi> + <Y,\nabla_{FX}\,\xi>$$
$$=< \nabla_Y FX - \overline{\nabla}_Y FX ,\xi > + < Y ,\nabla_{FX}\,\xi - \overline{\nabla}_{FX}\,\xi >$$
$$= h_{YFX\xi} + h_{FX\xi Y} = h_{YFX\xi} - h_{FXY\xi} .$$

**QED.**

**LEMMA 7.10.** $N^{(3)}(X) = 2(h_{\bar{\xi}} FX - Fh_{\bar{\xi}} X),\quad X \in \mathcal{X}(M).$

**Proof.** Using (6.19) we obtain

$$(L_\xi F)(X)=[\xi,FX]-F[\xi,X]=\nabla_\xi FX\xi-\nabla_{FX}\xi-F(\nabla_\xi X-\nabla_X\xi)$$
$$=\nabla_\xi FX-F\nabla_\xi X-(\nabla_{FX}\xi-\overline{\nabla}_{FX}\xi)+F(\nabla_X\xi-\overline{\nabla}_X\xi)$$
$$=(\nabla_\xi F)X-h_{FX}\xi+Fh_X\xi=h_\xi FX-Fh_\xi X-h_{FX}\xi+Fh_X\xi$$
$$=2h_\xi^- FX-2Fh_\xi^- X.$$

**QED.**

**LEMMA 7.11.** $N^{(4)}(X)=-h_{\xi X\xi},\ X\in\mathfrak{X}(M).$

**Proof.**

$$(L_\xi\eta)(X)=\xi<X,\xi>-<[\xi,X],\xi>=<\nabla_\xi X,\xi>-<[\xi,X],\xi>+<X,\nabla_\xi\xi>$$
$$=<\nabla_X\xi,\xi>+<X,\nabla_\xi\xi-\overline{\nabla}_\xi\xi>=h_{\xi\xi X}=-h_{\xi X\xi}.$$

**QED.**

**THEOREM 7.12.** Let a.c.m.s. be a quasi homogeneous structure, i.e., $\overline{\nabla}h=0$. Then $\overline{\nabla}N^{(1)}=\overline{\nabla}N^{(2)}=\overline{\nabla}N^{(3)}=\overline{\nabla}N^{(4)}=\overline{\nabla}N(F)=0$.

**Proof.** Let $Z$ be a vector fied and $\gamma$ a curve segment in $M$ defined by $Z$ or more precisely by a local $1$–parameter group of transformations induced by $Z$. We denote by $X, Y$ the vector fields defined on some neighbourhood of $\gamma$ which are obtained by the parallel translation of $X_p, Y_p\in T_p(M)$ along $\gamma$, $p\in\gamma$, in the connection $\overline{\nabla}$. So, we have $(\overline{\nabla}_Z X)_p=(\overline{\nabla}_Z Y)_p=0$. In this case

$$[(\overline{\nabla}_Z N^{(i)})(X,Y)]_p=[\overline{\nabla}_Z N^{(i)}(X,Y)]p \text{ and } [(\overline{\nabla}_Z h)(X,Y)]_p=[\overline{\nabla}_Z h_X Y]_p=0.$$

We know that $\overline{\nabla}\xi=\overline{\nabla}F=\overline{\nabla}g=0$. From *Lemmas 7.8 - 7.11* and (7.1) it follows that $[\overline{\nabla}_Z N^{(i)}(X,Y)]_p=0$, where $N^{(i)}$ denotes one of the tensor fields above.

**QED.**

**DEFINITION 7.3** [8],[73] . A.c.m.s. is called normal if $N^{(1)}=0$ on $M$.

**THEOREM 7.13**. A.c.m.s. is normal if and only if

$$h\in\overline{\boldsymbol{T}}_N=\overline{\boldsymbol{T}}_5\oplus\overline{\boldsymbol{T}}_8\oplus\overline{\boldsymbol{T}}_7\oplus\overline{\boldsymbol{T}}_{10}.$$

**Proof.** A.c.m.s. is normal if and only if

(7.3) $N(F)(X,Y)=(h_{XY\xi}-h_{YX\xi})\xi,\ \ N(F)(X,\xi)=-h_{\xi X\xi}\xi,\ \ X,Y\in V$.

Using (7.1) and (7.3) we obtain

(7.4) $h_{FX}FY-h_{FY}FX=(h_{XY\xi}-h_{YX\xi})\xi$,

(7.5) $F^2(h_XY-h_YX)+F(h_YFX-h_{FX}Y+h_{FY}X-h_XFY)=0$,

(7.6) $N(F)(X,\xi)=F^2(h_X\xi-h_\xi X)+F(h_\xi FX-h_{FX}\xi)=-h_{\xi X\xi}=0$.

It follows from (7.4) that

(7.7) $h_{FXFY\xi}-h_{FYFX\xi}=h_{XY\xi}-h_{YX\xi}$ or $h_{\bar{F}XFY\xi}=h_{\bar{X}Y\xi}$,

(7.8) $h_{FXFYZ}-h_{FYFXZ}=0$ or $h_{XYZ}=h_{YXZ}$, $X,Y,Z\in V$.

The condition (7.8) is equivalent to the folowing one

$h_{XYZ}=h_{YXZ}=-h_{YZX}=-h_{ZYX}=h_{ZXY}=h_{XZY}=-h_{XYZ}$ or

(7.9) $h_{XYZ}=0$.

We remark that (7.9) implies (7.5) and $h_{\xi X\xi}=0$ from (7.6). It also follows from (7.6) that

$$-h_{X\xi Y}+h_{\xi XY}-h_{\xi FXFY}+h_{FX\xi FY}=h_{XY\xi}-h_{FXFY\xi}+2h_{\xi XY}=0$$

By analogy $h_{YX\xi}-h_{FYFX\xi}+2h_{\xi YX}=0$

Adding the last two equalities we have

$$h_{FXFY\xi}+h_{FYFX\xi}=h_{XY\xi}+h_{YX\xi}$$

and adding obtained one with (7.7) we get

(7.10) $h_{FXFY\xi}=h_{XY\xi}$ and $h_{\xi XY}=0$.

Thus, $h\in\overline{\boldsymbol{T}}_N$.

Conversely, if $h\in\overline{\boldsymbol{T}}_N$, then (7.10), (7.9) are fulfilled and $h_{\xi X\xi}=0$ therefore

(7.4), (7.5), (7.6) are realized and (7.3) follows.

**QED.**

$3^0$. Let $X,Y,Z\in\mathfrak{X}(M)$, then a.c.m.s. $(F, \xi, \eta, g)$ is said to be:

Almost cosymplectic if $d\Phi=0$ and $d\eta=0$.

Quasi Sasakian if $d\Phi=0$ and a.c.m.s. is normal.

Nearly-K-cosymplectic if $(\nabla_X F)Y+(\nabla_Y F)X=0$ and $\nabla_X\xi=0$.

Quasi-K-cosymplectic if $(\nabla_X F)Y+(\nabla_{FX} F)FY=\eta(Y)\nabla_{FX}\xi$

Semi cosymplectic if $\delta\Phi=0$ and $\delta\eta=0$.

Trans-Sasakian if
$$(\nabla_X\Phi)(Y,Z)=-\frac{1}{2n}\{(<X,Y>\eta(Z)-<X,Z>\eta(Y))\delta\Phi(\xi)$$
$$+(<X,FY>\eta(Z)-<X,FZ>\eta(Y))\delta\eta\}.$$

Nearly-trans-Sasakian if

a) $(\nabla_X\Phi)(X,Y)=-\frac{1}{2n}\{<X,X>\delta\Phi(Y)-<X,Y>\delta\Phi(X)+<FX,Y>\eta(X)\delta\eta\}$,

b) $(\nabla_X\eta)Y=-\frac{1}{2n}\{<FX,FY>\delta\eta+<FX,Y>\delta\Phi(\xi)\}$.

Almost-K-contact if $\nabla_\xi F=0$.

In [14], it is explained how these classes above studied by various authors coincide with those introduced in *Tables 6.2, 6.3*.

$\boldsymbol{C}_2\oplus\boldsymbol{C}_9\cong\overline{\boldsymbol{T}}_2\oplus\overline{\boldsymbol{T}}_9$ = the class of almost cosymplectic manifolds.

$\boldsymbol{C}_5\oplus\boldsymbol{C}_6\cong\overline{\boldsymbol{T}}_7\oplus\overline{\boldsymbol{T}}_{10}$ = the class of trans-Sasakian manifolds.

$\boldsymbol{C}_6\oplus\boldsymbol{C}_7\cong\overline{\boldsymbol{T}}_5\oplus\overline{\boldsymbol{T}}_7$ = the class of quasi-Sasakian manifolds.

$\boldsymbol{C}_3\oplus\boldsymbol{C}_7\oplus\boldsymbol{C}_8\cong\overline{\boldsymbol{T}}_3\oplus\overline{\boldsymbol{T}}_5\oplus\overline{\boldsymbol{T}}_8$ = the class of semi-cosymplectic and normal manifolds.

$\boldsymbol{C}_1\oplus\boldsymbol{C}_5\oplus\boldsymbol{C}_6\cong\overline{\boldsymbol{T}}_1\oplus\overline{\boldsymbol{T}}_7\oplus\overline{\boldsymbol{T}}_{10}$ = the class of nearly-trans-Sasakian manifolds.

$\boldsymbol{C_1} \oplus \boldsymbol{C_2} \oplus \boldsymbol{C_9} \oplus \boldsymbol{C_{10}} \cong \overline{\boldsymbol{T}}_1 \oplus \overline{\boldsymbol{T}}_2 \oplus \overline{\boldsymbol{T}}_6 \oplus \overline{\boldsymbol{T}}_9$ = the class of quasi-K-cosymplectic manifolds.

$\bigoplus_{i \neq 11,12} \boldsymbol{C_i} \cong \bigoplus_{i \neq 11,12} \overline{\boldsymbol{T}}_i$ = the class of almost-K-contact manifolds.

$\bigoplus_{i \neq 4,5,6} \boldsymbol{C_i} \cong \bigoplus_{i \neq 4,7,10} \overline{\boldsymbol{T}}_i$ = the class of semi-cosymplectic manifolds.

**DEFINITION 7.4.** The structure affinor $F$ is said to be:

a) of $V$-invariant type if for all $X,Y \in V$, $(\nabla_X F)Y \in V$;

b) of $V$-antiinvariant type if for all $X,Y \in V$, $(\nabla_X F)Y \perp V$;

c) of $\xi$-antiinvariant type if for every $X \in V$, $(\nabla_\xi F)X \perp V$;

d) $V$-parallel if for all $X,Y \in V$, $(\nabla_X F)Y = 0$.

**PROPOSITION 7.14** [2] . Let $M$ be an almost contact metric manifold. Then

a) $F$ is of $V$-invariant type if and only if $h \in \overline{\boldsymbol{T}}_1 \oplus \overline{\boldsymbol{T}}_2 \oplus \overline{\boldsymbol{T}}_3 \oplus \overline{\boldsymbol{T}}_4 \oplus \overline{\boldsymbol{T}}_{11} \oplus \overline{\boldsymbol{T}}_{12}$;

b) $F$ is of $V$-antiinvariant type if and only if $h \in \bigoplus_{i \neq 1,2,3,4} \overline{\boldsymbol{T}}_i$

c) $F$ is of $\xi$-antiinvariant type if and only if $h \in \bigoplus_{i \neq 12} \overline{\boldsymbol{T}}_i$;

d) $F$ is $V$-parallel if and only if $h \in \overline{\boldsymbol{T}}_{11} \oplus \overline{\boldsymbol{T}}_{12}$.

**PROPOSITION 7.15** [2] . Let M be an almost contact metric manifold. Then

a) $h \in \overline{\boldsymbol{T}}_{11}$ if and only if $F$ is of $\xi$-antiinvariant type and $F$ is $V$-parallel.

b) $h \in \overline{\boldsymbol{T}}_5 \oplus \overline{\boldsymbol{T}}_6 \oplus \overline{\boldsymbol{T}}_7 \oplus \overline{\boldsymbol{T}}_8 \oplus \overline{\boldsymbol{T}}_9 \oplus \overline{\boldsymbol{T}}_{10} \oplus \overline{\boldsymbol{T}}_{12}$ if and only if $F$ is of $V$-antiinvariant type and the integral curves of $\xi$ are geodesics of $\nabla$.

c) $h \in \overline{\boldsymbol{T}}_1 \oplus \overline{\boldsymbol{T}}_2 \oplus \overline{\boldsymbol{T}}_3 \oplus \overline{\boldsymbol{T}}_4$ if and only if $F$ is of $V$-invariant type, F is of $\xi$-antiinvariant type and the integral curves of $\xi$ are geodesics.

The class $\overline{\boldsymbol{T}}_4$ is the analog of the class of conformally Kaehlerian manifolds

(see *Chapter 5, §3, $2^0$*).

**PROPOSITION 7.16** [2]. The class $\overline{\boldsymbol{T}}_4$ is characterized by the conditions:

a) $F$ is of $V$-invariant type;

b) $F$ is of $\xi$-invariant type;

c) the integral curves of $\xi$ are geodesics;

d) $(\nabla_X F)Y$ is in span $\{X, FX, Y, FY\}$ whenever $X,Y\in V$ and $X\perp Y,FY$.

In problems concerning the conformal change of the metric $g$, the class $\overline{\boldsymbol{T}}_4 \oplus \overline{\boldsymbol{T}}_7 \oplus \overline{\boldsymbol{T}}_{10} \oplus \overline{\boldsymbol{T}}_{11}$ plays the same role as the class of conformally Kaehlerian manifolds in the case of almost Hermitian manifolds.

**REMARK.** The text above shows that we can develop geometry of a.c.m.s. using the classification given in *Table 6.2*. It is isomorphic to one considered in *Table 6.3*. Our classification has a preference because it is given in terms of the tensor field $h=\nabla-\overline{\nabla}$, i. e., the canonical connection $\overline{\nabla}=\nabla-h$ has been constructed for every class. To study these classes of almost contact metric manifolds the various curvature tensor fields related to $\overline{\nabla}$ can be applied, see *Chapter 2, §1*. For example, with help of $\overline{\nabla}$ one can obtain characteristic classes.

## *§3. ALMOST CONTACT METRIC STRUCTURES ON RIEMANNIAN REGULAR σ–MANIFOLDS*

**$1^0$**. Let $(M,\{s_x\})$ be a R.l.r. σ-m. with one-dimensional disribution of mirrors $T^1(M)=[\xi]$, where $\xi\in\boldsymbol{X}(M)$, $\|\xi\|=1$, $T^2(M)=T^1(M)^\perp$ and $\pi_i$ $i=1,2$, are projections of $T(M)$ on $T^i(M)$.

For every point $p\in M$ we can choose such an open ball $B_p(R)$ of the radius $R$ that $B_p(R)\times B_p(R)\subset U$, where $U$ was considered in *Lemmas 4.57, 4.58*, therefore $s_x(y)=\mu(x,y)$ is defined for any $x,y\in B_p(R)$. Taking a ball $B_p(R/2)$ we see that $s_x(\mu(y,z))=x\cdot(y\cdot z)$ is defined for any $x,y,z\in B_p(R/2)$. So, for some concrete $k$ there exists such an open ball $B_p^0(R/2^k)$ of $p$ that all the theory developed in *Chapter*

*4, §3,* $1^0$, $2^0$, $3^0$ is true in $B_p^0$ and we obtain a localizaion of these results. Using localizations of *Propositions 4.26–4.32* we get a local infinitesimal automorphism $L(X)$ on $B_p^0$, where $X\in T_p(M)$, which is defined by

$$(7.11)\quad L(X)(x)=(I_p-S_p)^{-1}\pi_2 X\cdot(p^{-1}\cdot x),\quad x\in B_p^0,\ Y\in \mathfrak{X}(B_p^0),$$

and $\mathcal{L}_{L(X)}S=\mathcal{L}_{L(X)}g=0$ on $B_p^0$ for every such a vector $X$. Moreover, the canonical connection $\tilde{\nabla}$ of R.l.r. $\sigma$–m. $(M,\{s_x\})$ is given by

$$(7.12)\quad \tilde{\nabla}_X Y=\nabla_{\pi_1 X}Y+[L(X),Y](p)=\nabla_X Y-(\nabla_{(I-S)^{-1}\pi_2 X}S)(S^{-1}Y),$$

$$X\in T_p(M),\quad Y\in \mathfrak{X}(B_p^0).$$

**PROPOSITION 7.17**. $\mathcal{L}_{L(X)}\xi=0$ on $B_p^0$ for every $p\in M$ and $\mathcal{L}_\xi S=0$ on $M$.

**Proof.** Using that $\mathcal{L}_{L(X)}S=0$ on $B_p^0$ we obtain

$$\mathcal{L}_{L(X)}S\xi=[L(X),S\xi]=[L(X),S\xi]-S[L(X),\xi]+S[L(X),\xi]=(\mathcal{L}_{L(X)}S)\xi+S(\mathcal{L}_{L(X)}\xi)$$

$$=S(\mathcal{L}_{L(X)}\xi)=[L(X),\xi]=\mathcal{L}_{L(X)}\xi$$

So, $S(\mathcal{L}_{L(X)}\xi)=\mathcal{L}_{L(X)}\xi$ and $\mathcal{L}_{L(X)}\xi\in T^1$ or $\mathcal{L}_{L(X)}\xi=\alpha\xi$. Further,

$$(\mathcal{L}_{L(X)}g)(\xi,\xi)=L(X)<\xi,\xi>-<[L(X),\xi],\xi>-<\xi,[L(X),\xi]>$$

$$=2<\alpha\xi,\xi>=-2\alpha=0\,.$$

Thus, $\mathcal{L}_{L(X)}\xi=[L(X),\xi]=-\mathcal{L}_\xi L(X)=0$ on $B_p^0$. It follows from *Proposition 4.27* that $SL(X)=L(SX)$ on $B_p^0$. Further, we have

$$((\mathcal{L}_\xi S)L(X))(p)=([\xi,SL(X)]-S[\xi,L(X)])(p)=([\xi,SL(X)])(p)$$

$$=([\xi,L(SX)])(p)=0.$$

Since $L(X)(p)=\pi_2 X$, where $X\in T_p(M)$, then $(\mathcal{L}_\xi S)X=0$ for any $X\in T_p^2(M)$.

$((L_{\xi} S)\xi)(p) = ([\xi, S\xi] - S[\xi, \xi])(p) = 0.$

So, $(L_{\xi} S)(p)=0$ for every point $p$ of $M$.

**QED.**

**$2^0$**. It is easy to see that the distribution $T^2$ is integrable if and only if the Nijenhuis tensor field $N(\pi_2)$ vanishes on $M$. In general case it is clear that $M^0=\{x\in M:(N(\pi_2))(x)=0\}$ is a close subset of $M$ and $M'=\{x\in M:(N(\pi_2))(x)\neq 0\}$ has the induced structure of an open submanifold of $M$ with $dim\, M'=dim\, M$. For $M'$ *Definition 4.10* holds, therefore $M'$ is an invariant submanifold of $(M,\{s_x\})$ and $(M',\{s_x\})$ is also a R.l.r. $\sigma$–m. according to *Theorem 4.60*, maybe non-connected.

**PROPOSITION 7.18.** We have $\tilde{\nabla} g = \tilde{\nabla} h = \tilde{\nabla}\tilde{R} = 0$ on $M'$, where $\tilde{\nabla}$ is defined by (7.12), $h = \nabla - \tilde{\nabla}$ and $\tilde{R}$ is the curvature tensor field of $\tilde{\nabla}$, i.e., for every point $p\in M$ a connected component $M'(p)$ in $M$ containing $p$ is a Riemannian locally homogeneous manifold.

**Proof.** For every point $p\in M'$ there exists such an open ball $B_p^0 \subset M'$, which has been considered in $1^0$. Let $X,Y \in T_p^2(M')$ be such vectors that $N(\pi_2)(X,Y)\neq 0$, then $\pi_1[L(X),L(Y)]\neq 0$ on some open ball $B_p' \subset B_p^0$, where $L(X)$, $L(Y)$ are defined on $B_p^0$ by (7.11). It is clear because $L(Z)(p)=\pi_2 Z$ for any $Z\in T_p(M')$.

It follows from Proposition 4.52 that $\tilde{\nabla}_W h = \tilde{\nabla}_W \tilde{T} = \tilde{\nabla}_W \tilde{R} = 0$ for every $W\in T^2$. Let $K$ denote $h$, $\tilde{T}$ or $\tilde{R}$. The vector fields $L(X)$, $L(Y)$ are local infinitesimal affine transformations of $(M,\tilde{\nabla})$, see *Theorem 4.56*, therefore from (4.20) $L_{L(X)}K=0$ and $L_{L(Y)}K=0$. As we can see in [46]

$$(\tilde{\nabla}_{[L(X),L(Y)]}K)(p) = [L_{L(X)}(\tilde{\nabla}_{L(Y)}K) - \tilde{\nabla}_{L(Y)}(L_{L(X)}K)](p)=0$$

and $(\pi_1[L(X), L(Y)])(p)\neq 0$. So, $\tilde{\nabla} g = \tilde{\nabla} h = \tilde{\nabla}\tilde{R} = 0$ for any point $p\in M'$. Using (2.13) and covering a segment of a curve between two arbitrary points of $M'(p)$ by a finite number of balls like $B_p'$ we get the rest.

**QED.**

Let $M'= M$, then we have obtained from *Corollary 4.53* and *Proposition 7.18*

$$\text{a) } \tilde{\nabla} g = \tilde{\nabla} S = \tilde{\nabla} h = \tilde{\nabla} \tilde{T} = \tilde{\nabla} \tilde{R} = 0;$$

(7.13)

$$\text{b) } S(h) = h, \quad S(\tilde{T}) = \tilde{T}, \quad S(\tilde{R}) = \tilde{R}, \quad S(g) = g.$$

**DEFINITION 7.5.** Let $(M,\{s_x\})$ be a R.l.r. σ–m. with one-dimensional distribution of mirrors $T^1(M)$. We call $M$ a R.l.r. σ–m. of maximal torsion if $N(\pi_2) \neq 0$ on $M$, i.e., $M' = M$, and of minimal torsion if $T^2$ is an integrable distribution on $M$, i.e., $N(\pi_2) = 0$ on $M$ and $M^0 = M$.

**THEOREM 7.19.** Let $(M,g)$ be a Riemannian manifold with an affinor $S$, $S(g)=g$. If $\tilde{\nabla}$ is such a connection on $M$ that (7.13) are realized, where $h = \nabla - \tilde{\nabla}$, $\tilde{T}$, $\tilde{R}$ are the tensor fields of torsion and curvature of $\tilde{\nabla}$ respectively, then there exists such a structure $\{s_x\}$ of R.l.r. σ–m. on $(M, g)$ that $S_x = (s_x)_{*x}$ and $\tilde{\nabla}$ is a canonical connection of $(M,\{s_x\})$, see *Definition 4.1*.

**Proof** is similar to that considered in [48] for the case of locally regular *s*-manifolds.

In this theorem and in the following one the distribution $T^1$ is not required to be one-dimensional.

**THEOREM 7.20.** Let $(M, g)$ be a Riemannian manifold with an *O*-deformable affinor $S$, $S(g) = g$, and let $T^2$ be an integrable distribution, where $T^2 = T^{1\perp}$ and $T^1 = \{X \in \mathfrak{X}(M) : SX = X\}$. If $\tilde{\nabla}$ is such a connection on $M$ that the following conditions hold

(7.14)

$$\text{a) } \tilde{\nabla} g = \tilde{\nabla} S = 0;$$

$$\text{b) } S(h) = h, \; S(\tilde{R}) = \tilde{R};$$

$$\text{c) } \tilde{\nabla}_X h = \tilde{\nabla}_X \tilde{R} = 0 \text{ for any } X \in T^2;$$

where $h = \nabla - \tilde{\nabla}$, $\tilde{T}$, $\tilde{R}$ are the tensor fields of torsion and curvature of $\tilde{\nabla}$ respectively, then there exists such a structure $\{s_x\}$ of R.l.r. σ–m. on $(M, g)$ that $S_x = (s_x)_{*x}$ and $\tilde{\nabla}$ is a canonical connection of $(M,\{s_x\})$.

**Proof** will be given step by step in the next paragraphs.

1) A proof of *Theorem 4.1* is true in our case because $\tilde{\nabla} g = \tilde{\nabla} S = 0$ and $S(h)=h$, hence the distribution $T^1$ defines a foliation of mirrors. So, the almost product structure $T(M)=T^1(M)\oplus T^2(M)$ is also integrable.

2) It follows from (2.1),(2.3) that if $S(h)=h$ and $\tilde{\nabla}_X h = 0,\ X \in T^2$, then $S(\tilde{T}) = \tilde{T}$ and $\tilde{\nabla}_X \tilde{T} = 0$ too.

3) *Lemma 4.57* and proof of *Lemma 4.58* imply that for any point $p\in M$ there exists such an open ball $B_p$ that $\mu(x,y) = (e\tilde{x}p_x \cdot S \cdot e\tilde{x}p_x^{-1})(y)$ is defined on $B_p\times B_p$. If we consider $s_x = e\tilde{x}p_x \cdot S \cdot e\tilde{x}p_x^{-1}$, then $s_x$ is a local affine transformation of $\tilde{\nabla}$ because $S_x(\tilde{R}_x) = \tilde{R}_x$ and $S_x(\tilde{T}_x) = \tilde{T}_x$, see [46], $S$, $T^1$, $T^2$ are also invariant with respect to every $s_x$, therefore the corresponding foliations are invariant too.

4) The affinor $I-S$ has an inverse one on $T^2$ and conditions (7.13) are fulfilled on any of the maximal integral manifolds of the distribution $T^2$, hence every such a manifold is a locally regular s-manifold, see [48].

5) For each point $p\in M$ $s_p$ is identical on the connected component $\Lambda_p \cap B_p$ of the mirror containing $p$ and $s_p$ transforms mirrors onto mirrors.

6) The rest is easy the modification of the case of locally regular $s$-manifolds considered in [48].

**QED.**

$3^0$. Let $(M,\{s_x\})$ be a R.l.r. $\sigma$–m. with one-dimensional distribuion of mirrors $T^1(M)=[\xi]$, where $\xi \in \mathfrak{X}(M)$, $\|\xi\| = 1$. We denote $T^1$, $T^2$ by $L$, $V$ respectively. Further, let $S_x = (s_x)_{*x}$ have only complex eigenvalues $a_1 \pm b_1 i, \ldots, a_r \pm b_r i$ on $V$ and $D_i = \ker(S^2 - 2a_i S + I),\ i = 1,\ldots,r$. It is clear that $V = \bigoplus_{i=1}^{r} D_i$ and for every $X\in V$, $X = X_1 + \ldots + X_r$, where $X_i \in D_i,\ i = 1,\ldots,r$.

An affinor $F$ on $M$ is defined by

$$(7.15)\quad FX = \sum_{i=1}^{r} \frac{1}{b_i}(S - a_i I)X_i,\ \text{for } X\in V;\ \ F\xi = 0.$$

By similar argumentts as in *Chapter 5, §3, 6*$^0$ we obtain that for $X,Y\in V$

$$\langle FX, FY\rangle = \langle X, Y\rangle \text{ and } F^2X = -X.$$

So, $(F, \xi, \eta, g)$ is an a.c.m.s. on $M$. Let $\overline{\nabla}$ be its canonical connection defined by (5.1). Since $S\xi = \xi$, $S = a_i I + b_i F$ on $D_i$ and $\overline{\nabla} F = 0$, then $\overline{\nabla}_Y SX_i = S\overline{\nabla}_Y X_i$, $i = 1,...,r$, and $\overline{\nabla}S = \overline{\nabla}g = 0$ on *M*. It is clear that $D_i$ and *F* are invariant under any $s_x$, therefore we get from (5.1) that $\overline{\nabla}$ is also invariant under $s_x$. If $\tilde{\nabla}$ is the canonical connection of R.l.r. σ–m. defined by (4.1), then it follows from *Proposition 4.2* that $\overline{\nabla}_X = \tilde{\nabla}_X$ for $X \in V$.

From (4.1) it is obvious that $(\nabla_\xi S) = (\tilde{\nabla}_\xi S) = 0$, hence $\nabla_\xi F = 0$ too and our a.c.m.s. is almost-K-contact, see $3^0$.

Using (5.1) and formulas of $\pi_1$, $\pi_2$ we have that $\overline{\nabla}_\xi = \nabla_\xi = \tilde{\nabla}_\xi$. So, we have got the following

**THEOREM 7.21.** Let $(F, \xi, \eta, g)$ be an a.c.m.s. on $(M,\{s_x\})$ induced by a structure of the R.l.r. σ–m. as it have been shown above. If $\overline{\nabla}$ is the canonical connection of a.c.m.s. defined by (5.1) and $\tilde{\nabla}$ is that of $(M,\{s_x\})$ given by (4.1), then both the connections coincide,

$$\overline{\nabla} = \tilde{\nabla} \text{ on } M.$$

If *M* is a manifold of maximal torsion, then (7.13) are realized.
If *M* is one of minimal torsion, then (7.14) are fulfilled.

**REMARK.** The converse situations are described in *Theorems 7.19, 7.20* but we want to note that for given an a.c.m.s. a searching for a suitable affinor *S*, which have to define a structure of R.l.r. σ–m., is a sufficiently difficult problem.

$4^0$. We consider now a R.l.r. σ–m.o.3 $(M,\{s_x\})$, see *Definition 4.3*. So, $S^3=I$ and *S* has only three eigenvalues $1,\ -\frac{1}{2} \pm \frac{\sqrt{3}}{2}i$. Let $L=T^1=[\xi],\ \xi \in X(M),\ \|\xi\| = 1$, be one-dimensional distribution of mirrors corresponding to the eigenvalue *1* and $V = T^2 = L^\perp$. An affinor *F* on *M* is defined by the formula

$$(7.16)\quad F\xi = 0,\quad FX = \frac{1}{\sqrt{3}}(2S + I)X \quad \text{for } X \in V.$$

By similar arguments as in *Chapter 5, §3,* $6^0$ we easily obtain that

$F^2X = -X$ and $\langle FX, FY\rangle = \langle X, Y\rangle$ for $X,Y \in V$.

So, $(F, \xi, \eta, g)$ is an a.c.m.s. on *M* induced by affinor $S = (s_x)_{*x}$ *Theorem*

*7.21* holds in our case and it follows from (4.3) that

$$\overline{\nabla}_X Y = \frac{1}{3}(\nabla_X Y + S\nabla_X S^2 Y + S^2 \nabla_X SY), \ X,Y \in \mathfrak{X}(M).$$

Conversely, let we have an a.c.m.s. $(F, \xi, \eta, g)$ on $M$ and $\overline{\nabla}$ is its canonical connection. An affinor $S$ can be defined on $M$ by

$$(7.17) \quad S\xi = \xi, \ SX = (-\frac{1}{2}I + \frac{\sqrt{3}}{2}F)X \ \text{ for } X \in V.$$

The conditions for $(M, g, S, \overline{\nabla})$ to be a R.l.r. σ–m.o.3 of maximal or minimal torsion are described in *Theorems 7.19, 7.20*.

**REMARK.** With help of *Theorem 4.43,* where $\Lambda$ is one-dimensional Riemannian manifold, using results from [33] we can construct various examples of R.l.r. σ–m.o.3 and, therefore, those of induced a.c.m.s.

$5^0$. Let $(M,\{s_x\})$ be a R.l.r. σ–m.o.4 . Thus, $S^4=I$ and eigenvalues of $S$ are $\pm 1$, $\pm i$. We consider a case when $L=T^1=[\xi]$, $\xi \in \mathfrak{X}(M)$, $\|\xi\|=1$, $L$ is one-dimensional distribution corresponding to the eigenvalue $1$ or $-1$ and $S^2=-I$ on $V=L^{\perp}$.

An affinor $F$ on $M$ is defined by the formula

$$(7.18) \quad F\xi = 0, \quad F = S \ \text{ on } V.$$

So, $(F, \xi, \eta, g)$ is an a.c.m.s. on $M$, and *Theorem 7.21* holds, where

$$\overline{\nabla}_X Y = \frac{1}{4}(\nabla_X Y + S\nabla_X S^3 Y + S^2 \nabla_X S^2 Y + S^3 \nabla_X SY), \ X,Y \in \mathfrak{X}(M).$$

**PROPOSITION 7.22.** Let $(F, \xi, \eta, g)$ be an a.c.m.s. on $M$ induced by $S$ as it was shown above and $h = \nabla - \overline{\nabla}$

1) If $S\xi = -\xi$, then $h \in \overline{\boldsymbol{T}}_6 \oplus \overline{\boldsymbol{T}}_9$ and if $(M,\{s_x\})$ is of minimal torsion, then $h \in \boldsymbol{T}_9$
2) If $S\xi = \xi$, then $h = \overline{\boldsymbol{T}}_5 \oplus \overline{\boldsymbol{T}}_7 \oplus \overline{\boldsymbol{T}}_8 \oplus \overline{\boldsymbol{T}}_{10}$ and if $(M,\{s_x\})$ is of minimal torsion, then $h \in \overline{\boldsymbol{T}}_8 \oplus \overline{\boldsymbol{T}}_{10}$.

**Proof.** It was shown that the a.c.m.s is almost-K-contact, see $3^0$, hence

$\nabla_\xi S = \nabla_\xi F = 0$ and $h \in \bigoplus_{i \neq 11,12} \overline{\boldsymbol{T}}_i$. *Proposition 5.36* is also true in this case, therefore

a) $h_{XYZ} = < \pi_1 \nabla_X Y, Z > = 0$ for $X,Y,Z \in V$;

b) $h_{\xi X \xi} = 0$ because $h_\xi X \in V$;

c) $h_{\xi XY} = \frac{1}{2} < \nabla_\xi X + S \nabla_\xi SX, Y > = \frac{1}{2} < \nabla_\xi X - \nabla_\xi X, Y > = 0$.

Since $Sh_X Y = h_{SX} SY$, hence

$h_{FXFY\xi} = < h_{SX} SY, \xi > = < Sh_X Y, \xi > = < h_X Y, S^3 \xi >$.

If $S\xi = -\xi$, then $h_{FXFY\xi} = -h_{XY\xi}$ and $h \in \overline{\boldsymbol{T}}_6 \oplus \overline{\boldsymbol{T}}_9$.

If $S\xi = \xi$, then $h_{FXFY\xi} = h_{XY\xi}$ and $h \in \overline{\boldsymbol{T}}_5 \oplus \overline{\boldsymbol{T}}_7 \oplus \overline{\boldsymbol{T}}_8 \oplus \overline{\boldsymbol{T}}_{10}$.

The rest follows from *Lemma 7.1.*

**QED.**

Conversely, let we have an a.c.m.s. $(F, \xi, \eta, g)$ on $M$ and $\overline{\nabla}$ is its canonical connection. An affinor $S$ can be defined on $M$ by

(7.19) $S = F$ on $V$, $S\xi = -\xi$ or $S\xi = \xi$.

The conditions for $(M, g, S, \overline{\nabla})$ to be a R.l.r. σ–m.o.4 of maximal or minimal torsion are considered in *Theorems 7.19, 7.20*.

**REMARK.** Using *Theorem 4.43*, where $\Lambda$ is an one-dimensional Riemannian manifold, we can construct interesing examples of R.l.r. σ–m.o.4 with induced a.c.m.s. of various classes. A base manifold can be taken, for instance, from classification in [40].

# BIBLIOGRAPHY

The following bibliography includes all those references that are quoted in the text and it also contains some books and papers on the topics closely related to the content of our work.

1) V.A.Alexiev-G.Ganchev, On a classification of almost contact metric manifolds, pp.155-161;

2) V.A.Alexiev-G.Ganchev, On some classes of almost contact metric manifolds, pp.186-191;

Proceedings of the 15th spring conference of the Union of Bulgarian Mathematicians, Math. and Educ. in Math., (1986).

3) W.Ambrose-I.M.Singer, On homogeneous Riemannian manifolds, Duke Math. J., 25 (1958) pp.647-669.

4) J.E.D'Atri-W.Ziller, Naturally reductive metrics and Einstein metrics on compact Lie groups, Mem. Amer. Math. Soc., 215 (1979).

5) V.V.Balashchenko, Riemannian geometry of canonical structures on regular $\Phi$–spaces, Preprint N 174, Ruhr-Universitat Bochum (1994).

6) D.Bernard, Sur la geometrie differentielle des $G-$ structures, Ann. Inst.Fourier (Grenoble), V.10 (196O) pp.151-27O.

7) R.Bishop-R.Crittenden, Geometry of Manifolds, N.Y.: Academic Press, 1964.

8) D.E.Blair, Contact Manifolds in Riemannian Geometry, Lecture Notes in Math. 509, Springer Verlag, 1976.

9) D.E.Blair, On the of metrics associated to a symplectic or contact form, Bull. Inst. Math. Acad. Sinica, 11 (1983) pp.297-308.

10) D.E.Blair, Critical associated metrics on contact manifolds III, J. Austral. Math. Soc. (Series A), 50 (1991) pp.189-196.

11) D.E.Blair-L.Vanhecke, Symmetries and $\varphi$–symmetric spaces, Tohoku Math. J., 39, N 3 (1987) pp.373-383.

12) E.Boeckx, O.Kowalski, L.Vanhecke, Riemannian manifolds of conullity two, Singapore: World Sci., 1996.

13) S.S.Chern, The geometry of $G$–tructures, Bull. Amer. Math. Soc., V. 72 (1966) pp.167-219.

14) D.Chinea-C.Gonzalez, A classification of almost contact metric manifolds, Ann. Mat. pura appl., 156 (1990) pp.15-36.

15) R.S.Clark -M.Bruckheimer, Tensor structures on a differentiable manifolds, Ann. Mat. pura appl., 54 (1961) pp.123-141.

16) M.Djoric- L.Vanhecke, Almost Hermitian geometry, geodesic spheres and

symmetries, Math.J. Okayama Univ., 32 (1990) pp.187-206.
17) A.A.Ermolitski, The second quadratic form of a *G*–structure and nearly particular structures, C.R. Bulgarian Acad. Sci., 34, N 7 (1981) pp.963-964.
18) A.A.Ermolitski, The second quadratic form of a *G*–structure, VINITI, N 4789 (1981) pp.1-15 (Russian).
19) A.A.Ermolitski, The second quadratic form, the curvature and torsion of almost complex structure, VINITI, N 5612 (1981) pp.1-26 (Russian).
20) A.A.Ermolitski, A strong Kahlerian structure on differential manifold, Izv. Acad. Nauk BSSR. Ser.Fiz.-Mat. Nauk, N 1 (1981) 41-43 (Russian).
21) A.A.Ermolitski, The second quadratic form of a *G*–structure and its application to geometry of affinor structures, Dissertation, Minsk (1982) pp.1--101 (Russian).
22) A.A.Ermolitski, On the Riemannian geometry of classical structures, Izv. Acad.Nauk BSSR (Fiz.-mat.ser.), 2 (1986) pp.107-109 (Russian).
23) A.A.Ermolitski, Riemannian *G*–structures of type $T_1$, Diff.Geom. Mnogoobrazij Figur, 19 (1988) pp.62-63 (Russian).
24) A.A.Ermolitski, On a classification of Riemannian G-structures, Izv. VUZov. Mat., 12 (1988) pp.62-63 (Russian).
25) A.A.Ermolitski, On a classification of metric structures, Proceedings of the 9th Soviet Union conference at the modern geometry, Kishinev (1988) pp.111-112 (Russian).
26) A.A.Ermolitski, Periodic affinors and 2k-symmetric spaces, Dokl. Acad. Nauk BSSR, 34, N 2 (1990) pp.109-111 (Russian).
27) A.A.Ermolitski, Riemannian regular $\sigma$–manifolds, Czech. Math. J., 44 (1994) pp.57-66.
28) A.A.Ermolitski, About a classification of almost contact metric structures, Vesti BSPU, 2 (1994) pp.98-101 (Russian).
29) A.A.Ermolitski, About almost product structure, Vesti BSPU, 1 (1995) pp.92-94 (Russian).
30) A.S.Fedenko, Spaces with symmetries, Minsk: Izd BGU, 1977 (Russian).
31) F.R.Gantmaher, Theory of matrices, Moscow: Nauka, 1988 (Russian).
32) A.Gray, Nearly Kaehler manifolds, J. Diff. Geom., 4 (197O) pp.283-3O9.
33) A.Gray, Riemannian manifolds with geodesic symmetries of order 3, J. Diff. Geom., 7 N 3 - 4 (1972) pp.343-369.
34) A.Gray. The structure of nearly Kaehlerian manifolds, Math. Ann., 223 (1976) pp.233-248.
35) A.Gray - L.M.Hervella, The sixteen classes of almost Hermitian manifolds and their linear invariants, Ann. Mat. pura appl., 123 (1980) pp.35-58.

36) S.Helgason, Differential geometry and symmetric spaces, N.Y.: Academic Press, 1962.
37) L.M.Hervella - E.Vidal, Nouvelles geometries pseudo-kahleriennes $G_1$ et $G_2$, C.R.Acad. Sci. Paris, 283 (1976) pp.115-118 (French).
38) D.Janssens - L.Vanhecke, Almost contact structures and curvature tensors, Kodai Math. J., 4 (1981) pp.1-27.
39) C.Jeffries, O- deformable (1,1) tensor fields, J. Diff. Geom., 3-4, N 8 (1972) pp.575-583.
40) J.A.Jimenez, Riemannian 4-symmetric spaces, Trans.Amer.Math.Soc., 306 N 2 (1988) pp.715-734.
41) V.F.Kirichenko, Quasi homogeneous manifolds and generalised almost Hermitian structures, Izv. AN SSSR, 47 N 6 (1983) pp.1208-1223 (Russian).
42) V.F.Kirichenko, The axiom of $\Phi$–holomorphic planes in contact metric geometry, Izv. AN SSSR, 48 N 4 (1984) pp.711-734 (Russian).
43) V.F.Kirichenko, Methods of generalized Hermitian geometry in the theory of almost contact manifolds, Problems of geometry. V.18. Itogi Nauki i Tekniki. VINITI. Moscow, 1986, pp.25-71 (Russian).
44) E.J.Kobayashi, A remark on Nijenhuis tensor, Pacific J. Math., 12 (1962) pp.963-977.
45) S.Kobayashi, Transformation groups in differential geometry, Springer-Verlag, 1972.
46) S.Kobayashi- K.Nomizu, Foundations of differential geometry, N.Y.: Wiley, V.1, 1963.
47) S.Kobayashi- K.Nomizu, Foundations of differential geometry, N.Y.: Wiley, V.2, 1969.
48) O.Kowalski, Generalized symmetric spaces, Lecture Notes in Math. 8O5, Springer-Verlag, 198O.
49) O.Kowalski-F.Tricerri, Riemannian manifolds of dimension $n \leq 4$ admitting a homogeneous structure of class $T_2$, Conferenze del seminario di matematica dell universita di Bari, 222 (1987) pp.1-24.
50) O.Kowalski, F.Pruffer, L.Vanhecke, D Atri spaces. Topics in Geometry, In Memory of Joseph D Atri, Birkhauser, 1996, pp.241-284.
51) A.J.Ledger-L.Vanhecke, Symmetries and locally s-regular manifolds, Ann. Global Anal.Geom., 5 N 2 (1987) pp.151-160.
52) A.J.Ledger-L.Vanhecke, Symmetries on Riemannian manifolds, Math. Nachr., 136 (1988) pp.81-90.
53) A.J.Ledger-L.Vanhecke, Locally homogeneous S-structures, Bull. Austral. Math.Soc., 37 N 2 (1988) pp.241-246.

54) A.J.Ledger-L.Vanhecke, Naturally reductive S-manifolds, Note di Matematica, 10 N 2 (1990) pp.363-370.
55) J.Lehmann-Lejeune, Integrabilite des G-structures definies par une 1-forme O-deformable a valeurs dans le fibre tangent, Ann. Inst. Fourier (Grenoble), 16 N 2 (1966) pp.329-387 (French).
56) O.Loos, Spieglungsraume und homogene symmetrische Raume, Math. Z., 99 N 2 (1967) pp.141-170 (German).
57) O.Loos, Reflexion spaces of minimal and maximal torsion, Math. Z., 106 N 1 (1968) pp.67-72.
58) O.Loos, Symmetric Spaces, N.Y.- Amsterdam : Benjamin, 1969.
59) A.M.Naveira, A classification of Riemannian almost-product manifolds, Rend. mat. e appl., 3 N 3 (1983) pp.577-592.
60) J.Oubina, New classes of almost contact metric structures,Publicationes Mathematicae, 32 (1985) pp.187-193.
61) R.Palais, A global formulation of the Lie theory of transformation groups, Mem. Amer. Math. Soc., V.22 (1957).
62) L.V.Sabinin, About geometry of subsymmetric spaces, Nauchn. Dokl. V. Shk. Phys. Mat., (1958) (Russian).
63) A.P.Shirokov, Sructures on differentiable manifolds, A survey article in: Algebra,Topology,Geometry, V. 2, VINITI AN SSSR, 1974, pp.153-207 (Russian).
64) S.Sternberg, Lectures on differential geometry, Prentice-Hall, 1964.
65) W.P.Thurston, Some examples of symplectic manifolds, Proc. Amer. Math. Soc., 55 (1976) pp.467-468.
66) I.Tamura, Topology of foliations, Moscow: Editions Mir, 1979 (Russian).
67) F.Tricerri - L.Vanhecke, Homogeneous structures on Riemannian manifolds, London Math.Soc.Lect.Note Ser., 83 (1983).
68) V.V.Trofimov - A.T.Fomenko, Algebra and geometry of integrable Hamiltonian differential equations, Moscow: Fizmatlit Co., 1992 (Russian).
69) V.I.Vedernikov - A.S.Fedenko, Symmetric spaces and their generalizations, A survey article in: Algebra, Topology, Geometry, V. 14, VINITI AN SSSR, 1976 (Russian).
70) J.A.Wolf, Spaces of constant curvature, N.Y.: McGraw-Hill, 1967.
71) J.A.Wolf, The geometry and structures of isotropy irreducible homogeneous spaces, Acta Math., V. 12O (1968) pp.59-148.
72) K.Yano, On a structure defined by a tensor field f of type (1,1) satisfying $f^3 + f = 0$, Tensor N.S., 14 (1963) pp.99-1O9.
73) K.Yano - M.Kon, Structures on manifolds, Singapore: World Sci., 1984.

# SUBJECT INDEX

## *A*

## *C*

*D*

*E*

*F*

*G*

*H*

## *I*

## *K*

## *L*

## *M*

## *N*

## O

## P

## Q

## R

## S

## T

# LIST OF STANDARD DENOTATIONS AND ABBREVATIONS

| | |
|---|---|
| $M$ | $M$ a differentiable manifold of dimension n and of class $C^{\infty}$ ("smooth manifold") |
| $T_p(M)$, $M_p$ | the tangent space of $M$ at the point $p \in M$ |
| $T_p^*(M)$, $M_p^*$ | the cotangent space of $M$ at $p$ |
| $\varphi_*$, $\varphi_{*p}$ | the tangent mapping of a smooth mapping (on a manifold, at a point) |
| $T(M)$ | the tangent bundle of $M$ |
| $L(M)$ | the principal frame bundle of $M$ |
| $\mathcal{X}(M)$ | the Lie algebra of all smooth vector fields on $M$ |
| $P(G)$ | a $G$-structure over $M$ |
| $\Gamma$ | a connection in $P(G)$ |
| $\omega$ | the connection form of the given connection $\Gamma$ |
| ***o***, ***g***, ***h***, ***k*** | Lie algebras of the Lie groups $O(n)$, $G$, $H$, $K$ |
| ***m*** | a vector subspace of a Lie algebra |
| [ , ] | the Lie bracket |
| $L_x$ | the Lie derivative with respect to $X \in \mathcal{X}(M)$ |
| $\nabla$, $\nabla_X$ | Riemannian connection on $M$, covariant derivative with respect to $X \in \mathcal{X}(M)$ |
| $\overline{\nabla}$ | the canonical connection of the structure $(P(G),g)$ |
| $\tilde{\nabla}$ | the canonical connection of a homogeneous Riemannian space or Riemannian (locally) regular $\sigma$(or $s$)-manifold |
| $F, J, P, S$ | $O$-deformable (1,1) tensor fields (affinors) |
| $N(F)$ | Nijenhuis tensor of an affinor $F$ |
| $R, \overline{R}$ | curvature tensor fields of the connections $\nabla, \overline{\nabla}$ |
| $\overline{T}, \tilde{T}$ | torsion tensor fields of the connections $\overline{\nabla}, \tilde{\nabla}$ |
| $P(u_0)$ | the holonomy subbundle (of a connection) containing the frame $u_0 \in L(M)$ |
| $Tr(\overline{\nabla})$ | transvection group |
| a.c.m.s. | an almost contact metric structure |
| a.H.s. | an almost Hermitian structure |
| a.p.R.s. | an almost product Riemannian structure |
| f-s. | An *f*-structure |
| h.R.s. | an homogeneous Riemannian structure |
| k-s.R.s. | a *k*-symmetric Riemannian space |
| k-s.l.R.s. | a locally *k*-symmetric Riemannian space |
| R.r. σ-m. | an Riemannian regular σ-manifold |
| R.l.r. σ-m. | an Riemannian locally regular σ-manifold |
| R.r. σ-m.o.k | an Riemannian regular σ-manifold of order *k* |
| R.l.r. σ-m.o.k | an Riemannian locally regular σ-manifold of order *k* |